\documentclass[10pt]{article}

\usepackage{xcolor}
\usepackage{amsmath}
\usepackage{relsize}
\usepackage{bbold}
\usepackage{enumitem}
\usepackage{amssymb} 
\usepackage{mathrsfs}
\usepackage{caption}
\usepackage{graphicx}
\usepackage{authblk}
\usepackage[hypertexnames=false,colorlinks=true,linkcolor=blue,citecolor=blue]{hyperref}
\usepackage[numbers,comma,square,sort&compress]{natbib}
\usepackage[letterpaper,text={6.5in,9in},centering]{geometry}

\graphicspath{{eps/}{pdf/}}

\makeatletter\@addtoreset{equation}{section}\makeatother

\newtheorem{Lemma}{Lemma}[section]
\newtheorem{Theorem}{Theorem}
\newtheorem{Proposition}[Lemma]{Proposition}

\newtheorem{Remark}[Lemma]{Remark}

\def\change{\textcolor{black}}

\newenvironment{Proof}%
 {\begin{trivlist} \item[]{\bf Proof. }}%
 {\hspace*{\fill}$\rule{.4\baselineskip}{.4\baselineskip}$\end{trivlist}}

\newenvironment{Acknowledgment}%
 {\begin{trivlist}\item[]\textbf{Acknowledgments.}}{\end{trivlist}}

\begin{document}

\title{Localised Radial Patterns on the Free Surface of a Ferrofluid}
\author{Dan J. Hill}
\author{David J.B. Lloyd}
\author{Matthew R. Turner}

\affil{\small Department of Mathematics, University of Surrey, Guildford, GU2 7XH, UK}
\date{\today}
\maketitle
\begin{abstract}
\noindent
This paper investigates the existence of localised axisymmetric (radial) patterns on the surface of a ferrofluid in the presence of a uniform vertical magnetic field. We formally investigate all possible small-amplitude solutions which remain bounded close to the pattern's centre (the core region) and decay exponentially away from the pattern's centre (the far-field region). The results are presented for a finite-depth, infinite {expanse} of ferrofluid equipped with a linear magnetisation law. These patterns bifurcate at the Rosensweig instability, where the applied magnetic field strength reaches a critical threshold. Techniques for finding localised solutions to a non-autonomous PDE system are established; solutions are decomposed onto a basis which is independent of the radius, reducing the problem to an infinite set of nonlinear, non-autonomous ODEs. Using radial centre manifold theory, local manifolds of small-amplitude solutions are constructed in the core and far-field regions, respectively. Finally, using geometric blow-up coordinates, we match the core and far-field manifolds; any solution that lies on this intersection is a localised radial pattern. {Three} distinct classes of stationary radial solutions are found: {spot A and spot B solutions, which are equipped with two different amplitude scaling laws and achieve their maximum amplitudes at the core, and ring solutions, which achieve their maximum amplitudes away from the core.} These solutions correspond exactly to the classes of localised radial solutions found for the Swift-Hohenberg equation. Different values of the linear magnetisation and depth of the ferrofluid are investigated and parameter regions in which the various localised radial solutions emerge are identified. The approach taken in this paper outlines a route to rigorously establishing the existence of axisymmetric localised patterns in the future.

\end{abstract}

\section{Introduction}
\begin{figure}[t]
    \centering
    \includegraphics[height=5.5cm]{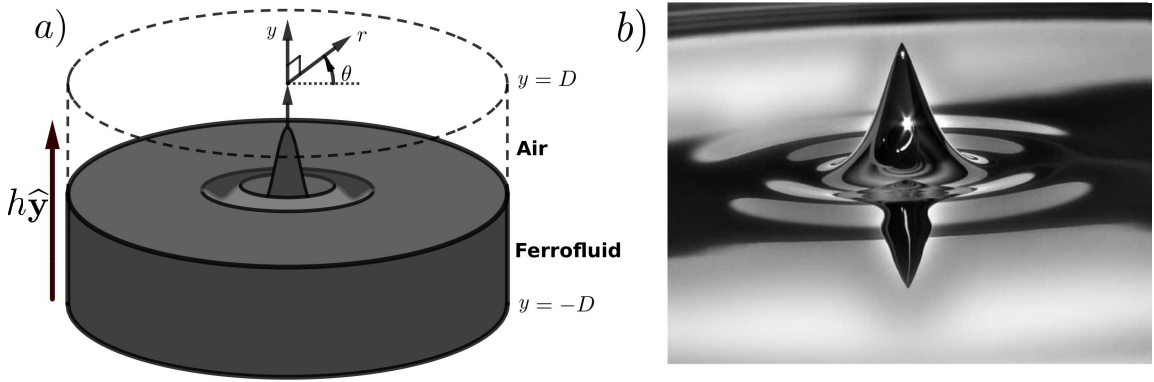}
    \caption{a) Static localised patterns appear on the free surface between a non-magnetisable fluid (transparent) and a ferrofluid (shaded) as the strength of a vertical applied magnetic field {$h$} is varied. b) Localised axisymmetric peaks, termed `spot A' solutions, have been observed experimentally (See Richter \cite{Richter2011Mountains}, reproduced with permission).}
    \label{fig:ffspot}
\end{figure}

{Ferrofluids are liquids with a colloidal suspension of magnetic nanoparticles; as a result, they can be manipulated by the application of an external magnetic field{; see} \cite{torres2014recent} \& \cite{rosensweig2013ferrohydrodynamics} for a review of ferrofluids and their associated experiments).} In an experiment, a quiescent layer of ferrofluid in a flat, horizontal container is subjected to a uniform vertical magnetic field of magnitude $h$, see Figure \ref{fig:ffspot}$a)$. As the applied magnetic induction exceeds a critical value, $h_{c}$, the homogeneous state becomes unstable to small perturbations and regular cellular spikes of patterns appear.
This phenomenon is known as the `Rosensweig' instability, and has been a subject of interest since the 1960s (see Cowley \& Rosensweig \cite{cowley1967interfacial}, Rosensweig \cite{rosensweig1987magnetic}). In 2005, Richter \& Barashenkov \cite{richter2005two} were the first to show that they could induce solitary axisymmetric spikes, as seen in Figure \ref{fig:ffspot}$b)$, using a local perturbation of the applied magnetic field. These radial `spots' persisted once the local magnetic perturbation was removed and exhibited some notable properties. The spikes emerged away from the centre of the container, were able to be moved around the container, and did not appear to feel the effects of the container boundary. This suggests the spots are well localised, possibly decaying exponentially fast to the flat state. This experimental evidence was further supported by a collection of numerical results (see Lloyd et al. \cite{lloyd2015homoclinic} for pseudo-spectral numerical methods and experimental results, and Lavrova et al. \cite{lavrova2008numerical}, Cao \& Ding \cite{cao2014formation} for finite-element simulations). Despite this, there is no analytical theory for proving the existence of exponentially localised axisymmetric solutions to the ferrohydrostatic problem.

{It is worth noting that static `ring-like' axisymmetric structures have been observed in { ferrofluid }experiments (see \cite{Knieling2007Growth,Reimann2003Oscillatory,spyropoulos2019spike}) and numerical studies (see \cite{spyropoulos2019spike,lavrova2016modeling}). However, in most cases, these solutions appear to be connected to magnetic effects occurring at the radial boundary (i.e. the wall of a finite size container). In contrast to this, we expect any localised solutions to be independent of the radial boundary, and so we complete our analysis on an infinitely wide layer of ferrofluid. These observed ring-like solutions also appear to maintain their axisymmetry relative to their container; hence, we expect these solutions to be more closely related to radial target patterns, as seen in \cite{KOPELL1981Target,Pomeau1985Axisymmetric}.}

\paragraph{Theoretical Approaches to Ferrofluids:}
The first attempt at an analytic study of domain-covering cellular patterns of the Rosensweig instability was by Gailitis \cite{gailitis1977formation}, who substituted a cellular free-surface ansatz into a hypothetical free energy  for the system (an infinite-depth ferrofluid with a linear magnetisation law). This resulted in finding regions of existence for various lattice patterns (squares, stripes and hexagons), and was later extended to finite-depth ferrofluids by Friedrichs \& Engel \cite{friedrichs2001pattern}.

Twombly \& Thomas \cite{twombly1983bifurcating} considered static, doubly-periodic solutions of the ferrohydrostatic equations near the onset of instability, extending the previous work of Zaitsev \& Shliomis \cite{Zaitsev1970Nature} to a finite-depth ferrofluid. These results, close to the bifurcation point, were supplemented by other studies of two-dimensional periodic free-surfaces; see Silber \& Knobloch \cite{silber1988pattern} for normal-form analysis, Bohlius et al. \cite{bohlius2011amplitude,Bohlius2007Adjoint} for deriving amplitude equations near onset and Groves \& Horn \cite{Horn2015Masters,Groves2018periodic} for a Dirichlet-Neumann formulation and local bifurcation theory. With the exception of the Groves \& Horn work, all these studies included a linear magnetisation law. 

The study of localised solutions to the ferrohydrostatic equations has had very little attention to date. Currently, the only rigorous result is a proof of existence for one-dimensional localised planar waves by Groves et al. \cite{groves2017pattern}. This work examined a finite-depth ferrofluid with a fully nonlinear magnetisation law, using variational methods to formulate a spatial Hamiltonian system. Using a combination of centre-manifold reduction techniques and normal-form theory, they were then able to prove the existence of spatially localised one-dimensional solutions. Notably, the nonlinearity of the magnetisation law does not qualitatively affect any of the results, and so it was concluded that a linear magnetisation law is a reasonable assumption for analytical studies.
\paragraph{Localised Radial Patterns:}
The ferrofluid problem strongly resembles that of gravity-capillary waves \cite{Iooss1992Small,buffoni1996plethora}; as a result, water wave methodologies provide a robust framework to study ferrofluids, notably in the case of ferrofluid jets \cite{blyth2014solitary,Groves2018Jet}. However, localised solutions in water waves are created via the introduction of a uniform flow, creating an inherent symmetry-breaking that makes localised radial solutions {unphysical}. Hence, new techniques must be developed in order to find localised radial solutions in the ferrofluid problem.

\begin{figure}[t!]
    \centering
 \includegraphics[width=\linewidth]{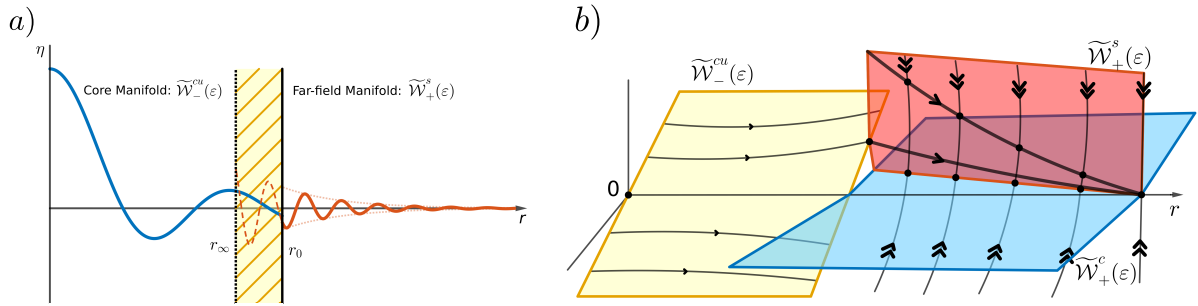}
    \caption{a) The free surface of a localised radial pattern lies on the intersection of solutions that are bounded as $r\to0$, called the `core manifold' $\widetilde{\mathcal{W}}^{cu}_{-}{(\varepsilon)}$, and solutions that decay exponentially as $r\to\infty$, called the `far-field manifold' $\change{\widetilde{\mathcal{W}}^{s}_{+}(\varepsilon)}$. b) The far-field manifold is constructed as perturbations from the decaying solutions that lie on the autonomous centre manifold $\change{\mathcal{W}^{c}_{+}(\varepsilon)}$. The mappings from $\change{\mathcal{W}^{c}_{+}(\varepsilon)}$ to $\change{\widetilde{\mathcal{W}}^{s}_{+}(\varepsilon)}$ are derived from foliation theory, as seen in \cite{mcquighan2014oscillons}.}
    \label{fig:match}
\end{figure}

Recently, there has been some progress in the study of localised radial patterns in reaction diffusion systems; for example, in 2003 Scheel \cite{scheel2003radially} introduced radial centre manifold and normal-form theory for stationary and time-periodic axisymmetric patterns in general reaction diffusion systems. Prototypical pattern forming systems have been analysed to prove the existence of localised radial solutions to the Swift-Hohenberg equation \cite{lloyd2009localized,mccalla2013spots,Mccalla2010snaking}, the complex Ginzburg-Landau equation \cite{mcquighan2014oscillons}, and others \cite{faye2013localized}. However, the non-autonomous nature of radial patterns complicates any theoretical understanding one may have for two- or three-dimensional patterns (as discussed in \cite{knobloch2008spatially}). In particular, the standard centre manifold reduction procedure \cite{Haragus2011Bifurcation,vanderbauwhede1992center,Mielke1986Reduction} for small amplitude patterns does not apply. One instead has to construct two sets of small amplitude solutions (one set where solutions remain bounded as $r\rightarrow0$ which we call the core manifold, $\widetilde{\mathcal{W}}_-^{cu}$, and the other set where solutions decay exponentially to the flat state as $r\rightarrow\infty$, which we call the far-field manifold, $\mathcal{W}_+^s$) and look for intersections of these two sets, see Figure~\ref{fig:match}(a), which we describe in more detail below.  For each problem, these sets must be constructed on a case-by-case basis.  

The present contribution is strongly motivated by the results found for the Swift-Hohenberg equation by Lloyd \& Sandstede \cite{lloyd2009localized}, who proved the existence of localised rings and elevated  spots using radial centre-manifold and normal-form theory of Scheel \cite{scheel2003radially}. Following this, McCalla \& Sandstede \cite{mccalla2013spots} employed geometric blow-up methods to prove the existence of a class of localised depressed spots, called spot B.

\paragraph{Overview:}
In the present paper we formally identify all possible small-amplitude localised radial patterns that bifurcate from the flat state in the ferrohydrostatic problem, providing a framework that guides future rigorous results on the existence of localised axisymmetric patterns. The expected properties of a localised radial solution $\mathbf{u}(\mathbf{x})$, where $\mathbf{x}=(r, \theta, y)$ are cylindrical polar coordinates, are the following:
\begin{enumerate}
    \item Axisymmetric: $\mathbf{u}(\mathbf{x})$ is independent of $\theta$, i.e. $\mathbf{u}(\mathbf{x})=\mathbf{u}(r,y)$.
    \item Non-singular: $\mathbf{u}(r,y)$ remains smooth, single-valued and bounded as $r\to0$.
    \item Localised: $\mathbf{u}(r,y)\to\mathbf{u}_{0}$ exponentially as $r\to\infty$, where $\mathbf{u}_{0}$ is the flat state of the system.
\end{enumerate}
Property 1 will inform our formulation of the radial ferrohydrostatic problem; we construct a free-surface problem with two immiscible fluids of equal arbitrary depth $D$ separated by an interface at $y=\eta(r)$, see Figure \ref{fig:ffspot}a), with a constant magnetic field of strength $h$ applied vertically. There exists a critical magnetic field strength, $h_{c}$, such that for values of $h>h_c$ the flat quiescent state destabilises and {domain-covering} patterns form, known as the Rosensweig instability. {We expect localised solutions to emerge for a sub-critical magnetic field strength, and so} we define $\widehat{\varepsilon}:={h_{c}-h}$ to be the bifurcation parameter of the system; the bifurcation point should occur at $\widehat{\varepsilon}=0$ and localised solutions should emerge for $0<{\widehat{\varepsilon}}\ll1$. The ferrofluid, taken to be the lower fluid, is assumed to have a linear magnetisation law, such that $\mathbf{M} = (\mu-1)\mathbf{H}$, where $\mu>1$ is the magnetic permeability of the ferrofluid. Using variational methods, we have formulated the radial ferrohydrostatic problem, for finite depth and a linear magnetisation law, as a quasilinear non-autonomous PDE of the form
\begin{align}
    \frac{\textnormal{d}}{\textnormal{d}r}\mathbf{u} &= \mathbf{L}\left(r\right)\mathbf{u} + \mathcal{F}(\mathbf{u},\varepsilon, r),\label{intro:eqn}
\end{align}
where $\mathbf{L}(r)$ is a non-autonomous linear differential operator containing $\partial_{yy}$ terms, $\varepsilon$ is a rescaled bifurcation parameter such that $\varepsilon \propto\widehat{\varepsilon}$, and $\mathcal{F}$ contains all the nonlinear terms of the system. The problem is non-autonomous due to the presence of $1/r$ terms that also create an apparent singularity at $r=0$. However, in the limit as $r\to\infty$, the system (\ref{intro:eqn}) reduces to the problem studied by Groves et al. \cite{groves2017pattern}. 

\begin{figure}[t!]
    \centering
    \includegraphics[height=6cm]{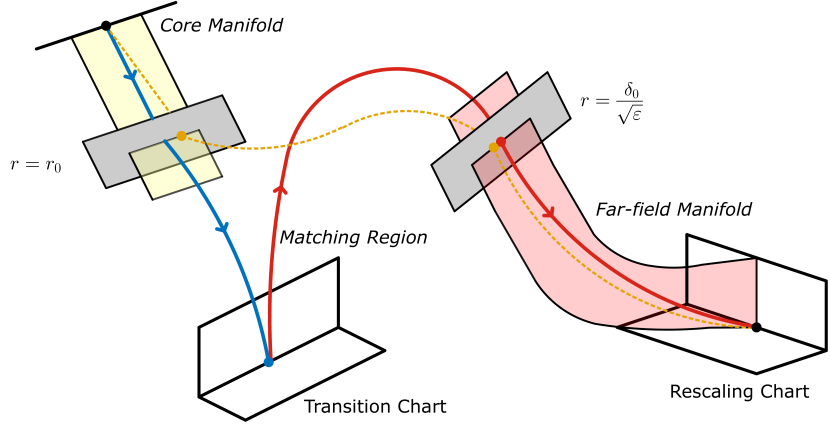}
    \caption{A sketch of the matching process in the geometric blow-up coordinates. The centre coordinates of the far-field manifold are written in `rescaling' coordinates, known as the `rescaling chart', in order to identify solutions that decay exponentially as $r\to\infty$. These solutions are then tracked back to the point $r=\delta_{0}\varepsilon^{-\frac{1}{2}}$, where $\delta_{0}>0$ is fixed. In order to retain control of the point $r=\delta_{0}\varepsilon^{-\frac{1}{2}}$ as $\varepsilon\to0$, solutions are tracked in `transition' coordinates, known as the `transition chart', to the point $r=r_{0}$, where they are matched with the core manifold.}
    \label{fig:geom}
\end{figure} 

Properties 2 \& 3 will be satisfied for local regions of $r$, specifically $r\in[0, r_{0}]$ and $r\in[r_{\infty}, \infty)$, called the `core' and `far-field', respectively. Here $r_{0},r_{\infty}>0$ are fixed constants, and $r_{0}$ is chosen to be larger than $r_{\infty}$ such that these regions overlap, see Figure \ref{fig:match}$a)$. In order to analyse \eqref{intro:eqn} in these regions, we decompose our solution onto a $y$-dependent eigenbasis, found from the spectrum of the linear operator $\mathbf{L}(r)$ in the limit $r\to\infty$. The eigenvalues of the matrix $\mathbf{L}_{\infty}:=\lim_{r\to\infty}\mathbf{L}(r)$ are the roots of the linear dispersion relation
\begin{align}
   \widetilde{\Delta}(\lambda) &:= \frac{\mu(\mu-1)^{2}}{(\mu+1)}\lambda \sin^2(\lambda D) - \left(\lambda^2 - \frac{\rho_{0}g\sigma}{\mu^{2} h^{4}}\right)\sin(\lambda D)\cos(\lambda D),\label{intro:disp}
\end{align}
where $\rho_{0}, g, \sigma$ are the relative density of the two fluids, gravitational constant, and surface tension, respectively. The autonomous linear problem $\mathbf{u}_{r} = \mathbf{L}_{\infty}\mathbf{u}$, related to \eqref{intro:eqn} as $r\to\infty$, exhibits a Hamiltonian-Hopf bifurcation at $h=h_{c}$. Then, at the bifurcation point $\varepsilon=0$, a complex number $\lambda\in\mathbb{C}$ is an eigenvalue of $\mathbf{L}_{\infty}$ if and only if $\lambda \in \left\{\pm\textnormal{i}k\right\} \cup\left\{\pm\lambda_{n}\right\}_{n\in\mathbb{N}}$, where $k,\lambda_{n}\in\mathbb{R}$, for $n\in\mathbb{N}$. The pair of imaginary eigenvalues $\{\pm\textnormal{i}k\}$ have double algebraic multiplicity, and we will treat the wavenumber $k$ as a parameter. By projecting onto the respective eigenvectors for each eigenvalue, the PDE system reduces to an infinite set of non-autonomous ordinary differential equations. Writing the full solution in this `spectral' decomposition, we can then solve the system in both the `core' and `far-field' regions, subject to Properties 2 and 3, respectively.

We first construct the set of all small-amplitude solutions that remain bounded as $r\to0$, which we call the `core' manifold and denote as $\widetilde{\mathcal{W}}^{cu}_{-}{(\varepsilon)}$, using notation from dynamical systems theory for a centre-unstable manifold. We note that this is a submanifold of the standard local centre-unstable manifold $\mathcal{W}^{cu}_{-}$, containing the set of all solutions that grow at most algebraically as $r\to0$. Then, we construct the set of all solutions that exponentially decay to the flat state as $r\to\infty$, called the `far-field' manifold and denoted as $\change{\widetilde{\mathcal{W}}^{s}_{+}(\varepsilon)}$, using the notation for a standard stable manifold. In order to perform our intended analysis, we parametrise the `far-field' manifold with respect to the decaying elements on an autonomous centre manifold $\change{\mathcal{W}^{c}_{+}(\varepsilon)}$, see Figure \ref{fig:match}$b)$, using the theory of foliations seen in \cite{mcquighan2014oscillons} for the Ginzburg-Landau equation.

Next we need to isolate the decaying solutions of $\mathcal{W}^{c}_{+}$, and so we introduce blow-up coordinates seen in \cite{mccalla2013spots} in order to restrict $\mathcal{W}^{c}_{+}$ to just its elements that decay as $r\to\infty$. More specifically, we will track solutions through two coordinate charts, called the `rescaling chart' and the `transition chart'. By rescaling $r=s/\sqrt{\varepsilon}$, where $\varepsilon$ is the bifurcation parameter of the system, we construct the `rescaling chart' which allows us to establish a parametrisation for decaying solutions that is valid for all $s\in[\delta_{0}, \infty)$ for fixed $\delta_{0}>0$. At the bifurcation point, i.e. for $\varepsilon=0$, we lose control of the point $r=\delta_{0}/\sqrt{\varepsilon}$; we construct the `transition chart' which allows us to maintain control of the solution as $\varepsilon\to0$, see Figure \ref{fig:geom}.

Once the `core' and `far-field' manifolds have both been constructed, we evaluate both at the large, but finite, value $r=r_{0}$ and apply asymptotic matching to identify any intersections between the two. Any solution that lives on the intersection of the `core' and `far-field' manifolds satisfies properties 1, 2 \& 3 and so is a small-amplitude localised radial pattern.

Through this procedure, we identify {three} types of localised radial patterns; {up- and down-spot A solutions, up- and down-spot B solutions, and up- and down-ring solutions,}. Here, {`A' and `B' refer to distinct scaling laws for the amplitude of the free-surface with respect to $\varepsilon$,} the prefix `up-' or `down-'  indicate whether the maximum point on the free surface is an elevation or a depression, and the class `spot' or `ring' indicates whether that maximum is at, or away from the core respectively. 

\paragraph{Main Results:}
The primary results of this work are related as follows; firstly, we formally show the existence of up-spot solutions, termed `spot A' in \cite{mccalla2013spots}, that were identified experimentally in \cite{richter2005two}:

\paragraph{Existence of Spot A:}
\textit{
Fix $\mu>1$, $\delta_{0}>0$, and {$D,k>0$} such that
\begin{align}
    {k D \neq \log(\sqrt{2} +1),}\label{cond:nu}
\end{align}
then for $0<\varepsilon\ll1$, the ferrohydrostatic equation \eqref{intro:eqn} has a stationary localised axisymmetric solution $\mathbf{u}(r,y)$. These solutions remain close to the trivial state $\mathbf{u}\equiv0$, and, for fixed $r_{0}>0$, the profile of the height of the interface is given by
\begin{align}
\eta_{A}(r) = \left\{ \begin{array}{cc}
    \displaystyle \varepsilon^{\frac{1}{2}}{\textnormal{sgn}\;(\nu)\frac{2\sqrt{c_{0}}}{m|\nu|}}\sqrt{\frac{k\pi}{2}} J_{0}(k r)+ \textnormal{O}(\varepsilon),& 0\leq r \leq r_{0},\\\,\\
     \displaystyle \varepsilon^{\frac{1}{2}}{\textnormal{sgn}\;(\nu)\frac{2\sqrt{c_{0}}}{m|\nu|}} \frac{1}{\sqrt{r}}\cos\left(k r - \frac{\pi}{4}\right)+ \textnormal{O}(\varepsilon),& r_{0} \leq r\leq  \delta_{0}\varepsilon^{-\frac{1}{2}}, \\\,\\
     \displaystyle \varepsilon^{\frac{1}{2}}{\textnormal{sgn}\;(\nu)\frac{2\sqrt{c_{0}}}{m|\nu|}} \textnormal{e}^{\sqrt{c_{0}}(\delta_{0} - \sqrt{\varepsilon} r)}  \frac{1}{\sqrt{r}}\cos\left(k r - \frac{\pi}{4}\right)+ \textnormal{O}(\varepsilon),& r \geq \delta_{0}\varepsilon^{-\frac{1}{2}},
\end{array}\right.  \nonumber
\end{align}
as $\varepsilon\to0$ uniformly for constants $c_{0}$, $\nu$, and $m$ defined in \eqref{c0}, \eqref{nu:defn}, \eqref{const:defn}, respectively.}

\begin{figure}[t!]
    \centering
   \includegraphics[width=0.6\linewidth]{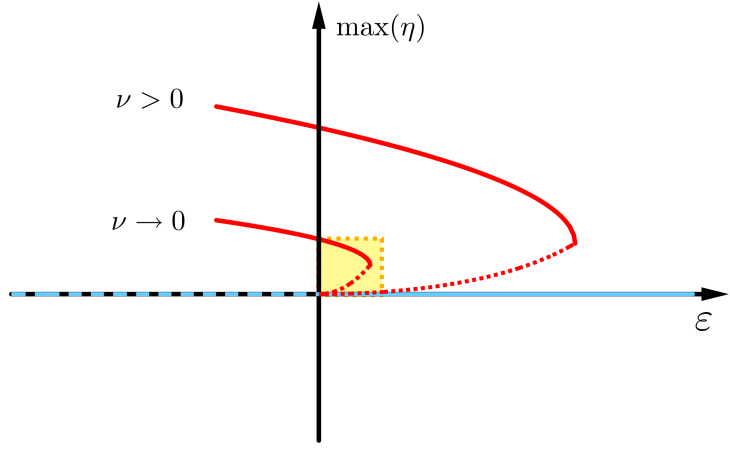}
    \caption{{Schematic of the bifurcation diagram for an up-spot A solution, where the maximum amplitude of the free surface $\eta(r)$ is plotted against values of the bifurcation parameter $\varepsilon$. Results found in this paper are only for small-amplitude solutions and small values of $\varepsilon$, illustrated by the yellow region. For $\nu>0$, our method only captures the spot A solution (red) bifurcating from the trivial state (light blue); for $\nu\to0$, however, we are able to capture the behaviour of a fold bifurcation. We postulate that this behaviour indicates a change in stability for the spot A solution; here, stability is indicated by dashed lines for unstable solutions and solid lines for stable solutions.}}
    \label{fig:fold}
\end{figure}

{We note that spot A solutions only require $k D \neq \log(\sqrt{2} +1)$ and do not require stripes to bifurcate subcritically; this is equivalent to the condition \eqref{cond:c3} and is required for localised rings to emerge}. {From the definition of $\nu$ in \eqref{nu:defn}, one can find that $\textnormal{sgn}(\nu) = \textnormal{sgn}(k D - \log(\sqrt{2} +1))$ and so the restriction {$k D \neq \log(\sqrt{2} +1)$ is equivalent to saying that $\nu\neq 0$};} see Figure \ref{fig:intro;spota}a) for a plot of {the sign of $\nu$ in $(k,D)$-space. When $\nu>0$, the} solution has its maximum at $r=0$ and behaves like the Bessel function $J_{0}(k r)$ near to the core, as depicted in Figure \ref{fig:intro;spota}b). {When $k D < \log(\sqrt{2} +1)$, the sign of the profile of spot A flips such that there is a depression at the core.} {Following the results found for the 2D radial Swift-Hohenberg equation seen in \cite{lloyd2009localized,Mccalla2010snaking}, we expect solution branches bifurcating from the trivial state to be unstable for $0<\varepsilon\ll1$. However, spot A solutions in the Swift-Hohenberg equation were found to restabilise at a fold bifurcation for moderate values of $\varepsilon$; this behaviour was captured for small $\varepsilon$ by incorporating the limit as $\nu\to0$ (see Figure \ref{fig:fold}). Hence, we also show the existence of a spot A solution when $\nu\approx0$.
\begin{figure}[t!]
    \centering
  \includegraphics[width=\linewidth]{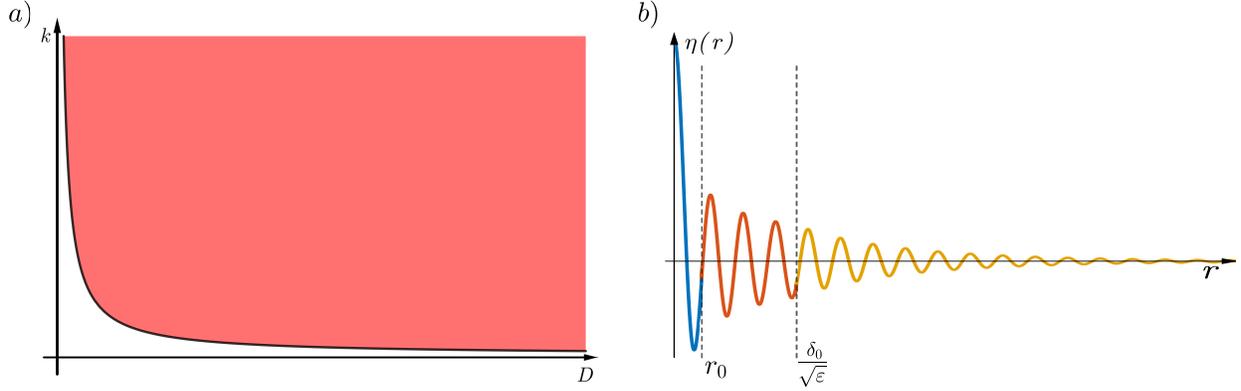}
    \caption{a) The sign of {$\nu$} is plotted for wavenumber $k>0$ and depth $D>0$, where {$\nu>0$} only in the shaded region. {Up-spot A} solutions emerge in the shaded region, i.e. when {$\nu>0$}, {and down-spot A solutions emerge in the white region}. b) The profile of the free surface {$\eta(r)$} for {up-}spot A solutions{; down-spot A solutions have the mirrored profile $-\eta(r)$}. Solutions have the following behaviour: they behave like the Bessel function $J_{0}(k r)$ in the core region (blue), maintain algebraic decay in the transition region (red), and exhibit exponential decay away from the core (yellow).}
    \label{fig:intro;spota}
\end{figure}
\paragraph{Existence of Spot A for $\nu\approx0$:}
\textit{
Fix $(D,M_{0})$ such that $\nu\approx0$ for $\nu$ defined in \eqref{nu:defn}. Then, spot A solutions undergo a fold bifurcation along the curve $\nu=\pm \varepsilon^{\frac{1}{4}}\left|\log\left(\varepsilon\right)\right|^{\frac{1}{2}}\sqrt{2c_{3}\sqrt{c_{0}}}+\textnormal{o}(1)$, where constants $c_{0}$ and $c_{3}$ are defined in \eqref{c0} and \eqref{c3}, respectively. These solutions remain close to the trivial state $\mathbf{u}\equiv0$, and, for fixed $r_{0}>0$, the profile of the height of the interface in the core region is given by
\begin{align}
\eta_{A}(r) = \textnormal{sgn}\;(\nu)\; \varepsilon^{\frac{1}{4}}|\log(\varepsilon)|^{-\frac{1}{2}}\frac{2}{m}\sqrt{\frac{2\sqrt{c_{0}}}{c_{3}}}\sqrt{\frac{k\pi}{2}}J_{0}(kr) + \textnormal{o}(1), \qquad \qquad 0\leq r\leq r_{0},\nonumber
\end{align}
as $\varepsilon\to0$ uniformly.}}

{We postulate that the existence of a fold bifurcation is indicative of a change in the stability of the spot A solution, thus allowing it to be observed experimentally.} We also show the existence of a class of up-ring and down-ring solutions, where the pattern's maximum exists away from their centre. These solutions have been found in the Swift-Hohenberg equation \cite{lloyd2009localized} {where the amplitude scales like} $\varepsilon^{\frac{3}{4}}$. To establish the existence of ring solutions, we first introduce re-scaled parameters $ k_{D}:=kD$ and $M_{0}:=\frac{\mu-1}{\mu+1}$; we find a restriction on $ k_{D}>0$ in terms of a fixed $0\leq M_{0}<1$ such that ring solutions emerge from the flat state.

\paragraph{Existence of Rings:}
\textit{
Fix $D>0$, $\delta_{0}>0$, and $( k_{D}, M_{0})$ such that
\begin{align}
    \widetilde{c}_{3} &:= -\left[\frac{ k_{D} \mathcal{M}^{2} M_{0}^{2}}{4}\left(\frac{(\cosh(4 k_{D}) - 4\cosh(2 k_{D}) - 3)(4 \,\textnormal{sech}(4 k_{D}) + \,\textnormal{sech}^{2}( k_{D}) - 2)}{(2 k_{D} \mathcal{M}\tanh(2 k_{D}) - \widetilde{\Upsilon}_{0} - 4 k_{D}^{2}) \cosh(2 k_{D})} + \frac{\,\textnormal{sech}^{2}( k_{D})}{\widetilde{\Upsilon}_{0}}\right)\,\textnormal{sech}^{2}( k_{D}) \right. \nonumber\\ & \qquad\qquad\qquad \left. + \frac{3  k_{D}}{2} + 4\mathcal{M}\left(M_{0}^{2} \,\textnormal{sech}(2 k_{D}) - \cosh^{2}( k_{D})\right)\textnormal{cosech}(2 k_{D})\right]<0,\label{cond:c3}
\end{align}
with
\begin{align}
    \mathcal{M} := \frac{2 k_{D}}{\tanh( k_{D}) +  k_{D}\,\textnormal{sech}^{2}( k_{D})}, \qquad \widetilde{\Upsilon}_{0} :=  k_{D}^{2}\frac{\tanh( k_{D}) -  k_{D}\,\textnormal{sech}^{2}( k_{D})}{\tanh( k_{D}) +  k_{D}\,\textnormal{sech}^{2}( k_{D})},\nonumber
\end{align}
then, for $0<\varepsilon\ll1$, the ferrohydrostatic equation \eqref{intro:eqn} has stationary localised axisymmetric solutions $\mathbf{u}^{\pm}(r,y)$. These solutions remain close to the trivial state $\mathbf{u}^{\pm}\equiv0$, and, for fixed $r_{0}>0$, the profile of the height of the interface is
\begin{align}
\eta^{\pm}_{R}(r) = \left\{ \begin{array}{cc}
    \displaystyle \pm \varepsilon^{\frac{3}{4}}\frac{2}{m} q_{0}\sqrt{\frac{k\pi}{2}}\left[ r J_{1}(k r) + \widetilde{b}_{D} J_{0}(k r)\right] + \textnormal{O}(\varepsilon),&   0\leq r \leq r_{0},\\\,\\
     \displaystyle \pm \varepsilon^{\frac{3}{4}}\frac{2}{m} q_{0}\left[ r^{\frac{1}{2}} \sin\left(k r-\frac{\pi}{4}\right) + \widetilde{b}_{D}\, r^{-\frac{1}{2}}\cos\left(k r - \frac{\pi}{4}\right)\right] + \textnormal{O}(\varepsilon),&   r_{0} \leq r\leq  \delta_{0}\varepsilon^{-\frac{1}{2}},\\\,\\
     \displaystyle \pm \varepsilon^{\frac{3}{4}} \frac{2}{m} \left[ \varepsilon^{-\frac{1}{4}}q(\varepsilon^{\frac{1}{2}}r) \sin\left(k r - \frac{\pi}{4}\right) + \widetilde{b}_{D}\, \varepsilon^{\frac{1}{4}} \,p\left(\varepsilon^{\frac{1}{2}} r\right)\cos\left(k r - \frac{\pi}{4}\right)\right] + \textnormal{O}(\varepsilon),&   r \geq \delta_{0}\varepsilon^{-\frac{1}{2}},
\end{array}\right.  \nonumber
\end{align}
as $\varepsilon\to0$ uniformly, where $\widetilde{b}_{D}:= b_{D} - D\tanh(k D) + k^{-1}$ and $m$, $b_{D}$, $q(\varepsilon^{\frac{1}{2}}r)$, and $p(\varepsilon^{\frac{1}{2}} r)$ are defined in \eqref{const:defn}, \eqref{bD:defn}, \eqref{q:defn}, and \eqref{pr:defn}, respectively.}

\begin{figure}[t!]
    \centering
   \includegraphics[width=\linewidth]{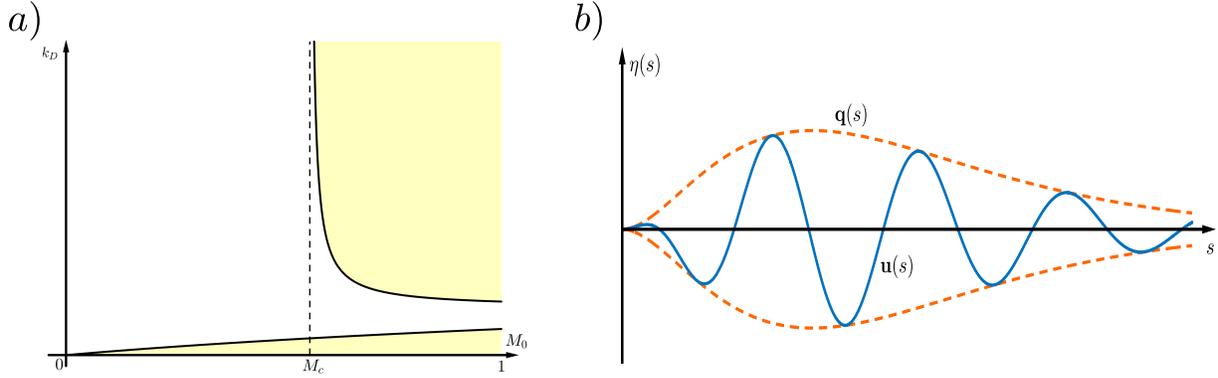}
    \caption{a) The sign of $\widetilde{c}_{3}$ is plotted for the scaled wavenumber $k_{D}:=kD$, and magnetic permeability $M_{0}:=\frac{\mu-1}{\mu+1}$, where $\widetilde{c}_{3}<0$ only in the shaded regions. {The lower shaded region corresponds to shallow-depth solutions, where our choice of a constant linear magnetisation law in a free surface is no longer physically relevant, and so we disregard this region in our analysis. This is discussed further at the end of \S\,\ref{subs:connect;P+}.} {Ring and spot B solutions emerge in the upper shaded region.} b) A sketch of the envelope function $q(s)$, which exhibits algebraic growth for small $s$ and exponential decay for large $s$. Any solution $\mathbf{u}(s)$ that lies close to the connecting orbit $\mathbf{q}_{+}(s)$ is bounded by $q(s)$.}
    \label{fig:intro;connectorq}
\end{figure}

\begin{figure}[t!]
    \centering
   \includegraphics[width=\linewidth]{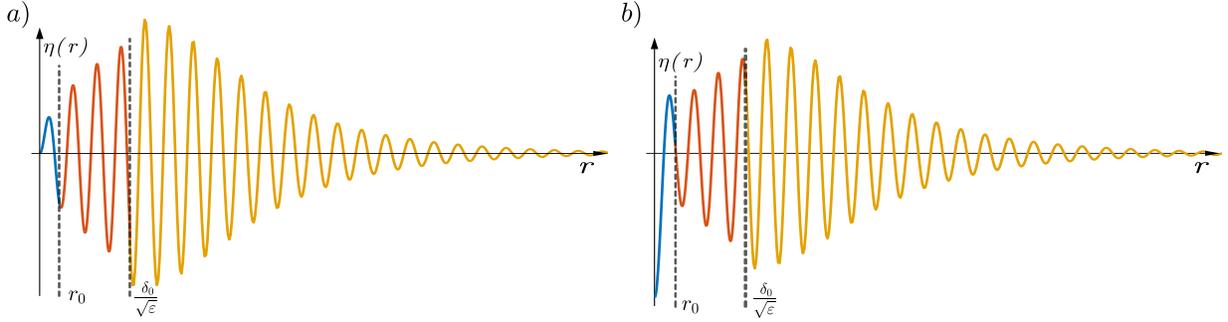}
    \caption{Profiles of the free surface are illustrated for a) up-ring, and b) {down-spot B} solutions. In both cases, solutions are bounded by $q(\varepsilon^{\frac{1}{2}}r)$ in the far-field region (yellow). In the core region (blue), solutions take the form a) $\eta(r)\thicksim r J_{1}(k r)$, b) $\eta(r)\thicksim - J_{0}(k r)$, and both solutions have algebraic growth in the transition region (red). Down-ring solutions are the reflection of up-ring solutions in the $r$-axis, {and up- or down-spot B} solutions have a much larger amplitude than {up- or down-spot A} solutions.}
    \label{fig:intro;spotbringprof}
\end{figure}

As before, we can numerically plot the parameter space in which ring solutions emerge from the flat state, as shown in Figure \ref{fig:intro;connectorq}a). The restriction $\widetilde{c}_{3}<0$ is identical to the one-dimensional problem, where $\widetilde{c}_{3}<0$ is the condition for homoclinic solutions to emerge. These solutions are restricted by a homoclinic envelope function $\mathbf{q}(\sqrt{\varepsilon}r)$, which exponentially decays as $r\to\infty$; an illustration of this function is depicted in Figure \ref{fig:intro;connectorq}b), while the free-surface profile of a ring solution is plotted in Figure \ref{fig:intro;spotbringprof}a). The ring solutions identified in this paper are closely related to the one-dimensional homoclinic solutions in \cite{groves2017pattern}, corresponding to planar waves in the three-dimensional ferrohydrostatic problem. 

The two previous classes of solutions, spot A and rings, were found for the prototypical Swift-Hohenberg {equation} by Lloyd \& Sandstede in 2009 \cite{lloyd2009localized}, and {their amplitudes were found to} have `standard' scaling laws $\varepsilon^{\frac{1}{2}}$ and $\varepsilon^{\frac{3}{4}}$, respectively. However, in 2013 McCalla \& Sandstede \cite{mccalla2013spots} also showed the existence of a down-spot solution, with a depression at its centre. {The amplitude of} this solution, termed `spot B', was found to have an unexpected scaling of $\varepsilon^{\frac{3}{8}}$ in the planar problem (and $\varepsilon^{\frac{1}{4}}$ for the three-dimensional equation). We show the existence of this spot B solution for the radial ferrohydrostatic problem and establish conditions such that spot B solutions emerge from the flat state.

\begin{figure}[t!]
    \centering
   \includegraphics[width=0.6\linewidth]{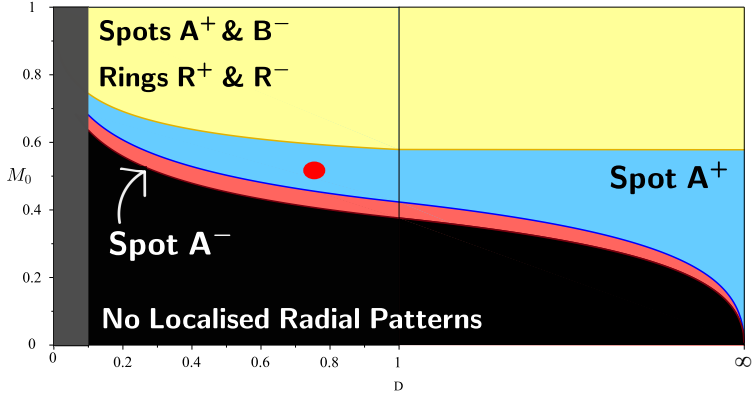}
    \caption{{Existence regions for different classes of localised radial patterns are plotted in terms of $D$, $M_{0}$. For $D>1$, the x-axis scales like $\frac{D-1}{D+1}$. The $\pm$ superscript notation refers to the sign of the pattern's maximum value, so $A^{+}$ and $R^{+}$ denote up-spot A and up-ring solutions, whereas $A^{-}$, $B^{-}$, and $R^{-}$ denote down-spot A, down-spot B, and down-ring solutions. The red circle represents the approximate parameter region in \cite{richter2005two} in which localised up-spot patterns were seen experimentally. The grey area represents small depth solutions where our choice of a constant linear magnetisation law in a free surface model is no longer reasonable.}}
    \label{fig:intro;paramspace}
\end{figure}

\paragraph{Existence of Spot B:}
\textit{
Fix $D,k, \delta_{0}>0$ and $0<M_{0}<1$ such that { \eqref{cond:nu} \& \eqref{cond:c3} are satisfied}, where $ k_{D}=kD$. Then, for $0<\varepsilon\ll1$, the ferrohydrostatic equation \eqref{intro:eqn} has a stationary localised axisymmetric \change{solution $\mathbf{u}(r,y)$}. This solution remains close to the trivial state \change{$\mathbf{u}\equiv0$}, exhibits exponential decay as $r\to\infty$, and, for fixed $r_{0}>0$, the profile of the height of the interface is
\begin{align}
\change{\eta_{B}}(r) = \left\{ \begin{array}{cc}
       \change{-\varepsilon^{\frac{3}{8}}\textnormal{sgn}\,\left(\nu\right)}\frac{2}{m} \sqrt{\frac{q_{0}}{|\nu|}} \sqrt{\frac{k\pi}{2}} J_{0} (kr) + \textnormal{O}(\varepsilon^{\frac{1}{2}}), &   0\leq r \leq r_{0}, \\\,\\  \change{-\varepsilon^{\frac{3}{8}}\textnormal{sgn}\,\left(\nu\right)}\frac{2}{m} \sqrt{\frac{q_{0}}{|\nu|}} r^{-\frac{1}{2}} \cos(kr - \frac{\pi}{4}) + \textnormal{O}(\varepsilon^{\frac{1}{2}}), &   r_{0}\leq r \leq \frac{\delta_{2}\varepsilon^{-\frac{3}{8}}}{\sqrt{q_{0}|\nu|}(1-\delta_{1})(1-\delta_{2})},\\\,\\
            \change{-\textnormal{sgn}\,\left(\nu\right)\frac{2}{m}}\bigg[\varepsilon^{\frac{3}{4}}q_{0} (1-\delta_{1})r^{\frac{1}{2}} + \varepsilon^{\frac{3}{8}}{\sqrt{\frac{q_{0}}{|\nu|}}} r^{-\frac{1}{2}}\bigg]\cos\left(k r - \frac{\pi}{4}\right) + \textnormal{O}(\varepsilon^{\frac{3}{4}}),& \frac{\delta_{2}\varepsilon^{-\frac{3}{8}}}{\sqrt{q_{0}|\nu|}(1-\delta_{1})(1-\delta_{2})} \leq r \leq \frac{\varepsilon^{-\frac{3}{8}}}{\delta_{1}\sqrt{q_{0}|\nu|}}, \\\,\\   \change{-\varepsilon^{\frac{3}{4}} \textnormal{sgn}\,\left(\nu\right)}\frac{2}{m} q_{0} r^{\frac{1}{2}} \cos\left(k r - \frac{\pi}{4}\right) + \textnormal{O}(\varepsilon),& \frac{\varepsilon^{-\frac{3}{8}}}{\delta_{1}\sqrt{q_{0}|\nu|}} \leq r \leq \delta_{0}\varepsilon^{-\frac{1}{2}}, \\\,\\
       \change{-\varepsilon^{\frac{1}{2}}\textnormal{sgn}\,\left(\nu\right)}\frac{2}{m} q(\varepsilon^{\frac{1}{2}}r) \cos\left(k r - \frac{\pi}{4}\right) + \textnormal{O}(\varepsilon),&   r \geq \delta_{0}\varepsilon^{-\frac{1}{2}},
\end{array}\right.  \nonumber
\end{align}
    as $\varepsilon\to0$ uniformly, where $m, b_{D}>0$, $\nu$, and $q(\varepsilon^{\frac{1}{2}}r)$ are defined in \eqref{const:defn}, \eqref{nu:defn}, and \eqref{q:defn}, respectively, and $\delta_{1},\delta_{2}>0$ are appropriate small constants.}
The free-surface profile of this spot B solution is plotted in Figure \ref{fig:intro;spotbringprof}b). By taking a fixed depth $D>0$, we can plot the restriction of the rescaled wavenumber $ k_{D}$, as seen in Figure \ref{fig:intro;connectorq}a), as a restriction of the wavenumber $k$ in terms of the rescaled magnetic permeability $M_{0}$. Then, we can combine Figures \ref{fig:intro;spota}a) \& \ref{fig:intro;connectorq}a) into one plot for the existence regions of all localised radial solutions in {$(D, M_{0})$-parameter space, as shown in Figure \ref{fig:intro;paramspace}. We note that restriction \eqref{cond:c3} also implies that $\nu>0$, and so there is no choice of parameters $(D,M_{0})$ such that a spot B solution emerges for $\nu\approx0$.
As seen in Figure \ref{fig:intro;paramspace}, one can numerically verify that there are three parameter curves $M_{1}(D)$, $M_{2}(D)$, and $M_{3}(D)$ such that: 
\begin{itemize}
    \item $0\leq M_{0} \leq M_{1}(D)$: no Hamiltonian-Hopf bifurcation occurs, and so no localised radial solutions emerge.
    \item $M_{1}(D)< M_{0}< M_{2}(D)$: $\nu<0$ and $c_{3}>0$, and so only depressed spot A solutions emerge.
    \item $M_{2}(D)< M_{0}< M_{3}(D)$: $\nu>0$ and $c_{3}>0$, and so only elevated spot A solutions emerge.
    \item $M_{3}(D)< M_{0}$: $\nu>0$ and $c_{3}<0$, and so elevated spot A, depressed spot B, elevated ring, and depressed ring solutions emerge.
\end{itemize}}
{We note that the fold curve $\nu\approx0$ is equivalent to $M_{0}\approx M_{2}(D)$, and so we predict there is a critical magnetisation at which spot A solutions have an alternate amplitude scaling $\varepsilon^{\frac{1}{4}}|\log(\varepsilon)|^{-\frac{1}{2}}$.}
\paragraph{Outline:}
The outline of the rest of the paper is as follows. 

\underline{Section \ref{s:setup}} (\textit{Formulation}) We formulate the radial ferrohydrostatic problem with two fluids of finite depth; the ferrofluid lies below a non-magnetisable fluid, where the ferrofluid obeys a linear magnetisation law. The key variables are axisymmetric magnetic potentials $\phi^{\pm}(r,y)$ of the upper and lower fluid respectively, as well as the height of the free surface $\eta(r)$, such that the two fluids form an interface at $y=\eta(r)$. Here, $r$ and $y$ are the radius and the height variables in cylindrical polar coordinates, respectively. The physical system is defined by the following equations: Maxwell's equations for both fluids with continuity conditions on the interface; and the steady Euler equations with a magnetic forcing term, solved on the interface in conjunction with a normal stress continuity condition. These equations can be found from variations of a `free energy', or Lagrangian, $\widehat{\mathscr{L}}$ of the form
\begin{align}
 \widehat{\mathscr{L}}(\phi^{-},\phi^{+}, \eta) =  \int_{0}^{\infty}\bigg[&\frac{\mu \mu_{0}}{2}\int_{-D}^{\eta(r)} \left|\nabla \phi^{-}\right|^2 \,\textnormal{d}y + \frac{\mu_{0}}{2}\int_{\eta(r)}^{D} \left|\nabla\phi^{+}\right|^2 \,\textnormal{d}y + \mu_{0}\mu h \left(\left.\phi^{-}\right|_{y=-D}-\left.\phi^{+}\right|_{y=D}\right)\nonumber\\
 & -\frac{\rho_{0}g \eta^2}{2} - \sigma\left(\sqrt{1+\left(\eta_{r}\right)^2}-1\right) - \frac{\mu_{0}}{2}\mu(\mu-1)h^2 \eta \bigg]r\;\textnormal{d}r,
\nonumber\end{align}
where variations are taken with respect to $\phi^{-}$, $\phi^{+}$, and $\eta$, and subject to the constraint $\phi^{-} = \phi^{+}$ at $y=\eta(r)$. Here, $\mu_{0}$ is the magnetic permeability of free space, $\nabla$ is the axisymmetric gradient operator, and the other parameters were defined previously for \eqref{intro:disp}.

\begin{figure}[ht!]
    \centering
    \includegraphics[height=4.2cm]{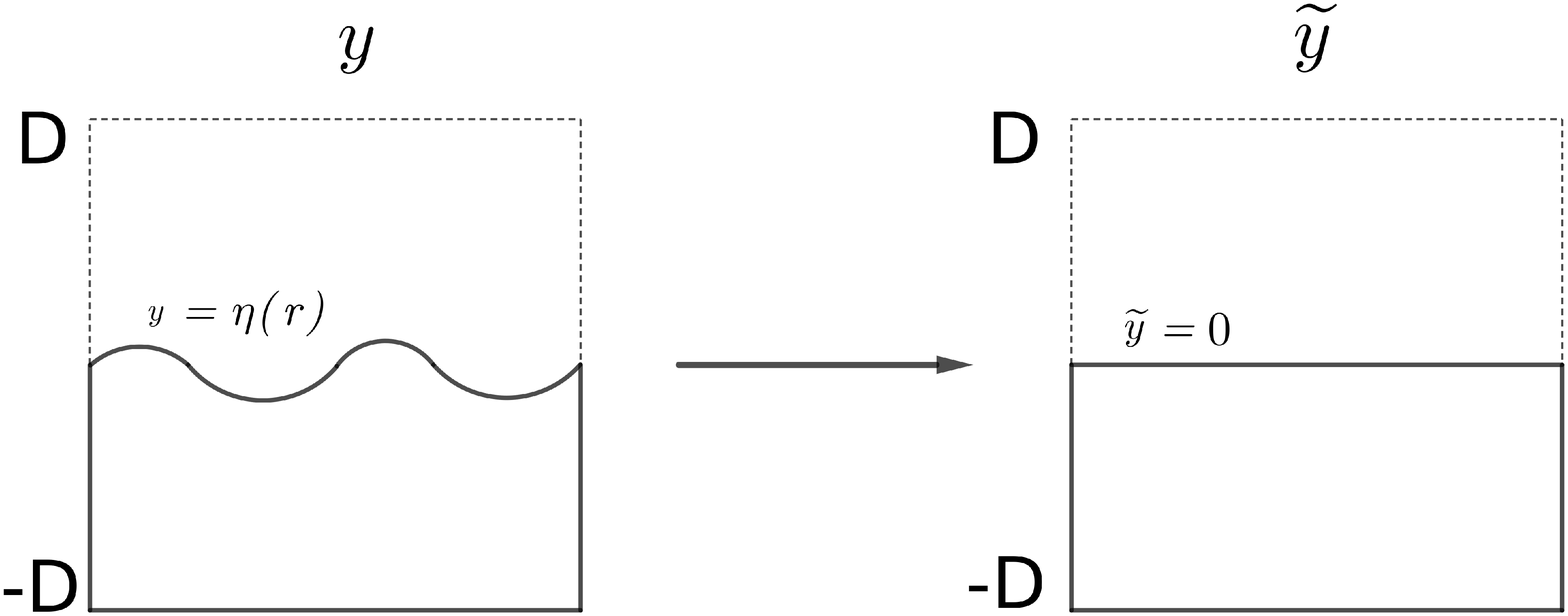}
    \caption{The transformation $y\to\widetilde{y}$ maps the domains $\mathscr{D}^{+}=\{(r,y)\,:\; \eta(r)\leq y \leq D\}$ and $\mathscr{D}^{-}=\{(r,y)\,:\; -D\leq y \leq \eta(r)\}$ to the infinite strips $\mathbb{R}\times[0,D]$ and $\mathbb{R}\times[-D,0]$, respectively. {T}he boundaries $y=\pm D$ are maintained as $\widetilde{y}=\pm D$, while the interface $y=\eta(r)$ is transformed into $\widetilde{y}=0$, which is now independent of $r$.}
    \label{fig:intro;flattening}
\end{figure}

There are various problems making analytical methods difficult to apply:
\begin{itemize}
    \item The domains of the system, i.e. the upper and lower fluids, are non-uniform in $r$ due to the interface at $y=\eta(r)$.
    \item The magnetic potentials have translational symmetry: if $(\phi^{\pm})$ are solutions, then $(\phi^{\pm}+c)$ are also solutions, for any arbitrary constant $c$. 
    \item The system has nonlinear boundary conditions, so standard invariant manifold theory cannot be applied in its current state.
\end{itemize}
To deal with these problems, multiple coordinate transformations are applied. We begin by non-dimensionalising the problem and applying a `flattening' transformation often seen in water wave literature. We rescale the height variable $y\to\widetilde{y}(r)$ with respect to the free surface $\eta(r)$, such that each domain defines a uniform rectangle in $\widetilde{y}$ coordinates, see Figure \ref{fig:intro;flattening}, at the expense of increasing the system's nonlinearity. We also isolate the constant applied magnetic field by defining $\phi^{+} = \psi^{+} + \mu h y$, $\phi^{-} = \psi^{-} + h y$ and perform a Legendre transformation to turn our second-order system for three variables into a first-order system for six variables. We remove the translational symmetry by defining the `average' of $\psi^{\pm}$ and removing it from the system; this imposes an extra constraint on the system and adds to the nonlinearity.

We deal with the nonlinear boundary conditions by employing a coordinate transformation seen in \cite[Section 4]{groves2017pattern}. This transformation leaves the linear problem invariant, and the new coordinates still satisfy the constraints of the system, now equipped with linear boundary conditions. Then, we define $\mathbf{u}=(\widetilde{\psi}^{-}, \widetilde{\psi}^{+}, \widetilde{\eta}, \widetilde{\alpha}^{-}, \widetilde{\alpha}^{+}, \widetilde{\gamma})^{\intercal}$, where $\widetilde{\psi}^{\pm}, \widetilde{\eta}$ are the transformed variables related to the respective original variables $\psi^{\pm}, \eta$ and the transformed variables $\widetilde{\alpha}^{\pm}, \widetilde{\gamma}$ are related to the Legendre conjugates of $\psi^{\pm}, \eta$, respectively. We can then write the full system in the form \eqref{intro:eqn}, with linear boundary conditions $\mathbf{B}_{1}\mathbf{u}=\mathbf{0}$.

\underline{Section \ref{s:spect}} (\textit{Decomposition}) We want to be able to construct local solutions in both the `core' and `far-field' regions and match them at a fixed $r=r_{0}$. Therefore, we would like to decompose solutions of \eqref{intro:eqn} into a basis of $y$-dependent vectors with $r$-dependent amplitudes. This means that, in both the `core' and `far-field' regions, we can match the amplitudes without having to worry about our basis changing across regions. Hence, in Section \ref{s:spect}, we construct a basis of eigenmodes for the matrix $\mathbf{L}_{\infty}:=\lim_{r\to\infty}\mathbf{L}(r)$, where $\mathbf{L}(r)$ is the non-autonomous linear operator for \eqref{intro:eqn}. The eigenvalues of $\mathbf{L}_{\infty}$ are found from the dispersion relation \eqref{intro:disp}, with corresponding eigenmodes that are found to form a complete $y$-dependent basis. Hence we can write,
\begin{align}
    \mathbf{u}= a(r)\;\mathbf{e}(y) + b(r)\;\mathbf{f}(y) + \overline{a}(r)\;\overline{\mathbf{e}}(y) + \overline{b}(r)\;\overline{\mathbf{f}}(y) + \sum_{n=1}^{\infty} \left\{ a_{n}(r)\mathbf{e}_{n}(y) + a_{-n}(r)\mathbf{e}_{-n}(y)\right\},\label{intro:u}
\end{align}
where $\overline{\,\cdot\,}$ indicates complex conjugation, $\{\mathbf{e}, \mathbf{f}\}$ correspond to the double mulitplicity eigenvalue $\textnormal{i}k$, and $\{\mathbf{e}_{\pm n}\}$ correspond to the respective real eigenvalues $\pm\lambda_{n}$, for $n\in\mathbb{N}$. Then, projecting onto each of these eigenmodes, we reduce our PDE system to infinitely-many nonlinear coupled ODEs, where the linear problem decouples into disjoint pairs,
\begin{align}
    &\frac{\textnormal{d}}{\textnormal{d}\,r} a = \textnormal{i}k\, a + b - \frac{1}{2r}\left(a-\overline{a}\right),
    &\quad
    &\frac{\textnormal{d}}{\textnormal{d}\,r} a_{n} = \lambda_{n}\, a_{n} - \frac{1}{2r}\left(a_{n}-a_{-n}\right),
    &\quad
    &\,
    &\nonumber\\
    &\frac{\textnormal{d}}{\textnormal{d}\,r} b = \textnormal{i}k\, b - \frac{1}{2r}\left(b+\overline{b}\right),
    &\quad
    &\frac{\textnormal{d}}{\textnormal{d}\,r} a_{-n} = -\lambda_{n}\, a_{-n} + \frac{1}{2r}\left(a_{n}-a_{-n}\right),
    &\quad
    & \forall n\in\mathbb{N}.
    &\label{intro:amp}
\end{align}

\underline{Section \ref{s:core}} (\textit{Core Parametrisation}) In this section, we discuss solutions near the core and construct a parametrisation for the core manifold, $\widetilde{\mathcal{W}}^{cu}_{-}{(\varepsilon)}$, which contains all small-amplitude solutions that are bounded as $r\to0$.

To do this, we first find explicit solutions to the linear problem \eqref{intro:amp}, noting that half the solutions are bounded as $r\to0$ and half are not. Taking a general solution to \eqref{intro:amp} with coefficients $d_{1}$, $d_{2}$, $d_{3}$, $d_{4}$, $c_{1,n}$, and $c_{2,n}$, we apply the method of variation of constants \cite{Walter1998ODEs} to find a general solution for the nonlinear problem; here, $d_{1}$, $d_{2}$, and $c_{1,n}$ are coefficients of the asymptotically stable solutions as $r\rightarrow\infty$ and $d_{3}$, $d_{4}$, and $c_{2,n}$ are coefficients of the unstable solutions as $r\rightarrow\infty$, for all $n\in\mathbb{N}$. Then, for $d_{3}=d_{4}=c_{2,n}=0$ and sufficiently small constants $d_{1}$, $d_{2}$, and $c_{1,n}$ for all $n\in\mathbb{N}$, we define a parametrisation that captures all small-amplitude solutions of \eqref{intro:eqn} with decomposition \eqref{intro:u} that are bounded as $r\to0$.

We explicitly find the coefficient $\nu$ in front of $d_{2}^{2}$ in the parametrisation of the amplitude $b(r_{0})$; we will require this to be {non-zero} for spot solutions to emerge from the flat state.

\underline{Section \ref{s:farfield}} (\textit{Far-field Decomposition}) In this section, we discuss solutions in the far-field, constructing a parametrisation for the far-field manifold {$\widetilde{\mathcal{W}}^{s}_{+}(\varepsilon)$}, which contains all small-amplitude solutions that decay exponentially as $r\to\infty$.

We first extend the system by defining $\sigma:=\frac{1}{r}$, where $0<\sigma\ll1$ in the far-field region. This makes the system \eqref{intro:amp} autonomous, and so standard invariant manifold theory can be applied; we note, for the linear problem, $\sigma=\sigma^{*}$ forms an invariant subspace, where $\sigma^{*}$ is {a small fixed value. For sufficiently small solutions, the centre and stable submanifolds $\mathcal{W}^{c}_{+}{(\varepsilon)|_{\sigma=\sigma^{*}}}, \mathcal{W}^{s}_{+}{(\varepsilon)|_{\sigma=\sigma^{*}}}\subset\mathcal{W}^{cs}_{+}{(\varepsilon)|_{\sigma=\sigma^{*}}}$ form a transverse intersection, i.e. they parametrise the entire manifold $\mathcal{W}^{cs}_{+}{(\varepsilon)|_{\sigma=\sigma^{*}}}$. Here, $\mathcal{W}^{c}_{+}{(\varepsilon)|_{\sigma=\sigma^{*}}}$ is an autonomous centre manifold, containing all small-amplitude solutions to the non-autonomous problem that are bounded for all large $r$, and $\change{\mathcal{W}^{s}_{+}(\varepsilon)|_{\sigma=\sigma^{*}}}$ is the stable manifold containing all small-amplitude solutions that exponentially decay to zero as $r\to\infty$. } {We evaluate the system at $\sigma=\frac{1}{r}$ in order to parametrise our solutions; this creates higher-order correction terms in our parametrisation, but leaves the leading-order expansion invariant..} We initially use the same variation of constants approach as in the core problem to construct a centre-stable manifold $\mathcal{W}^{cs}_{+}{(\varepsilon)|_{\sigma=\frac{1}{r}}}$ for a fixed $r\gg 1$, containing all small-amplitude solutions that {are bounded} as $r\to\infty$. 

In order to be able to match the far-field manifold ${\widetilde{\mathcal{W}}^{s}_{+}(\varepsilon):=\mathcal{W}^{s}_{+}(\varepsilon)|_{\sigma=\frac{1}{r}}}$ to the core manifold $\widetilde{\mathcal{W}}^{cu}_{-}{(\varepsilon)}$, we parametrise $\mathbf{q}\in{\widetilde{\mathcal{W}}^{s}_{+}(\varepsilon)}$ as a perturbation from a base point $\mathbf{p}\in{\mathcal{W}^{s}_{c}(\varepsilon)|_{\sigma=\frac{1}{r}}}$, where ${\mathcal{W}^{s}_{c}(\varepsilon)|_{\sigma=\frac{1}{r}}}\subset{\mathcal{W}^{c}_{+}(\varepsilon)|_{\sigma=\frac{1}{r}}}$ is the set of all small-amplitude solutions on the centre manifold that exhibit exponential decay as $r\to\infty$, such that $\mathbf{q}\to \mathbf{p}$ as $r\to\infty$. Then, to fully parametrise the far-field manifold, we need to identify initial conditions for a point $\mathbf{p}\in{\mathcal{W}^{c}_{+}(\varepsilon)|_{\sigma=\frac{1}{r}}}$ on the centre manifold such that $\mathbf{p}\in{\mathcal{W}^{s}_{c}(\varepsilon)|_{\sigma=\frac{1}{r}}}$, i.e. they decay exponentially to zero as $r\to\infty$. We restrict the far-field system onto the centre manifold, applying a normal-form transformation to make the next analytical step tractable.

\underline{Section \ref{s:geoblowup}} (\textit{Blow-up Coordinates}) We formulate the geometric blow-up coordinates introduced in \cite{mccalla2013spots}. These coordinates are required to parametrise decaying solutions on the centre manifold. However, we also need to control the far-field parametrisation as $r$ approaches $r_{0}$ at the bifurcation point $\varepsilon=0$, reconciling the algebraic behaviour in the core with the exponential decay in the far-field. 

Firstly, we introduce the `rescaling' chart from \cite[Section 2.3]{mccalla2013spots}: we define $s:= \varepsilon^{\frac{1}{2}} r$, $\sigma_{2}:=\sigma/ \varepsilon^{\frac{1}{2}}$, $A_{2}(s):=A(r)/ \varepsilon^{\frac{1}{2}}$, $\varepsilon_{2}:= \varepsilon^{\frac{1}{2}}$, and introduce $z_{2}(s) \approx \frac{\textnormal{d}}{\textnormal{d}s}A_{2}/A_{2}$. The rescaled system has equilibria at the points $ Q_{\pm}:(A_{2}, z_{2}, \sigma_{2}, \varepsilon_{2}) = (0, \pm\sqrt{c_{0}}, 0 ,  \varepsilon^{\frac{1}{2}})$ and, for $A_{2}$ restricted to $\mathbb{R}$, can be reduced to the non-autonomous real Ginzburg-Landau equation \eqref{Ginzburg-Landau,c3}. This has solutions that grow and decay exponentially with rate $\sqrt{c_{0}}$ as $s\to\infty$, respectively. Then, in order to find solutions that decay exponentially, we look for solutions that tend to $ Q_{-}$ as $s\to\infty$; these solutions are found to be valid for $s\in[\delta_{0},\infty)$, for some fixed $\delta_{0}>0$. However, as $\varepsilon\to0$, we lose control of the point $r=\delta_{0}/ \varepsilon^{\frac{1}{2}}$ and so a second transformation is required.

To maintain control of the solutions as $\varepsilon\to0$, we introduce a second transformation called the `transition' chart \cite[Section 2.4]{mccalla2013spots}. By scaling our variables by $\sigma$ and introducing `exponential time' $\rho={\log}\; r$, the system can be seen to have two equilibria $P_{\pm}$. Finding explicit solutions to \eqref{Ginzburg-Landau,c3}, we establish heteroclinic orbits $\mathbf{q}_{\pm}(s)$ that connect $ Q_{-}$ in the rescaling chart to $P_{\pm}$ in the transition chart, respectively. In the transition chart, the dynamics of the system is well-defined close to these equilibria, and so we will use $P_{+}$ and $P_{-}$ as a guide to track our solutions back to the point $\rho = {\log}\; r_{0}$. Then, substituting the parametrisation of the centre coordinates back into the foliation parametrisation, we will be able to match the far-field manifold with the core manifold at the point, as seen in the next section.

\underline{Section \ref{s:match}} (\textit{Matching}) In the final section, we follow the matching procedure for each localised radial pattern: spot A, spot B and rings. For spot A solutions, we track solutions from $ Q_{-}$ in the rescaling chart, staying close to the orbit $\mathbf{q}_{-}(s)$, to some fixed point ${s} = \delta_{0}>0$. Then, in the transition chart, solutions are tracked backwards in $\rho$, close to the equilibrium $P_{-}$ up to the point $\rho = {\log}\; r_{0}$. {Results are found in the cases when $\nu\neq0$ and $\nu\approx0$, separately.} For spot B solutions, we track solutions from $ Q_{-}$ in the rescaling chart, staying close to the orbit $\mathbf{q}_{+}(s)$, to ${s} = \delta_{0}$. Then, in the transition chart, solutions are tracked close to the equilibrium $P_{+}$ up to some point $\rho = \rho_{1}>{\log}\; r_{0}$. Following this, solutions transfer from the neighbourhood of $P_{+}$ to the neighbourhood of $P_{-}$, and then solutions are tracked close to $P_{-}$ up to the point $\rho={\log}\; r_{0}$. For ring solutions, we track solutions from $ Q_{-}$ in the rescaling chart, staying close to the orbit $\mathbf{q}_{+}(s)$, to ${s} = \delta_{0}$. Then, in the transition chart, solutions are tracked close to the equilibrium $P_{+}$ up to the point $\rho = {\log}\; r_{0}$. In each case, the full far-field parametrisation is then matched with the core manifold at the point $r=r_{0}$, to find the full profile of each localised radial pattern.

\paragraph{Novelty:} 
The novelty of this work can be identified as the following:
\begin{enumerate}
    \item We establish the formal existence of localised radial patterns for a free-boundary problem.
    \item We provide an approach to reduce the ferrofluid problem to an infinite system of ODEs.
    \item Using methods from ODE theory, we find {three} classes of localised radial patterns; {spot A, spot B and ring solutions. Each class is equipped with up- and down- variants, which indicate the sign of the maximum amplitude of the free surface}.
\end{enumerate}
Notably, this analytically verifies the existence of {up-spot A} solutions which have been observed experimentally \cite{richter2005two}, but also predicts the existence of {down-spot A}, {up-spot B}, {down-spot B}, up-ring, and down-ring solutions. These {latter} solutions have not been observed in the Rosensweig instability experiment, to the author's knowledge, but have been proven to exist in the Swift-Hohenberg equation \cite{lloyd2009localized,mccalla2013spots}.

{We wish to emphasise that the results presented here only serve to demonstrate the theoretical existence of localised radial patterns. In order to predict whether one could find these solutions experimentally, we would first need to investigate the stability of such solutions, which is not covered by this work. For example, in \cite{friedrichs2001pattern} one-dimensional stripes were found to emerge on the surface of a ferrofluid in the normal-field experiment. However, after investigating the energy of the system, these solutions were determined to be unstable everywhere. Hence, although down-spot A, up-spot B, down-spot B, up-ring, and down-ring solutions may exist, they may be unstable for all parameter values. {We note that, in the Swift-Hohenberg equation, spot A solutions are initially unstable at bifurcation but then restabilise at a fold bifurcation; this may also be the case in the ferrofluid experiment.}}

\paragraph{Future Work:}
This work serves as a framework for a rigorous proof of existence of each localised solution, which is in progress. Many parts of this work have been proven rigorously, however there are still some areas that prove to be more difficult. For example, the existence of infinitely many foliations and the regularity of our solutions are both non-trivial exercises to prove, in part due to the quasilinearity of the system, and the non-autonomous nature of the problem. There has recently been some progress made for quasilinear systems (see Chen et al. \cite{chen2019center}) that may prove fruitful in this problem.

As mentioned in the previous subsection, we find an appropriate basis decomposition that allows us to reduce the problem to an ODE system. However, this is an approach that is tailored for our specific problem, whereas we would prefer to find a generalisation of these methods for a general quasilinear PDE of this form. The recent work of Beck et al. \cite{beck2019exponential} for semilinear systems in radial domains may be useful for future progress in this area of study.

Some other areas of interest include: applying recent studies of radial invasion fronts (see Castillo-Pinto \cite{castillo2019extended}) to the ferrofluid problem, where there is an opportunity for experimental validation; looking at the question of homoclinic snaking (see \cite{lloyd2015homoclinic}) for localised radial patterns on a ferrofluid; and studying stability of our solutions, which is non-trivial due to the form of the Euler equations in our formulation of the problem.

\begin{Acknowledgment}
The authors would like to thank Mark Groves \& Reinhard Richter for illuminating discussions on this problem{, as well as the two anonymous referees for their thorough review and many helpful comments.} DJH acknowledges the EPSRC whose Doctoral Training Partnership Grant (EP/N509772/1) helped fund his PhD. DJBL acknowledges support from an EPSRC grant (Nucleation of Ferrosolitons and Ferropatterns) EP/H05040X/1. No new data was created during this study. 
\end{Acknowledgment}

\section{The Radial Ferrohydrostatic Problem}
\label{s:setup} 
 \begin{figure}[t]
    \centering
    \includegraphics[height=6cm]{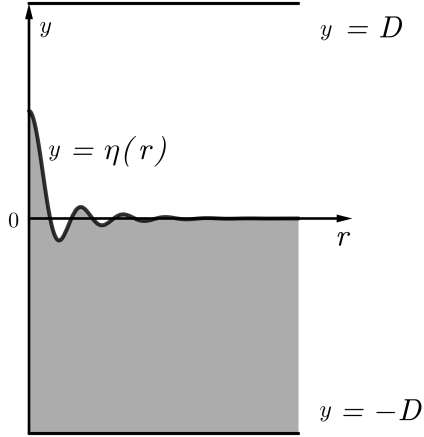}
    \caption{A diagram of the two-dimensional domains $\mathscr{D}^{\pm}$ (white and shaded, respectively). The three dimensional axisymmetric problem is in essence an adaptation of the 2 dimensional problem, where the coordinate frame is centred at the middle of the pattern.}
    \label{fig:radial}
\end{figure}
In order to formulate the radial ferrohydrostatic problem, we will assume that all of the relevant variables are axisymmetric. Therefore, for any function $f(r,\theta,y)$, $\frac{\textnormal{d}}{\textnormal{d}\theta} f(r,\theta, y) = 0$, where $(r,\theta, y)$ denote cylindrical polar coordinates. Then, the full 3-dimensional, finite-depth cylinder illustrated in Figure \ref{fig:ffspot} can be reduced to a 2-dimensional infinite strip in $(r, y)$ coordinates. Hence, we consider a domain consisting of two quiescent immiscible inviscid incompressible fluids
\begin{align}
\mathscr{D}^{+} &= \left\{\left( r , y\right) :\; r\in\left(0, \infty\right), \quad \eta(r)< y <D\right\},\nonumber\\
\mathscr{D}^{-} &= \left\{\left( r , y\right) :\; r\in\left(0, \infty\right), \quad -D < y < \eta(r)\right\},
\nonumber\end{align}
separated by a free surface $y=\eta\left(r\right)$, as shown in Figure \ref{fig:radial}. The upper fluid has density $\rho^{+}$ and unit relative permeability, while the lower fluid is a ferrofluid with density $\rho^{-}$. Denoting $\mathbf{H^{+}}$, $\mathbf{H^{-}}$ and $\mathbf{B^{+}}$, $\mathbf{B^{-}}$ as the respective magnetising and induction fields of the fluids, we assume that these fields satisfy the relation
\begin{align}
\mathbf{B^{+}} = \mu_{0}\mathbf{H^{+}}, \qquad \mathbf{B^{-}} = \mu_{0} \mu \mathbf{H^{-}}, \label{magn:lin}
\end{align}
where $\mu_{0}$ is the magnetic permeability of free space and the ferrofluid is assumed to have a linear magnetisation $\mathbf{M} = \left(\mu-1\right)\mathbf{H^{-}}$, where $\mu>1$. For a static fluid with no electric field, Maxwell's equations state that the magnetising and induction fields are irrotational and solenoidal, respectively. Defining magnetic potential functions $\phi^{-}, \phi^{+}$, we can write 
\begin{align}
    \left\{\begin{array}{rl}
           \nabla\cdot\mathbf{B}^{+}& = 0,\\
         \nabla\cdot\mathbf{B}^{-}& = 0, 
    \end{array}\right. \quad 
    \left\{\begin{array}{rl}
         \mathbf{H}^{+} & = - \nabla \phi^{+},  \\
         \mathbf{H}^{-} & = - \nabla \phi^{-}, 
    \end{array}\right. \qquad \textnormal{where}\;
    \begin{array}{rl}
         \nabla f&:= \mathbf{\widehat{r}}\frac{\partial f }{\partial r} + \mathbf{\widehat{y}}\frac{\partial f }{\partial y},\\
         \nabla\cdot\mathbf{f} &:= \frac{1}{r}\frac{\partial }{\partial r}\left(r\mathbf{f}\cdot\mathbf{\widehat{r}}\right) +  \frac{\partial }{\partial y} \left(\mathbf{f}\cdot\mathbf{\widehat{y}}\right)
    \end{array}\label{Max}
\end{align}
are the axisymmetric gradient and divergence operators for axisymmetric cylindrical polar coordinates $(r, y)$ with related unit vectors $(\mathbf{\widehat{r}}, \mathbf{\widehat{y}})$. Substituting in the linear relation \eqref{magn:lin}, \eqref{Max} reduces to Laplace's equation in both domains
\begin{align}
\Delta \phi^{\pm} = 0, \qquad \textnormal{where} \;\; \Delta f:= \frac{1}{r}\frac{\partial}{\partial r}\left(r\frac{\partial f}{\partial r}\right) + \frac{\partial^{2} f}{\partial y^{2}},
\label{Laplace:phi;xi}\end{align}
with the following magnetic continuity equations at $y=\eta(r)$
\begin{align}
\left(\mathbf{H}^{+}-\mathbf{H}^{-}\right){\times\mathbf{\widehat{n}} = \mathbf{0},} \qquad \left(\mathbf{B}^{+}-\mathbf{B}^{-}\right){\cdot\mathbf{\widehat{n}} = 0.}
\nonumber\end{align}
Here $\mathbf{\widehat{n}}=\frac{\nabla\left(y-\eta(r)\right)}{\left|\nabla\left(y-\eta(r)\right)\right|}$ is the unit normal vector to the interface, and it follows that
\begin{align}
\phi^{-}-\phi^{+} = 0, \qquad \mu\left(\phi^{-}_{r}\eta_{r} - \phi^{-}_{y}\right) - \left(\phi^{+}_{r}\eta_{r} - \phi^{+}_{y}\right) = 0,
\label{continuity:phi;xi}\end{align}
at $y=\eta(r)$, where subscripts denote partial derivatives. The ferrohydrostatic Euler equations for an incompressible fluid are given by 
\begin{align}
 \nabla\left(p^{+} + \rho^{+}g y\right)&=0, \qquad 
 \nabla\left(p^{-} + \rho^{-}g y - \frac{\mu_{0}}{2} (\mu-1)\left|\mathbf{H}^{-}\right|^2\right)=0,
\nonumber\end{align}
(Rosensweig \cite[Section 5.1]{rosensweig2013ferrohydrodynamics}), where $g$ is the acceleration due to gravity, and $p^{+}$, $p^{-}$ are the respective pressures of the upper and lower fluid. Integrating over their respective domains and equating on the interface, these equations are equivalent to 
\begin{align}
\frac{\mu_{0}}{2}(\mu-1)\left|\nabla \phi^{-}\right|^2-\left(p_{0} + \rho_{0}g\eta\right) = b_{0},
\label{ferrohydrostaticEuler:phi;eta}\end{align}
where $b_{0}:=b^{-}-b^{+}$, $p_{0}:=p^{-}-p^{+}$, $\rho_{0}:=\rho^{-}-\rho^{+}$, and $b^{\pm}$ are the respective Bernoulli constants for each domain $\mathscr{D}^{\pm}$. The ferrohydrostatic boundary condition at $y=\eta(r)$ is the continuity of the normal stress
\begin{align}
 p_{0} + \frac{\mu_{0}}{2}(\mu-1)^2\left(\mathbf{H}^{-}\cdot\mathbf{\widehat{n}}\right)^2 = 2\sigma \kappa, \label{p0}
\end{align}
(Rosensweig \cite[Section 5.2]{rosensweig2013ferrohydrodynamics}), where $\sigma>0$ is the coefficient of surface tension and 
\begin{align}
 2\kappa = \nabla\cdot\mathbf{\widehat{n}} = \frac{1}{r}\frac{\partial}{\partial r}\left(\frac{-r\,\eta_{r}}{\sqrt{1+\left(\eta_{r}\right)^2}}\right),
\nonumber\end{align}
is the mean curvature of the free surface. By substituting \eqref{p0} into \eqref{ferrohydrostaticEuler:phi;eta} to eliminate the relative pressure $p_{0}$, we find
\begin{align}
 \frac{\mu_{0}}{2}(\mu-1)\left|\nabla \phi^{-}\right|^2- \rho_{0}g\eta - b_{0} + \frac{\mu_{0}}{2}(\mu-1)^2\left(\mathbf{H}^{-}\cdot\mathbf{\widehat{n}}\right)^2 = 2\sigma \kappa.
\nonumber\end{align}
That is, $\phi^{-}$, $\phi^{+}$, and $\eta$ must satisfy
\begin{align}
 \frac{\mu_{0}\mu}{2}\left|\nabla \phi^{-}\right|^2 - \frac{\mu_{0}}{2}\left|\nabla\phi^{+}\right|^2 + \mu_{0}\mu\left(\phi^{-}_{r}\eta_{r} - \phi^{-}_{y}\right)\left(\phi^{-}_{y} - \phi^{+}_{y}\right) - \rho_{0}g\eta - 2\sigma \kappa - b_{0} =0,
\label{Ferrohydrostatic:phi;xi,eta}\end{align}
at $y=\eta(r)$. By choosing $b_{0} = \frac{\mu_{0}}{2}\mu(\mu-1)h^2$, the trivial state $(\phi^{-},\phi^{+},\eta) = (hy,\mu hy, 0)$ satisfies \eqref{Laplace:phi;xi}, \eqref{continuity:phi;xi} and \eqref{Ferrohydrostatic:phi;xi,eta}, and corresponds to a flat free surface and a uniform magnetic field with applied field strength $h$. In the standard Rosensweig instability experiment, $h$ has a critical value at which domain-covering peaks begin to form; we call this value $h_{c}$ and define the bifurcation parameter of the problem to be $\widehat{\varepsilon}:=h-h_{c}$ such that the Rosensweig bifurcation occurs at $\widehat{\varepsilon}=0$. Also, choosing inhomogeneous Neumann boundary conditions
\begin{align}
\left.\phi^{-}_{y}\right|_{y=-D} = h, \qquad \left.\phi^{+}_{y}\right|_{y=D} = \mu h,
\nonumber\end{align}
guarantees that the system of equations \eqref{Laplace:phi;xi}, \eqref{continuity:phi;xi}, and \eqref{Ferrohydrostatic:phi;xi,eta} has a variational structure: they are obtained from variations of the Lagrangian
\begin{align}
 \widehat{\mathscr{L}}(\phi^{-},\phi^{+}, \eta) =  \int_{0}^{\infty}\bigg[&\mu_{0}\frac{\mu}{2}\int_{-D}^{\eta(r)} \left|\nabla \phi^{-}\right|^2 \,\textnormal{d}y + \mu_{0}\frac{1}{2}\int_{\eta(r)}^{D} \left|\nabla \phi^{+}\right|^2 \,\textnormal{d}y + \mu_{0}\mu h \left(\left.\phi^{-}\right|_{y=-D}-\left.\phi^{+}\right|_{y=D}\right)\nonumber\\
 & -\frac{\rho_{0}g \eta^2}{2} - \sigma\left(\sqrt{1+\left(\eta_{r}\right)^2}-1\right) - \frac{\mu_{0}}{2}\mu(\mu-1)h^2 \eta \bigg]r\;\textnormal{d}r,
\nonumber\end{align}
where variations are taken with respect to $\phi^{-}$, $\phi^{+}$, and $\eta$, and subject to the constraint $\phi^{-} = \phi^{+}$ at $y=\eta(r)$. As we are investigating smooth solutions, we will include the restriction that, for an arbitrary function $f(r,y)$,
\begin{align}
    f_{r}(0,y) = 0 \;\implies \left(r\,f_{r}\right)_{r}(0,y)=0. \label{r0:cond}
\end{align}
This condition corresponds to radial symmetry at the origin, such that the $\frac{1}{r} f_{r}$ terms in Laplace's equation remain bounded as $r\to0$. We now introduce dimensionless variables
\begin{align}
 (\widehat{r}, \widehat{y}, \widehat{D}, \widehat{\eta}) :&= \frac{\mu_{0}h^2}{\sigma}(r,y, D, \eta), \qquad \widehat{\phi}^{\pm} := \frac{\mu_{0}h}{\sigma}\phi^{\pm}, \nonumber
\end{align}
and write $(\widehat{\phi}^{-},\; \widehat{\phi}^{+}, \;\widehat{\eta}) = (\widehat{\psi}^{-} + \mu\, \widehat{y}, \;\widehat{\psi}^{+} + \widehat{y}, \;\widehat{\eta})$, such that $(\widehat{\psi}^{-},\widehat{\psi}^{+},\widehat{\eta})=(0,\;0,\;0)$ is a solution to the full problem. Then, we find that
\begin{align}
\Delta \psi^{\pm} = 0, \qquad \left.\psi^{\pm}_{y}\right|_{y= \pm D} = 0,
\label{Laplace:psi;zeta}\end{align}
and
\begin{align}
\psi^{-} - \psi^{+} - (\mu-1)\eta =0,
\label{boundary:dyn}\end{align}
\begin{align}
\mu\left(\psi^{-}_{r}\eta_{r}-\psi^{-}_{y}\right) - \left(\psi^{+}_{r}\eta_{r} - \psi^{+}_{y}\right) = 0,
\label{boundary:kin}\end{align}
\begin{align}
    \frac{\mu}{2}\left|\nabla \psi^{-}\right|^2 - \frac{1}{2}\left|\nabla \psi^{+}\right|^2 + \mu\left(\psi^{-}_{r}\eta_{r} - \psi^{-}_{y}\right)\left(\psi^{-}_{y} - \psi^{+}_{y} - (\mu-1)\right) - \Upsilon\eta - 2\kappa =0,
\label{boundary:ferrohyd}\end{align}
at $y = \eta(r)$, where $\Upsilon:=\rho_{0}g\sigma \mu^{-2}h^{-4}$
and the hat notation has been dropped for notational simplicity. Note, equations \eqref{Laplace:psi;zeta}-\eqref{boundary:ferrohyd} follow from the variational principle $\delta \mathscr{L} =0$, where $\mathscr{L}$ is the non-dimensionalised `free energy' of the system:
\begin{align}
 \mathscr{L}(\psi^{-}, \psi^{+}, \eta) =  \int_{0}^{\infty}\bigg[&\frac{\mu}{2}\int_{-D}^{\eta(r)} \left|\nabla \psi^{-}\right|^2 \,\textnormal{d}y + \frac{1}{2}\int_{\eta(r)}^{D} \left|\nabla \psi^{+}\right|^2 \,\textnormal{d}y \nonumber\\
 & -\frac{\Upsilon \eta^2}{2} - \left(\sqrt{1+\left(\eta_{r}\right)^2}-1\right) + \mu \left(\left.\psi^{-}\right|_{y=\eta(r)}-\left.\psi^{+}\right|_{y=\eta(r)} - (\mu-1)\eta\right) \bigg]r\;\textnormal{d}r,
\nonumber\end{align}
and the variations are taken with respect to $\psi^{-}$, $\psi^{+}$, $\eta$, with $\psi^{-} - \psi^{+} - (\mu-1)\eta = 0$ at $y=\eta(r)$.
\paragraph{Hamiltonian formulation:}
One of the first difficulties encountered in this problem is the $r$-dependent nature of the domains $\mathscr{D}^{\pm}$, due to the boundary $y=\eta(r)$. In order to apply analytical techniques to solve \eqref{Laplace:psi;zeta}, we need the domains to be fixed for all $r$. This problem is often overcome in the water waves literature via a `flattening' transformation of the variable $y$, as seen for the ferrofluid problem in Twombly \& Thomas \cite{twombly1983bifurcating}, to
\begin{align}
 \widetilde{y}:=\left\{\begin{array}{cc}
    \frac{y-\eta(r)}{D-\eta(r)}D,  & y\geq \eta(r), \\
     \frac{y-\eta(r)}{D+\eta(r)}D,  & y < \eta(r), 
 \end{array}\right.
\nonumber\end{align}
which maps the respective domains $\mathscr{D}^{+}$, $\mathscr{D}^{-}$ to the infinite strips $\mathbb{R}\times\left(0,D\right)$ and $\mathbb{R}\times\left(-D,0\right)$, and the free surface $ y=\eta(r)$ to $\widetilde{y}=0$. The `flattened' variables $\widetilde{\psi}^{\pm}(r, \widetilde{y}) = \psi^{\pm}(r,y)$ satisfy the equations
\begin{align}
\psi^{\pm}_{rr} + \frac{1}{r}\psi^{\pm}_{r} - 2\eta_{r}K_{1}^{\pm}\psi^{\pm}_{ry} - \left[\eta_{rr}K_{1}^{\pm} \pm \frac{2}{D}\eta_{r}^2 K_{1}^{\pm}K_{2}^{\pm} + \frac{1}{r}\eta_{r}K_{1}^{\pm}\right]\psi^{\pm}_{y} + \left[\eta_{r}^2 \left(K_{1}^{\pm}\right)^2 + \left(K_{2}^{\pm}\right)^2\right]\psi^{\pm}_{yy} &= 0,
\label{Laplace:zeta;flat}
\end{align}
with boundary conditions
\begin{align}
\left.\psi^{\pm}_{y}\right|_{y=\pm D} = 0,
\label{bc:psi:zeta;flat}\end{align}
and
\begin{align}
&\psi^{-} - \psi^{+} - (\mu-1)\eta =0,
\label{boundary:dyn;flat}\\
&\mu\left(\psi^{-}_{r}\eta_{r}-\left(1+\eta_{r}^2\right)K_{2}^{-}\psi^{-}_{y}\right) - \left(\psi^{+}_{r}\eta_{r} - \left(1+\eta_{r}^2\right)K_{2}^{+}\psi^{+}_{y}\right) = 0,
\label{boundary:kin;flat}\\
    &\frac{\mu}{2}\left[\left(\psi^{-}_{r} - \eta_{r}K_{2}^{-}\psi^{-}_{y}\right)^2 + \left(K_{2}^{-}\psi^{-}_{y}\right)^2\right] - \frac{1}{2}\left[\left(\psi^{+}_{r} - \eta_{r}K_{2}^{+}\psi^{+}_{y}\right)^2 + \left(K_{2}^{+}\psi^{+}_{y}\right)^2\right] \quad \nonumber\\ & \qquad\qquad + \mu\left(\psi^{-}_{r}\eta_{r}-\left(1+\eta_{r}^2\right)K_{2}^{-}\psi^{-}_{y}\right)\left(K_{2}^{-}\psi^{-}_{y} - K_{2}^{+}\psi^{+}_{y} - (\mu-1)\right) - \Upsilon\eta - 2\kappa =0
\label{boundary:ferrohyd;flat}\end{align}
at $y=0$, where
\begin{align}
 K_{1}^{\pm} := \frac{D\mp y}{D\mp\eta}, \qquad K_{2}^{\pm} := \frac{D}{D\mp\eta},
\nonumber\end{align}
and the tildes have been dropped for notational simplicity.
Note that equations \eqref{Laplace:zeta;flat}-\eqref{boundary:ferrohyd;flat} can be obtained from the new variational principle $\delta \widetilde{\mathscr{L}}=0$, where 
\begin{align}
 \widetilde{\mathscr{L}}(\psi^{-}, \psi^{+}, \eta) =  \int_{0}^{\infty}\bigg[&\frac{\mu}{2}\int_{-D}^{0} \left[\left(\psi^{-}_{r} - \eta_{r}K_{1}^{-}\psi^{-}_{y}\right)^2 + \left(K_{2}^{-}\psi^{-}_{y}\right)^2\right]\left(K_{2}^{-}\right)^{-1} \,\textnormal{d}y\nonumber \\ & + \frac{1}{2}\int_{0}^{D} \left[\left(\psi^{+}_{r} - \eta_{r}K_{1}^{+}\psi^{+}_{y}\right)^2 + \left(K_{2}^{+}\psi^{+}_{y}\right)^2\right]\left(K_{2}^{+}\right)^{-1} \,\textnormal{d}y\label{L,flat} \\
 & -\frac{\Upsilon \eta^2}{2} - \left(\sqrt{1+\left(\eta_{r}\right)^2}-1\right) \bigg]r\;\textnormal{d}r,
\nonumber\end{align}
and the variations are taken with respect to $\psi^{-}$, $\psi^{+}$, $\eta$ satisfying the constraint \eqref{boundary:dyn;flat} at $y=0$. We consider $\widetilde{\mathscr{L}}$ as an action functional
\begin{align}
 \widetilde{\mathscr{L}}(\psi^{-}, \psi^{+}, \eta) = \int_{0}^{\infty} \mathcal{L}\left(\psi^{-}, \psi^{+}, \eta, \psi^{-}_{r}, \psi^{+}_{r}, \eta_{r}\right) \; \textnormal{d} r,
\nonumber\end{align}
where $\mathcal{L}$ is the integrand on the right-hand side of equation
\eqref{L,flat}. We carry out a Legendre transformation and introduce conjugate variables $\alpha^{-}$, $\alpha^{+}$, $\gamma$ defined by
\begin{align}
 \alpha^{-} = \frac{\delta \mathcal{L}}{\delta \psi^{-}_{r}} &= r\mu\left(\psi^{-}_{r} - \eta_{r}K_{1}^{-}\psi^{-}_{y}\right) \left(K_{2}^{-}\right)^{-1} ,\nonumber\\
 \alpha^{+} = \frac{\delta \mathcal{L}}{\delta \psi^{+}_{r}} &= r\left(\psi^{+}_{r} - \eta_{r}K_{1}^{+}\psi^{+}_{y}\right) \left(K_{2}^{+}\right)^{-1},\label{conj}\\
 \gamma = \frac{\delta \mathcal{L}}{\delta \eta_{r}} &= -r \mu\int_{-D}^{0}\left(\psi^{-}_{r} - \eta_{r}K_{1}^{-}\psi^{-}_{y}\right) \left(K_{2}^{-}\right)^{-1} K_{1}^{-} \psi^{-}_{y} \;\textnormal{d}y - r\int_{0}^{D}\left(\psi^{+}_{r} - \eta_{r}K_{1}^{+}\psi^{+}_{y}\right) \left(K_{2}^{+}\right)^{-1} K_{1}^{+} \psi^{+}_{y} \;\textnormal{d}y \nonumber\\ 
 & \;\qquad - \frac{r\eta_{r}}{\sqrt{1+\left(\eta_{r}\right)^2}},\nonumber
\end{align}
and thus the Hamiltonian function $\mathcal{H}\left(\psi^{-}, \psi^{+}, \eta, \alpha^{-}, \alpha^{+}, \gamma\right)$ is given by
\begin{align}
 \mathcal{H} &= \int_{-D}^{0} \alpha^{-} \psi^{-}_{r} \;\textnormal{d}y + \int_{0}^{D} \alpha^{+} \psi^{+}_{r} \;\textnormal{d}y + \gamma \eta_{r} - \mathcal{L}\left(\psi^{-}, \psi^{+}, \eta, \alpha^{-}, \alpha^{+}, \gamma\right),\label{H:defn}\\
 &=  r\frac{\mu}{2}\int_{-D}^{0} \left[\left(\frac{\alpha^{-}}{\mu r}\right)^2 - (\psi^{-}_{y})^2\right]K_{2}^{-} \,\textnormal{d}y + r\frac{1}{2}\int_{0}^{D} \left[\left(\frac{\alpha^{+}}{r}\right)^2 - (\psi^{+}_{y})^2\right]K_{2}^{+} \,\textnormal{d}y  +\frac{r\Upsilon \eta^2}{2} + r\left(\sqrt{1-W^2}-1\right),
\nonumber\end{align}
where
\begin{align}
 W = -\frac{\gamma}{r} - \int_{-D}^{0} \frac{\alpha^{-}}{r} K_{1}^{-}\psi^{-}_{y} \; \textnormal{d}y - \int_{0}^{D} \frac{\alpha^{+}}{r} K_{1}^{+}\psi^{+}_{y} \; \textnormal{d}y.
\nonumber\end{align}
Hamilton's equations are given explicitly by
\begin{align}
 \psi^{-}_{r} = \frac{\delta \mathcal{H}}{\delta \alpha^{-}} &= K_{2}^{-}\frac{\alpha^{-}}{\mu r} + \frac{W}{\sqrt{1-W^2}}K_{1}^{-}\psi^{-}_{y},
\label{Ham:psi}\\
 \psi^{+}_{r} = \frac{\delta \mathcal{H}}{\delta \alpha^{+}} &= K_{2}^{+}\frac{\alpha^{+}}{r} + \frac{W}{\sqrt{1-W^2}}K_{1}^{+}\psi^{+}_{y},
\label{Ham:zeta}\\
 \eta_{r} = \frac{\delta \mathcal{H}}{\delta \gamma} &= \frac{W}{\sqrt{1-W^2}},
\label{Ham:eta}\\
 \alpha^{-}_{r} = -\frac{\delta \mathcal{H}}{\delta \psi^{-}} &= -K_{2}^{-}\mu r \psi^{-}_{yy} + \frac{W}{\sqrt{1-W^2}}\left(K_{1}^{-}\alpha^{-}\right)_{y},
\label{Ham:alpha}\\
 \alpha^{+}_{r} = -\frac{\delta \mathcal{H}}{\delta \psi^{+}} &= -K_{2}^{+} r \psi^{+}_{yy} + \frac{W}{\sqrt{1-W^2}}\left(K_{1}^{+}\alpha^{+}\right)_{y},
\label{Ham:beta}\\
 \gamma_{r} = -\frac{\delta \mathcal{H}}{\delta \eta} &= -r\Upsilon \eta + \mu r (\mu-1) K_{2}^{-}\left[\psi^{-}_{y} - \frac{W}{\sqrt{1-W^2}} \frac{\alpha^{-}}{\mu r}\right]_{y=0}\nonumber\\ &\quad + \frac{1}{2D}\left( \mu r \left(K_{2}^{-}\right)^2 \int_{-D}^{0} \left[\left(\frac{\alpha^{-}}{\mu r}\right)^2 - (\psi^{-}_{y})^2\right]\;\textnormal{d} y -  r \left(K_{2}^{+}\right)^2 \int_{0}^{D} \left[\left(\frac{\alpha^{+}}{ r}\right)^2 - (\psi^{+}_{y})^2\right]\;\textnormal{d} y\right)\nonumber\\
 &\quad + \frac{1}{D}\frac{W}{\sqrt{1-W^2}}\left[ \mu r K_{2}^{-}\int_{-D}^{0} K_{1}^{-} \frac{\alpha^{-}}{\mu r} \psi^{-}_{y} \; \textnormal{d} y - r K_{2}^{+}\int_{0}^{D} K_{1}^{+} \frac{\alpha^{+}}{ r} \psi^{+}_{y} \; \textnormal{d} y\right],
\label{Ham:gamma}\end{align}
and are complemented by the constraint \eqref{boundary:dyn;flat} and boundary conditions
\begin{align}
    &\left.\psi^{\pm}_{y}\right|_{y=\pm D} = 0,\label{bd:Ham;upper}\\
    &\left[\mu K_{2}^{-}\left(\psi^{-}_{y} - \frac{W}{\sqrt{1-W^2}}\frac{\alpha^{-}}{\mu r}\right) - K_{2}^{+}\left(\psi^{+}_{y} - \frac{W}{\sqrt{1-W^2}}\frac{\alpha^{+}}{ r}\right)\right]_{y=0} = 0,\label{bd:Ham;kin}\\
    &\left[\mu K_{2}^{-}\left( \frac{\alpha^{-}}{\mu r} + \frac{W}{\sqrt{1-W^2}}\psi^{-}_{y}\right) - K_{2}^{+}\left( \frac{\alpha^{+}}{r} + \frac{W}{\sqrt{1-W^2}}\psi^{+}_{y}\right)\right]_{y=0} - (\mu-1)\frac{W}{\sqrt{1-W^2}} = 0.\label{bd:Ham;cont}
\end{align}
The boundary conditions \eqref{bd:Ham;upper} and \eqref{bd:Ham;kin} are found as a result of taking variations of \eqref{H:defn} with respect to $\psi^{-}$, $\psi^{+}$ and $\eta$, subject to \eqref{boundary:dyn;flat}. The condition \eqref{bd:Ham;cont} is equivalent to saying $\left[\psi^{-}_{r} - \psi^{+}_{r}\right]_{y=0} - (\mu-1)\eta_{r} = 0$ in the original, non-`flattened' coordinates. We note that our equations are invariant under the transformation $\psi^{-}\mapsto \psi^{-} + c$,\; $\psi^{+}\mapsto \psi^{+} + c$ for any constant $c$. To eliminate the translational symmetry of the problem, we replace $\left(\psi^{-},\psi^{+},\alpha^{-},\alpha^{+}\right)$ with $\left(\widehat{\psi}^{-}, \widehat{\psi}^{+}, \widehat{\alpha}^{-}, \widehat{\alpha}^{+},  A_{\psi}, A_{\alpha}\right)$, where
\begin{align}
    \widehat{\psi}^{\pm} := \psi^{\pm} -  A_{\psi}(r) , \qquad \widehat{\alpha}^{\pm} := \alpha^{\pm} - A_{\alpha}(r),\nonumber\end{align}
and
    \begin{align}
     A_{\psi} := \frac{1}{(\mu+1)D}\left[\mu\int_{-D}^{0} \psi^{-} \;\textnormal{d} y + \int_{0}^{D} \psi^{+} \;\textnormal{d} y\right], \qquad A_{\alpha} := \frac{1}{2D}\left[\int_{-D}^{0} \alpha^{-} \;\textnormal{d} y + \int_{0}^{D} \alpha^{+} \;\textnormal{d} y\right].\nonumber\end{align}
    Then, we find that
    \begin{align}\mu\int_{-D}^{0} \widehat{\psi}^{-} \;\textnormal{d} y + \int_{0}^{D} \widehat{\psi}^{+} \;\textnormal{d} y = 0, \qquad \int_{-D}^{0} \widehat{\alpha}^{-} \;\textnormal{d} y + \int_{0}^{D} \widehat{\alpha}^{+} \;\textnormal{d} y = 0.\nonumber
\end{align}
Note that $A_{\alpha}$ is a conserved quantity, and so we set $A_{\alpha}=0$, without loss of generality. Dropping the hats for notational simplicity, we observe that \eqref{Ham:psi} and \eqref{Ham:zeta} from Hamilton's equations become
\begin{align}
 \psi^{-}_{r} &= K_{2}^{-}\frac{\alpha^{-}}{\mu r} + \frac{W}{\sqrt{1-W^2}}K_{1}^{-}\psi^{-}_{y} - \frac{\textnormal{d}}{\textnormal{d} r}  A_{\psi},
\label{Ham:psi;average}\\
 \psi^{+}_{r} &= K_{2}^{+}\frac{\alpha^{+}}{r} + \frac{W}{\sqrt{1-W^2}}K_{1}^{+}\psi^{+}_{y}-\frac{\textnormal{d}}{\textnormal{d} r}  A_{\psi},
\label{Ham:zeta;average}
\end{align}
while equations \eqref{Ham:eta}-\eqref{Ham:gamma}, boundary conditions \eqref{bd:Ham;upper}-\eqref{bd:Ham;cont}, and constraint \eqref{boundary:dyn;flat} remain unchanged; the quantity $ A_{\psi}$ is found by quadrature from the equation
\begin{align}
    \frac{\textnormal{d}}{\textnormal{d} r}  A_{\psi} &= \frac{1}{(\mu+1)D}\left[\mu\int_{-D}^{0} \left(K_{2}^{-}\frac{\alpha^{-}}{\mu r} + \frac{W}{\sqrt{1-W^2}}K_{1}^{-}\psi^{-}_{y}\right) \;\textnormal{d} y + \int_{0}^{D} \left(K_{2}^{+}\frac{\alpha^{+}}{r} + \frac{W}{\sqrt{1-W^2}}K_{1}^{+}\psi^{+}_{y}\right)\;\textnormal{d} y\right].\nonumber
\end{align}
Note that equations \eqref{Ham:psi;average},\eqref{Ham:zeta;average} and \eqref{Ham:eta} appear to be singular as $r\to0$, while equations \eqref{Ham:alpha}-\eqref{Ham:gamma} appear unbounded as $r\to\infty$. In order to resolve the singular behaviour as $r\to\infty$, we introduce the new variables
\begin{align}
    \widetilde{\alpha}^{-} :&= \frac{\alpha^{-}}{\mu r}, \qquad \widetilde{\alpha}^{+} := \frac{\alpha^{+}}{r}, \qquad \widetilde{\gamma} := -\frac{\gamma}{r}.\label{conj:tilde}
\end{align}
This transformation again appears to be singular as $r\to0$, however we recall that the condition \eqref{r0:cond} implies that $\frac{1}{r} f_{r}$ remains bounded as $r\to0$ for any arbitrary function $f(r,y)$, and $\alpha^{\pm}$, $\gamma$ are Legendre conjugates associated with $\psi^{\pm}_{r}$, $\eta_{r}$, respectively. In particular, from \eqref{conj}, we see that 
\begin{align}
    \widetilde{\alpha}^{\pm} = \left(\psi^{\pm}_{r} - \frac{W}{\sqrt{1-W^{2}}} K_{1}^{\pm} \psi^{\pm}_{y}\right)\left(K_{2}^{\pm}\right)^{-1}, \qquad \widetilde{\gamma} = W + \mu\int_{-D}^{0} K_{1}^{-}\widetilde{\alpha}^{-}\psi^{-}_{y} \textnormal{d} y + \int_{0}^{D} K_{1}^{+}\widetilde{\alpha}^{+}\psi^{+}_{y} \textnormal{d} y,\nonumber
\end{align}
which remain bounded as $r\to0$. Applying the transformation \eqref{conj:tilde}, our equations become
\begin{align}
 \psi^{-}_{r} &= K_{2}^{-}\alpha^{-} + \frac{W}{\sqrt{1-W^2}}K_{1}^{-}\psi^{-}_{y}-\nu_{0},
\label{Ham:psi;auto}\\
 \psi^{+}_{r} &= K_{2}^{+}\alpha^{+} + \frac{W}{\sqrt{1-W^2}}K_{1}^{+}\psi^{+}_{y}-\nu_{0},
\label{Ham:zeta;auto}\\
 \eta_{r} &= \frac{W}{\sqrt{1-W^2}},
\label{Ham:eta;auto}\\
 \alpha^{-}_{r}  &= -\frac{1}{r}\alpha^{-} -K_{2}^{-} \psi^{-}_{yy} + \frac{W}{\sqrt{1-W^2}}\left(K_{1}^{-}\alpha^{-}\right)_{y},
\label{Ham:alpha;auto}\\
 \alpha^{+}_{r}  &= -\frac{1}{r}\alpha^{+} -K_{2}^{+}  \psi^{+}_{yy} + \frac{W}{\sqrt{1-W^2}}\left(K_{1}^{+}\alpha^{+}\right)_{y},
\label{Ham:beta;auto}\\
 \gamma_{r} &= -\frac{1}{r}\gamma + \Upsilon \eta - \mu  (\mu-1) K_{2}^{-}\left[\psi^{-}_{y} - \frac{W}{\sqrt{1-W^2}} \alpha^{-}\right]_{y=0}\nonumber\\ &\quad - \frac{1}{2D}\left( \mu \left(K_{2}^{-}\right)^2 \int_{-D}^{0} \left[(\alpha^{-})^2 - (\psi^{-}_{y})^2\right]\;\textnormal{d} y -   \left(K_{2}^{+}\right)^2 \int_{0}^{D} \left[(\alpha^{+})^2 - (\psi^{+}_{y})^2\right]\;\textnormal{d} y\right)\nonumber\\
 &\quad - \frac{1}{D}\frac{W}{\sqrt{1-W^2}}\left[ \mu  K_{2}^{-}\int_{-D}^{0} K_{1}^{-} \alpha^{-} \psi^{-}_{y} \; \textnormal{d} y -  K_{2}^{+}\int_{0}^{D} K_{1}^{+} \alpha^{+} \psi^{+}_{y} \; \textnormal{d} y\right],
\label{Ham:gamma;auto}\end{align}
where we have removed the tildes for notational simplicity, with constraints \eqref{boundary:dyn;flat} and
\begin{align}&\mu\int_{-D}^{0} \psi^{-} \;\textnormal{d} y + \int_{0}^{D} \psi^{+} \;\textnormal{d} y = 0, \qquad\qquad \mu\int_{-D}^{0} \alpha^{-} \;\textnormal{d} y + \int_{0}^{D} \alpha^{+} \;\textnormal{d} y = 0,\label{Ham:constraint;average}\\
\intertext{and boundary conditions \eqref{bd:Ham;upper} and}
&\left[\mu K_{2}^{-}\left(\psi^{-}_{y} - \frac{W}{\sqrt{1-W^2}}\alpha^{-}\right) - K_{2}^{+}\left(\psi^{+}_{y} - \frac{W}{\sqrt{1-W^2}}\alpha^{+}\right)\right]_{y=0} = 0,\label{Ham:bndry;kin}\\
&\left[\mu K_{2}^{-}\left( \alpha^{-} + \frac{W}{\sqrt{1-W^2}}\psi^{-}_{y}\right) - K_{2}^{+}\left( \alpha^{+} + \frac{W}{\sqrt{1-W^2}}\psi^{+}_{y}\right)\right]_{y=0} - (\mu-1)\frac{W}{\sqrt{1-W^2}} = 0.\label{Ham:bndry;cont}
\end{align}
Here, we have defined $\nu_{0}:= \frac{\textnormal{d}}{\textnormal{d}r}\,A_{\psi}$, so
\begin{align}
\nu_{0} &= \frac{1}{(\mu+1)D}\left[\mu\int_{-D}^{0} \left(K_{2}^{-}\alpha^{-} + \frac{W}{\sqrt{1-W^2}}K_{1}^{-}\psi^{-}_{y}\right) \;\textnormal{d} y + \int_{0}^{D} \left(K_{2}^{+}\alpha^{+} + \frac{W}{\sqrt{1-W^2}}K_{1}^{+}\psi^{+}_{y}\right)\;\textnormal{d} y\right],\nonumber\\
 W &= \gamma - \mu\int_{-D}^{0} \alpha^{-} K_{1}^{-}\psi^{-}_{y} \; \textnormal{d}y - \int_{0}^{D} \alpha^{+} K_{1}^{+}\psi^{+}_{y} \; \textnormal{d}y.
\nonumber\end{align}
In the form of $\Upsilon$ defined following \eqref{boundary:ferrohyd}, we consider $h=h_{c} + \widehat{\varepsilon}$ as a perturbation from the critical field strength for the Rosensweig instability. Thus, $\Upsilon = \Upsilon_{0} - \varepsilon$, where
\begin{align}
\Upsilon_{0}:= \frac{\rho_{0}g\sigma}{\mu^2 h_{c}^{4}}, \qquad \qquad \varepsilon =  \frac{4\Upsilon_{0}}{h_{c}} \widehat{\varepsilon} + \textnormal{O}(\widehat{\varepsilon}^{2}).\nonumber
\end{align}
{However, in the one-dimensional problem \cite{groves2017pattern}, localised patterns were found to bifurcate for super-critical values of $\Upsilon$, i.e. when $\varepsilon<0$. Hence, for convenience we define $\widetilde{\varepsilon}:=-\varepsilon$, such that localised patterns emerge in the region $0<\widetilde{\varepsilon}\ll1$, and then we remove the tilde for notational simplicity such that $\Upsilon=\Upsilon_{0}+\varepsilon$.} The ferrohydrostatic problem can then be written as
\begin{align}
    \mathbf{u}_{r} = \mathbf{g}\left(\mathbf{u}, \varepsilon, r\right) = \mathbf{L}(r)\mathbf{u} + \mathscr{F}(\mathbf{u},\varepsilon,r),\nonumber
\end{align}
where $\mathbf{u}=\left( \psi^{-}, \psi^{+}, \eta, \alpha^{-}, \alpha^{+}, \gamma\right)^{\intercal}$ and $\mathbf{g}(\mathbf{u},\varepsilon,r)$ is a nonlinear function containing (up to and including) second-order derivatives in $y$. The function $\mathbf{g}$ is given by the right-hand side of \eqref{Ham:psi;auto}-\eqref{Ham:gamma;auto}, and can be written as a non-autonomous linear part, $\mathbf{L}(r)\mathbf{u}$, and a collection of nonlinear terms $\mathscr{F}(\mathbf{u},\varepsilon,r)$. In addition, we have the linear constraints \eqref{boundary:dyn;flat} and \eqref{Ham:constraint;average} and nonlinear boundary conditions which we write as $\mathbf{B}(\mathbf{u})=\mathbf{0}$, where $\mathbf{B}$ is given by the left-hand sides of \eqref{bd:Ham;upper}, \eqref{Ham:bndry;kin}, and \eqref{Ham:bndry;cont}. In order to apply centre-manifold theory to our problem, and to keep consistent with the one dimensional problem in \cite{groves2017pattern}, we require the problem to have linear boundary conditions. We can overcome this problem via the transformation $\widetilde{\mathbf{u}} = \mathcal{G}\mathbf{u}$, where
\begin{align}
    \widetilde{\psi}^{\pm} &= \psi^{\pm} + \int_{0}^{y}\left\{K_{2}^{\pm}\psi^{\pm}_{y} - K_{1}^{\pm}\frac{W}{\sqrt{1-W^2}}\alpha^{\pm} - \psi^{\pm}_{y}\right\} \;\textnormal{d}t - \frac{\tau}{(\mu+1)D}\nonumber,\\
    \widetilde{\alpha}^{-} &=K_{2}^{-}\alpha^{-} + K_{1}^{-}\frac{W}{\sqrt{1-W^2}}\psi^{-}_{y} - \frac{3y^2 }{(\mu+1)D^3}\nu_{1}\nonumber,\\
    \widetilde{\alpha}^{+} &=K_{2}^{+}\alpha^{+} + K_{1}^{+}\frac{W}{\sqrt{1-W^2}}\psi^{+}_{y} - (\mu-1)\left[\gamma - \frac{W}{\sqrt{1-W^2}}\right] - \frac{3y^2 }{(\mu+1)D^3}\nu_{1},\nonumber
\end{align}
with $\eta$ and $\gamma$ left unchanged, and $\tau$ and $\nu_{1}$ are defined so that the constraint \eqref{Ham:constraint;average} is satisfied (see \cite{groves2017pattern} for an almost identical transformation). Then, $\widetilde{\mathbf{u}} := \mathcal{G}\mathbf{u}$ satisfies the equations
\begin{align}
     {\mathbf{u}}_{r} =  {\mathbf{g}}\left( {\mathbf{u}}, \varepsilon, r\right) &= \textnormal{d}\mathcal{G}\left[\mathcal{G}^{-1}( {\mathbf{u}})\right]\mathbf{g}\left(\mathcal{G}^{-1}( {\mathbf{u}}), \varepsilon, r\right),\nonumber\\
     &= \mathbf{L}(r)\mathbf{u} + \mathcal{F}(\mathbf{u},\varepsilon,r),\label{full:syst;linbd}
\end{align}
with boundary conditions $\mathbf{B}_{1}\mathbf{u}=\mathbf{0}$, where $\mathbf{B}_{1}\mathbf{u}$ is the linearisation of $\mathbf{B}(\mathbf{u})$ and we have dropped the tildes for notational simplicity. The transformation $\mathcal{G}$ is a near-identity operator, leaving the linear part of $\mathbf{g}$ invariant; therefore the related linearised system is 
\begin{align}
    \mathbf{u}_{r} = \mathbf{L}(r)\mathbf{u},\qquad \qquad \mathbf{B}_{1}\mathbf{u}=\mathbf{0},\nonumber
\end{align}
where 
\begin{align}
    \mathbf{L}(r)\mathbf{u} &= \begin{pmatrix} \alpha^{-} \\ \alpha^{+} \\ \gamma \\ -\frac{1}{r}\alpha^{-} - \psi^{-}_{yy} \\ -\frac{1}{r}\alpha^{+} - \psi^{+}_{yy} \\ -\frac{1}{r}\gamma + \Upsilon_{0} \eta -(\mu-1)\mu\left.\psi^{-}_{y}\right|_{y=0}\end{pmatrix},\qquad \mathbf{B}_{1}\mathbf{u} = \begin{pmatrix} \left.\psi^{-}_{y}\right|_{y=-D}\\\left.\psi^{+}_{y}\right|_{y=D}\\
    \left[\mu \psi^{-}_{y} - \psi^{+}_{y}\right]_{y=0}\\
    \left[ \alpha^{-} - \alpha^{+} - (\mu-1)\gamma \right]_{y=0}\\\end{pmatrix},\label{L:lin}
\end{align}
with linear constraints \eqref{boundary:dyn;flat} and \eqref{Ham:constraint;average}. We define the autonomous linear operator $\mathbf{L}_{\infty}:=\lim_{r\to\infty}\mathbf{L}(r)$, which is equivalent to the linear operator of the one-dimensional ferrohydrostatic problem found in \cite{groves2017pattern}.

\section{`Spectral' Decomposition}\label{s:spect}
\begin{figure}[t!]
    \centering
\includegraphics[width=\linewidth]{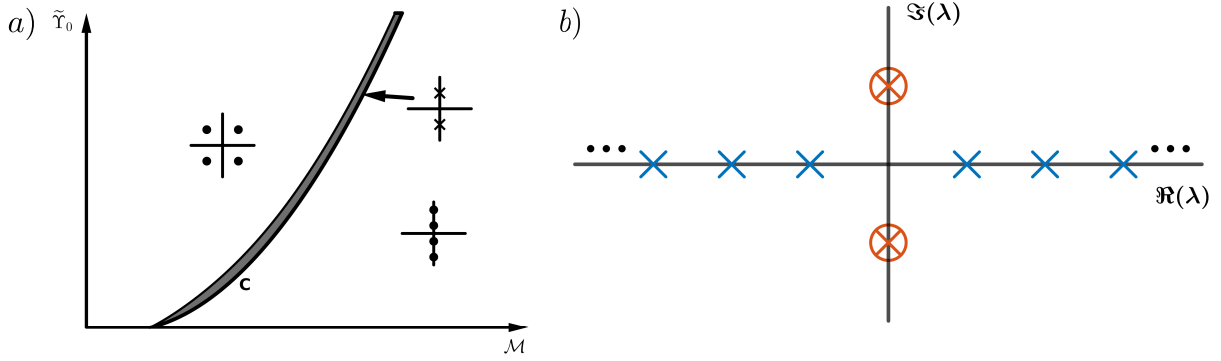}
    \caption{{a) Purely imaginary eigenvalues plotted as functions of the parameters $\mathcal{M}$ and $\widetilde{\Upsilon}_{0}$. A Hamiltonian Hopf bifurcation occurs at points on the curve $C$, and localised patterns emerge in the shaded region. b) The spectrum of $\mathbf{L}_{\infty}$ at the bifurcation point $\varepsilon=0$. There are two eigenvalues with double algebraic multiplicity at $\lambda = \pm \textnormal{i}k$, and infinitely-many real eigenvalues $\lambda=\pm\lambda_{n}$, with $\frac{\pi}{2D}<\lambda_{n}<\lambda_{n+1}$ for all $n\in\mathbb{N}$.}}
    \label{fig:hamhopf}
\end{figure}

We wish to decompose $\mathbf{u}$ onto an $r$-independent basis. We recall that $\mathbf{L}(r)\to\mathbf{L}_{\infty}$ as $r\to\infty$, and we will show that the eigenmodes of $\mathbf{L}_{\infty}$, which are independent of $r$, are a good choice of basis. Here $\mathbf{L}_{\infty}$ is a linear operator acting on a subset of the space $\mathcal{X}$, where
\begin{align}
    \mathcal{X} :&= \left\{\mathbf{u}\in\textnormal{H}^{1}(-D,0)\times\textnormal{H}^{1}(0,D)\times\mathbb{R}\times\textnormal{L}^{2}(-D,0)\times\textnormal{L}^{2}(0,D)\times\mathbb{R} \;\; \textnormal{s.t. $\mathbf{u}$ satisfies \eqref{boundary:dyn;flat} and \eqref{Ham:constraint;average}}\right\}.\nonumber
\end{align}
Through explicit calculations, one can show that a non-zero complex number $\lambda$ is an eigenvalue of the linear problem $\mathbf{u}_{r}=\mathbf{L}_{\infty}\mathbf{u}$ if and only if $\lambda = \frac{\varsigma}{D}$, and $\varsigma$ satisfies the dispersion relation
\begin{align}
    \Delta_{0} (\varsigma) :&= \mathcal{M}\;\varsigma \sin^2(\varsigma) - \left(\varsigma^2 - \widetilde{\Upsilon}_{0}\right)\sin(\varsigma)\cos(\varsigma) = 0,\label{disp}\\
    \intertext{where}
    \mathcal{M}&= \frac{\mu(\mu-1)^2 D}{(\mu +1)}, \qquad\qquad \widetilde{\Upsilon}_{0} = \Upsilon_{0}D^2; \nonumber
\end{align}
furthermore, even though $\varsigma=0$ satisfies \eqref{disp}, $\lambda=0$ is not an eigenvalue of $\mathbf{L}_{\infty}$ when restricted to $\mathcal{X}$. A purely imaginary number $\lambda=\textnormal{i} k$ is an eigenvalue of \eqref{disp} if and only if $k = \frac{ k_{D}}{D}$, where $ k_{D}$ satisfies
\begin{align}
    \mathcal{M}\;  k_{D} \tanh( k_{D}) - \left( k_{D}^2 + \widetilde{\Upsilon}_{0}\right)= 0.\nonumber
\end{align}
This equation has either zero, one or two pairs $\pm  k_{D}$ of solutions, as seen in Figure \ref{fig:hamhopf}{a)}. A Hamiltonian Hopf bifurcation takes place for values $(\widetilde{\Upsilon}_{0},\mathcal{M}) = (\widetilde{\Upsilon}_{H},\mathcal{M}_{H})( k_{D})$, where
\begin{align}
    \widetilde{\Upsilon}_{H}( k_{D}) &=  k_{D}^2 \left(\frac{\tanh( k_{D}) -  k_{D} \,\textnormal{sech}^2( k_{D})}{\tanh( k_{D}) + k_{D}  \,\textnormal{sech}^2( k_{D})}\right), \qquad\qquad \mathcal{M}_{H}( k_{D}) =\frac{2 k_{D}}{\tanh( k_{D}) +  k_{D} \,\textnormal{sech}^2( k_{D})} \label{M:Upsilon}
\end{align}
for any $ k_{D}\in(0,\infty)$; at this point, two pairs of simple, purely imaginary eigenvalues become complex by colliding on the imaginary axis at $\pm\textnormal{i} k$, and forming two Jordan chains of length 2. We therefore set $(\widetilde{\Upsilon}_{0},\mathcal{M}) = (\widetilde{\Upsilon}_{H},\mathcal{M}_{H})( k_{D})$ for some $ k_{D}\in(0,\infty)$. 

The spectrum of $\mathbf{L}_{\infty}$ now consists of two eigenvalues $\pm\textnormal{i}k$ with double algebraic multiplicity and a countably infinite family $\{\pm\lambda_{n}\}_{n\in\mathbb{N}}$, see Figure \ref{fig:hamhopf}{b)}, where $\{\lambda_{n}\}_{n\in\mathbb{N}}$ is the set of all positive real solutions of $\Delta_{0}(\lambda_{n}D)=0$. Furthermore,  $\{\lambda_{n}\}_{n\in\mathbb{N}}$ is an increasing sequence, and $\lambda_{1}>\frac{\pi}{2 D}$ (see Appendix \ref{app:spect}). The limiting linear system $\mathbf{u}_{r} = \mathbf{L}_{\infty}\mathbf{u}\;$ has a Hamiltonian structure, equipped with a symplectic two-form $\Omega$, defined as
\begin{align}
    \Omega\left(\mathbf{u}_{1},\mathbf{u}_{2}\right) = \mu\int_{-D}^{0} \left[\psi^{-}_{1}\alpha^{-}_{2} - \alpha^{-}_{1}\psi^{-}_{2}\right]\;\textnormal{d}y + \int_{0}^{D} \left[\psi^{+}_{1}\alpha^{+}_{2} - \alpha^{+}_{1}\psi^{+}_{2}\right]\;\textnormal{d}y - \left[\eta_{1}\gamma_{2} - \gamma_{1}\eta_{2}\right], \nonumber
\end{align}
where $\mathbf{u}_{i}=\left( \psi^{-}_{i}, \psi^{+}_{i}, \eta_{i}, \alpha^{-}_{i}, \alpha^{+}_{i}, \gamma_{i}\right)$. We choose basis vectors $\mathbf{e},\mathbf{f}, \mathbf{e}_{\pm n}$ such that
\begin{align}
    \left[\mathbf{L}_{\infty} - \textnormal{i}k\mathbb{1}\right] \mathbf{e}=\mathbf{0}, \qquad 
    \left[\mathbf{L}_{\infty} - \textnormal{i}k\mathbb{1}\right] \mathbf{f}=\mathbf{e}, \qquad  \left[\mathbf{L}_{\infty} \mp \lambda_{n}\mathbb{1}\right] \mathbf{e}_{\pm n}=\mathbf{0},\label{defn:eigenmode}
\end{align}
for $n\in\mathbb{N}$, where 
\begin{align}
    \Omega(\mathbf{f},\overline{\mathbf{e}}) =\Omega(\mathbf{e}_{-n}, \mathbf{e}_{n})= 1, \qquad  \Omega(\mathbf{e},\overline{\mathbf{f}}) =\Omega(\mathbf{e}_{n},\mathbf{e}_{-n})= -1,\nonumber
\end{align}
and the `symplectic products' of all other distinct combinations are zero. Here, $\mathbb{1}$ denotes the identity matrix and overbars denote the complex conjugate. The set $\{\mathbf{e},\overline{\mathbf{e}},\mathbf{f},\overline{\mathbf{f}}\}\cup\{\mathbf{e}_{n}\}_{n\in\mathbb{Z}\backslash\{0\}}$ can be proven to be a Riesz basis in $\mathcal{X}$ \cite{Groves2019Complete}, such that the set is a complete basis that converges independently of its ordering (called \textit{unconditional convergence}). This means that we can decompose our solution $\mathbf{u}$ in terms of a `spectral decomposition'
\begin{align}
    \mathbf{u}= \widetilde{\mathbf{u}} := a(r)\;\mathbf{e} + b(r)\;\mathbf{f} + \overline{a}(r)\;\overline{\mathbf{e}} + \overline{b}(r)\;\overline{\mathbf{f}} + \sum_{n=1}^{\infty} \left\{ a_{n}(r)\mathbf{e}_{n} + a_{-n}(r)\mathbf{e}_{-n}\right\},\label{u:defn}
\end{align}
where $a,b\in\mathbb{C}$, $a_{\pm n}\in\mathbb{R}$ for all $n\in\mathbb{N}$, and the dynamics of $\mathbf{u}$ is tied directly to the dynamics of the amplitudes $a,b$, and $a_{\pm n}$. We introduce the projections
\begin{align}
    \mathcal{P}_{0}\mathbf{v} := \Omega\left(\mathbf{v}, \overline{\mathbf{e}}\right)\mathbf{f} + \Omega\left(\mathbf{v}, \mathbf{e}\right)\overline{\mathbf{f}} -\Omega\left(\mathbf{v}, \overline{\mathbf{f}}\right)\mathbf{e} - \Omega\left(\mathbf{v}, \mathbf{f}\right)\overline{\mathbf{e}} , \qquad \mathcal{P}_{n}\mathbf{v} := \Omega\left(\mathbf{v}, \mathbf{e}_{n}\right)\mathbf{e}_{-n} -\Omega\left(\mathbf{v}, \mathbf{e}_{-n}\right)\mathbf{e}_{n},  \nonumber
\end{align}
with respective complements $ Q_{j}:= \mathbb{1} - \mathcal{P}_{j}$. By applying $\mathcal{P}_{0}$ to the full system \eqref{full:syst;linbd} and using the fact that $\{\mathbf{e},\overline{\mathbf{e}},\mathbf{f},\overline{\mathbf{f}}\}\cup\{\mathbf{e}_{n}\}_{n\in\mathbb{Z}\backslash\{0\}}$ is a Riesz basis, we find that $a(r), b(r)$ satisfy the complex amplitude equations
\begin{align}
    \frac{\textnormal{d}}{\textnormal{d}\,r} a &= \textnormal{i}k\, a + b - \frac{1}{2r}\left(a-\overline{a}\right) - \Omega\left(\mathcal{F}, \overline{\mathbf{f}}\right),\label{amp:a}\\
    \frac{\textnormal{d}}{\textnormal{d}\,r} b &= \textnormal{i}k\, b - \frac{1}{2r}\left(b+\overline{b}\right) + \Omega\left(\mathcal{F}, \overline{\mathbf{e}}\right),\label{amp:b}
\end{align}
where $\mathcal{F}(\varepsilon, r, a,b,a_{1},a_{-1},a_{2},\dots)= \mathcal{F}(\mathbf{u},\varepsilon, r)$ for $\mathbf{u}$ defined in \eqref{u:defn}. Similarly, by applying $\mathcal{P}_{n}$ for each $n\in\mathbb{N}$ respectively, we find that 
\begin{align}
    \frac{\textnormal{d}}{\textnormal{d}\,r} a_{n} &= \lambda_{n}\, a_{n} - \frac{1}{2r}\left(a_{n}-a_{-n}\right) - \Omega\left(\mathcal{F}, \mathbf{e}_{-n}\right),\label{amp:a;n}\\
    \frac{\textnormal{d}}{\textnormal{d}\,r} a_{-n} &= -\lambda_{n}\, a_{-n} + \frac{1}{2r}\left(a_{n}-a_{-n}\right) + \Omega\left(\mathcal{F}, \mathbf{e}_{n}\right).\label{amp:a;-n}
\end{align}
\section{Core Solutions}
\label{s:core}
We first wish to construct the `core' manifold, i.e. the set of all small-amplitude solutions to \eqref{full:syst;linbd} that remain bounded as $r\to0$. We begin by investigating the linear problem of \eqref{amp:a} and \eqref{amp:b}, i.e. with $\mathcal{F}\equiv\mathbf{0}$, thus
\begin{align}
    \frac{\textnormal{d}}{\textnormal{d}\,r} a &= \textnormal{i}k\, a + b - \frac{1}{2r}\left(a-\overline{a}\right),\qquad 
    \frac{\textnormal{d}}{\textnormal{d}\,r} b = \textnormal{i}k\, b - \frac{1}{2r}\left(b+\overline{b}\right).\nonumber
\end{align}
By performing a change of coordinates into the real and imaginary parts of both $a$ and $b$ respectively, a simple calculation provides the general linear solution
\begin{align}
 \mathbf{a} = \widetilde{\mathbf{a}} &:= \sum_{j=1}^{4} d_{j} \mathbf{V}_{j}(r),\label{soln:lin;a}\end{align}
for $\mathbf{a}:=(a,b)\in\mathbb{C}^2$, where
\begin{align}
 \mathbf{V}_{1}(r) &= \sqrt{\frac{k \pi}{2}}\Big(r J_{1}(k r) + \textnormal{i}\left[k^{-1} J_{1}(k r) - rJ_{0}(k r)\right], \;J_{1}(k r) - \textnormal{i} J_{0}(k r) \Big), \quad
 \mathbf{V}_{2}(r) = \sqrt{\frac{k \pi}{2}}\Big(J_{0}(k r) + \textnormal{i} J_{1}(k r), \;0 \Big),\nonumber\\
 \mathbf{V}_{3}(r) &= \sqrt{\frac{k \pi}{2}}\Big(r Y_{1}(k r) + \textnormal{i}\left[k^{-1} Y_{1}(k r) - rY_{0}(k r)\right], \;Y_{1}(k r) - \textnormal{i} Y_{0}(k r) \Big),\quad
 \mathbf{V}_{4}(r) = \sqrt{\frac{k \pi}{2}}\Big(Y_{0}(k r) + \textnormal{i} Y_{1}(k r),\; 0 \Big),\nonumber
\end{align}
where $J_{\nu}, Y_{\nu}$ are $\nu$-th order Bessel functions of the first and second kind, respectively. Similarly, for the linear problem related to \eqref{amp:a;n} and \eqref{amp:a;-n}
\begin{align}
    \frac{\textnormal{d}}{\textnormal{d}\,r} a_{n} &= \lambda_{n}\, a_{n} - \frac{1}{2r}\left(a_{n}-a_{-n}\right),\qquad
    \frac{\textnormal{d}}{\textnormal{d}\,r} a_{-n} = -\lambda_{n}\, a_{-n} + \frac{1}{2r}\left(a_{n}-a_{-n}\right),\nonumber
\end{align}
we find, by defining $b_{n}:=\frac{a_{n}+a_{-n}}{2}$ and $b_{-n}:=\frac{a_{n}-a_{-n}}{2}$ and solving the resulting equations, the general linear solution
\begin{align}
 \mathbf{a}_{n} = \widetilde{\mathbf{a}}_{n} &:= \sum_{j=1}^{2} c_{j,n} \mathbf{W}_{j,n}(r),\label{soln:lin;an}\end{align}
for $\mathbf{a}_{n}:=(a_{n}, a_{-n})\in\mathbb{R}^2$, where
\begin{align}
 \mathbf{W}_{1,n}(r) &= \sqrt{\frac{\lambda_{n}}{2}}\Big(I_{0}(\lambda_{n} r) + I_{1}(\lambda_{n} r), \; I_{0}(\lambda_{n} r) - I_{1}(\lambda_{n} r) \Big), \nonumber\\
 \mathbf{W}_{2,n}(r) &= \sqrt{\frac{\lambda_{n}}{2}}\Big(K_{0}(\lambda_{n} r) - K_{1}(\lambda_{n} r), \; K_{0}(\lambda_{n} r) + K_{1}(\lambda_{n} r) \Big),\nonumber
\end{align}
for $n\in\mathbb{N}$. Here, $I_{\nu}$ and $K_{\nu}$ are $\nu$-th order modified Bessel functions of the first and second kind, respectively. We note that, in both cases, the respective linear solutions $\mathbf{V}_{j}$ and $\mathbf{W}_{j,n}$ split into two types; solutions ($\mathbf{V}_{1}$, $\mathbf{V}_{2}$ and $\mathbf{W}_{1,n}$) that are bounded as $r\to0$, and solutions ($\mathbf{V}_{3}$, $\mathbf{V}_{4}$ and $\mathbf{W}_{2,n}$) that are unbounded as $r\to0$; see Table \ref{tab:asymptote}. By denoting $\langle\cdot,\cdot\rangle$ as the complex dot product
\begin{align}
    \left\langle\cdot, \cdot\right\rangle: \mathbb{C}^{2}\times\mathbb{C}^{2}\to\mathbb{R}, \quad \left\langle\mathbf{x}, \mathbf{y}\right\rangle = \frac{1}{2}\sum_{n=1}^{2} \left\{ \overline{x}_{n} y_{n} + x_{n}\overline{y}_{n}\right\},\nonumber
\end{align}
then the respective eigenmodes $V_{i}(r)$ and $W_{i, n}(r)$ satisfy
\begin{align}
    \left\langle \mathbf{V}^{*}_{i}(r), \mathbf{V}_{j}(r)\right\rangle = \delta_{ij}, \qquad \left\langle \mathbf{W}^{*}_{i,n}(r), \mathbf{W}_{j,n}(r)\right\rangle = \delta_{ij},\nonumber
\end{align}
where $\delta_{ij}$ is the Kronecker delta, and $\mathbf{V}_{j}^{*}(r)$, $\mathbf{W}^{*}_{i,n}(r)$ are the respective adjoint vectors of $\mathbf{V}_{j}(r)$, $\mathbf{W}_{i,n}(r)$, defined by
\begin{align}
 \mathbf{V}^{*}_{1}(r) &= \sqrt{\frac{k \pi}{2}}\Big(0, \;r\left[Y_{0}(k r) + \textnormal{i} Y_{1}(k r) \right]\Big), \quad
 \mathbf{W}^{*}_{1,n}(r) = \sqrt{\frac{\lambda_{n}}{2}}\Big(r\left[K_{1}(\lambda_{n} r) + K_{0}(\lambda_{n} r)\right], \; r\left[K_{1}(\lambda_{n} r) - K_{0}(\lambda_{n} r)\right] \Big),\nonumber\\
 \mathbf{V}^{*}_{3}(r) &= \sqrt{\frac{k \pi}{2}}\Big(0, \;-r\left[J_{0}(k r) + \textnormal{i} J_{1}(k r) \right]\Big),\quad
  \mathbf{W}^{*}_{2,n}(r) = \sqrt{\frac{\lambda_{n}}{2}}\Big(r\left[I_{1}(\lambda_{n} r) - I_{0}(\lambda_{n} r)\right], \; r\left[I_{1}(\lambda_{n} r) + I_{0}(\lambda_{n} r)\right] \Big),\nonumber\\
 \mathbf{V}^{*}_{2}(r) &= \sqrt{\frac{k \pi}{2}}\Big(-r\left[Y_{1}(k r) - \textnormal{i}Y_{0}(k r)\right], \;r^2\left[Y_{1}(k r) - \frac{1}{k r}Y_{0}(k r) - \textnormal{i}Y_{0}(k r)\right] \Big), \nonumber\\
 \mathbf{V}^{*}_{4}(r) &= \sqrt{\frac{k \pi}{2}}\Big(r\left[J_{1}(k r) - \textnormal{i}J_{0}(k r)\right], \;-r^2\left[J_{1}(k r) - \frac{1}{k r}J_{0}(k r) - \textnormal{i}J_{0}(k r)\right] \Big).\nonumber
\end{align}
\begin{table}[hb!]
    \centering
    \begin{tabular}{|c|c|c|}
         \hline
         $ $ & $r\to0$ & $r\to\infty$  \\
             \hline
    $J_{0}(r)$ & $1 + \textnormal{O}\left(r^2\right)$ &  $\sqrt{\frac{2}{\pi r}}\cos\left(r-\frac{\pi}{4}\right) + \textnormal{O}\left(r^{-\frac{3}{2}}\right)$\\ 
    $J_{1}(r)$ & $\frac{r}{2} + \textnormal{O}\left(r^3\right)$ &  $\sqrt{\frac{2}{\pi r}}\sin\left(r-\frac{\pi}{4}\right) + \textnormal{O}\left(r^{-\frac{3}{2}}\right)$\\ 
    $Y_{0}(r)$ & $\frac{2}{\pi}\left(1 + \textnormal{O}\left(r^2\right)\right)\log(r) + \textnormal{O}(1)$ &  $\sqrt{\frac{2}{\pi r}}\sin\left(r-\frac{\pi}{4}\right) + \textnormal{O}\left(r^{-\frac{3}{2}}\right)$\\ 
    $Y_{1}(r)$ & $\frac{1}{\pi}\left(1 + \textnormal{O}\left(r^2\right)\right)r\log(r) - \frac{2}{\pi r} + \textnormal{O}(1)$ &  $-\sqrt{\frac{2}{\pi r}}\cos\left(r-\frac{\pi}{4}\right) + \textnormal{O}\left(r^{-\frac{3}{2}}\right)$\\ 
    $I_{0}(r)$ & $1 + \textnormal{O}\left(r^2\right)$ & $\frac{1}{\sqrt{2 \pi r}}\textnormal{e}^{r} + \textnormal{O}\left(r^{-\frac{3}{2}}\right)$\\ 
    $I_{1}(r)$ &  $\frac{r}{2} + \textnormal{O}\left(r^3\right)$ & $\frac{1}{\sqrt{2 \pi r}}\textnormal{e}^{r} + \textnormal{O}\left(r^{-\frac{3}{2}}\right)$\\
    $K_{0}(r)$ & $-\left(1 + \textnormal{O}\left(r^2\right)\right)\log(r) + \textnormal{O}(1)$ &  $\sqrt{\frac{\pi}{2 r}}\textnormal{e}^{-r} + \textnormal{O}\left(r^{-\frac{3}{2}}\right)$\\ 
    $K_{1}(r)$ & $\frac{1}{2}\left(1 + \textnormal{O}\left(r^2\right)\right)r\log(r) + \frac{1}{r} + \textnormal{O}(1)$ & $\sqrt{\frac{\pi}{2 r}}\textnormal{e}^{-r} + \textnormal{O}\left(r^{-\frac{3}{2}}\right)$\\
    \hline
    \end{tabular}
    \caption{Expansions of the Bessel functions $J_{\nu}(r), Y_{\nu}(r), I_{\nu}(r), K_{\nu}(r)$ for $r\to0$ and $r\to\infty$ (see \cite{abramowitz1972handbook}, (9.1.10-11), (9.6.10-11), $\S$9.2 and $\S$9.7).}
    \label{tab:asymptote}
\end{table}
We apply the method of variation of constants \cite[\S2.II]{Walter1998ODEs} in conjunction with the general linear solution \eqref{soln:lin;a}, to find a general solution to the full nonlinear system \eqref{amp:a}-\eqref{amp:b} for the complex amplitudes $a,b$. This solution has the form
\begin{align}
    \mathbf{a}(r) &= \sum_{j=1}^{4} \langle \mathbf{V}^{*}_{j},\widetilde{\mathbf{a}}\rangle\mathbf{V}_{j}(r) + \sum_{j=1}^{2} \mathbf{V}_{j}(r) \int_{r_{0}}^{r}\langle \mathbf{V}^{*}_{j},\mathbf{F}\rangle \;\textnormal{d}s + \sum_{j=3}^{4} \mathbf{V}_{j}(r)\int_{0}^{r} \langle \mathbf{V}^{*}_{j},\mathbf{F}\rangle\;\textnormal{d}s, \label{fixed:a}\\
    &= \sum_{j=1}^{4} d_{j}\,\mathbf{V}_{j}(r) + \sum_{j=1}^{2} \mathbf{V}_{j}(r) \int_{r_{0}}^{r}\langle \mathbf{V}^{*}_{j},\mathbf{F}\rangle \;\textnormal{d}s + \sum_{j=3}^{4} \mathbf{V}_{j}(r)\int_{0}^{r} \langle \mathbf{V}^{*}_{j},\mathbf{F}\rangle\;\textnormal{d}s,\nonumber
\end{align}
where $\mathbf{F} := \left(-\Omega\left(\mathcal{F}, \overline{\mathbf{f}}\right), \Omega\left(\mathcal{F}, \overline{\mathbf{e}}\right)\right)$ and $\mathbf{e}(y),\mathbf{f}(y)$ are the eigenmodes defined in \eqref{defn:eigenmode} for the imaginary eigenvalues $\textnormal{i}k$ and a large fixed radial value $r_{0}>0$. Similarly, for the pairs of real-valued amplitudes $\mathbf{a}_{n}=(a_{n}, a_{-n})$, we take the linear solution \eqref{soln:lin;an} and apply the method of variation of constants to derive a general solution for the full nonlinear system \eqref{amp:a;n}-\eqref{amp:a;-n}, which has the form
\begin{align}
    \mathbf{a}_{n}(r) &= \sum_{j=1}^{2} \langle \mathbf{W}^{*}_{j,n},\widetilde{\mathbf{a}}_{n}\rangle\mathbf{W}_{j,n}(r) + \mathbf{W}_{1,n}(r) \int_{r_{0}}^{r}\langle \mathbf{W}^{*}_{1,n},\mathbf{F}_{n}\rangle \;\textnormal{d}s  + \mathbf{W}_{2,n}(r)\int_{0}^{r}\langle \mathbf{W}^{*}_{2,n},\mathbf{F}_{n}\rangle\;\textnormal{d}s, \label{fixed:a;n}\\
    &= \sum_{j=1}^{2} c_{j,n}\,\mathbf{W}_{j,n}(r) + \mathbf{W}_{1,n}(r) \int_{r_{0}}^{r}\langle \mathbf{W}^{*}_{1,n},\mathbf{F}_{n}\rangle \;\textnormal{d}s  + \mathbf{W}_{2,n}(r)\int_{0}^{r}\langle \mathbf{W}^{*}_{2,n},\mathbf{F}_{n}\rangle\;\textnormal{d}s,\nonumber
\end{align}
where $\mathbf{F}_{n} := \left(-\Omega\left(\mathcal{F}, \mathbf{e}_{-n}\right), \Omega\left(\mathcal{F}, \mathbf{e}_{n}\right)\right)$ and $\mathbf{e}_{n}(y), \mathbf{e}_{- n}(y)$ are the eigenmodes defined in \eqref{defn:eigenmode} for the respective real eigenvalues $\lambda_{n}, -\lambda_{n}$. Using the $r\to0$ asymptotic forms for the Bessel functions detailed in Table \ref{tab:asymptote}, we note that
\begin{align}
    &\left|\mathbf{V}_{1}(r)\right| = \textnormal{O}(1), \quad & &\left|\mathbf{V}_{2}(r)\right| = \textnormal{O}(1), \quad & &\left|\mathbf{V}_{3}(r)\right| = \textnormal{O}(r^{-1}), \quad & &\left|\mathbf{V}_{4}(r)\right| = \textnormal{O}(r^{-1}), \nonumber\\
    &\left|\mathbf{V}^{*}_{1}(r)\right| = \textnormal{O}(1), \quad & &\left|\mathbf{V}^{*}_{2}(r)\right| = \textnormal{O}(1), \quad & &\left|\mathbf{V}^{*}_{3}(r)\right| = \textnormal{O}(r), \quad & &\left|\mathbf{V}^{*}_{4}(r)\right| = \textnormal{O}(r), \nonumber\\
    &\left|\mathbf{W}_{1,n}(r)\right| = \textnormal{O}(1), \quad & &\left|\mathbf{W}_{2,n}(r)\right| = \textnormal{O}(r^{-1}), \quad & &\left|\mathbf{W}^{*}_{1,n}(r)\right| = \textnormal{O}(1), \quad & &\left|\mathbf{W}^{*}_{2,n}(r)\right| = \textnormal{O}(r) \nonumber
\end{align}
as $r\to0$. We see that the linear terms in \eqref{fixed:a} are only bounded as $r\to0$ if $d_{3}=d_{4}=0$, and similarly the linear terms in \eqref{fixed:a;n} are only bounded as $r\to0$ if $c_{2,n}=0$ for each $n\in\mathbb{N}$. {We note that there are no linear terms that grow exponentially as $r\to0$, and so} solutions of \eqref{fixed:a} and \eqref{fixed:a;n} form a parametrisation for the full {local} centre-unstable manifold $\mathcal{W}^{cu}_{-}{(\varepsilon)}$, containing all solutions that have at most algebraic growth as $r\to0$. To parametrise the core manifold $\widetilde{\mathcal{W}}^{cu}_{-}{(\varepsilon)}$, we want to restrict each respective solution of \eqref{fixed:a} and \eqref{fixed:a;n} to be solutions of equations \eqref{amp:a}-\eqref{amp:b} and \eqref{amp:a;n}-\eqref{amp:a;-n} which are bounded on $r\in[0, r_{0}]$. Setting $d_{3}=d_{4}=c_{2,n}=0$ for all $n\in\mathbb{N}$, the integral terms on the right-hand side of \eqref{fixed:a} and \eqref{fixed:a;n} are bounded as $r\to0$, therefore $\mathbf{a}(r)$ and $\mathbf{a}_{n}(r)$ are bounded for all $r\in[0,r_{0}]$. For notational purposes, we define 
\begin{align}
    {\mathbf{c}_{i}:=(c_{i,1},\, c_{i,2},\, \dots), \qquad |\mathbf{c}_{i}|_{1}:= \sum_{j=1}^{\infty}|c_{i, j}|.} \label{ci:defn}
\end{align}
{The} respective equations \eqref{fixed:a} and \eqref{fixed:a;n} can be solved for sufficiently small $(d_{1}, d_{2},\mathbf{c}_{1})$ and $\varepsilon$. These solutions take the form
\begin{align}
    \mathbf{a}(r) &= \sum_{j=1}^{4} \widetilde{d}_{j}(d_{1}, d_{2}, \mathbf{c}_{1}{; \varepsilon}) \mathbf{V}_{j}(r), \qquad \mathbf{a}_{n}(r) = \sum_{j=1}^{2} \widetilde{c}_{j,n}(d_{1}, d_{2}, \mathbf{c}_{1}{; \varepsilon}) \mathbf{W}_{j,n}(r)\nonumber
\end{align}
on $[0,r_{0}]$, with respective nonlinear functions $\widetilde{d}_{j}(d_{1}, d_{2}, \mathbf{c}_{1}{; \varepsilon})$ and $\widetilde{c}_{j,n}(d_{1}, d_{2}, \mathbf{c}_{1}{; \varepsilon})$
\begin{align}
        \widetilde{c}_{1,n} &= c_{1,n} + \textnormal{O}_{r}\left(\left(|d_{1}| + |d_{2}| + |\mathbf{c}_{1}|_{1}\right)\left(|\varepsilon| + |d_{1}| + |d_{2}| + |\mathbf{c}_{1}|_{1}\right)\right),\nonumber\\
        \widetilde{c}_{2,n} &= \textnormal{O}_{r}\left(\left(|d_{1}| + |d_{2}| + |\mathbf{c}_{1}|_{1}\right)\left(|\varepsilon| + |d_{1}| + |d_{2}| + |\mathbf{c}_{1}|_{1}\right)\right),\label{amp:core;param}\\
        \widetilde{d}_{j} &= d_{j} + \textnormal{O}_{r}\left(\left(|d_{1}| + |d_{2}| + |\mathbf{c}_{1}|_{1}\right)\left(|\varepsilon| + |d_{1}| + |d_{2}| + |\mathbf{c}_{1}|_{1}\right)\right),\nonumber \\
        \intertext{for $j=1,2$, and}
        \widetilde{d}_{j} &= \textnormal{O}_{r}\left(\left(|d_{1}| + |d_{2}| + |\mathbf{c}_{1}|_{1}\right)\left(|\varepsilon| + |d_{1}| + |d_{2}| + |\mathbf{c}_{1}|_{1}\right)\right),\nonumber 
\end{align}
for $j=3,4$. Here, $\textnormal{O}_{r}(\dots)$ denotes the standard Landau symbol where the bounding constants may depend on a fixed $r\in[0,r_{0}]$. We evaluate \eqref{fixed:a} and \eqref{fixed:a;n} at $r=r_{0}$, and determine the quadratic coefficient for $d_{2}$ in front of $\mathbf{V}_{3}(r_{0})$, 
\begin{align}
    \widetilde{d}_{3} = \left[\nu + \textnormal{O}(r_{0}^{-\frac{1}{2}})\right] d_{2}^{2} + \textnormal{O}_{r}\left(\left(|d_{1}| + |d_{2}| + |\mathbf{c}_{1}|_{1}\right)\left(|\varepsilon| + |d_{1}| + |\mathbf{c}_{1}|_{1}\right) + |d_{2}|^{3}\right),\nonumber 
\end{align}
where $\nu$ is found in Appendix \ref{app:core;quad} to be
\begin{align}
    \nu &={M_{0}\sqrt{\frac{3 k \pi}{2}}\frac{k \mathcal{M}}{ D m^{3}}\left[1- 2\;\textnormal{sech}^{2}(k D)\right]},\label{nu:defn}
\end{align}
with $M_{0}:= \frac{\mu-1}{\mu+1}$ and $m:=\sqrt{1+(kD\tanh(kD)-1)\mathcal{M}\,\textnormal{sech}^{2}(kD)}$. Then, we can parametrise the core manifold $\widetilde{\mathcal{W}}^{cu}_{-}{(\varepsilon)}$, containing all small-amplitude solutions of \eqref{full:syst;linbd} which are bounded as $r\to0$, as the following
\begin{align}
    \widetilde{\mathcal{W}}^{cu}_{-}{(\varepsilon)} = \left\{\mathbf{u} =\widetilde{\mathbf{u}}\,:\; (a,b) = \sum_{j=1}^4 \widetilde{d}_{j}\left(d_{1}, d_{2}, \mathbf{c}_{1}{; \varepsilon}\right) \mathbf{V}_{j}(r), \quad (a_{n}, a_{-n}) = \sum_{j=1}^2 \widetilde{c}_{j,n}\left(d_{1}, d_{2}, \mathbf{c}_{1}{; \varepsilon}\right) \mathbf{W}_{j,n}(r)\right\},\nonumber
\end{align}
where $\widetilde{\mathbf{u}}$ is the decomposition defined in \eqref{u:defn}. It is convenient to express the core manifold up to leading order for large $r_{0}$
\begin{align}
    a(r_{0}) &= \textnormal{e}^{\textnormal{i}\left(k r_{0} - \frac{\pi}{4}\right)}r_{0}^{-\frac{1}{2}}\left([-\textnormal{i}+\textnormal{O}(r_{0}^{-1})]r_{0} d_{1} + [1+\textnormal{O}(r_{0}^{-1})] d_{2}\right) + \textnormal{O}_{r_{0}}\left(|\mathbf{d}|_{1}\left[|\varepsilon| + |\mathbf{d}|_{1}\right]\right),\nonumber\\
    b(r_{0}) &= \textnormal{e}^{\textnormal{i}\left(k r_{0} - \frac{\pi}{4}\right)}r_{0}^{-\frac{1}{2}}\left(\left[-\textnormal{i} +\textnormal{O}(r_{0}^{-1})\right] d_{1} - \left[\nu + \textnormal{O}(r_{0}^{-\frac{1}{2}})\right] d_{2}^{2}\right)+ \textnormal{O}_{r_{0}}\left(|\mathbf{d}|_{1}\left[|\varepsilon| + |\mathbf{d}_{2}|_{1}\right] + |d_{2}|^{3}\right),\nonumber\\
    a_{n}(r_{0}) &= \textnormal{e}^{\lambda_{n}r_{0}}r_{0}^{-\frac{1}{2}}\left[\frac{1}{\sqrt{\pi}} + \textnormal{O}(r_{0}^{-1})\right] c_{1,n} + \textnormal{O}_{r_{0}}\left(|\mathbf{d}|_{1}\left[|\varepsilon| + |\mathbf{d}|_{1}\right]\right),\label{match:core}\\
    a_{-n}(r_{0}) &= \textnormal{O}_{r_{0}}\left(|\mathbf{d}|_{1}\left[|\varepsilon| + |\mathbf{d}|_{1}\right]\right),\nonumber
\end{align}
where we have used the large-argument expansions of Bessel functions noted in Table \ref{tab:asymptote} and defined $\mathbf{d}, \mathbf{d}_{2}$ such that { $\mathbf{d}=(d_{1}, d_{2}, c_{1,1}, c_{1,2}, \dots)$ and $\mathbf{d}_{2}=(d_{1}, 0, c_{1,1}, c_{1,2}, \dots)$, equipped with the norm $|\mathbf{d}|_{1}:= |d_{1}| + |d_{2}| + |\mathbf{c}_{1}|_{1}$.}
\section{Far-field Solutions}
\label{s:farfield}
In this section, we now turn our attention to the far-field region, where $r$ is larger than some fixed {$r_{\infty}\gg1$}, and construct the `far-field' manifold; the set of all small-amplitude solutions that decay to zero exponentially fast as $r\to\infty$. We augment \eqref{full:syst;linbd} by introducing the variable $\sigma:=\frac{1}{r}$ so that the extended system becomes
\begin{align}
    \mathbf{u}_{r} &= \mathbf{L}\left(\frac{1}{\sigma}\right)\mathbf{u} + \mathcal{F}(\mathbf{u},\varepsilon,\sigma), \nonumber\\
    \sigma_{r} &= -\sigma^2.\label{sigma:r}
\end{align}
{By introducing $\sigma$ and treating it as a {small} unconstrained function of $r$, (\ref{sigma:r}) is now considered to be an autonomous system. This is because the $1/r$ terms in $\mathbf{L}(r)$ in \eqref{full:syst;linbd} are now of `quadratic' order, and hence they are absorbed into the nonlinearity, so the remaining linear operator $\mathbf{L}_{\infty}$ is independent of $r$. We will formulate our parametrisations for the far-field manifolds in the extended system \eqref{sigma:r}; however, we will then return to the non-autonomous system \eqref{full:syst;linbd} by evaluating at $\sigma=\frac{1}{r}$, in order to complete the matching process with the core manifold. We note that \eqref{sigma:r} recovers $\sigma(r)=\frac{1}{r+c}$ for some constant $c$; for $r\gg1$, we see that $\sigma(r)=\frac{1}{r} + \textnormal{O}\left(\frac{1}{r^{2}}\right)$, and so any correction terms arising from setting $\sigma=\frac{1}{r}$ will not appear in the leading order expansion of our parametrisations, as they will be $\textnormal{O}\left(\frac{1}{r^{2}}\right)$ (or, equivalently, $\textnormal{O}\left(\sigma^{2}\right)$ ).} By applying the same `spectral decomposition' as detailed in Section \ref{s:spect}, we arrive at the extended versions of \eqref{amp:a}-\eqref{amp:a;-n}, namely
\begin{align}
    &\frac{\textnormal{d}}{\textnormal{d}r}a = \textnormal{i}k a + b - \frac{\sigma}{2}(a-\overline{a}) - \Omega(\mathcal{F}, \overline{\mathbf{f}}),  & \quad &\frac{\textnormal{d}}{\textnormal{d}r}a_{n} = \lambda_{n} a_{n} - \frac{\sigma}{2}(a_{n}-a_{-n}) - \Omega(\mathcal{F}, \mathbf{e}_{-n}), & \nonumber\\
    &\frac{\textnormal{d}}{\textnormal{d}r}b = \textnormal{i}k b - \frac{\sigma}{2}(b+\overline{b}) + \Omega(\mathcal{F}, \overline{\mathbf{e}}),  & \quad & \frac{\textnormal{d}}{\textnormal{d}r}a_{-n} = -\lambda_{n} a_{-n} + \frac{\sigma}{2}(a_{n}-a_{-n}) + \Omega(\mathcal{F}, \mathbf{e}_{n}), & \nonumber\\
    &{\frac{\textnormal{d}}{\textnormal{d}r}\sigma = -\sigma^{2}}, & & &\label{syst:farf;amp}
\end{align}
for $n\in\mathbb{N}$. {We define $U(r):=(a,b,a_{1},a_{-1},a_{2},a_{-2},\dots)(r)$ and denote the full solution to \eqref{syst:farf;amp} as $(U,\sigma)(r)$. We note that \{$U(r)=\mathbf{0}$\} and \{$\sigma(r)=\frac{1}{r}$\} are invariant subspaces, and any $\sigma(r)>0$ that satisfies $\frac{\textnormal{d}}{\textnormal{d}r}\sigma=-\sigma^{2}$ is a forward-bounded solution in $r$. We begin by investigating the linear dynamics of solutions close to the trivial state.}
{\paragraph{Linear Dynamics} At the bifurcation point $\varepsilon=0$, we linearise \eqref{syst:farf;amp}} about the trivial {equilibrium {$(U,\sigma)=(\mathbf{0},0)$}; then,} the system has the form
\begin{align}
    &\frac{\textnormal{d}}{\textnormal{d}r}a = \textnormal{i}k a + b,  & \quad &\frac{\textnormal{d}}{\textnormal{d}r}a_{n} = \lambda_{n} a_{n}, & \quad &\frac{\textnormal{d}}{\textnormal{d}r}\sigma = 0,\nonumber\\
    &\frac{\textnormal{d}}{\textnormal{d}r}b = \textnormal{i}k b,  & \quad & \frac{\textnormal{d}}{\textnormal{d}r}a_{-n} = -\lambda_{n} a_{-n}. & \label{lin:far;ext}
\end{align}

\begin{figure}[t!]
    \centering
    \includegraphics[width=\linewidth]{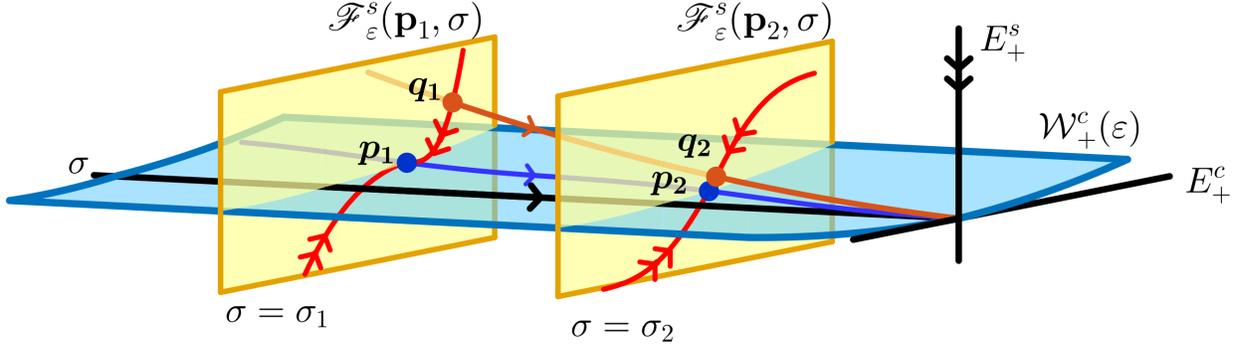}
   \caption{The centre-stable manifold $\mathcal{W}^{cs}_{+}$ is parametrised by the foliation $\mathscr{F}^{s}(\mathbf{p},\sigma)$ over base points $\mathbf{p}$ in the centre manifold $\mathcal{W}^{c}_{+}{|_{\sigma=\frac{1}{r}}}$. Stable fibres are displayed in respective Poincar\'e sections $\sigma = \sigma_{1}, \sigma_{2}$, with base points $\mathbf{p}_{1}, \mathbf{p}_{2}\in\mathcal{W}^{c}_{+}{|_{\sigma=\frac{1}{r}}}$. Solutions related to initial data ${(\mathbf{q}_{1},\sigma_{1})}\in\mathscr{F}^{s}(\mathbf{p}_{1},\sigma_{1})$ approach the solution related to the initial {data $(\mathbf{p}_{1},\sigma_{1})$, where $\mathbf{p}_{1}\in\mathcal{W}^{c}_{+}{|_{\sigma=\frac{1}{r}}}$}, exponentially as $r\to\infty$.}
    \label{fig:foliation}
\end{figure}

\begin{figure}[!t]
    \centering
\includegraphics[width=\linewidth]{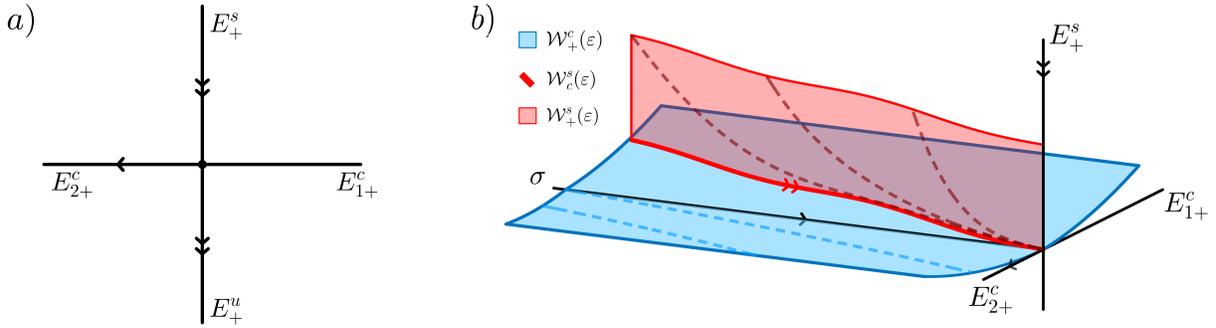}
    \caption{ a) The linear eigenspace at $\sigma={\sigma^{*}}${, for $0\leq\sigma^{*}\ll1$ fixed, close to $U=0$}; there are infinite-dimensional stable and unstable spaces $E^{s}_{+}$ and $E^{u}_{+}$, respectively, as well as a four-dimensional centre space. The double and single arrows indicate exponential and algebraic behaviour, respectively. b) The {extended} far-field manifold $\mathcal{W}^{s}_{+}{(\varepsilon)}$ (red) is seen as a perturbation away from $\mathcal{W}^{s}_{c}{(\varepsilon)}$ (bold red line), which contains solutions {$(U,\sigma)(r)$ which are} bounded for all $r${, and $U(r)$ decays exponentially} as $r\to\infty$. In turn, the manifold $\mathcal{W}^{s}_{c}{(\varepsilon)}$ lies within the centre manifold $\mathcal{W}^{c}_{+}{(\varepsilon)}$ (blue).}
    \label{fig:farfman}
\end{figure}

We note that the equation for $\sigma$ decouples from the rest of the system and so we can solve this to find $\sigma=\sigma^{*}$, where $\sigma^{*}$ is an arbitrary {small} constant; {then, $\{\sigma=\sigma^{*}\}$ is an invariant set of \eqref{lin:far;ext}}. Therefore, for solutions {sufficiently close to the trivial equilibrium}, we can look at the spatial dynamics of {$U(r)$ for} the linear problem \eqref{lin:far;ext} on the Poincar\'e section $\sigma = \sigma^{*}$, where $\sigma^{*}$ is constant, see Figure \ref{fig:foliation}. For a fixed $\sigma$, the rest of the linear problem \eqref{lin:far;ext} has {a four dimensional centre eigenspace $E^{c}_{+}:=\{a_{\pm n}=0, \forall n\in\mathbb{N}\}$, a countably infinite dimensional stable eigenspace $E^{s}_{+}:=\{ a=b=a_{ n}=0, \forall n\in\mathbb{N}\}$, and a countably infinite dimensional unstable eigenspace $E^{u}_{+}:=\{ a=b=a_{- n}=0, \forall n\in\mathbb{N}\}$}; see Figure \ref{fig:farfman}a). { We then define the centre eigenspace of the extended system to be $\widetilde{E}^{c}_{+}:=E^{c}_{+}\times\mathbb{R}$, while the stable and unstable eigenspaces are unchanged. We {expect} that, for sufficiently small fixed $\varepsilon$, a local centre manifold $\mathcal{W}^{c}_{0}$ exists and, for a neighbourhood $\mathcal{O}$ of the equilibrium $(U,\sigma)=(\mathbf{0},0)$, has the following properties (see \cite[Thm. 2.1]{Haragus2011Bifurcation}, \cite[Thm. 3]{vanderbauwhede1992center} \& \cite[Thm. 2.9]{mielke2006hamiltonian}):
\begin{enumerate}[label=\roman*.]
    \item If a solution $(U,\sigma)(r)$ is globally-bounded and $(U,\sigma)(r)\in\mathcal{O}$ for all $r\in\mathbb{R}^{+}$, then $(U,\sigma)(r)\in\mathcal{W}^{c}_{0}$.
    \item If, for $r=r_{1}$, a solution $(U,\sigma)(r_{1})\in\mathcal{W}^{c}_{0}\cap\mathcal{O}$ and $(U,\sigma)(r)\in\mathcal{O}$ for all $r$ in a neighbourhood of $r=r_{1}$, then $(U,\sigma)(r)\in\mathcal{W}^{c}_{0}$ for all $r$ in a neighbourhood of $r=r_{1}$.
    \item $\mathcal{W}^{c}_{0}$ is tangential to $\widetilde{E}^{c}_{+}$ at the equilibrium $(U,\sigma)=(\mathbf{0},0)$.
    \item $\textnormal{dim}\,[\mathcal{W}^{c}_{0}] = \textnormal{dim}\,[\widetilde{E}^{c}_{+}]$.
    \item $\mathcal{W}^{c}_{0}$ can be expressed as the graph of a smooth function $\tau^{c}$ that maps the centre eigenspace $\widetilde{E}^{c}_{+}$ to the full solution space.
\end{enumerate}
{We note that the resolvent estimates necessary in order to prove the existence of $\mathcal{W}^{c}_{0}$ have not been computed, and so this remains a key assumption.} Hence, taking the equilibrium $(U,\sigma)=(\mathbf{0},0)$ to be the limit as $r\to\infty$, solutions contained in $\mathcal{W}^{c}_{+}(\varepsilon)$ remain sufficiently small for all large $r$, i.e. for $r\in[r_{\infty},\infty)$, for some large fixed $r_{\infty}$. Similarly, a centre-stable manifold $\mathcal{W}^{cs}_{+}(\varepsilon)$ contains all sufficiently small solutions that remain sufficiently small as $r$ increases, for $r\in[r_{\infty},\infty)$. We refer to these solutions as being "locally bounded" and "locally forward-bounded", respectively. These manifolds are smooth in $0\leq\varepsilon\ll1$, and hence can be parametrised by $\textnormal{O}(\varepsilon)$-perturbations of the manifolds $\mathcal{W}^{c}_{+}(0)$ and $\mathcal{W}^{cs}_{+}(0)$. {We also note that the respective centre and centre-stable manifolds are locally bounded and forward-bounded in the sense of exponentially-weighted spaces, and thus contain solutions with algebraic growth.} {Hence, it is worth remarking that, for fixed $\sigma=\sigma^{*}\ll1$,} the four-dimensional centre {eigenspace $E^{c}_{+}$} can be decomposed further: two directions define solutions that remain bounded for all $r$, while the other two directions define solutions that grow algebraically as $r\to\infty$. Then, we expect that the {stable part of the} expected centre manifold $\mathcal{W}^{c}_{+}{(\varepsilon)}|_{\sigma=\sigma^{*}}$, containing all solutions { in $\mathcal{W}^{c}_{+}(\varepsilon)|_{\sigma=\sigma^{*}}$ that decay exponentially as $r\to\infty$}, is {a two-dimensional submanifold of the entire four-dimensional centre manifold.
\paragraph{Constructing the far-field manifold:} Having investigated the linear dynamics of \eqref{syst:farf;amp} for sufficiently small $|U|$, $|\sigma|$, \& $|\varepsilon|$, we will now use this to motivate the construction of the far-field manifold; the rest of this section is dedicated to achieving the following aim:}\\
\underline{Goal:} We wish to isolate solutions of \eqref{amp:a}-\eqref{amp:a;-n} that exhibit exponential decay as $r\to\infty$. Equivalently, we wish to find functions of the form $(U,\sigma)(r)=(U^{*}(r),\frac{1}{r})$ that solve \eqref{syst:farf;amp} such that $U^{*}(r)\in\mathcal{W}^{s}_{+}(\varepsilon)|_{\sigma=\frac{1}{r}}$ for $r\in[r_{\infty},\infty)$ and sufficiently small $|U^{*}|$, $|\varepsilon|$.} 

The following theory is well-defined in infinite dimensions; however, in order to discuss the dimensions of each manifold, we will refer to a countably infinite set of dimensions as being $n$-dimensional. Conceptually, we think of this as an order $n$-truncation as $n$ becomes very large. Formally, we would expect to find a $({5}+n)$-dimensional centre-stable manifold $\mathcal{W}^{cs}_{+}{(\varepsilon)}$, containing all {small} solutions that {remain bounded as} $r\to\infty$; an $n$-dimensional stable manifold $\mathcal{W}^{s}_{+}{(\varepsilon)}$, containing all {small} solutions that exponentially decay as $r\to\infty$; a {five}-dimensional centre manifold $\mathcal{W}^{c}_{+}{(\varepsilon)}$, containing all {small, locally bounded} solutions; {and a three-dimensional submanifold $\mathcal{W}^{{s}}_{{c}}{(\varepsilon)}\subseteq\mathcal{W}^{c}_{+}{(\varepsilon)}$, that contains { all small} solutions ${(U,\sigma)\in}\mathcal{W}^{c}_{+}{(\varepsilon)}$ {such that $U(r)$ decays exponentially} as $r\to\infty$}; see Figure \ref{fig:farfman}b) for a sketch of this idea.

We want to parametrise the far-field stable manifold $\mathcal{W}^{s}_{+}{(\varepsilon)}{|_{\sigma=\frac{1}{r}}}$ {in the invariant subspace $\{\sigma=\frac{1}{r}\}$}, containing all small-amplitude solutions that exponentially decay to zero as $r\to\infty$, which we would expect to be $n$-dimensional. However, there are a few problems:
\begin{itemize}
    \item The core manifold $\widetilde{\mathcal{W}}^{cu}_{-}{(\varepsilon)}$ parametrised in $\S\ref{s:core}$ is $(2+n)$-dimensional, so there are not currently enough parameters to match {$\mathcal{W}^{s}_{+}(\varepsilon)|_{\sigma=\frac{1}{r}}$ with $\widetilde{\mathcal{W}}^{cu}_{-}(\varepsilon)$}.
    \item In order to perform standard analytical techniques, we wish to apply normal-form transformations to make our system more tractable.
    However, these are hard to formulate unless working on a centre manifold \cite{scheel2003radially}.
\end{itemize}
Hence, we need to find a way to parametrise the far-field manifold in $(2+n)$ dimensions, such that any required analysis can be applied on the centre manifold. In \cite{mcquighan2014oscillons}, McQuighan \& Sandstede faced a similar problem for the case of a finite-dimensional ODE system, which they overcame by utilising the theory of stable foliations; see \cite{chow1990smooth,aulbach2003foliation}. 

The motivating theory behind this can be seen as follows; we already established that we expect to find a $({5}+n)$-dimensional centre-stable manifold. By definition, the centre and stable manifolds should both be contained within the centre-stable manifold, such that we could parametrise the stable manifold as an $n$-dimensional manifold contained in the larger centre-stable manifold. This parametrisation would isolate all of the solutions that decay exponentially as $r\to\infty$, but there would only be $n$ parameters, causing the matching problem.

We could parametrise the stable manifold with respect to the coordinates of the centre-stable manifold, that is, displacements in all of the centre and stable directions. By performing analysis on each parameter to ensure they decay to zero as $r\to\infty$ {and evaluating at $\sigma=\frac{1}{r}$}, we could reduce this manifold to $(2+n)$-dimensions, and so the matching problem would be solved. However, formulating the required analysis remains a difficult task.

We will instead choose to parametrise the stable manifold in relation to the {five}-dimensional centre manifold. By taking a point in the centre manifold, parametrised {by $E^{c}_{+}\cup\{(U,\sigma):\,U=0\}$}, and then taking displacements from that point in each of the stable directions, we can construct a $({5}+n)$-dimensional parametrisation that covers the entire centre-stable manifold close to the trivial state $\mathbf{u}\equiv\mathbf{0}$. As before, in order to isolate the far-field stable manifold and reduce the parametrisation to $(2+n)$ dimensions, we would need to {evaluate $\sigma=\frac{1}{r}$ and} ensure that every parameter decays to zero as $r\to\infty$. We know that the stable displacements tend to zero, and so the final analysis can be conducted on the remaining centre manifold parameters. This allows us to perform normal-form transformations on the four-dimensional centre manifold parameters in order to find the set of solutions that lie on the centre manifold and decay exponentially as $r\to\infty$, which we call the stable part of the centre manifold and denote as ${\mathcal{W}^{s}_{c}(\varepsilon)|_{\sigma=\frac{1}{r}}}\subset\mathcal{W}^{c}_{+}{(\varepsilon)|_{\sigma=\frac{1}{r}}}$. Since, {by definition,} any solution in $\mathcal{W}^{s}_{c}{(\varepsilon)|_{\sigma=\frac{1}{r}}}$ {does not grow algebraically} as $r\to\infty$, we { expect that} $\mathcal{W}^{s}_{c}{(\varepsilon)|_{\sigma=\frac{1}{r}}}$ can be parametrised as a two-dimensional manifold.

In particular, as illustrated in Figure \ref{fig:foliation}, the centre-stable manifold {$\mathcal{W}^{cs}_{+}{(\varepsilon)}|_{\sigma=\frac{1}{r}}$} can be written as a foliation, i.e. a union of stable fibres, $\bigcup_{\mathbf{p}}\{\mathscr{F}^{s}{_{\varepsilon}}(\mathbf{p},{\sigma})\}$ over base points $\mathbf{p}$ in the full centre manifold $\mathcal{W}^{c}_{+}{(\varepsilon)}|_{{\sigma=\frac{1}{r}}}\subset \mathcal{W}^{cs}_{+}{(\varepsilon)}|_{{\sigma=\frac{1}{r}}}$. {That is,
\begin{align}
    {\mathcal{W}^{cs}_{+}{(\varepsilon)}|_{\sigma=\frac{1}{r}}} = \bigcup_{\mathbf{p}\in\mathcal{W}^{c}_{+}{(\varepsilon)}|_{\sigma={\frac{1}{r}}}} \left\{\mathscr{F}^{s}{_{\varepsilon}}(\mathbf{p}, {\sigma})\right\},\nonumber
\end{align}
where $\mathscr{F}^{s}_{\varepsilon}(\cdot,\sigma)$ is smooth in $\varepsilon$ and $\sigma$.} Then, a trajectory beginning at $(\mathbf{q},{\sigma}){\in\mathcal{W}^{cs}_{+}}{(\varepsilon)}$ will converge to zero as $r\to\infty$ if and only if its related base point $(\mathbf{p},{\sigma}){\in\mathcal{W}^{c}_{+}}{(\varepsilon)}$ on the centre manifold does. The stable manifold, defined as the set of all solutions that exponentially tend to zero as $r\to\infty$, is hence given by the union of all stable fibres associated with solutions on the centre manifold that exponentially decay as $r\to\infty$,
\begin{align}
    {\mathcal{W}^{s}_{+}{(\varepsilon)}|_{\sigma=\frac{1}{r}}} := \bigcup_{\mathbf{p}\in\mathcal{W}^{s}_{c}{(\varepsilon)}|_{\sigma={\frac{1}{r}}}} \left\{\mathscr{F}^{s}{_{\varepsilon}}(\mathbf{p}, {\sigma})\right\}.\nonumber
\end{align}

We begin by finding a parametrisation of all stable fibres $\mathscr{F}^{s}{_{\varepsilon}}(\mathbf{p},{\sigma})$ with respect to their base points $\mathbf{p}$ on the full centre manifold $\mathcal{W}^{c}_{+}{(\varepsilon)}|_{\sigma={\frac{1}{r}}}$ {(see \S\ref{subs:far-param})}. Following this, {we perform the following procedure to find localised radial patterns:
\begin{itemize}
    \item We reduce the extended system \eqref{syst:farf;amp} to the five-dimensional centre-manifold $\mathcal{W}^{c}_{+}(\varepsilon)$ and perform normal-form transformations to make the problem analytically tractable (see \S\ref{subs:far-cent}).
    \item In later sections, we use geometric blow-up coordinates to find solutions $(U,\sigma)=(\widetilde{U}(r),\widetilde{\sigma}(r))\in\mathcal{W}^{c}_{+}(\varepsilon)$, where  $\widetilde{U}(r)$ decays exponentially as $r\to\infty$ {(see \S\ref{s:geoblowup}, \S\ref{s:match})}.
    \item Evaluating this solution at $\sigma=\frac{1}{r}$, we have $\left(\widehat{U}(r),\frac{1}{r}\right)$, where $\widehat{U}(r)$ solves \eqref{amp:a}-\eqref{amp:a;-n} and $\widehat{U}(r)\in\mathcal{W}^{s}_{c}(\varepsilon)|_{\sigma=\frac{1}{r}}$ {(see \S\ref{subs:spotA;track},  \S\ref{s:fold}, \S\ref{subs:spot;B,3} \& \S\ref{subs:rings;track})}.
    \item We substitute our solution $\widehat{U}(r)\in\mathcal{W}^{s}_{c}(\varepsilon)|_{\sigma=\frac{1}{r}}$ into the foliation parametrisation for $\mathcal{W}^{s}_{+}(\varepsilon)|_{\sigma=\frac{1}{r}}$, resulting in a $(2+n)$-dimensional parametrisation {(see \S\ref{subs:spotA;match}, \S\ref{s:fold;match}, \S\ref{subs:spot;B,match} \& \S\ref{subs:rings;match})}.
    \item Finally, the core manifold $\widetilde{\mathcal{W}}^{cu}_{-}(\varepsilon)$ and far-field manifold $\widetilde{\mathcal{W}}^{s}_{+}(\varepsilon):=\mathcal{W}^{s}_{+}(\varepsilon)|_{\sigma=\frac{1}{r}}$ are both evaluated at $r=r_{0}$, and the respective coordinates of each parametrisation are matched {(see \S\ref{subs:spotA;match}, \S\ref{s:fold;match}, \S\ref{subs:spot;B,match} \& \S\ref{subs:rings;match})}. Any solution that lies on $\widetilde{\mathcal{W}}^{cu}_{-}(\varepsilon)|_{r=r_{0}}\cap\widetilde{\mathcal{W}}^{s}_{+}(\varepsilon)|_{r=r_{0}}$ is, by definition, a localised radial pattern.
\end{itemize}
}

\subsection{Parametrisation of the far-field manifolds}\label{subs:far-param}
{We begin our parametrisation by restricting to the invariant subspace $\{\sigma(r)=\frac{1}{r}\}$, for some fixed $r\geq r_{\infty}$.} By the same methods as outlined in Section \ref{s:core}, we can construct a local `centre-stable' manifold
\begin{align}
    \mathcal{W}^{cs}_{+}{(\varepsilon)|_{\sigma=\frac{1}{r}}} = \left\{\mathbf{u}=\widetilde{\mathbf{u}}\,:\; (a,b) = \sum_{j=1}^4 d^{cs}_{j}\left(\mathbf{d}, \mathbf{c}_{2}{; \varepsilon}\right) \widetilde{\mathbf{V}}_{j}(r), \quad (a_{n}, a_{-n}) = \sum_{j=1}^2 c^{cs}_{j,n}\left(\mathbf{d}, \mathbf{c}_{2}{; \varepsilon}\right) \widetilde{\mathbf{W}}_{j,n}(r)\right\},\nonumber
\end{align}
which contains all {small locally forward-bounded solutions to the full far-field system}. Here $\mathbf{d}:=(d_{1}, d_{2}, d_{3}, d_{4})$, $\mathbf{c}_{2}$ is of the form defined in \eqref{ci:defn}, and $\widetilde{\mathbf{V}}_{j}(r)$ and $\widetilde{\mathbf{W}}_{j,n}(r)$ are solutions to the linear far-field system 
\begin{align}
 \widetilde{\mathbf{V}}_{1} &= \Big(r\textnormal{e}^{\textnormal{i}k r}, \; \textnormal{e}^{\textnormal{i}k r}\Big), \quad
 \widetilde{\mathbf{V}}_{2} = \Big(-\textnormal{i} r\textnormal{e}^{\textnormal{i}k r}, \; -\textnormal{i}\textnormal{e}^{\textnormal{i}k r}\Big),\quad \widetilde{\mathbf{V}}_{3} = \Big(\textnormal{e}^{\textnormal{i}k r} , \; 0\Big),\quad
 \widetilde{\mathbf{V}}_{4} = \Big(-\textnormal{i}\textnormal{e}^{\textnormal{i}k r} , \; 0\Big),\nonumber\\
  \widetilde{\mathbf{W}}_{1,n} &= \Big(\textnormal{e}^{\lambda_{n} r}, \; 0\Big), \quad\,
 \widetilde{\mathbf{W}}_{2,n} = \Big(0, \; \textnormal{e}^{-\lambda_{n} r}\Big),\nonumber
\end{align}
with respective adjoint solutions
\begin{align}
  \widetilde{\mathbf{V}}^{*}_{1} &= \Big(0, \; \textnormal{e}^{\textnormal{i}k r}\Big),\quad
 \widetilde{\mathbf{V}}^{*}_{2} = \Big(0, \;-\textnormal{i}\textnormal{e}^{\textnormal{i}k r} \Big),\quad \widetilde{\mathbf{V}}^{*}_{3} = \Big(\textnormal{e}^{\textnormal{i}k r}, \; -r\textnormal{e}^{\textnormal{i}k r}\Big), \quad
 \widetilde{\mathbf{V}}^{*}_{4} = \Big(-\textnormal{i}\textnormal{e}^{\textnormal{i}k r}, \; \textnormal{i} r\textnormal{e}^{\textnormal{i}k r}\Big),\nonumber\\
 \widetilde{\mathbf{W}}^{*}_{1,n} &= \Big(\textnormal{e}^{-\lambda_{n} r}, \; 0\Big), \quad
 \widetilde{\mathbf{W}}^{*}_{2,n} = \Big(0, \; \textnormal{e}^{\lambda_{n} r}\Big),\nonumber
\end{align}
where $\langle \widetilde{\mathbf{V}}^{*}_{i}, \widetilde{\mathbf{V}}_{j}\rangle = \delta_{i,j}$, $\langle \widetilde{\mathbf{W}}^{*}_{i,n}, \widetilde{\mathbf{W}}_{j,n}\rangle = \delta_{i,j}$. We note that $\{\widetilde{\mathbf{V}}_{3}, \widetilde{\mathbf{V}}_{4}\}$ are bounded as $r\to\infty$, $\{\widetilde{\mathbf{V}}_{1}, \widetilde{\mathbf{V}}_{2}\}$ grow algebraically as $r\to\infty$, $\widetilde{\mathbf{W}}_{2,n}(r)$ decays exponentially as $r\to\infty$, and $\widetilde{\mathbf{W}}_{1,n}(r)$ grows exponentially as $r\to\infty$, for any $n\in\mathbb{N}$. Applying the method of variation of parameters and solving the respective equations similar to \eqref{fixed:a}, \eqref{fixed:a;n} for small enough $(\mathbf{d}, \mathbf{c}_{2})$, the parametrisation of the centre-stable manifold takes a form similar to \eqref{amp:core;param},
\begin{align}
    d^{cs}_{j}(\mathbf{d}, \mathbf{c}_{2}{; \varepsilon}) &= \textstyle d_{j} + \textnormal{O}_{\sigma}\left(\big[|\mathbf{d}| + |\mathbf{c}_{2}|_{1}\big]\big[|\varepsilon| + |\mathbf{d}| + |\mathbf{c}_{2}|_{1}\big]\right),\nonumber\\
    c^{cs}_{1,n}(\mathbf{d}, \mathbf{c}_{2}{; \varepsilon}) &=\textstyle \textnormal{O}_{\sigma}\left(\big[|\mathbf{d}| + |\mathbf{c}_{2}|_{1}\big]\big[|\varepsilon| + |\mathbf{d}| + |\mathbf{c}_{2}|_{1}\big]\right), \nonumber\\
    c^{cs}_{2,n}(\mathbf{d}, \mathbf{c}_{2}{; \varepsilon}) &=\textstyle c_{2,n} + \textnormal{O}_{\sigma}\left(\big[|\mathbf{d}| + |\mathbf{c}_{2}|_{1}\big]\big[|\varepsilon| + |\mathbf{d}| + |\mathbf{c}_{2}|_{1}\big]\right), \nonumber
\end{align}
where $\textnormal{O}_{\sigma}(\cdot)$ denotes the standard Landau symbol where the bounding constants may depend on the value of $\sigma$. We introduce the centre manifold $\mathcal{W}^{c}_{+}{(\varepsilon)_{\sigma=\frac{1}{r}}}$, which contains all {small locally bounded solutions,}
\begin{align}
    \mathcal{W}^{c}_{+}{(\varepsilon)|_{\sigma=\frac{1}{r}}} = \left\{\mathbf{u}=\widetilde{\mathbf{u}}\,:\; (a,b) = \sum_{j=1}^4 d^{c}_{j}\left(\mathbf{d}{; \varepsilon}\right) \widetilde{\mathbf{V}}_{j}(r), \quad (a_{n}, a_{-n}) = \sum_{j=1}^2 c^{c}_{j,n}\left(\mathbf{d}{; \varepsilon}\right) \widetilde{\mathbf{W}}_{j,n}(r)\right\}.\nonumber
\end{align}
We note that any solution that remains bounded for all $r$, i.e. is on the centre manifold $\mathcal{W}^{c}_{+}{(\varepsilon)|_{\sigma=\frac{1}{r}}}$, must remain bounded as $r\to\infty$, by definition. Hence, the set of all solutions bounded for all $r$ is contained within the set of all solutions bounded as $r\to\infty$, and so $\mathcal{W}^{c}_{+}{(\varepsilon)|_{\sigma=\frac{1}{r}}}\subset\mathcal{W}^{cs}_{+}{(\varepsilon)|_{\sigma=\frac{1}{r}}}$ for any {$r\gg1$}. Thus, the centre and centre-stable parametrisations should coincide when $c_{2,n}$ is slaved to the centre coordinates $\mathbf{d}$ for all $n\in\mathbb{N}$; this is a key assumption which we do not prove in this work. That is, when $c_{2,n} = \textnormal{O}_{\sigma}\left(|\mathbf{d}|\big[|\varepsilon| + |\mathbf{d}|\big]\right)$, we find
\begin{align}
    c^{c}_{2,n}(\mathbf{d}{; \varepsilon}) &=\textstyle \textnormal{O}_{\sigma}\left(|\mathbf{d}|\big[|\varepsilon| + |\mathbf{d}|\big]\right), \nonumber\\
    c^{c}_{1,n}(\mathbf{d}{; \varepsilon}) &=c^{cs}_{1,n}(\mathbf{d}, c^{c}_{(2,i)}(\mathbf{d}{; \varepsilon}){; \varepsilon}), \nonumber\\
    \; &=\textstyle \textnormal{O}_{\sigma}\left(|\mathbf{d}|\big[|\varepsilon| + |\mathbf{d}|\big]\right), \label{farf:cent;param}\\
    d^{c}_{j}(\mathbf{d}{; \varepsilon}) &= d^{cs}_{j}(\mathbf{d}, c^{c}_{(2,i)}(\mathbf{d}{; \varepsilon}){; \varepsilon}),  \nonumber\\
    \; &= \textstyle d_{j} + \textnormal{O}_{\sigma}\left(|\mathbf{d}|\big[|\varepsilon| + |\mathbf{d}|\big]\right).\nonumber
\end{align}
\begin{figure}[t!]
    \centering
    \includegraphics[height=7cm]{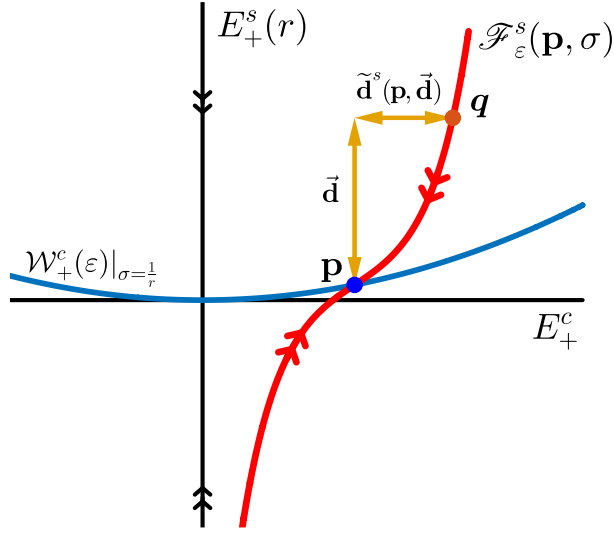}
    \caption{Parametrisation of a strong stable fibre $\mathscr{F}^{s}{_{\varepsilon}}(\mathbf{p},{\sigma})$ with base point $\mathbf{p}\in\mathcal{W}^{c}_{+}{(\varepsilon)}|_{\sigma={\frac{1}{r}}}$ for some {$r\geq r_{\infty}$} fixed. Any solution with initial {data $(\mathbf{q},\sigma)$, where} $\mathbf{q}\in\mathscr{F}^{s}_{{\varepsilon}}(\mathbf{p},{\sigma})$, will tend to solutions with initial {data $(\mathbf{p},\sigma)$, where}  $\mathbf{p}\in\mathcal{W}^{c}_{+}{(\varepsilon)}|_{\sigma={\frac{1}{r}}}$ as $r\to\infty$.}
    \label{fig:foliationparam}
\end{figure}

We now wish to parametrise $\mathcal{W}^{cs}_{+}{(\varepsilon)|_{\sigma=\frac{1}{r}}}$ as the union of strong stable fibres $\bigcup_{\mathbf{p}}\mathscr{F}^{s}{_{\varepsilon}}(\mathbf{p},{\sigma})$ over base points $\mathbf{p}\in\mathcal{W}^{c}_{+}{(\varepsilon)|_{\sigma=\frac{1}{r}}}$. We refer to \cite{chow1990smooth} for background on stable foliations, \cite{mcquighan2014oscillons} for a finite dimensional example, and Figure \ref{fig:foliationparam} for an illustration. We introduce a base point $\mathbf{p}\in\mathcal{W}^{c}_{+}{(\varepsilon)|_{\sigma=\frac{1}{r}}}$ with the decomposition
\begin{align}
    \mathbf{p} :&= A_{0}(r)\;\mathbf{e} + B_{0}(r)\;\mathbf{f} + \overline{A}_{0}(r)\;\overline{\mathbf{e}} + \overline{B}_{0}(r)\;\overline{\mathbf{f}} + \sum_{n=1}^{{\infty}} \left\{ A_{n}(r)\mathbf{e}_{n} + A_{-n}(r)\mathbf{e}_{-n}\right\} \label{p:defn}
    \end{align}
such that, since $\mathbf{p}\in\mathcal{W}^{c}_{+}{(\varepsilon)|_{\sigma=\frac{1}{r}}}$, we can use the parametrisation \eqref{farf:cent;param} to write
    \begin{align}
    A_{\pm n}(r) &= \textnormal{O}_{\sigma}\left(\left[|A_{0}| + |B_{0}|\right]\left[|\varepsilon| + |A_{0}| + |B_{0}|\right]\right),\label{p:defn;n}
\end{align}
for all $n\in\mathbb{N}$. Then, we can write the stable fibre of trajectories related to the point $\mathbf{p}$ as
\begin{align}
    \mathscr{F}^{s}{_{\varepsilon}}(\mathbf{p}, {\sigma}) = \left\{\mathbf{u} = \mathbf{p} + \widetilde{\mathbf{u}}\,:\; (a,b) = \sum_{j=1}^4 \widetilde{d}^{s}_{j}\left(\mathbf{p}, \vec{\mathbf{d}}{; \varepsilon}\right) \widetilde{\mathbf{V}}_{j}(r), \quad (a_{n}, a_{-n}) = \sum_{j=1}^2 \widetilde{c}^{s}_{j,n}\left(\mathbf{p}, \vec{\mathbf{d}}{; \varepsilon}\right) \widetilde{\mathbf{W}}_{j,n}(r)\right\},\nonumber
\end{align}
where ${\vec{\mathbf{d}} := (\,\overline{d}_{1}, \overline{d}_{2}, \dots)}$, and any trajectory $(\mathbf{q}, {\sigma})$ with $\mathbf{q}\in\mathscr{F}^{s}{_{\varepsilon}}(\mathbf{p}, {\sigma})$ approaches the solution associated with initial conditions $(\mathbf{p}, {\sigma})$, as $r\to\infty$. We define $\vec{d}_{n}$ as the distance between a point $\mathbf{q}\in\mathscr{F}^{s}{_{\varepsilon}}(\mathbf{p},{\sigma})$ and $\mathbf{p}\in\mathcal{W}^{c}_{+}{(\varepsilon)|_{\sigma=\frac{1}{r}}}$ in each respective $\widetilde{\mathbf{W}}_{2,n}(r)$ direction, see Figure \ref{fig:foliationparam}. If $\vec{d}_{n} =0$ for all $n\in\mathbb{N}$, then $\mathbf{q} = \mathbf{p}$ and so we write
\begin{align}
    \widetilde{d}^{s}_{j}(\mathbf{p}, \vec{\mathbf{d}}{; \varepsilon}) &= \textstyle\textnormal{O}_{\sigma}\left(|\vec{\mathbf{d}}|_{1}\left[|\varepsilon| + |\mathbf{p}| + |\vec{\mathbf{d}}|_{1}\right]\right), \label{farf:stab;param}\\
    \widetilde{c}^{s}_{1,n}(\mathbf{p}, \vec{\mathbf{d}}{; \varepsilon}) &=\textstyle \textnormal{O}_{\sigma}\left(|\vec{\mathbf{d}}|_{1}\left[|\varepsilon| + |\mathbf{p}| + |\vec{\mathbf{d}}|_{1}\right]\right), \nonumber\\
    \widetilde{c}^{s}_{2,n}(\mathbf{p}, \vec{\mathbf{d}}{; \varepsilon}) &=\textstyle \vec{d}_{n}. \nonumber
\end{align}
Substituting the parametrisation \eqref{p:defn}-\eqref{p:defn;n} into \eqref{farf:stab;param}, we can parametrise the far-field manifold as,
\begin{align}
    a(r) &= A_{0}(r) +\textnormal{O}_{\sigma}\left(|\vec{\mathbf{d}}|_{1}\left[|\varepsilon| + |A_{0}| + |B_{0}| + |\vec{\mathbf{d}}|_{1}\right]\right),\nonumber\\
    b(r) &= B_{0}(r) +\textnormal{O}_{\sigma}\left( |\vec{\mathbf{d}}|_{1}\left[|\varepsilon| + |A_{0}| + |B_{0}| + |\vec{\mathbf{d}}|_{1}\right]\right),\nonumber\\
    a_{n}(r) &= \textnormal{O}_{\sigma}\left(\left[|A_{0}| + |B_{0}| + |\vec{\mathbf{d}}|_{1}\right]\left[|\varepsilon| + |A_{0}| + |B_{0}| + |\vec{\mathbf{d}}|_{1}\right]\right),\label{match:farf:folia}\\
    a_{-n}(r) &= \vec{d}_{n}\textnormal{e}^{-\lambda_{n}r}  + \textnormal{O}_{\sigma}\left(\left[|A_{0}| + |B_{0}|\right]\left[|\varepsilon| + |A_{0}| + |B_{0}| + |\vec{\mathbf{d}}|_{1}\right]\right) ,\nonumber
\end{align}
for the complex amplitudes $A_{0}(r)$, $B_{0}(r)$ defined such that the base point $\mathbf{p}\in\mathcal{W}^{c}_{+}{(\varepsilon)|_{\sigma=\frac{1}{r}}}$ defined in \eqref{p:defn} decays to zero as $r\to\infty$.\\

We have now parametrised the centre-stable manifold $\mathcal{W}^{cs}_{+}{(\varepsilon)|_{\sigma=\frac{1}{r}}}$ as the union of the stable foliations $\bigcup_{\mathbf{p}}\mathscr{F}^{s}{_{\varepsilon}}(\mathbf{p},{\sigma})$ over every point $\mathbf{p}\in\mathcal{W}^{c}_{+}{(\varepsilon)|_{\sigma=\frac{1}{r}}}$ as defined in \eqref{p:defn}. In order to identify the far-field manifold {$\widetilde{\mathcal{W}}^{s}_{+}(\varepsilon):=\mathcal{W}^{s}_{+}(\varepsilon)|_{\sigma=\frac{1}{r}}$} containing all solutions that decay exponentially as $r\to\infty$, we must isolate the values of the centre coordinates $(A_{0},B_{0})$ such that $\mathbf{p}$ exponentially decays to zero as $r\to\infty$. Hence, we reduce the problem to the centre manifold{, and look to find exponentially decaying bounded solutions. Following this, we track our solutions back in $r$ to the matching point $r=r_{0}$, where we obtain the leading-order parametrisation of $\widetilde{\mathcal{W}}^{s}_{+}(\varepsilon)|_{r=r_{0}}$ from \eqref{match:farf:folia}, which we can then match to our core manifold $\widetilde{\mathcal{W}}^{cu}_{-}(\varepsilon)|_{r=r_{0}}$.}
\subsection{Reduction to the centre manifold}\label{subs:far-cent}
In order to identify solutions that are exponentially decaying as $r\to\infty$, we need to analyse the dynamics of the centre manifold $\mathcal{W}^{c}_{+}{(\varepsilon)}$. To do this, we first derive a normal-form for the extended  system
\begin{align}
    \mathbf{u}_{r} &= \mathbf{L}\left(\frac{1}{\sigma}\right)\mathbf{u} + \mathcal{F}(\mathbf{u},\varepsilon,\sigma), \nonumber\\
    \sigma_{r} &= -\sigma^2,\nonumber
\end{align}
restricted to the centre manifold. We need that the hyperbolic part of the system is slaved to the centre modes; a key assumption not proven here. We recall the centre manifold parametrisation introduced in \eqref{p:defn}-\eqref{p:defn;n}, and define a centre manifold reduction
\begin{align}
    \mathbf{u}^{c}:&= A_{0}\mathbf{e} + B_{0}\mathbf{f} + \overline{A}_{0}\overline{\mathbf{e}} + \overline{B}_{0}\overline{\mathbf{f}},\qquad 
    \textnormal{and} \quad
    Q(A_{0},B_{0},\sigma;\varepsilon) = \textnormal{O}\left(\big[ |A_{0}| + |B_{0}|\big]\big[|\varepsilon| + |A_{0}| + |B_{0}|\big]\right).\nonumber
\end{align}
 Then, the reduced vector field on $\mathcal{W}^{c}_{+}{(\varepsilon)}$ projected onto $\mathcal{P}_{0}$ is 
\begin{align}
    \frac{\textnormal{d}}{\textnormal{d} r}\mathbf{u}^{c} &= \widetilde{\mathbf{L}}(\sigma)\mathbf{u}^{c} + \mathcal{P}_{0}\mathcal{F}\left(\mathbf{u}^{c}+Q(A_{0},B_{0},\sigma;\varepsilon), \sigma; \varepsilon\right),\nonumber
\end{align}
where $\mathcal{P}_{0}$ is defined in Section \ref{s:spect}, and $\widetilde{\mathbf{L}}(\sigma)$ is the projection of linear operator $\mathbf{L}\left(\frac{1}{\sigma}\right)$ restricted to the centre manifold $\mathcal{W}^{c}_{+}{(\varepsilon)}$. Then, $A_{0},B_{0},\sigma$ satisfy a complex ODE system of the form
\begin{align}
    \frac{\textnormal{d}}{\textnormal{d}r} A_{0} &= \textnormal{i}k A_{0} + B_{0} - \frac{\sigma}{2}\left(A_{0}-\overline{A}_{0}\right) + \textnormal{O}\left(\left[|\varepsilon| + \big[1 + |\sigma|\big]\big(|A_{0}| + |B_{0}|\big)\right]\big[|A_{0}| + |B_{0}|\big] \right), \nonumber\\
    \frac{\textnormal{d}}{\textnormal{d}r} B_{0} &= \textnormal{i}k B_{0} - \frac{\sigma}{2}\left(B_{0}+\overline{B}_{0}\right) + \textnormal{O}\left(\left[|\varepsilon| + \big[1 + |\sigma|\big]\big(|A_{0}| + |B_{0}|\big)\right]\big[|A_{0}| + |B_{0}|\big] \right), \label{amp:a0;b0}\\
    \frac{\textnormal{d}}{\textnormal{d}r} \sigma &= -\sigma^2. \nonumber
\end{align}
In order to transform \eqref{amp:a0;b0} into a normal-form that we can analyse, we utilise the following result,
\begin{Lemma}{\cite{lloyd2009localized}}
There exists a change of coordinates \label{Res:normal;form}
\begin{align}
    \begin{pmatrix} A_{0} \\ B_{0}\end{pmatrix} &\mapsto \begin{pmatrix} A \\ B\end{pmatrix} := \textnormal{e}^{-\textnormal{i}\phi(r)}\left[ 1 + \mathcal{T}(\sigma)\right]\begin{pmatrix} A_{0} \\ B_{0}\end{pmatrix} + \textnormal{O}\left(\left[|\varepsilon| + |A_{0}| + |B_{0}|\right]\big[|A_{0}| + |B_{0}|\big] \right), \label{A0:A;trans}\\
    \intertext{such that \eqref{amp:a0;b0} is transformed into}
    \frac{\textnormal{d}}{\textnormal{d}r} A &= -\frac{\sigma}{2}A + B + \mathcal{R}_{A}(A,B, \sigma, \varepsilon), \nonumber\\
    \frac{\textnormal{d}}{\textnormal{d}r} B &= - \frac{\sigma}{2}B + c_{0}\varepsilon A + c_{3}|A|^2 A + \mathcal{R}_{B}(A,B, \sigma, \varepsilon), \label{amp:A;B}\\
    \frac{\textnormal{d}}{\textnormal{d}r} \sigma &= -\sigma^2, \nonumber
\end{align}
where we define the remainder terms $\mathcal{R}_{A}$, $\mathcal{R}_{B}$ below. This coordinate change is polynomial in $(A, B, \sigma)$ and smooth in $\varepsilon$, and $\mathcal{T}(\sigma) = \textnormal{O}(\sigma)$ is linear and upper triangular for each $\sigma$, while $\phi(r)$ satisfies 
\begin{align}
    \phi_{r} &= k + \textnormal{O}(\varepsilon + |\sigma|^3 + |A|^2), \qquad \phi(0) = 0. \label{phase:defn}
\end{align}
The remainder terms satisfy
\begin{align}
\mathcal{R}_{A}\left(A, B, \sigma, \varepsilon \right) &= \textnormal{O}\left(\sum_{j=0}^{2} |A^{j}B^{3-j}| + (|A|+|B|)^{5} + |\sigma|^2|B| + |\varepsilon||\sigma|^{m}(|A|+|B|)\right),\label{normal:remainder}\\
\mathcal{R}_{B}\left(A, B, \sigma, \varepsilon \right) &= \textnormal{O}\left(\sum_{j=0}^{1} |A^{j}B^{3-j}|+ |\varepsilon|(|\varepsilon| + |\sigma|^3 + |A|^2)|A| + (|A|+|B|)^{5} + |\varepsilon||\sigma|^{m}|B|\right),\nonumber
\end{align}
for fixed $0<m<\infty$, $c_{0}\in\mathbb{R}$, and $c_{3}\in\mathbb{R}$.
\end{Lemma}
This result, with different values of $c_{0},c_{3}$, is proven in \cite[Lemma 2]{lloyd2009localized}, and so we will omit the details. This change of coordinates is called a normal-form transformation, where \eqref{amp:A;B} is called the normal-form of the system. As seen in \cite{lloyd2009localized} and \cite{scheel2003radially}, \eqref{amp:A;B} is equivalent to the normal-form of the radial Swift-Hohenberg equation. 

\begin{Remark} The radial problem is an $\textnormal{O}(r^{-1})$ perturbation from the one-dimensional problem and so, with the added non-autonomous terms, we expect $c_{0} = \widetilde{c}_{0} + \textnormal{O}(r^{-1})$ and $c_{3} = \widetilde{c}_{3} + \textnormal{O}(r^{-1})$, where $\widetilde{c}_{0}$, $\widetilde{c}_{3}$ are the respective coefficients for the one-dimensional problem. Then, extending the system with respect to $\sigma=r^{-1}$, these $\textnormal{O}(r^{-1})$ perturbations of $\widetilde{c}_{0}$, $\widetilde{c}_{3}$ become $\textnormal{O}(\sigma)$ terms and so are lifted to the coefficients of the higher order terms $\sigma\varepsilon A$ and $\sigma |A|^{2} A$, respectively. Therefore, the values of $c_{0}$, $c_{3}$ are identical to those found in \cite[\S5]{groves2017pattern} for the one-dimensional ferrohydrostatic problem subject to a linear magnetisation law. That is, for fixed $ k_{D}\in(0,\infty)$ and $M_{0}\in(0,1)$, where $ k_{D}:=k D$ and $M_{0}:=\frac{\mu-1}{\mu+1}$, we find
\begin{align}
    c_{0} &=  k_{D}^{2}\left[\widetilde{\Upsilon}_{0} +  k_{D}^{3}\mathcal{M}\tanh( k_{D})\,\textnormal{sech}^{2}( k_{D})\right]^{-1},\label{c0}
\end{align}
and 
\begin{align}
    c_{3} &= -\frac{c_{0} k_{D}^{3}}{2D^{4}}\left[\frac{ k_{D} \mathcal{M}^{2} M_{0}^{2}}{4}\left(\frac{(\cosh(4 k_{D}) - 4\cosh(2 k_{D}) - 3)(4 \,\textnormal{sech}(4 k_{D}) + \,\textnormal{sech}^{2}( k_{D}) - 2)}{(2 k_{D} \mathcal{M}\tanh(2 k_{D}) - \widetilde{\Upsilon}_{0} - 4 k_{D}^{2}) \cosh(2 k_{D})} + \frac{\,\textnormal{sech}^{2}( k_{D})}{\widetilde{\Upsilon}_{0}}\right)\,\textnormal{sech}^{2}( k_{D}) \right. \nonumber\\ & \qquad\qquad\qquad \left. + \frac{3  k_{D}}{2} + 4\mathcal{M}\left(M_{0}^{2} \,\textnormal{sech}(2 k_{D}) - \cosh^{2}( k_{D})\right)\textnormal{cosech}(2 k_{D})\right],\label{c3}
\end{align}
where $(\mathcal{M}, \widetilde{\Upsilon}_{0})=(\mathcal{M}_{H}, \widetilde{\Upsilon}_{H})( k_{D})$, as defined in \eqref{M:Upsilon}.
\label{rem:c0c3}\end{Remark}
It is straightforward to show that $c_{0}$ is positive for all $ k_{D}$ and $M_{0}$, one can find a restriction on $ k_{D}$ for a given $M_{0}$ such that $c_{3}<0\,$; this is required for the existence of homoclinic solutions in \cite{groves2017pattern} and will be discussed in greater detail in Section \ref{subs:orbit}.

It is convenient to extend the transformation in Lemma \ref{Res:normal;form} and apply it to our matching coordinates $(a,b)$
\begin{align}
    \begin{pmatrix} a \\ b\end{pmatrix} \mapsto \begin{pmatrix} \widetilde{a} \\ \widetilde{b}\end{pmatrix} :&= \textnormal{e}^{-\textnormal{i}\phi(r)}\left[ 1 + \mathcal{T}(\sigma)\right]\begin{pmatrix} a \\ b\end{pmatrix} + \textnormal{O}\left(\left[|\varepsilon| + |a| + |b|\right]\big[|a| + |b|\big] \right), \label{ab:tilde}\\
    \intertext{such that, in the far-field}
    \widetilde{a}(r) &= A(r) +\textnormal{O}_{\sigma}\left(|\vec{\mathbf{d}}|_{1}\left[|\varepsilon| + |A| + |B| + |\vec{\mathbf{d}}|_{1}\right]\right),\nonumber\\
    \widetilde{b}(r) &= B(r) +\textnormal{O}_{\sigma}\left(|\vec{\mathbf{d}}|_{1}\left[|\varepsilon| + |A| + |B| + |\vec{\mathbf{d}}|_{1}\right]\right).\nonumber
\end{align}
Using these new coordinates, we can write our core parametrisation \eqref{match:core} as
\begin{align}
    \widetilde{a}(r_{0}) &= \textnormal{e}^{-\textnormal{i}\left(\frac{\pi}{4} + \textnormal{O}(r_{0}^{-2}) + \textnormal{O}_{r_{0}}(|\varepsilon| + |\mathbf{d}|_{1})\right)}r_{0}^{-\frac{1}{2}}\left([-\textnormal{i}+\textnormal{O}(r_{0}^{-1})]r_{0} d_{1} + [1+\textnormal{O}(r_{0}^{-1})]d_{2}\right) + \textnormal{O}_{r_{0}}\left(|\mathbf{d}|_{1}\left[|\varepsilon| + |\mathbf{d}|_{1}\right]\right),\nonumber\\
    \widetilde{b}(r_{0}) &= \textnormal{e}^{-\textnormal{i}\left(\frac{\pi}{4} + \textnormal{O}(r_{0}^{-2}) + \textnormal{O}_{r_{0}}(|\varepsilon| + |\mathbf{d}|_{1})\right)}r_{0}^{-\frac{1}{2}}\left(\left[-\textnormal{i} +\textnormal{O}(r_{0}^{-1})\right]d_{1} - \left[\nu + \textnormal{O}(r_{0}^{-\frac{1}{2}})\right] d_{2}^{2}\right)+ \textnormal{O}_{r_{0}}\left(|\mathbf{d}|_{1}\left[|\varepsilon| + |\mathbf{d}_{2}|_{1}\right] + |d_{2}|^{3}\right),\nonumber\\
    a_{n}(r_{0}) &= \textnormal{e}^{\lambda_{n}r_{0}}r_{0}^{-\frac{1}{2}}\left[\frac{1}{\sqrt{\pi}} + \textnormal{O}(r_{0}^{-1})\right]\widetilde{c}_{1,n} + \textnormal{O}_{r_{0}}\left(|\mathbf{d}|_{1}\left[|\varepsilon| + |\mathbf{d}|_{1}\right]\right),\label{match:core;normal}\\
    a_{-n}(r_{0}) &= \textnormal{O}_{r_{0}}\left(|\mathbf{d}|_{1}\left[|\varepsilon| + |\mathbf{d}|_{1}\right]\right),\nonumber
\end{align}
where we have used 
\begin{align}
    \phi(r_{0}) = kr_{0} + \textnormal{O}(r_{0}^{-2}) + \textnormal{O}_{r_{0}}(|\varepsilon| + |\mathbf{d}|_{1}),\nonumber
\end{align}
from \eqref{phase:defn}. We have also exploited the fact that $\mathcal{T}(\sigma)$ is upper triangular and the coefficient in front of $d_{2}^{2}$ for $\widetilde{b}(r_{0})$ scales with $\sqrt{r_{0}}$ so that it is only affected by the transformation \eqref{ab:tilde} at higher order.
\section{Geometric Blow-up}\label{s:geoblowup}

So far, we have constructed the core manifold in Section \ref{s:core} on bounded intervals $r\in[0, r_{0}]$, and the far-field stable manifold for $r\geq r_{\infty}$. We choose $r_{0}$ large enough such that $r_{\infty}< r_{0}$, ensuring that $\widetilde{\mathcal{W}}^{cu}_{-}{(\varepsilon)}$ and ${\widetilde{\mathcal{W}}^{s}_{+}(\varepsilon)}$ both exist at the matching point $r=r_{0}$. {We have also taken the hyperbolic modes to be parametrised by the coordinates on the centre manifold $\mathcal{W}^{c}_{+}(\varepsilon)$. As a result, we have restricted our analysis to the centre-coordinates $A(r),B(r),\sigma(r)$, as defined in \eqref{A0:A;trans},\eqref{amp:A;B}.} 
{
In order to characterise the behaviour of exponentially decaying solutions on the far-field centre manifold, 
we first present some formal analysis as motivation for the upcoming procedures.
\paragraph{{Motivation}}
\begin{figure}[t!]
    \centering
    \includegraphics[height=5cm]{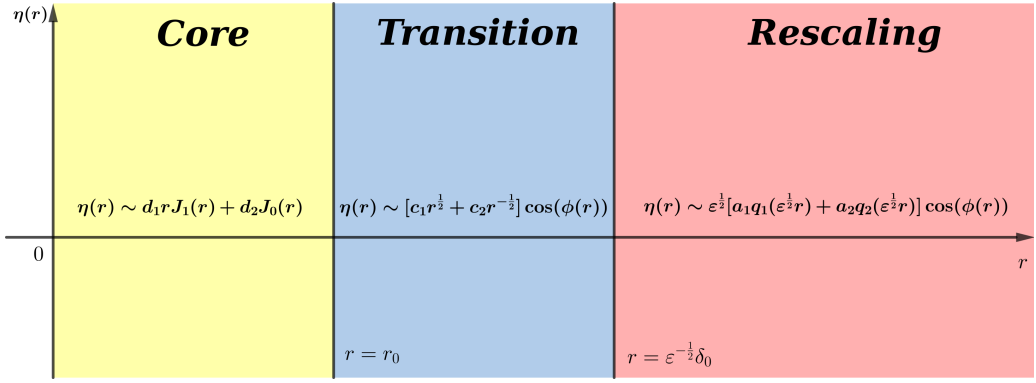}
    \caption{{The radial coordinate $r$ is formally divided into three regions: the core region $r \in[0, r_{0}]$, where solutions look like a combination of Bessel functions; the transition region $r\in[r_{0},\varepsilon^{-\frac{1}{2}}\delta_{0}]$, where solutions have an amplitude that behaves like $r^{\frac{1}{2}}$ or $r^{-\frac{1}{2}}$; and the rescaling region $r\in[\varepsilon^{-\frac{1}{2}}\delta_{0},\infty)$, where solutions have an amplitude defined by solutions of \eqref{G-L:rescal;formal}. Here, $d_{1}$, $d_{2}$, $c_{1}$, $c_{2}$, $a_{1}$, and $a_{2}$ are arbitrary real constants and $\phi(r)$ is some $r$-dependent phase. }}
    \label{fig:formalbehav}
\end{figure}
We begin by formally setting the remainder terms in the normal form \eqref{amp:A;B} to zero and evaluating at $\sigma=\frac{1}{r}$; then, we find that $A(r)$ satisfies,
\begin{align}
 \left[\left(\frac{\textnormal{d}}{\textnormal{d}r} + \frac{1}{2r}\right)^{2} - c_{0}\varepsilon - c_{3}|A|^{2}\right]A = 0.\label{G-L:formal}
\end{align}
We assume that any solution emerging for $0<\varepsilon\ll1$ must have $A=\textnormal{O}(\varepsilon^{\frac{1}{2}})$, such that $c_{0}\varepsilon, c_{3}|A|^{2}$ are of the same order. {In order to balance the differential operator and the $c_{0}\varepsilon, c_{3}|A|^{2}$ terms, we would require $r=\textnormal{O}\left(\varepsilon^{-\frac{1}{2}}\right)$. Then, for $r_{0}\leq r\leq \delta_{0} \varepsilon^{-\frac{1}{2}}$, where $0<\delta_{0}\ll1$ is fixed such that $r$ is "not too large"}, \eqref{amp:A;B} can be approximated as,
\begin{align}
 \begin{array}{c}
     \left(\frac{\textnormal{d}}{\textnormal{d}r} + \frac{1}{2r}\right)A = B,   \\
       \left(\frac{\textnormal{d}}{\textnormal{d}r} + \frac{1}{2r}\right)B = 0,
 \end{array} \qquad \implies \qquad  \begin{array}{c}
        A(r) \sim r^{-\frac{1}{2}}, \\
       B(r) \sim 0,
 \end{array} \quad \textnormal{or} \quad  \begin{array}{c}
        A(r) \sim r^{\frac{1}{2}}, \\
       B(r) \sim r^{-\frac{1}{2}}.
 \end{array} \label{transition:formal}
\end{align}
Similarly, for $\delta_{0}\varepsilon^{-\frac{1}{2}}\leq r$, we can introduce  $A_{2}=\varepsilon^{-\frac{1}{2}}A$ and $s=\varepsilon^{\frac{1}{2}}r$ such that \eqref{G-L:formal} becomes
\begin{align}
 \left[\left(\frac{\textnormal{d}}{\textnormal{d}s} + \frac{1}{2s}\right)^{2} - c_{0} - c_{3}|A_{2}|^{2}\right]A_{2} = 0,\label{G-L:rescal;formal}
\end{align}
{which can be found to have solutions $A_{2}(s)=q_{1}(s)$, $A_{2}(s)=q_{2}(s)$; these are discussed further in \S\ref{subs:connect;P+}.} Therefore, we have two regions in the far-field with differing spatial dynamics: the transition region, where $r_{0}\leq r\leq \delta_{0}\varepsilon^{-\frac{1}{2}}$, in which solutions behave like $r^{\frac{1}{2}}$ or $r^{-\frac{1}{2}}$; and the rescaling region, where $r\geq \delta_{0}\varepsilon^{-\frac{1}{2}}$, in which solutions behave like functions $q_{1}(s)$, $q_{2}(s)$ that solve the nonlinear equation \eqref{G-L:rescal;formal}. This partitioning of the radial coordinate $r$ is illustrated in Figure \ref{fig:formalbehav}, along with the leading order behaviour of solutions. 
\begin{figure}[t!]
    \centering
    \includegraphics[height=7cm]{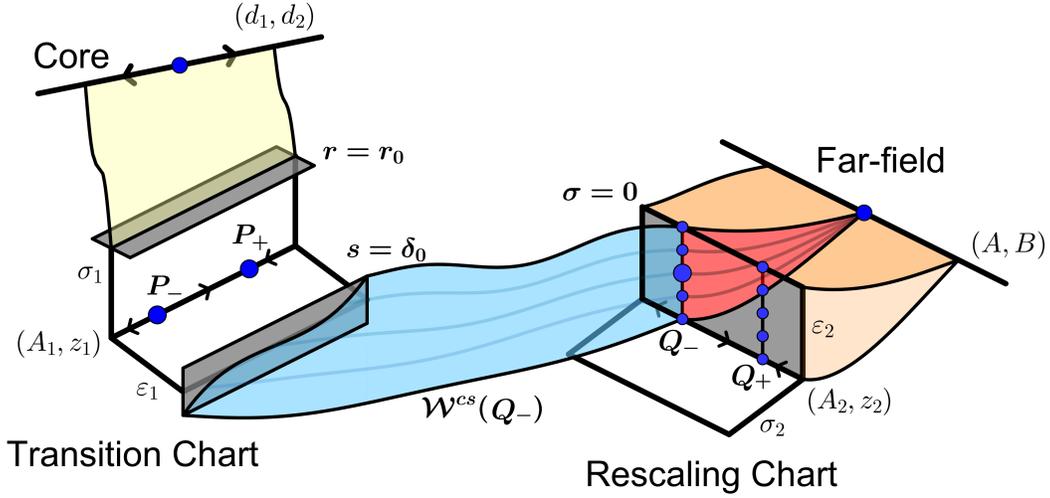}
    \caption{A schematic overview of the various coordinates involved in the matching process. The far-field centre coordinates $(A,B)$ are transformed into $(A_{2}, z_{2}, \sigma_{2}, \varepsilon_{2})$, with equilibria $ Q_{-}(\varepsilon_{2})$, $ Q_{+}(\varepsilon_{2})$. Then, the manifold $\mathcal{W}^{cs}( Q_{-}{[\varepsilon_{2}]})$ is tracked back from $ Q_{-}{(\varepsilon_{2})}$ to the section $s=\delta_{0}$, where another coordinate transformation is applied. The $(A_{1}, z_{1}, \sigma_{1}, \varepsilon_{1})$ coordinates then travel through the transition chart to the section $r=r_{0}$, where solutions are matched with centre coordinates $(d_{1},d_{2})$ on the core manifold $\widetilde{\mathcal{W}}^{cu}_{-}\change{(\varepsilon)}$.}
    \label{fig:GBUManifolds}
\end{figure}
 We will introduce the same coordinate charts used by McCalla \& Sandstede \cite{mccalla2013spots} for the two-dimensional Swift-Hohenberg equation, illustrated in Figure \ref{fig:GBUManifolds}, since the normal-form \eqref{amp:A;B} is equivalent up to leading order. }
\subsection{Rescaling chart}\label{subs:rescale}
Following the procedure implemented in \cite{mccalla2013spots}, we introduce the variable $z\thicksim \left(\frac{\textnormal{d}}{\textnormal{d}r}A\right)/A$. From \eqref{amp:A;B}, we define
\begin{align}
    z:&= -\frac{\sigma}{2} + \frac{B}{A}\in\mathbb{C}, \label{z:def}\\
    \intertext{and introduce rescaling coordinates}
    A_{2} = \frac{A}{ \varepsilon^{\frac{1}{2}}}, \qquad z_{2} = \frac{z}{ \varepsilon^{\frac{1}{2}}}, \qquad \varepsilon_{2} &=  \varepsilon^{\frac{1}{2}}, \qquad \sigma_{2} = \frac{\sigma}{ \varepsilon^{\frac{1}{2}}}, \qquad s =  \varepsilon^{\frac{1}{2}} r.\label{scale:2}
\end{align}
In these new coordinates \eqref{amp:A;B} becomes
\begin{align}
    &\frac{\textnormal{d}}{\textnormal{d} s}A_{2} = A_{2} z_{2} + \frac{1}{\varepsilon_{2}^2}\mathcal{R}_{A}, & \qquad 
    &\frac{\textnormal{d}}{\textnormal{d} s}\sigma_{2} = -\sigma_{2}^{2},
    &\nonumber\\
    &\frac{\textnormal{d}}{\textnormal{d} s}z_{2} = c_{0} + \frac{\sigma_{2}^{2}}{4} + c_{3}|A_{2}|^2 - \sigma_{2}z_{2}- z_{2}^{2} - \frac{z_{2}+ \sigma_{2}/2}{\varepsilon_{2}^2 A_{2}}\mathcal{R}_{A} + \frac{1}{\varepsilon_{2}^{3}A_{2}}\mathcal{R}_{B},
    &\qquad
    &\frac{\textnormal{d}}{\textnormal{d} s}\varepsilon_{2} = 0,
    &\nonumber
\end{align}
where it can be found from \eqref{normal:remainder} that 
\begin{align}
\frac{1}{\varepsilon_{2}^{2}}\mathcal{R}_{A}\left(\varepsilon_{2}A_{2}, \varepsilon_{2}^{2}A_{2}(z_{2} + \sigma_{2}/2), \varepsilon_{2}\sigma_{2}, {\varepsilon_{2}^{2}} \right) &= |A_{2}|\textnormal{O}\left(|\varepsilon_{2}|^2\right),\nonumber\\
\frac{1}{\varepsilon_{2}^{3} A_{2}}\mathcal{R}_{B}\left(\varepsilon_{2}A_{2}, \varepsilon_{2}^{2}A_{2}(z_{2} + \sigma_{2}/2), \varepsilon_{2}\sigma_{2}, {\varepsilon_{2}^{2}} \right) &= \textnormal{O}\left(|\varepsilon_{2}|^2\right).\nonumber
\end{align}
Hence, we obtain the system
\begin{align}
    &\frac{\textnormal{d}}{\textnormal{d} s}A_{2} = A_{2} \left( z_{2} + \textnormal{O} (|\varepsilon_{2}|^{2})\right),
    & \qquad 
    &\frac{\textnormal{d}}{\textnormal{d} s}\sigma_{2} = -\sigma_{2}^{2},\nonumber\\
    &\frac{\textnormal{d}}{\textnormal{d} s}z_{2} = c_{0} + \frac{\sigma_{2}^{2}}{4} + c_{3}|A_{2}|^2 - \sigma_{2}z_{2}- z_{2}^{2}  + \textnormal{O}(|\varepsilon_{2}|^{2}),
    &\qquad
    &\frac{\textnormal{d}}{\textnormal{d} s}\varepsilon_{2} = 0.\label{amp:a2;z2,eps}
\end{align}
where $(A_{2}, z_{2}, \sigma_{2}, \varepsilon_{2})\in\mathbb{C}\times\mathbb{C}\times\mathbb{R}\times\mathbb{R}$. Note that the last equation implies that the variable $\varepsilon_{2}$ can be arbitrarily fixed; setting $\varepsilon_{2}=0$, we arrive at the system
\begin{align}
    &\frac{\textnormal{d}}{\textnormal{d} s}A_{2} = A_{2} z_{2},
    &\qquad
    &\frac{\textnormal{d}}{\textnormal{d} s}\sigma_{2} = -\sigma_{2}^{2},
    &\nonumber\\
    &\frac{\textnormal{d}}{\textnormal{d} s}z_{2} = c_{0} + \frac{\sigma_{2}^{2}}{4} + c_{3}|A_{2}|^2 - \sigma_{2}z_{2}- z_{2}^{2},
    &\qquad
    &\frac{\textnormal{d}}{\textnormal{d} s}\varepsilon_{2} = 0,
    &\label{amp:a2;z2,0}
\end{align}
which, for a small fixed $\varepsilon_{2}$, has two families of equilibria; $ Q_{\pm}(\varepsilon_{2}) = \left(0, \pm\sqrt{c_{0}} { + \textnormal{O}\left(|\varepsilon_{2}|^{2}\right)}, 0, {\textnormal{O}(|\varepsilon_{2}|})\right)$. {Inverting \eqref{scale:2} and recalling $ \frac{\textnormal{d}}{\textnormal{d}r}A\sim z A$, the equilibria $Q_{\pm}(\varepsilon_{2})$ can be thought of as solutions of the form $A(r) = \varepsilon^{\frac{1}{2}}a\left(\varepsilon^{\frac{1}{2}}r\right)\textnormal{e}^{\pm\sqrt{c_{0}}\varepsilon^{\frac{1}{2}}r} + \textnormal{O}\left(\varepsilon\right)$, in the limit as $|a|\to0$.} {The linearisation of \eqref{amp:a2;z2,eps} has {four} eigenvalues $\big\{\pm\sqrt{c_{0}} { + \textnormal{O}\left(|\varepsilon_{2}|^{2}\right)}, \mp2\sqrt{c_{0}} { + \textnormal{O}\left(|\varepsilon_{2}|^{2}\right)}, 0, 0\big\}$ close to $Q_{\pm}(\varepsilon_{2})$, with related eigenvectors $\{\mathbf{v}_{1},\mathbf{v}_{2},\mathbf{v}_{3}(\varepsilon_{2}),\mathbf{v}_{4}\}$, where
\begin{align}
\mathbf{v}_{1} = \left(1,0,0,0\right)^{\intercal}, \qquad \mathbf{v}_{2} = \left(0,1,0,0\right)^{\intercal}, \qquad \mathbf{v}_{3}(\varepsilon_{2}) = \left(0,-\frac{1}{2}+\textnormal{O}(|\varepsilon_{2}|^{2}),1,0\right)^{\intercal}, \qquad \mathbf{v}_{4} = \left(0,0,0,1\right)^{\intercal}. \nonumber
\end{align}
We introduce the centre-stable manifold $\mathcal{W}^{cs}(Q_{\pm}[\varepsilon_{2}])$ which contains all small solutions $(A_{2},z_{2},\sigma_{2},\varepsilon_{2})(s)$ such that $(A_{2},z_{2},\sigma_{2},\varepsilon_{2})(s)\to Q_{\pm}(\varepsilon_{2})$ as $s\to\infty$. We note that $\varepsilon_{2}(s) = \varepsilon_{2} = \varepsilon^{\frac{1}{2}}$ is constant for all $s$, and so the $\mathbf{v}_{4}$ direction is neutral and not contained in any leading order expansion of $\mathcal{W}^{cs}(Q_{\pm}[\varepsilon_{2}])$. For sufficiently small $\varepsilon_{2}$, the equilibria $Q_{\pm}(\varepsilon_{2})$ and invariant manifolds $\mathcal{W}^{cs}(Q_{\pm}[\varepsilon_{2}])$ persist and depend smoothly on $\varepsilon_{2}$.} 
For the subspace $A_{2},z_{2}\in\mathbb{R}$, \eqref{amp:a2;z2,0} can be reduced to the non-autonomous real Ginzburg-Landau equation,
\begin{align}
    \frac{\textnormal{d}^2}{\textnormal{d}s^2}A_{2} &= -\frac{1}{s}\frac{\textnormal{d}}{\textnormal{d}s}A_{2} + \frac{1}{4s^2}A_{2} + c_{0}A_{2} + c_{3}A_{2}^{3}, \qquad A_{2}\in\mathbb{R}.\label{Ginzburg-Landau,c3}
\end{align}
By analysing \eqref{Ginzburg-Landau,c3} linearised about the trivial state for large values of $s$, one can see that solutions to this equation decay or grow exponentially with rate $\pm\sqrt{c_{0}}$ as $s\to\infty$. We see from \eqref{amp:a2;z2,0} that $z_{2} = A'_{2}/A_{2}$ and so, if $A_{2}\thicksim \exp(\pm\sqrt{c_{0}}s)$, then $z_{2}=\pm\sqrt{c_{0}}$. Therefore, to find solutions that are exponentially decaying, we are interested in the centre-stable manifold $\mathcal{W}^{cs}( Q_{-})$ of the family of equililbria $ Q_{-}(\varepsilon_{2})$. Then, we will capture all solutions of \eqref{amp:A;B} of size $ \varepsilon^{\frac{1}{2}}$ that decay as $r\to\infty$, for $\varepsilon>0$.
\subsection{Transition chart}
We rewrite \eqref{amp:A;B} in terms of coordinates $(A_{1}, z_{1})$ that are obtained by rescaling $(A,z)$ by $\sigma=1/r$, where $z$ is defined in \eqref{z:def}. In particular, we define
\begin{align}
    A_{1} &= \frac{A}{\sigma}, \qquad z_{1} = \frac{z}{\sigma} = -\frac{1}{2} + \frac{B}{\sigma A}, \qquad \varepsilon_{1} = \frac{ \varepsilon^{\frac{1}{2}}}{\sigma}, \qquad \sigma_{1} = \sigma, \qquad \rho = {\log}\;(r). \label{scale:1}
\end{align}
Then, we can write \eqref{amp:A;B} as 
\begin{align}
    &\frac{\textnormal{d}}{\textnormal{d} \rho}A_{1} = A_{1}[1+z_{1}] + \frac{1}{\sigma_{1}^{2}}\mathcal{R}_{A},
    &\qquad
    &\frac{\textnormal{d}}{\textnormal{d} \rho}\sigma_{1} = -\sigma_{1},
    &\nonumber\\
    &\frac{\textnormal{d}}{\textnormal{d} \rho}z_{1} = c_{0}\varepsilon^{2}_{1} + \frac{1}{4} + c_{3}|A_{1}|^2 - z_{1}^{2} - \frac{z_{1} + 1/2}{\sigma_{1}^{2} A_{1}}\mathcal{R}_{A} + \frac{1}{\sigma_{1}^{3} A_{1}}\mathcal{R}_{B},
    &\qquad
    &\frac{\textnormal{d}}{\textnormal{d} \rho}\varepsilon_{1} =\varepsilon_{1},
    &\nonumber
\end{align}
where $\mathcal{R}_{A}$, $\mathcal{R}_{B}$ are evaluated at $(A, B, \sigma, {\varepsilon}) = \left(A_{1}, \sigma_{1}^{2}A_{1}\left(z_{1}+\frac{1}{2}\right), \sigma_{1}, {\sigma_{1}^{2}\varepsilon_{1}^{2}}\right)$. From \eqref{normal:remainder}, one can find that 
\begin{align}
\frac{1}{\sigma_{1}^{2}A_{1}}\mathcal{R}_{A}\left(A_{1}, \sigma_{1}^{2}A_{1}(z_{1}+1/2), \sigma_{1}, {\sigma_{1}^{2}\varepsilon_{1}^{2}}\right) &= \textnormal{O}\left(|\sigma_{1}|^2\right),\nonumber\\
\frac{1}{\sigma_{1}^{3} A_{1}}\mathcal{R}_{B}\left(A_{1}, \sigma_{1}^{2}A_{1}(z_{1}+1/2), \sigma_{1}, {\sigma_{1}^{2}\varepsilon_{1}^{2}}\right) &= |\sigma_{1}|^2\textnormal{O}\left(|A_{1}|^{4} + |z_{1} + 1/2| + |\varepsilon_{1}|^{2}\right),\nonumber
\end{align}
as shown in \cite[\S2.4]{mccalla2013spots}. Thus, the flow in the transition chart is defined by
\begin{align}
    &\frac{\textnormal{d}}{\textnormal{d} \rho}A_{1} = A_{1}\left[1+z_{1} + \textnormal{O}(|\sigma_{1}|^{2})\right],
    &\qquad
    &\frac{\textnormal{d}}{\textnormal{d} \rho}\sigma_{1} = -\sigma_{1},
    &\nonumber\\
    &\frac{\textnormal{d}}{\textnormal{d} \rho}z_{1} = c_{0}\varepsilon_{1}^{2} + \frac{1}{4} + c_{3}|A_{1}|^2 - z_{1}^{2} + |\sigma_{1}|^{2}\textnormal{O}\left(|A_{1}|^{4} + \left|z_{1} + \frac{1}{2}\right| + |\varepsilon_{1}|^{2}\right),
    &\qquad
    &\frac{\textnormal{d}}{\textnormal{d} \rho}\varepsilon_{1} =\varepsilon_{1},    
    &\label{amp:a1;b1}
\end{align}
where $(A_{1}, z_{1}, \sigma_{1}, \varepsilon_{1})\in \mathbb{C}\times\mathbb{C}\times\mathbb{R}\times\mathbb{R}$. For $c_{3}<0$, \eqref{amp:a1;b1} has precisely two equilibria $P_{\pm}:= \left(0, \pm \frac{1}{2}, 0, 0\right)$ {which, by inverting \eqref{scale:1} and recalling $\frac{\textnormal{d}}{\textnormal{d}r}A \sim z A$, correspond to solutions where $A(r)\sim r^{\pm\frac{1}{2}}$, respectively, as seen in \eqref{transition:formal}}. We require a good understanding of the dynamics close to each equilibria, and so we introduce two sets of coordinates that move each respective equilibria to the origin. First, near $P_{+}$, we use the new variables $(A_{+}, z_{+}, \sigma_{+}, \varepsilon_{+}) = (A_{1}, z_{1}-\frac{1}{2}, \sigma_{1}, \varepsilon_{1})$ to get
\begin{align}
    &\frac{\textnormal{d}}{\textnormal{d} \rho}A_{+} = A_{+}\left[\frac{3}{2} +z_{+} + \textnormal{O}(|\sigma_{+}|^{2})\right],
    &\qquad
    &\frac{\textnormal{d}}{\textnormal{d} \rho}\sigma_{+} = -\sigma_{+},
    &\nonumber\\
    &\frac{\textnormal{d}}{\textnormal{d} \rho}z_{+} = c_{0}\varepsilon_{+}^{2} - z_{+} - z_{+}^{2} + c_{3}|A_{+}|^2 + \textnormal{O}\left(|\sigma_{+}|^{2}\right),
    &\qquad
    &\frac{\textnormal{d}}{\textnormal{d} \rho}\varepsilon_{+} =\varepsilon_{+}.
    &\label{amp:a1;b1,+}
\end{align}
The linearisation of \eqref{amp:a1;b1,+} about the origin has the eigenvalues $\{\frac{3}{2}, -1, -1, 1\}$ { with related eigenvectors $\{\mathbf{w}_{1},\mathbf{w}_{2},\mathbf{w}_{3},\mathbf{w}_{4}\}$, where
\begin{align}
    \mathbf{w}_{1} = \left(1,0,0,0\right)^{\intercal}, \qquad  \mathbf{w}_{2} = \left(0,1,0,0\right)^{\intercal}, \qquad  \mathbf{w}_{3} = \left(0,0,1,0\right)^{\intercal}, \qquad  \mathbf{w}_{4} = \left(0,0,0,1\right)^{\intercal}. \nonumber
\end{align}
We introduce the unstable manifold $\mathcal{W}^{u}(P_{\pm})$, containing all small solutions $(A_{1},z_{1},\sigma_{1},\varepsilon_{1})(\rho)$ such that $(A_{1},z_{1},\sigma_{1},\varepsilon_{1})(\rho)\to P_{\pm}$ as $\rho\to-\infty$. Then, the two-dimensional unstable manifold $\mathcal{W}^{u}(P_{+})$ close to the equilibrium $P_{+}$ is defined to leading order by the $(A_{+},\varepsilon_{+})$-directions.} Near $P_{-}$, we choose the variables
$(A_{-}, z_{-}, \sigma_{-}, \varepsilon_{-}) = (A_{1}, z_{1}+\frac{1}{2}, \sigma_{1}, \varepsilon_{1})$ to get
\begin{align}
    &\frac{\textnormal{d}}{\textnormal{d} \rho}A_{-} = A_{-}\left[\frac{1}{2} +z_{-} + \textnormal{O}(|\sigma_{-}|^{2})\right],
    &\qquad
    &\frac{\textnormal{d}}{\textnormal{d} \rho}\sigma_{-} = -\sigma_{-},
    &\nonumber\\
    &\frac{\textnormal{d}}{\textnormal{d} \rho}z_{-} = c_{0}\varepsilon_{-}^{2} + z_{-} - z_{-}^{2} + c_{3}|A_{-}|^2 + |\sigma_{-}|^{2}\textnormal{O}\left(|A_{-}|^{4} + \left|z_{-}\right| + |\varepsilon_{-}|\right),
    &\qquad 
    &\frac{\textnormal{d}}{\textnormal{d} \rho}\varepsilon_{-} =\varepsilon_{-},
    &\label{amp:a1;b1,-}
\end{align}
where the linearisation about the origin has the eigenvalues $\{\frac{1}{2}, 1, -1, 1\}$ { with eigenvectors $\{\mathbf{w}_{1},\mathbf{w}_{2},\mathbf{w}_{3},\mathbf{w}_{4}\}$, as defined above. Hence, the three-dimensional unstable manifold $\mathcal{W}^{u}(P_{-})$ close to the origin is defined to leading order by the $(A_{-}, z_{-},\varepsilon_{-})$-directions. Thus, we have seen that both equilibria $P_{\pm}$ and unstable manifolds $\mathcal{W}^{u}(P_{\pm})$ are invariant under the restriction $\sigma_{1}=\sigma_{+}=\sigma_{-}=0$.} Finally, we note that the coordinates in the transition and rescaling charts are related by the transformation
{\begin{align}
    A_{1} &= \frac{A_{2}}{\sigma_{2}}, \qquad z_{1} = \frac{z_{2}}{\sigma_{2}}, \qquad \varepsilon_{1} = \frac{1}{\sigma_{2}}, \qquad \sigma_{1} = \varepsilon_{2}\sigma_{2}, \qquad \rho = \log\left(\varepsilon^{-\frac{1}{2}}s\right), \label{scale:rel}\\
    A_{2} &= \frac{A_{1}}{\varepsilon_{1}}, \qquad z_{2} = \frac{z_{1}}{\varepsilon_{1}}, \qquad \sigma_{2} = \frac{1}{\varepsilon_{1}}, \qquad \varepsilon_{2} = \sigma_{1}\varepsilon_{1},  \qquad s = \varepsilon^{\frac{1}{2}}\textnormal{e}^{\rho}, \nonumber
\end{align}}

and we can therefore transform from one chart to the other in the transverse section {$\varepsilon_{1} = \sigma^{-1}_{2} = d$, where $d>0$ is some fixed constant, for example.}
\subsection{Connecting orbits}\label{subs:orbit}
\begin{figure}[t!]
    \centering
    \includegraphics[height=6cm]{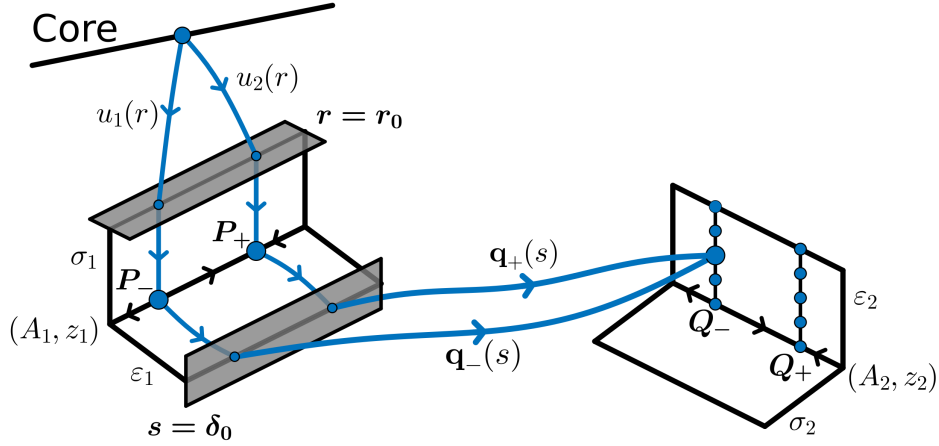}
    \caption{A sketch of the potential connecting orbits between the core and the equilibrium $ Q_{-}(\varepsilon_{2})$ in the transition chart. There exist two connecting orbits: $\mathbf{q}_{\pm}(s)$, connecting $P_{\pm}$ and $ Q_{-}$, defined in \eqref{soln:exp;spotA} and \eqref{soln:chart}, respectively. The connections between the core and $P_{+}$ and $P_{-}$ are $u_{2}(r)$, which corresponds to rings, and $u_{1}(r)$, which corresponds to spots, respectively.}
    \label{fig:GBUTrajectories}
\end{figure}
We now investigate the dynamics of the system in the transition and rescaling charts. Our goal is to show the existence of heteroclinic orbits $\mathbf{q}_{\pm}(s)$ that connect the respective equilibria $P_{\pm}$ in the transition chart to the equilibrium $ Q_{-}(\varepsilon_{2})$ in the rescaling chart in the limit ${\varepsilon_{2}}=0$; see Figure \ref{fig:GBUTrajectories} for a schematic diagram. {For small fixed $\varepsilon_{2}$, the equilibria $Q_{-}(\varepsilon_{2}), P_{\pm}$ and manifolds $\mathcal{W}^{cs}(Q_{-}[\varepsilon_{2}]), \mathcal{W}^{u}(P_{\pm})$ persist and depend smoothly on $\varepsilon_{2}$, and so we take the centre-stable manifold $\mathcal{W}^{cs}(Q_{-}[\varepsilon_{2}])$ to be an $\textnormal{O}(|\varepsilon_{2}|^{2})$-perturbation from the manifold $\mathcal{W}^{cs}(Q_{-}[0])$. Hence, in Section \ref{s:match}, we parametrise our exponentially decaying solutions in the rescaling chart as $\textnormal{O}(|\varepsilon_{2}|^{2})$-perturbations of these heteroclinic orbits $\mathbf{q}_{\pm}(s)$.} First we note that the subspace $\sigma_{1}=\varepsilon_{2}=0$ is invariant under the flow of the equations in the transition and rescaling charts. Setting $\sigma_{1}=0$ in \eqref{amp:a1;b1}, we have
\begin{align}
    &\frac{\textnormal{d}}{\textnormal{d} \rho}A_{1} = A_{1}\left[1+z_{1}\right],
    &\qquad
    &\frac{\textnormal{d}}{\textnormal{d} \rho}z_{1} = c_{0}\varepsilon_{1}^{2} + \frac{1}{4} + c_{3}|A_{1}|^2 - z_{1}^{2},
    &\qquad
    &\frac{\textnormal{d}}{\textnormal{d} \rho}\varepsilon_{1} = \varepsilon_{1},
    &\label{amp:a1;b1,eps0}
\end{align}
with $(A_{1},z_{1}, \varepsilon_{1})\in\mathbb{C}\times\mathbb{C}\times\mathbb{R}$. Setting $\varepsilon_{2}=0$ in \eqref{amp:a2;z2,0}, we obtain
\begin{align}
    &\frac{\textnormal{d}}{\textnormal{d} s}A_{2} = A_{2} z_{2},
    &\qquad
    &\frac{\textnormal{d}}{\textnormal{d} s}z_{2} = c_{0} + \frac{\sigma_{2}^{2}}{4} + c_{3}|A_{2}|^2 - \sigma_{2}z_{2}- z_{2}^{2},
    &\qquad
    &\frac{\textnormal{d}}{\textnormal{d} s}\sigma_{2} = -\sigma_{2}^{2},
    &\label{amp:a2;z2,eps0}
\end{align}
with $(A_{2},z_{2},\sigma_{2})\in\mathbb{C}\times\mathbb{C}\times\mathbb{R}$. We note that the transformation \eqref{scale:rel} between transition and rescaling charts maps $\sigma_{1}=0$ into $\varepsilon_{2}=0$. Also, the real subspaces $(A_{1},z_{1},\varepsilon_{1})\in\mathbb{R}^3$ and $(A_{2},z_{2},\sigma_{2})\in\mathbb{R}^3$ are invariant for \eqref{amp:a1;b1,eps0} and \eqref{amp:a2;z2,eps0}, respectively, and so we will initially restrict our attention to these two subspaces. 
\subsubsection{Connecting $P_{+}$ and $Q_{-}$} \label{subs:connect;P+}
{We start by analysing the non-autonomous real Ginzburg-Landau equation \eqref{Ginzburg-Landau,c3} in the rescaling chart}. We have the following result,
\begin{Lemma}{\cite{mccalla2013spots}}
for $c_{3}<0$, \eqref{Ginzburg-Landau,c3} has a bounded nontrivial solution $A_{2}(s) = q(s)$, and there are constants $q_{0}>0$ and $q_{+}\neq0$ so that
\begin{align}
    q(s) &= \left\{\begin{array}{cc}
         q_{0}s^{1/2} + \textnormal{O}(s^{3/2}),& s\to0,  \\
         (q_{+} + \textnormal{O}(\textnormal{e}^{-\sqrt{c_{0}}s}))\frac{\textnormal{e}^{-\sqrt{c_{0}}s}}{\sqrt{s}},& s\to\infty.
    \end{array}\right. \label{q:defn}
\end{align}
In addition, the linearisation of \eqref{Ginzburg-Landau,c3} about $q(s)$ does not have a nontrivial solution that is bounded uniformly on $\mathbb{R}^{+}$. If $c_{3}>0$, then the only bounded solution of \eqref{Ginzburg-Landau,c3} on $\mathbb{R}^{+}$ is $A_{2}(s)\equiv 0$.
\end{Lemma} 
This Lemma is proven in \cite[Lemma 2.3]{mccalla2013spots}, and relies on a hypothesis \cite[Hypothesis 1]{mccalla2013spots} that the general non-autonomous real Ginzburg-Landau equation, where $c_{0}=1$, and $c_{3} = -1$, has a nontrivial bounded solution $q(s)$ on $\mathbb{R}^{+}$ and the linearisation around $q(s)$ does not have a nontrivial bounded solution. This hypothesis was later proven with rigorous numerical methods by van den Berg et al. \cite[Theorem 1.1]{vandenberg2015Rigorous}. Next we write the solution $q(s)$ as
\begin{align}
    z_{2} = \frac{1}{A_{2}}\frac{\textnormal{d}}{\textnormal{d}s}A_{2}, \qquad \sigma_{2} = \frac{1}{s},\nonumber
\end{align}
{and transform between the coordinates of the transition and rescaling charts via \eqref{scale:rel}; then, we conclude that the real functions} 
\begin{align}
    \left(A_{1}^{+}, z_{1}^{+}, \varepsilon_{1}^{+}\right)(\rho) = {\left( \varepsilon^{\frac{1}{2}}\textnormal{e}^{\rho} q\left( \varepsilon^{\frac{1}{2}}\textnormal{e}^{\rho}\right), \frac{\frac{\textnormal{d}}{\textnormal{d}\rho} q\left( \varepsilon^{\frac{1}{2}}\textnormal{e}^{\rho}\right)}{q\left( \varepsilon^{\frac{1}{2}}\textnormal{e}^{\rho}\right)},  \varepsilon^{\frac{1}{2}}\textnormal{e}^{\rho}\right)} \quad \overset{\eqref{scale:rel}}{\mathlarger{\thicksim}} \quad  
    \left(A_{2}^{+}, z_{2}^{+}, \sigma_{2}^{+}\right)(s) = \left(q\left(s\right), \frac{{\frac{\textnormal{d}}{\textnormal{d}s} q}\left(s\right)}{q\left(s\right)}, \frac{1}{s}\right),\label{soln:chart}
\end{align}
satisfy \eqref{amp:a1;b1,eps0} and \eqref{amp:a2;z2,eps0}, respectively. {We note that, after inverting \eqref{scale:2} \eqref{scale:1}, and \eqref{z:def}, we can write \eqref{soln:chart} as
\begin{align}
    A^{+}(\rho) &= \varepsilon^{\frac{1}{2}} q\left(\varepsilon^{\frac{1}{2}}\textnormal{e}^{\rho}\right), \qquad B^{+}(\rho) = \varepsilon^{\frac{1}{2}}\textnormal{e}^{-\rho} \left(\frac{\textnormal{d}}{\textnormal{d}\rho}+ \frac{1}{2}\right)\; q\left(\varepsilon^{\frac{1}{2}}\textnormal{e}^{\rho}\right), \qquad \sigma^{+}(\rho) = \textnormal{e}^{-\rho}, \qquad\textnormal{where}\quad \rho=\log(r)\in\mathbb{R}.\nonumber
\end{align}
This solves \eqref{amp:A;B} as $\varepsilon\to0$ and has the property
\begin{align}
    \left(\frac{A^{+}(\rho)}{\sigma^{+}(\rho)},\; \frac{B^{+}(\rho)}{\sigma^{+}(\rho)A^{+}(\rho)} - \frac{1}{2},\; \frac{\varepsilon^{\frac{1}{2}}}{\sigma^{+}(\rho)}\right)=\left(A_{1}^{+}, z_{1}^{+}, \varepsilon_{1}^{+}\right)(\rho)&\to\left(0, \frac{1}{2}, 0\right), \qquad \textnormal{as}\; \rho\to-\infty,\nonumber\\
    \left(\frac{A^{+}(\rho)}{\varepsilon^{\frac{1}{2}}},\; \frac{B^{+}(\rho)}{\varepsilon^{\frac{1}{2}}A^{+}(\rho)} - \frac{\sigma^{+}(\rho)}{2\varepsilon^{\frac{1}{2}}},\; \frac{\sigma^{+}(\rho)}{\varepsilon^{\frac{1}{2}}}\right)=\left(A_{2}^{+}, z_{2}^{+}, \sigma_{2}^{+}\right)(\rho)&\to\left(0, -\sqrt{c_{0}}, 0\right), \qquad \textnormal{as}\; \rho\to\infty.\nonumber
\end{align}
Hence, we see that this solution lies in both $\mathcal{W}^{u}(P_{+})$ in the transition chart and $\mathcal{W}^{cs}(Q_{-}[0])$ in the rescaling chart, and so we call the solution represented by \eqref{soln:chart} a singular {(as $\rho\to-\infty$)} heteroclinic orbit between equilibria $P_{+}$ and $Q_{-}(0)$.} In \cite[Lemma 2.4]{mccalla2013spots}, McCalla \& Sandstede prove that \eqref{soln:chart} is a connecting orbit that forms a transverse intersection {between $\mathcal{W}^{u}(P_{+})$ and $\mathcal{W}^{cs}(Q_{-}[0])$ in $\mathbb{R}^{3}$}. We note that the real-valued heteroclinic orbit given in \eqref{soln:chart} generates a one-parameter family
\begin{align}
    \mathbf{q}_{+}(s): \quad \left(A_{1}, z_{1}, \varepsilon_{1}\right)(\rho) = \left(\textnormal{e}^{\textnormal{i}Y}A^{+}_{1}, z_{1}^{+}, \varepsilon_{1}^{+}\right) \quad \overset{\eqref{scale:rel}}{\mathlarger{\thicksim}} \quad 
    \left(A_{2}, z_{2}, \sigma_{2}\right)(s) = \left(\textnormal{e}^{\textnormal{i}Y}A^{+}_{2}, z_{2}^{+}, \sigma_{2}^{+}\right),\label{soln:chart;comp}
\end{align}
of heteroclinic orbits of \eqref{amp:a1;b1,eps0} and \eqref{amp:a2;z2,eps0} between equilibria $P_{+}$, $ Q_{-}(0)$ that are parametrised by $Y\in\mathbb{R}$, where orbits related to the choice of $Y=0$ and $Y=\pi$ both lie in the invariant subspace $\mathbb{R}^{3}$. 

The existence of the connecting orbit \eqref{soln:chart;comp} is reliant on the sign of $c_{3}$ being negative, which we will now discuss in more detail. As we alluded to in Remark \ref{rem:c0c3}, the value of $c_{3}$ in the radial problem is equivalent to leading order to that found for the one-dimensional problem in \cite{groves2017pattern}, and can be written down explicitly in terms of $ k_{D}$ and $M_{0}$, as defined in \eqref{c3}. We see from \eqref{c3} that one can find a restriction on $ k_{D}$ for a given $M_{0}$ such that $c_{3} <0$, as seen in Figure \ref{fig:intro;connectorq}a). For large depth fluids one has a critical magnetic permeability $ M_{c}\approx 0.56$, or equivalently $\mu_{c} \approx 3.54$, such that $c_{3}<0$ for all $M_{0}>M_{c}$; for experimentally relevant ferrofluids (where $2<\mu<6$) one requires $ k_{D}> k_{D}^{*} \approx 1.8$. For sufficiently small depths, Figure \ref{fig:intro;connectorq}a) appears to suggest that $\mathbf{q}_{+}(s)$ exists for all values of $\mu$. However, our free-surface model becomes unphysical for these shallow depths, since one should then also include a magnetic potential below the ferrofluid. Therefore, we will ignore the small-$ k_{D}$ shaded region seen in Figure \ref{fig:intro;connectorq}a).

\subsubsection{Connecting $P_{-}$ and $ Q_{-}$}
We start by analyzing \eqref{amp:a1;b1,eps0}: rewriting the system in $(A_{-},z_{-},\sigma_{-},\varepsilon_{-})$ coordinates, where $\sigma_{-}=0$, we find
\begin{align}
    \frac{\textnormal{d}}{\textnormal{d} \rho}A_{-} = A_{-}\left[1/2 +z_{-}\right],\qquad
    \frac{\textnormal{d}}{\textnormal{d} \rho}z_{-} = z_{-} + c_{0}\varepsilon_{-}^{2} - z_{-}^{2} + c_{3}|A_{-}|^2, \qquad 
    \frac{\textnormal{d}}{\textnormal{d} \rho}\varepsilon_{-} =\varepsilon_{-},\nonumber
\end{align}
and we obtain an explicit solution given by
\begin{align}
    \left(A^{-}_{-}, z^{-}_{-}, \varepsilon^{-}_{-}\right) = \left(0, -\sqrt{c_{0}}\,\textnormal{e}^{\rho} , \textnormal{e}^{\rho}\right)  \quad \overset{\eqref{scale:rel}}{\mathlarger{\thicksim}} \quad \left(A^{-}_{2}, z^{-}_{2}, \sigma^{-}_{2}\right) = \left(0, -\sqrt{c_{0}}-\frac{1}{2s}, \frac{1}{s}\right).\label{soln:exp;spotA}
\end{align}
This solution \eqref{soln:exp;spotA} satisfies $z_{2}(s) \to -\sqrt{c_{0}}$ as $s\to\infty$ {(which is equivalent to $\rho\to\infty$)} and $z_{-}(\rho)\to 0$ as $\rho\to-\infty$, and therefore lies in the intersection of $\mathcal{W}^{u}(P_{-})$ and $\mathcal{W}^{cs}( Q_{-})$. {For sufficiently small solutions, the unstable manifold $\mathcal{W}^{u}(P_{-})$ encompasses the entire invariant subspace $\mathbb{R}^{3}$, and so the transversality of this intersection is trivial.} We note the real-valued heteroclinic orbit given in \eqref{soln:exp;spotA} generates a one-parameter family
\begin{align}
\mathbf{q}_{-}: \quad (A_{-}, z_{-}, \varepsilon_{-}) = \left(\textnormal{e}^{\textnormal{i}Y}A^{-}_{-}, z^{-}_{-}, \varepsilon^{-}_{-}\right) \quad \overset{\eqref{scale:rel}}{\mathlarger{\thicksim}} \quad \left(A_{2}, z_{2}, \sigma_{2}\right) = \left(\textnormal{e}^{\textnormal{i}Y}A^{-}_{2}, z^{-}_{2}, \sigma^{-}_{2}\right),\label{qtilde:defn}
\end{align}
of heteroclinic orbits of \eqref{amp:a1;b1,eps0} and \eqref{amp:a2;z2,eps0} between the equilibria $P_{-}$, $ Q_{-}(0)$ that are parametrised by $Y\in\mathbb{R}$, where orbits related to the choice of $Y=0$ and $Y=\pi$ both lie in the invariant subspace $\mathbb{R}^{3}$. 
\section{Radial Solutions: Construction and Matching}\label{s:match}
{This section is devoted to showing the existence of three separate classes of localised radial patterns: spot A solutions are found in \S\ref{s:spotA}, spot B solutions in \S\ref{s:spotB} and ring solutions in \S\ref{s:rings}. For each class of pattern, we track solutions through the coordinate charts defined in Section \ref{s:geoblowup} in order to find a parametrisation for exponentially decaying solutions on the centre-manifold $\mathcal{W}^{c}_{+}(\varepsilon)|_{\sigma=\frac{1}{r}}$ at the matching point $r=r_{0}$. Then, by applying \eqref{match:farf:folia}, we return to the full infinite-dimensional problem and arrive at a parametrisation for the far-field stable manifold $\widetilde{\mathcal{W}}^{s}_{+}(\varepsilon)$ at $r=r_{0}$. Finally, for each class of pattern, we match our parametrisation for $\widetilde{\mathcal{W}}^{s}_{+}(\varepsilon)$ with the core manifold $\widetilde{\mathcal{W}}^{cu}_{-}(\varepsilon)$ in order to write down the leading-order profile of each pattern across the entire range $r\in(0,\infty)$.}
\subsection{Spot A}\label{s:spotA}
\begin{figure}[t!]
    \centering
    \includegraphics[height=6cm]{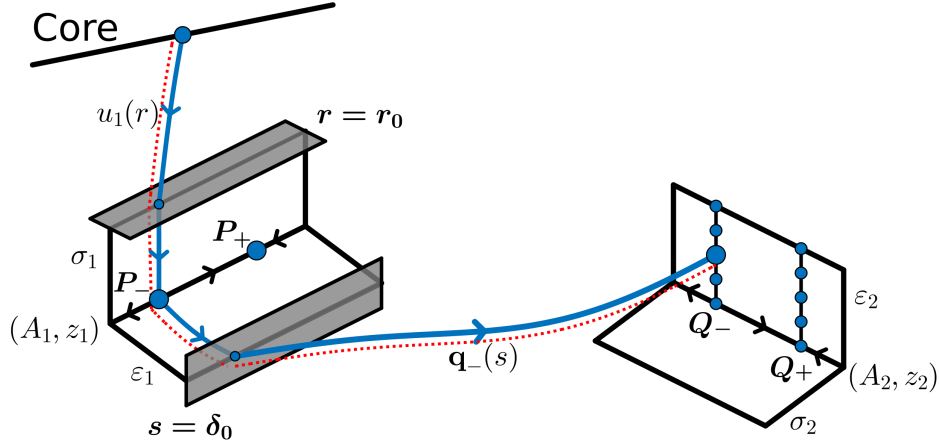}
    \caption{The path followed by the spot A solution; it begins by following $u_{1}(r)$ to the equilibrium $P_{-}$, and then tracking the connecting orbit $\mathbf{q}_{-}(s)$, as defined in \eqref{soln:exp;spotA}, to the equilibrium $ Q_{-}$.}
    \label{fig:GBUSpotA}
\end{figure}
The construction of the up-spot, termed spot A in \cite{Mccalla2010snaking,mccalla2013spots}, proceeds as follows; we follow the centre-stable manifold $\mathcal{W}^{cs}( Q_{-})$ near the heteroclinic orbits $\mathbf{q}_{-}(s)$ defined in \eqref{qtilde:defn} for $Y=0, \pi$ backward in `spatial time' $\rho$ past the equilibria $P_{-}$, as shown in Figure \ref{fig:GBUSpotA}. This will define the initial conditions required for solutions to decay to zero as $r\to\infty$, parametrising the centre coordinates of the far-field manifold at $r=r_{0}$. These centre coordinates are then substituted into the foliation parametrisation defined in \eqref{match:farf:folia} and the final matching between the core and far-field manifolds will take place at $r=r_{0}$.
{\paragraph{{Formal Scaling Analysis}: Spot A}
For this formal analysis, we have set all parameters equal to 1, other than the coordinates we wish to match. In the rescaling chart, solutions follow close to the connecting orbit $\mathbf{q}_{-}(s)$ and so, at the point $\varepsilon_{-}=1$, the manifold $\mathcal{W}^{cs}( Q_{-})$ can be written as 
\begin{align}
    \mathcal{W}^{cs}( Q_{-}): A_{-}(0)= \widetilde{a}\textnormal{e}^{\textnormal{i}Y},\qquad z_{-}(0)=-1, \qquad \sigma_{-}(0)=\varepsilon^{\frac{1}{2}}, \qquad \varepsilon_{-}(0)=1.\nonumber
\end{align}
We recall that spot A solutions lie close to the linear flow of the equilibrium $P_{-}$, and so will be subject to the equations
\begin{align}
    \frac{\textnormal{d}}{\textnormal{d}\rho}A_{-}=\frac{1}{2}A_{-}, \qquad \frac{\textnormal{d}}{\textnormal{d}\rho}z_{-}=z_{-}, \qquad \frac{\textnormal{d}}{\textnormal{d}\rho}\sigma_{-}=-\sigma_{-}, \qquad \frac{\textnormal{d}}{\textnormal{d}\rho}\varepsilon_{-}=\varepsilon_{-}.\nonumber
\end{align}
Hence, evolving our solutions with respect to the linear flow near $P_{-}$, we find
\begin{align}
    A_{-}(\rho)= \widetilde{a}\textnormal{e}^{\textnormal{i}Y}\textnormal{e}^{\frac{\rho}{2}},\qquad z_{-}(\rho)=-\textnormal{e}^{\rho}, \qquad \sigma_{-}(\rho)=\varepsilon^{\frac{1}{2}}\textnormal{e}^{-\rho}, \qquad \varepsilon_{-}(\rho)=\textnormal{e}^{\rho},\nonumber
\end{align}
which we evolve backwards in $\rho$, until the point $\rho^{*}:=\log\left(\varepsilon^{\frac{1}{2}}\right)$, such that $\sigma_{-}(\rho^{*})=1$. Then, at $\rho=\rho^{*}$, our solutions take the form,
\begin{align}
    A_{-}(\rho^{*})= \varepsilon^{\frac{1}{4}}\widetilde{a}\textnormal{e}^{\textnormal{i}Y},\qquad z_{-}(\rho^{*})=-\varepsilon^{\frac{1}{2}}, \qquad \sigma_{-}(\rho^{*})=1, \qquad \varepsilon_{-}(\rho^{*})=\varepsilon^{\frac{1}{2}},\nonumber
\end{align}
which we convert back into $A,B$ coordinates via the transformation $A =\sigma_{-}A_{-}$, $B=\sigma_{-}^{2}A_{-}z_{-}$; so,
\begin{align}
    A(r_{0})= \varepsilon^{\frac{1}{4}}\widetilde{a}\textnormal{e}^{\textnormal{i}Y},\qquad B(r_{0})=-\varepsilon^{\frac{3}{4}}\widetilde{a}\textnormal{e}^{\textnormal{i}Y},\nonumber
\end{align}
Our core parametrisation can be written as 
\begin{align}
    \widetilde{\mathcal{W}}^{cu}_{-}: A(r_{0})= -\textnormal{i}d_{1} + d_{2},\qquad B(r_{0})=-\textnormal{i}d_{1} - \nu d_{2}^{2},\nonumber
\end{align}
where we have set all other parameters equal to 1, and absorbed the complex phase into the definition of $Y$. Matching our far-field solutions to the core parametrisation and taking real and imaginary parts, we find
\begin{align}
    d_{1} &= -\varepsilon^{\frac{1}{4}}\widetilde{a}\sin(Y),\qquad
    d_{2} = \varepsilon^{\frac{1}{4}}\widetilde{a}\cos(Y),\qquad 
    d_{1} = \varepsilon^{\frac{3}{4}}\widetilde{a}\sin(Y),\qquad
    \nu d_{2}^{2} = \varepsilon^{\frac{3}{4}}\widetilde{a}\cos(Y),\nonumber
\end{align}
and so, we see that $d_{1}=Y=0$, and
\begin{align}
    d_{2}=\widetilde{a}\varepsilon^{\frac{1}{4}}, \qquad \widetilde{a} = \frac{\varepsilon^{\frac{1}{4}}}{\nu}.\nonumber
\end{align}
Hence, we expect to find $\widetilde{a}=\textnormal{O}\left(\varepsilon^{\frac{1}{4}}\right)$, $d_{2} = \textnormal{O}\left(\varepsilon^{\frac{1}{2}}\right)$, and $d_{1}$ to be of higher order in $\varepsilon$. Also, we note that this analysis breaks down when $\nu=0$; as such, we will consider the cases when $\nu\neq0$ and $\nu\approx0$ independently.
}
{\subsubsection{Tracking Solutions: Spot A for $\nu\neq0$}\label{subs:spotA;track}} 
We recall some of the results established in Section \ref{s:geoblowup}; the $(A_{-}, z_{-}, \sigma_{-}, \varepsilon_{-})$ coordinates satisfy the system \eqref{amp:a1;b1,-} with equilibrium point $P_{-}:(A_{-}, z_{-}, \sigma_{-}, \varepsilon_{-})=(0, 0, 0, 0)$. The associated linear problem to \eqref{amp:a1;b1,-} is hyperbolic at $P_{-}$, with a three dimensional unstable eigenspace and a one dimensional stable eigenspace. Furthermore, the tangent space of $\mathcal{W}^{u}(P_{-})$ at $P_{-}$ is spanned by the $A_{-}$, $z_{-}$, and $\varepsilon_{-}$ directions. Similarly, the $(A_{2}, z_{2}, \sigma_{2}, \varepsilon_{2})$ coordinates satisfy the system \eqref{amp:a2;z2,eps} with equilibrium point $ Q_{-}(\varepsilon_{2}):(A_{2}, z_{2}, \sigma_{2}, \varepsilon_{2})=\left(0, -\sqrt{c_{0}}, 0 , \textnormal{O}(|\varepsilon_{2}|^{2})\right)$, where $\varepsilon_{2}$ plays the role of a parameter. The associated linear problem of \eqref{amp:a2;z2,eps} about the solution $ Q_{-}(\varepsilon_{2})$ has a three dimensional centre-stable eigenspace and a one dimensional unstable eigenspace. In particular, the tangent space of $\mathcal{W}^{cs}( Q_{-})$ at $ Q_{-}(\varepsilon_{2})$ is spanned by the $A_{2}$, $\varepsilon_{2}$, and $\sigma_{2}$ directions.

\begin{figure}[t!]
    \centering
    \includegraphics[width=\linewidth]{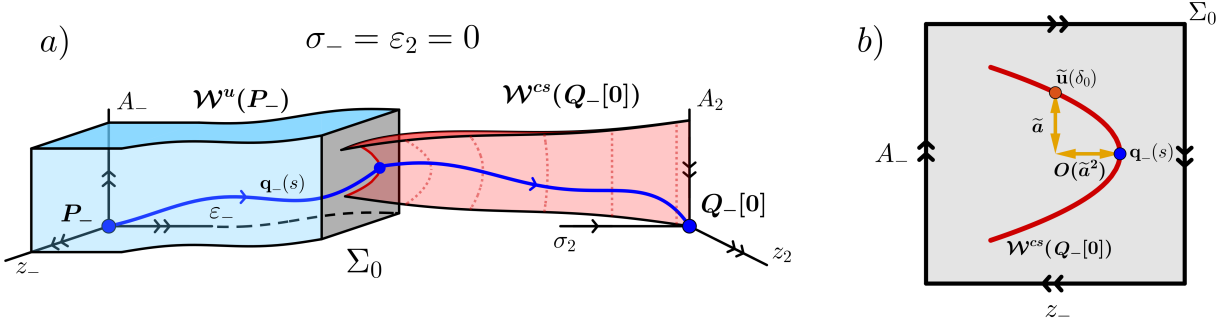}
    \caption{a) An illustration of the parametrisation of the intersection of $\mathcal{W}^{u}(P_{-})$ and $\mathcal{W}^{cs}( Q_{-})$ at the Poincar\'e section $\Sigma_{0}$, displayed in the invariant subspace $\sigma_{-}=\varepsilon_{2}=0$ for simplicity. b) The quadratic tangency between $\mathcal{W}^{u}(P_{-})$ and $\mathcal{W}^{cs}( Q_{-})$ is shown on the section $\Sigma_{0}$; any small $\textnormal{O}(\widetilde{a})$ perturbation in the $A_{-}$-direction yields an $\textnormal{O}(\widetilde{a}^{2})$ perturbation in the $z_{-}$-direction. }
    \label{fig:GBUMatching}
\end{figure}

We begin by tracking the solution in the invariant subspace $\varepsilon_{2}=0$, i.e. at the bifurcation point. We know that the solution $\mathbf{q}_{-}(s)$ lies on the intersection $\mathcal{W}^{u}(P_{-})\cap\mathcal{W}^{cs}( Q_{-}[0])$. Therefore, we will first parametrise the spot A solution as a small perturbation away from the solution $\mathbf{q}_{-}(s)$ at a fixed small value $s=\delta_{0}>0$. By taking a small perturbation of $\mathbf{q}_{-}(s)$, we define the coordinates $A_{2} = \widetilde{A}_{2}$, $z_{2} = -\sqrt{c_{0}} - \frac{1}{2s} + {\widetilde{z}_{2}}$, and $\sigma_{2}=\frac{1}{s}$, where $\widetilde{A}_{2}$ must satisfy the flow \eqref{amp:a2;z2,0},
\begin{align}
    \frac{\textnormal{d}}{\textnormal{d} s}\widetilde{A}_{2} &= \left[-\sqrt{c_{0}}-\frac{1}{2s}\right]\widetilde{A}_{2} + \textnormal{O}(|\widetilde{A}_{2}||\widetilde{z}_{2}|),\nonumber\\
    \intertext{and so $A_{2}(s)$, for small perturbations from $\mathbf{q}_{-}(s)$, has the form}
    A_{2}(s)\approx \widetilde{A}_{2}(s) &= \widetilde{A}_{2}(\delta_{0})\sqrt{\frac{\delta_{0}}{s}}\textnormal{e}^{\sqrt{c_{0}}[\delta_{0}-s]}, \qquad \textnormal{for} \, s\in[\delta_{0},\infty).\label{A2:s}
\end{align}
We define the section $\Sigma_{0}:=\left\{(A_{2}, z_{2}, \sigma_{2}, \varepsilon_{2})(s) \,: \;s=\delta_{0}\right\}$ and investigate the parametrisation of $\mathcal{W}^{cs}( Q_{-})\cap \Sigma_{0}$. We define $\widetilde{a}:=\widetilde{A}_{2}(\delta_{0})$ and use the $S^{1}$-symmetry in the normal-form to introduce the parameter $Y\in\mathbb{R}$, such that  $(A_{2},z_{2},\sigma_{2},\epsilon_{2})(\delta_{0})=(\widetilde{a}\textnormal{e}^{\textnormal{i}Y}, -\sqrt{c_{0}} - \frac{1}{2\delta_{0}} {+ \widetilde{z}_{2}(\delta_{0})}, \frac{1}{\delta_{0}},0)$. We perform the coordinate change \eqref{subs:rescale} to write this as our initial parametrisation for $\mathcal{W}^{u}(P_{-})$ at the section $\Sigma_{0}=\left\{(A_{-}, z_{-}, \sigma_{-}, \varepsilon_{-}) \,: \;{\varepsilon_{-}}=\delta_{0}\right\}$. Linearising \eqref{amp:a1;b1,eps0} about the solution $\mathbf{q}_{-}(s)$ at $\varepsilon_{-}=\delta_{0}$, we find that, defining $A_{-} = \widetilde{A}_{-}$ and $z_{-} = -\sqrt{c_{0}}\delta_{0} + \widetilde{z}_{-}$, the small variables $(\widetilde{A}_{-}, \widetilde{z}_{-})$ must satisfy
\begin{align}
    &\frac{\textnormal{d}}{\textnormal{d} \rho}\widetilde{A}_{-} = \frac{1}{2}\widetilde{A}_{-} + \textnormal{O}(|\widetilde{A}_{-}|(|\delta_{0}| + |\widetilde{z}_{-}|)),&
    \qquad 
    &\frac{\textnormal{d}}{\textnormal{d} \rho}\widetilde{z}_{-} = z_{-} + \textnormal{O}((|\widetilde{A}_{-}| + |\widetilde{z}_{-}| + |\delta_{0}|)^{2}),&\nonumber
\end{align}
and so, for small-amplitude solutions $\widetilde{z}_{-}\propto\textnormal{e}^{\rho}\propto|\widetilde{A}_{-}|^{2}$. This explains the quadratic dependency illustrated in Figure \ref{fig:GBUMatching}. Therefore, we can paramatrise $\mathcal{W}^{cs}( Q_{-})\cap\Sigma_{0}$ in $(A_{-},z_{-},\sigma_{-})$ coordinates, when $\sigma_{-}=0$, as
\begin{align}
\left.\mathcal{W}^{cs}( Q_{-})\cap\Sigma_{0}\right|_{\sigma_{-}=0} &= \left\{ (A_{-}, z_{-}, \sigma_{-}) = (\widetilde{a}\textnormal{e}^{\textnormal{i}Y}, -\sqrt{c_{0}}\delta_{0} + \textnormal{O}(\widetilde{a}^{2}), 0) \,:\; |\widetilde{a}|< a_{0}\right\},\label{Wcs:eps0}
\end{align}
for each small $\delta_{0}$ and constant $a_{0}$. Having parametrised solutions travelling from $\varepsilon_{-}=\delta_{0}$ to $ Q_{-}(0)$, we now want to perturb our parametrisation by a small value of $\varepsilon_{2}$, such that solutions tend to $ Q_{-}(\varepsilon_{2})$ as $s\to\infty$. Setting $\varepsilon_{2}= \varepsilon^{\frac{1}{2}}$, we define 
\begin{align}
    \left(A_{-}, z_{-}, \sigma_{-}, \varepsilon_{-}\right) = \left( \textnormal{O}( {\varepsilon}), -\sqrt{c_{0}}\textnormal{e}^{\rho} + \textnormal{O}( {\varepsilon}),  \varepsilon^{\frac{1}{2}}\textnormal{e}^{-\rho}, \textnormal{e}^{\rho} \right) \; \overset{\eqref{scale:rel}}{\mathlarger{\thicksim}} \; \left(A_{2}, z_{2}, \sigma_{2}, \varepsilon_{2}\right) = \left( \textnormal{O}( {\varepsilon}), -\sqrt{c_{0}}- \frac{1}{2s} + \textnormal{O}( {\varepsilon}), \frac{1}{s},  \varepsilon^{\frac{1}{2}} \right),\nonumber
\end{align}
which, for small values of $ {\varepsilon}$, is asymptotically close to $\mathbf{q}_{-}(\delta_{0})$. We note that adding the parameter $ {\varepsilon}$ yields additional $\textnormal{O}( {\varepsilon})$ terms, since the remainder terms in the rescaling chart \eqref{amp:a2;z2,eps} are of order $\textnormal{O}(|\varepsilon_{2}|^{2})$.\\
Therefore, for each small $\delta_{0}>0$, there are constants $a_{0}, \varepsilon_{0}>0$ such that
\begin{align}
\mathcal{W}^{cs}( Q_{-})\cap\Sigma_{0} &= \left\{ (A_{-}, z_{-}, \sigma_{-}) = \left(\widetilde{a}\textnormal{e}^{\textnormal{i}Y} + \textnormal{O}( {\varepsilon}), -\sqrt{c_{0}}\delta_{0} + \textnormal{O}(\widetilde{a}^{2} +  {\varepsilon}), \frac{ \varepsilon^{\frac{1}{2}}}{\delta_{0}}\right) \,:\; |\widetilde{a}|< a_{0}, \; | {\varepsilon}|<\varepsilon_{0}\right\}.\label{Wcs;param;delta0}
\end{align}
Again, we see that \eqref{Wcs;param;delta0} is asymptotically close to \eqref{Wcs:eps0} as $ {\varepsilon}\to0$. {From our formal analysis above, we set $\widetilde{a}:=\delta^{\frac{1}{2}}_{0}  {\varepsilon^{\frac{1}{4}}} a$, where $a$ is a fixed constant}. Then, we obtain the initial data
{\begin{align}
    A_{-}(\rho_{0}) = a\delta_{0}^{\frac{1}{2}} {\varepsilon^{\frac{1}{4}}}\textnormal{e}^{\textnormal{i}Y} + \textnormal{O}( {\varepsilon}), \qquad z_{-}(\rho_{0}) = -\sqrt{c_{0}}\delta_{0} + \textnormal{O}( {\varepsilon}), \qquad \sigma_{-}(\rho_{0}) = \frac{ \varepsilon^{\frac{1}{2}}}{\delta_{0}}, \qquad \varepsilon_{-}(\rho_{0}) = \delta_{0},\label{inital:spota}
\end{align}
for which we need to solve \eqref{amp:a1;b1,-} for $\rho^{*}\leq\rho\leq\rho_{0}$, where 
\begin{align}
    \rho_{0} = {\log}\;\,\frac{\delta_{0}}{ \varepsilon^{\frac{1}{2}}}, \qquad \rho^{*} = {\log}\;\,\frac{1}{\sigma_{0}}.\nonumber
\end{align}
In particular, we find
\begin{align}
    \varepsilon_{-}(\rho) = \varepsilon_{-}(\rho_{0})\textnormal{e}^{\rho-\rho_{0}}, \qquad \varepsilon_{-}(\rho) =  \varepsilon^{\frac{1}{2}} \textnormal{e}^{\rho}, \qquad \sigma_{-}(\rho) = \sigma_{0}\textnormal{e}^{\rho^{*}-\rho},\nonumber
\end{align}
for $\rho^{*}\leq\rho\leq\rho_{0}$.} Next, we define $A_{-}(\rho) = \widetilde{A}_{-}(\rho)\textnormal{e}^{\rho/2}$, $z_{-}(\rho) = \widetilde{z}_{-}(\rho)\textnormal{e}^{\rho}$
and obtain the system
{\begin{align}
    \frac{\textnormal{d}}{\textnormal{d}\rho} \widetilde{A}_{-} &= \widetilde{A}_{-}\left[\widetilde{z}_{-}\textnormal{e}^{\rho} + \textnormal{e}^{2(\rho^{*}-\rho)}\textnormal{O}\left(|\sigma_{0}|^2\right)\right],\label{Az:tilde;syst}\\
    \frac{\textnormal{d}}{\textnormal{d}\rho} \widetilde{z}_{-} &= -\textnormal{e}^{\rho}\widetilde{z}_{-}^{2} + c_{0}|\varepsilon_{-}(\rho_{0})|^2 \textnormal{e}^{\rho} + c_{3}|\widetilde{A}_{-}|^2 + \textnormal{e}^{2(\rho^{*}-\rho)}|\sigma_{0}|^2 \textnormal{O}\left(|\widetilde{A}_{-}|^{4} + |\widetilde{z}_{-}| + |\varepsilon_{-}(\rho_{0})|\right),\nonumber
\end{align}}
which we consider with initial conditions \eqref{inital:spota}, written in our new coordinates as
{\begin{align}
    \widetilde{A}_{-}(\rho_{0}) = \widetilde{A}_{-}^{0}:= a \varepsilon^{\frac{1}{2}}\textnormal{e}^{\textnormal{i}Y} + \textnormal{O}( {\varepsilon}), \qquad \widetilde{z}_{-}(\rho_{0}) = \widetilde{z}_{-}^{0}:= -\sqrt{c_{0}} \varepsilon^{\frac{1}{2}} + \textnormal{O}( {\varepsilon}). \label{initial:spota;tilde}
\end{align}
We write the solution to \eqref{Az:tilde;syst} as the variation of constants equation
\begin{align}
    \widetilde{A}_{-}(\rho) &= \widetilde{A}^{0}_{-} + \int_{\rho_{0}}^{\rho} \widetilde{A}_{-}(u)\left[\widetilde{z}_{-}(u)\textnormal{e}^{u} + \textnormal{e}^{2(\rho^{*}-u)}\textnormal{O}\left(|\sigma_{0}|^2\right)\right] \textnormal{d}u, \label{fixedpoint:spota}\\
    \widetilde{z}_{-}(\rho) &= \widetilde{z}^{0}_{-} + \int_{\rho_{0}}^{\rho} -\textnormal{e}^{u}\widetilde{z}_{-}^{2}(u) + c_{0}|\varepsilon_{-}(\rho_{0})|^2 \textnormal{e}^{u} + c_{3}|\widetilde{A}_{-}(u)|^2 + \textnormal{e}^{2(\rho^{*}-u)}|\sigma_{0}|^2 \textnormal{O}\left(|\widetilde{A}_{-}(u)|^{4} + |\widetilde{z}_{-}(u)| + |\varepsilon_{-}(\rho_{0})|\right) \textnormal{d}u, \nonumber
\end{align}}
which, for sufficiently small variables, has a unique solution \cite{Sandstede1997Convergence}. Furthermore, there is a uniform constant $C$ with
\begin{align}
    \|\widetilde{A}_{-}\| \leq C|\widetilde{A}_{-}^{0}|, \qquad \|\widetilde{z}_{-}\| \leq C\left[|\widetilde{z}_{-}^{0}| + {|\varepsilon_{-}(\rho_{0})| + |\rho^{*}-\rho_{0}|}|\widetilde{A}_{-}^{0}|^2\right].\nonumber
\end{align}
Using these estimates, along with \eqref{fixedpoint:spota} and \eqref{initial:spota;tilde}, we obtain
{\begin{align}
&A_{-}(\rho) =a \varepsilon^{\frac{1}{2}} \textnormal{e}^{\rho/2}\textnormal{e}^{\textnormal{i}Y}(1 + \Delta_{A}),& \quad 
    &z_{-}(\rho) =-\sqrt{c_{0}} \varepsilon^{\frac{1}{2}}\textnormal{e}^{\rho}(1+\Delta_{A}),& \quad&\sigma_{-}(\rho) = \sigma_{0}\textnormal{e}^{\rho^{*}-\rho},& \quad &\varepsilon_{-}(\rho) =  \varepsilon^{\frac{1}{2}}\textnormal{e}^{\rho},&\nonumber\\
&A_{-}(\rho^{*}) =\frac{a}{\sqrt{\sigma_{0}}} \varepsilon^{\frac{1}{2}}\textnormal{e}^{\textnormal{i}Y}(1 + \Delta_{A}),& \quad 
    &z_{-}(\rho^{*}) =-\sqrt{c_{0}}\frac{ \varepsilon^{\frac{1}{2}}}{\sigma_{0}}(1+\Delta_{A}),& \quad&\sigma_{-}(\rho^{*}) = \sigma_{0},& \quad& \varepsilon_{-}(\rho^{*}) = \frac{ \varepsilon^{\frac{1}{2}}}{\sigma_{0}},&\label{initial:spot;rho0}
\end{align}}
with remainder terms
\begin{align}
    \Delta_{A} := \textnormal{O}(\sigma_{0} + \delta_{0} +  {\varepsilon^{\frac{1}{4}}}).\nonumber
\end{align}
Now, we return to the original $(A,B)$ coordinates by converting \eqref{initial:spot;rho0} into $(A_{1}, z_{1})$ coordinates and inverting the coordinate change \eqref{scale:1}. That is, we write $A = \sigma_{-}A_{-}$ and $B = \sigma_{-}^{2}A_{-}z_{-}$, and {evaluate at $\sigma=\frac{1}{r}$}; then, we obtain 
{\begin{align}
    A(r_{0}) &=\frac{a \varepsilon^{\frac{1}{2}}}{\sqrt{r_{0}}}\textnormal{e}^{\textnormal{i}Y}(1 + \Delta_{A}), \qquad 
    B(r_{0}) =-\sqrt{c_{0}} \frac{a {\varepsilon}}{\sqrt{r_{0}}} \textnormal{e}^{\textnormal{i}Y}(1+\Delta_{A}). \label{initial:spot;AB}
\end{align}}
{\subsubsection{Matching Core and Far-Field Solutions: Spot A for $\nu\neq0$}\label{subs:spotA;match}}
We are now fully equipped to find a nontrivial solution contained in the intersection of the core manifold $\widetilde{\mathcal{W}}^{cu}_{-}\change{(\varepsilon)}$ and the far-field manifold $\change{\widetilde{\mathcal{W}}^{s}_{+}(\varepsilon)}$. We begin by substituting the initial conditions \eqref{initial:spot;AB} into the far-field parametrisation \eqref{match:farf:folia} following the normal-form transformation \eqref{ab:tilde} to obtain
\begin{align}
    \widetilde{a}(r_{0}) &= A(r_{0}) +\textnormal{O}_{r_{0}}\left(|\vec{\mathbf{d}}|_{1}\left[|\varepsilon| + |A| + |B| + |\vec{\mathbf{d}}|_{1}\right]\right),\nonumber\\
    &= \textnormal{e}^{\textnormal{i}Y} a r_{0}^{-\frac{1}{2}}\varepsilon^{\frac{1}{2}}(1 + \Delta_{A}) +\textnormal{O}_{r_{0}}\left(|\vec{\mathbf{d}}|_{1}\left[{|\varepsilon|^{\frac{1}{2}}} + |\vec{\mathbf{d}}|_{1}\right]\right),\nonumber\\
    \widetilde{b}(r_{0}) &= B(r_{0}) +\textnormal{O}_{r_{0}}\left(|\vec{\mathbf{d}}|_{1}\left[|\varepsilon| + |A| + |B| + |\vec{\mathbf{d}}|_{1}\right]\right),\nonumber\\
    &= -\textnormal{e}^{\textnormal{i}Y} a \sqrt{c_{0}}r_{0}^{-\frac{1}{2}} \varepsilon (1+\Delta_{A}) +\textnormal{O}_{r_{0}}\left(|\vec{\mathbf{d}}|_{1}\left[{|\varepsilon|^{\frac{1}{2}}} + |\vec{\mathbf{d}}|_{1}\right]\right),\nonumber\\
    a_{n}(r_{0}) &= \textnormal{O}_{r_{0}}\left(\left[|A| + |B| + |\vec{\mathbf{d}}|_{1}\right]\left[|\varepsilon| + |A| + |B| + |\vec{\mathbf{d}}|_{1}\right]\right),\nonumber\\
    &= \textnormal{O}_{r_{0}}\left(\left[{|\varepsilon|^{\frac{1}{2}}} + |\vec{\mathbf{d}}|_{1}\right]^{2}\right),\nonumber\\
    a_{-n}(r_{0}) &= \vec{d}_{n}\textnormal{e}^{-\lambda_{n}r_{0}}  + \textnormal{O}_{r_{0}}\left(\left[|A| + |B| \right]\left[|\varepsilon| + |A| + |B| + |\vec{\mathbf{d}}|_{1}\right]\right) ,\nonumber\\
     &= \vec{d}_{n}\textnormal{e}^{-\lambda_{n}r_{0}}  + \textnormal{O}_{r_{0}}\left({|\varepsilon|^{\frac{1}{2}}}\left[ {|\varepsilon|^{\frac{1}{2}}} + |\vec{\mathbf{d}}|_{1}\right]\right). \nonumber
\end{align}
We also recall the transformed core parametrisation \eqref{match:core;normal}
\begin{align}
    \widetilde{a}(r_{0}) &= \textnormal{e}^{-\textnormal{i}\left(\frac{\pi}{4} + \textnormal{O}(r_{0}^{-2}) + \textnormal{O}_{r_{0}}(|\varepsilon| + |\mathbf{d}|_{1})\right)}r_{0}^{-\frac{1}{2}}\left([-\textnormal{i}+\textnormal{O}(r_{0}^{-1})]r_{0} d_{1} + [1+\textnormal{O}(r_{0}^{-1})]d_{2}\right) + \textnormal{O}_{r_{0}}\left(|\mathbf{d}|_{1}\left[|\varepsilon| + |\mathbf{d}|_{1}\right]\right),\nonumber\\
    \widetilde{b}(r_{0}) &= \textnormal{e}^{-\textnormal{i}\left(\frac{\pi}{4} + \textnormal{O}(r_{0}^{-2}) + \textnormal{O}_{r_{0}}(|\varepsilon| + |\mathbf{d}|_{1})\right)}r_{0}^{-\frac{1}{2}}\left(\left[-\textnormal{i} +\textnormal{O}(r_{0}^{-1})\right]d_{1} - \left[\nu + \textnormal{O}(r_{0}^{-\frac{1}{2}})\right] d_{2}^{2}\right)+ \textnormal{O}_{r_{0}}\left(|\mathbf{d}|_{1}\left[|\varepsilon| + |\mathbf{d}_{2}|_{1}\right] + |d_{2}|^{3}\right),\nonumber\\
    a_{n}(r_{0}) &= \textnormal{e}^{\lambda_{n}r_{0}}r_{0}^{-\frac{1}{2}}\left[\frac{1}{\sqrt{\pi}} + \textnormal{O}(r_{0}^{-1})\right]\widetilde{c}_{1,n} + \textnormal{O}_{r_{0}}\left(|\mathbf{d}|_{1}\left[|\varepsilon| + |\mathbf{d}|_{1}\right]\right),\nonumber\\
    a_{-n}(r_{0}) &= \textnormal{O}_{r_{0}}\left(|\mathbf{d}|_{1}\left[|\varepsilon| + |\mathbf{d}|_{1}\right]\right),\nonumber
\end{align}
where $|\mathbf{d}|_{1} = |d_{1}| + |d_{2}| + |\mathbf{c}_{1}|_{1}$, and $|\mathbf{d}_{2}|_{1} = |d_{1}| + |\mathbf{c}_{1}|_{1}$. Setting these parametrisations equivalent to each other for each respective coordinate $(\widetilde{a}, \widetilde{b}, a_{1}, a_{-1}, a_{2}, a_{-2}, \dots)(r_{0})$ is the same as finding the zeros of the functional 
\begin{align}
    G:&\left(d_{1}, d_{2}, \widetilde{Y},  \mathbf{c}_{1}, \vec{\mathbf{d}}; a, \varepsilon\right) \mapsto \left( G^{C}_{1}, G^{C}_{2}, G^{H}_{1}, G^{H}_{-1}, G^{H}_{2}, G^{H}_{-2}, \dots \right), \nonumber
\end{align}
where 
\begin{align}
    G^{C}_{1} &= [-\textnormal{i}+\textnormal{O}(r_{0}^{-1})]r_{0} d_{1} + [1+\textnormal{O}(r_{0}^{-1})]d_{2} - \textnormal{e}^{\textnormal{i}\widetilde{Y}} a \varepsilon^{\frac{1}{2}}(1 + \Delta_{A}) + \mathscr{R}^{C}_{1},\label{G:SpotA;C1}\\
    G^{C}_{2} &= \left[-\textnormal{i} +\textnormal{O}(r_{0}^{-1})\right] d_{1} - \left[\nu + \textnormal{O}(r_{0}^{-\frac{1}{2}})\right] d_{2}^{2} +\textnormal{e}^{\textnormal{i}\widetilde{Y}} a \sqrt{c_{0}} \varepsilon (1+\Delta_{A}) + \mathscr{R}^{C}_{2},\label{G:SpotA;C2}\\
    G^{H}_{n} &= \textnormal{e}^{\lambda_{n}r_{0}}r_{0}^{-\frac{1}{2}}\left[\frac{1}{\sqrt{\pi}} + \textnormal{O}(r_{0}^{-1})\right]c_{1,n} + \mathscr{R}^{H}_{1},\label{G:SpotA;Hn}\\
    G^{H}_{-n} &= - \vec{d}_{n}\textnormal{e}^{-\lambda_{n}r_{0}} + \mathscr{R}^{H}_{2},\label{G:SpotA;-Hn}
\end{align}
and we have defined
\begin{align}
    \widetilde{Y} :&= Y + \frac{\pi}{4} + \textnormal{O}(r_{0}^{-2}) + \textnormal{O}_{r_{0}}(|\varepsilon| + |d_{1}| + |d_{2}| + |\mathbf{c}_{1}|_{1}),\nonumber\\
    \mathscr{R}^{C}_{1} :&= \textnormal{O}_{r_{0}}\left(\left[|d_{1}| + |d_{2}| + |\mathbf{c}_{1}|_{1}\right]\left[|\varepsilon| + |d_{1}| + |d_{2}| + |\mathbf{c}_{1}|_{1}\right] + |\vec{\mathbf{d}}|_{1}\left[{|\varepsilon|^{\frac{1}{2}}} + |\vec{\mathbf{d}}|_{1}\right]\right),\nonumber\\
    \mathscr{R}^{C}_{2} :&= \textnormal{O}_{r_{0}}\left(\left[|d_{1}| + |d_{2}| + |\mathbf{c}_{1}|_{1}\right]\left[|\varepsilon| + |d_{1}| + |\mathbf{c}_{1}|_{1}\right] + |d_{2}|^{3} + |\vec{\mathbf{d}}|_{1}\left[{|\varepsilon|^{\frac{1}{2}}} + |\vec{\mathbf{d}}|_{1}\right]\right),\nonumber\\
    \mathscr{R}^{H}_{1} :&= \textnormal{O}_{r_{0}}\left(\left[|d_{1}| + |d_{2}| + |\mathbf{c}_{1}|_{1}\right]\left[|\varepsilon| + |d_{1}| + |d_{2}| + |\mathbf{c}_{1}|_{1}\right] + \left[{|\varepsilon|^{\frac{1}{2}}} + |\vec{\mathbf{d}}|_{1}\right]^{2}\right),\nonumber\\
    \mathscr{R}^{H}_{2} :&= \textnormal{O}_{r_{0}}\left(\left[|d_{1}| + |d_{2}| + |\mathbf{c}_{1}|_{1}\right]\left[|\varepsilon| + |d_{1}| + |d_{2}| + |\mathbf{c}_{1}|_{1}\right] + {|\varepsilon|^{\frac{1}{2}}}\left[{|\varepsilon|^{\frac{1}{2}}} + |\vec{\mathbf{d}}|_{1}\right]\right).\nonumber
\end{align}
We introduce the scaling $(d_{1}, d_{2}) = (\varepsilon \widetilde{d}_{1}, \varepsilon^{\frac{1}{2}} \widetilde{d}_{2})$ so that terms in \eqref{G:SpotA;C1} and \eqref{G:SpotA;C2} scale with the same order in $\varepsilon$. Initially, setting $\varepsilon=0$, we investigate \eqref{G:SpotA;Hn} and \eqref{G:SpotA;-Hn}
\begin{align}
    F_{n} &= \textnormal{e}^{\lambda_{n}r_{0}}r_{0}^{-\frac{1}{2}}\left[\frac{1}{\sqrt{\pi}} + \textnormal{O}(r_{0}^{-1})\right]c_{1,n} + \textnormal{O}_{r_{0}}\left( |\mathbf{c}_{1}|^{2}_{1} + |\vec{\mathbf{d}}|_{1}^{2}\right),\nonumber\\
    F_{-n} &= - \vec{d}_{n}\textnormal{e}^{-\lambda_{n}r_{0}} + \textnormal{O}_{r_{0}}\left( |\mathbf{c}_{1}|^{2}_{1}\right).\nonumber
\end{align}
Defining the functional $F:(\mathbf{c}_{1}, \vec{\mathbf{d}})\mapsto \left(F_{j}\right)_{j\in\mathbb{Z}\backslash\{0\}}$, it clear that $F(\mathbf{0},\mathbf{0})=\mathbf{0}$, where ${\mathbf{0} = (0,\dots)}$. Furthermore, the {Jacobian} $DF(\mathbf{0},\mathbf{0})$ is invertible, and so, using the implicit function theorem we can solve \eqref{G:SpotA;Hn} and \eqref{G:SpotA;-Hn} for all values of $n\in\mathbb{N}$ uniquely for sufficiently small $0<\varepsilon\ll1$. Matching orders of $\varepsilon$ in \eqref{G:SpotA;Hn} and \eqref{G:SpotA;-Hn}, we find that 
\begin{align}
    c_{1,n} = \textnormal{O}_{r_{0}}\left({|\varepsilon|}\right), \qquad \vec{d}_{n}  =\textnormal{O}_{r_{0}}\left({|\varepsilon|}\right), \qquad \forall n\in\mathbb{N}.\nonumber
\end{align}
Returning to \eqref{G:SpotA;C1} and \eqref{G:SpotA;C2}, we have
\begin{align}
    \widetilde{G}_{1} &= [1+\textnormal{O}(r_{0}^{-1})]\widetilde{d}_{2} - \textnormal{e}^{\textnormal{i}\widetilde{Y}} a (1 + \Delta_{A}) + [-\textnormal{i}+\textnormal{O}(r_{0}^{-1})]\varepsilon^{\frac{1}{2}}r_{0}\widetilde{d}_{1} + \varepsilon^{-\frac{1}{2}}\mathscr{R}^{C}_{1},\label{Match:SpotA;C}\\
    \widetilde{G}_{2} &= \left[-\textnormal{i} +\textnormal{O}(r_{0}^{-1})\right] \widetilde{d}_{1} - \left[\nu + \textnormal{O}(r_{0}^{-\frac{1}{2}})\right]\widetilde{d}_{2}^{2} +\textnormal{e}^{\textnormal{i}\widetilde{Y}} a \sqrt{c_{0}} (1+\Delta_{A}) + \varepsilon^{-1}\mathscr{R}^{C}_{2},\nonumber\\
    \intertext{where}
    \mathscr{R}^{C}_{1} &= |\varepsilon|\textnormal{O}_{r_{0}}\left(\left[|\varepsilon|^{\frac{1}{2}}|\widetilde{d}_{1}| + |\widetilde{d}_{2}| + {|\varepsilon|^{\frac{1}{2}}}\right]^{2} + {|\varepsilon|^{\frac{1}{2}}}\right),\nonumber\\
    \mathscr{R}^{C}_{2} &=  |\varepsilon|\textnormal{O}_{r_{0}}\left(\left[|\varepsilon|^{\frac{1}{2}}|\widetilde{d}_{1}| + |\widetilde{d}_{2}| + {|\varepsilon|^{\frac{1}{2}}}\right]\left[|\varepsilon|^{\frac{1}{2}}|\widetilde{d}_{1}| + {|\varepsilon|^{\frac{1}{2}}}\right] + |\varepsilon|^{\frac{1}{2}}|\widetilde{d}_{2}|^{3} + {|\varepsilon|^{\frac{1}{2}}}\right).\nonumber
\end{align}
Initially setting $\varepsilon=0$, we obtain the system
\begin{align}
    \widetilde{G}_{1} &= [1+\Delta_{A}]\widetilde{d}_{2} - \textnormal{e}^{\textnormal{i}\widetilde{Y}} a (1 + \Delta_{A}),\label{Match:SpotA;C,eps0}\\
    \widetilde{G}_{2} &= \left[-\textnormal{i} +\Delta_{A}\right] \widetilde{d}_{1} - \left[\nu + \Delta_{A}\right]\widetilde{d}_{2}^{2} +\textnormal{e}^{\textnormal{i}\widetilde{Y}} a \sqrt{c_{0}} (1+\Delta_{A}).\nonumber
\end{align}
We formally set $\Delta_{A}=0$ and separate \eqref{Match:SpotA;C,eps0} into real and imaginary parts: this is equivalent to finding zeros of the functional
\begin{align}
    \widetilde{G}(\widetilde{d}_{1}, \widetilde{d}_{2}, \widetilde{Y}, a) = \begin{pmatrix} \widetilde{d}_{2} - a \cos(\widetilde{Y}) \\ -a \sin(\widetilde{Y}) \\ -\nu \widetilde{d}_{2}^{2} + a \sqrt{c_{0}}\cos(\widetilde{Y}) \\ - \widetilde{d}_{1} + a \sqrt{c_{0}}\sin(\widetilde{Y})
    \end{pmatrix}.\nonumber
\end{align}
It is apparent that the vector
\begin{align}
    \left(\widetilde{d}_{1}^{*}, \widetilde{d}_{2}^{*}, \widetilde{Y}^{*}, a^{*}\right) = \left(0, {\textnormal{sgn}\;(\nu)\frac{\sqrt{c_{0}}}{|\nu|}}, {\arccos[\textnormal{sgn}\;(\nu)]}, {\frac{\sqrt{c_{0}}}{|\nu|}}\right)\nonumber
\end{align}
is a root of $\widetilde{G}$ with Jacobian
\begin{align}
    D\widetilde{G}\left(\widetilde{d}_{1}^{*}, \widetilde{d}_{2}^{*}, \widetilde{Y}^{*}, a^{*}\right) = \begin{pmatrix} 0 & 1 & 0 & -{\textnormal{sgn}\;(\nu)} \\ 0 & 0 & -{\textnormal{sgn}\;(\nu)\frac{\sqrt{c_{0}}}{|\nu|}} & 0 \\ 0 & -2\sqrt{c_{0}} & 0 & {\textnormal{sgn}\;(\nu)\sqrt{c_{0}}} \\ -1 & 0 & {\textnormal{sgn}\;(\nu)\frac{c_{0}}{|\nu|}} & 0\end{pmatrix}.\nonumber
\end{align}
Since $c_{0}>0$, the Jacobian is invertible and we can, using the implicit function theorem, solve \eqref{Match:SpotA;C,eps0} uniquely for all sufficiently small $\Delta$, that is, for $r_{0}$ large enough and $\delta_{0}$ small enough, and subsequently \eqref{Match:SpotA;C} for all $0<\varepsilon\ll1$. Reversing the scaling for $d$, we find that 
\begin{align}
    d_{1} &= \varepsilon\,\textnormal{O}(r_{0}^{-1} + \delta_{0} + \varepsilon^{\frac{1}{4}}),\label{SpotA:d1}\\
    d_{2} &= \varepsilon^{\frac{1}{2}}{\textnormal{sgn}\;(\nu)\frac{\sqrt{c_{0}}}{|\nu|}}\left(1 + \textnormal{O}(r_{0}^{-1} + \delta_{0} + \varepsilon^{\frac{1}{4}})\right).\label{SpotA:d2}
\end{align}
Hence, we have found a spot A solution. We recall that our solution $\mathbf{u}(r,y)$ takes the form
\begin{align}
    \mathbf{u}(r,y) &= a(r)\mathbf{e}(y) + b(r)\mathbf{f}(y) + \overline{a}(r)\overline{\mathbf{e}}(y) + \overline{b}(r)\overline{\mathbf{f}}(y) + \Sigma_{n=1}^{\infty} \left\{a_{n}(r)\mathbf{e}_{n}(y) + a_{-n}(r)\mathbf{e}_{-n}(y)\right\}, \nonumber
\end{align}
and 
\begin{align}
    \begin{pmatrix}a \\ b \end{pmatrix}(r) & = \sum_{i=1}^{4} \widetilde{d}_{i}\mathbf{V}_{i}(r), \qquad \qquad \begin{pmatrix}a_{n} \\ a_{-n} \end{pmatrix}(r)  = \sum_{i=1}^{2} \widetilde{c}_{i,n}\mathbf{W}_{i,n}(r).\nonumber
\end{align}
Substituting \eqref{SpotA:d1} and \eqref{SpotA:d2} into this form, we can write the spot A solution $\mathbf{u}_{A}$ as
\begin{align}
    \mathbf{u}_{A}(r,y) &= \varepsilon^{\frac{1}{2}}{\textnormal{sgn}\;(\nu)\frac{\sqrt{c_{0}}}{|\nu|}}\sqrt{\frac{k\pi}{2}}\left[ J_{0}(k r) (\mathbf{e} + \overline{\mathbf{e}}) + \textnormal{i}J_{1}(k r)(\mathbf{e} - \overline{\mathbf{e}})\right] + \textnormal{O}(\varepsilon)\nonumber
\end{align}
for all $r\in(0,r_{0})$. In particular, using the explicit forms of $\mathbf{e},\mathbf{f}$ \eqref{eigvec;c:defn} defined in Appendix \ref{app:basis}, the height of the free surface $\eta_{A}(r)$ has the form
\begin{align}
    \eta_{A}(r) &= \varepsilon^{\frac{1}{2}}{\textnormal{sgn}\;(\nu)\frac{2\sqrt{c_{0}}}{m|\nu|}}\sqrt{\frac{k\pi}{2}} J_{0}(k r) + \textnormal{O}(\varepsilon) \nonumber
\end{align}
for all $r\in(0, r_{0})$. Similarly, for $(A_{-},z_{-})(\rho)$ in \eqref{initial:spot;rho0}, where $\rho={\log}\; r $, we can invert the coordinate transformations \eqref{scale:1}, \eqref{z:def}, and \eqref{A0:A;trans} to find the free surface profile in the transition chart,
\begin{align}
    \eta_{A}(r) &= \varepsilon^{\frac{1}{2}}{\textnormal{sgn}\;(\nu)\frac{2\sqrt{c_{0}}}{m|\nu|}}\frac{1}{\sqrt{r}}\cos\left(k r - \frac{\pi}{4}\right)+ \textnormal{O}(\varepsilon), \nonumber
\end{align}
for $r\in[r_{0}, \delta_{0}\varepsilon^{-\frac{1}{2}}]$. Finally, for $(A_{2},z_{2})(s)$ in \eqref{A2:s}, we can invert transformations \eqref{scale:2}, \eqref{z:def}, and \eqref{A0:A;trans} to find the free surface profile in the far-field,
\begin{align}
    \eta_{A}(r) &= \varepsilon^{\frac{1}{2}}{\textnormal{sgn}\;(\nu)\frac{2\sqrt{c_{0}}}{m|\nu|}}\textnormal{e}^{\sqrt{c_{0}}(\delta_{0} - \sqrt{\varepsilon} r)}  \frac{1}{\sqrt{r}}\cos\left(k r - \frac{\pi}{4}\right)+ \textnormal{O}(\varepsilon), \nonumber
\end{align}
for $r\in[\delta_{0}\varepsilon^{-\frac{1}{2}}, \infty)$. This completes the result for spot A.
{
\subsubsection{Tracking Solutions: Spot A for $\nu\approx 0$}\label{s:fold}
As discussed at the beginning of Section \ref{s:spotA}, we expect spot A solutions to {undergo a fold in the region $\nu\approx0$}. In this case, we no longer assume that $\widetilde{a}=\textnormal{O}\left( {\varepsilon^{\frac{1}{4}}}\right)$, and we will adapt the results seen previously for solutions where $ {\varepsilon^{\frac{1}{4}}} \ll |\widetilde{a}|\ll 1$. As seen in \eqref{Wcs;param;delta0}, we can parametrise solutions close to the connecting orbit $\mathbf{q}_{-}(s)$ evaluated at $s=\delta_{0}$ such that,
\begin{align}
\mathcal{W}^{cs}( Q_{-})\cap\Sigma_{0} &= \bigg\{ (A_{-}, z_{-}, \sigma_{-}) = \left(\widetilde{a}\textnormal{e}^{\textnormal{i}Y} + \textnormal{O}( {\varepsilon}), -\sqrt{c_{0}}\delta_{0} + \textnormal{O}(\widetilde{a}^{2} +  {\varepsilon}), \frac{ \varepsilon^{\frac{1}{2}}}{\delta_{0}}\right) \,:\; |\widetilde{a}|< a_{0}, \; | {\varepsilon}|<\varepsilon_{0}\bigg\}.\label{FoldWcs;param;delta0}
\end{align}
Then, we obtain the initial data
\begin{align}
    A_{-}(\rho_{0}) = \widetilde{a}\textnormal{e}^{\textnormal{i}Y} + \textnormal{O}( {\varepsilon}), \qquad z_{-}(\rho_{0}) = -\sqrt{c_{0}}\delta_{0} + \textnormal{O}(\widetilde{a}^{2} +  {\varepsilon}), \qquad \sigma_{-}(\rho_{0}) = \frac{ \varepsilon^{\frac{1}{2}}}{\delta_{0}}, \qquad \varepsilon_{-}(\rho_{0}) = \delta_{0},\label{FoldInital:spota}
\end{align}
for which we need to solve \eqref{amp:a1;b1,-}
for $\rho^{*}\leq\rho\leq\rho_{0}$, where
\begin{align}
    \rho_{0} = {\log}\;\,\frac{\delta_{0}}{ \varepsilon^{\frac{1}{2}}},\qquad \rho^{*} = {\log}\;\,\frac{1}{\sigma_{0}}.\nonumber
\end{align}
In particular, we find
\begin{align}
    \varepsilon_{-}(\rho) = \varepsilon_{-}(\rho_{0})\textnormal{e}^{\rho-\rho_{0}}, \qquad \varepsilon_{-}(\rho) =  \varepsilon^{\frac{1}{2}}\textnormal{e}^{\rho}, \qquad \sigma_{-}(\rho) = \sigma_{0}\textnormal{e}^{\rho^{*}-\rho},\nonumber
\end{align}
for $\rho_{0}\leq\rho\leq0$. Next, we define $A_{-}(\rho) = \widetilde{A}_{-}(\rho)\textnormal{e}^{\rho/2}$, $z_{-}(\rho) = \widetilde{z}_{-}(\rho)\textnormal{e}^{\rho}$
and obtain the system
\begin{align}
    \frac{\textnormal{d}}{\textnormal{d}\rho} \widetilde{A}_{-} &= \widetilde{A}_{-}\left[\widetilde{z}_{-}\textnormal{e}^{\rho} + \textnormal{e}^{2(\rho^{*}-\rho)}\textnormal{O}\left(|\sigma_{0}|^2\right)\right],\label{FoldAz:tilde;syst}\\
    \frac{\textnormal{d}}{\textnormal{d}\rho} \widetilde{z}_{-} &= -\textnormal{e}^{\rho}\widetilde{z}_{-}^{2} + c_{0}|\varepsilon_{-}(\rho_{0})|^2 \textnormal{e}^{\rho} + c_{3}|\widetilde{A}_{-}|^2 + \textnormal{e}^{2(\rho^{*}-\rho)}|\sigma_{0}|^2 \textnormal{O}\left(|\widetilde{A}_{-}|^{4} + |\widetilde{z}_{-}| + |\varepsilon_{-}(\rho_{0})|\right),\nonumber
\end{align}
which we consider with initial conditions \eqref{FoldInital:spota}, written in our new coordinates as
\begin{align}
    \widetilde{A}_{-}(\rho_{0}) = \widetilde{A}_{-}^{0}:= \widetilde{a} {\varepsilon^{\frac{1}{4}}}\delta_{0}^{-\frac{1}{2}}\textnormal{e}^{\textnormal{i}Y} + \textnormal{O}( {\varepsilon}), \qquad \widetilde{z}_{-}(\rho_{0}) = \widetilde{z}_{-}^{0}:= -\sqrt{c_{0}} \varepsilon^{\frac{1}{2}} + \textnormal{O}(\widetilde{a}^{2} +  {\varepsilon}). \label{Foldinitial:spota;tilde}
\end{align}
We write the solution to \eqref{FoldAz:tilde;syst} as the variation of constants equation
\begin{align}
    \widetilde{A}_{-}(\rho) &= \widetilde{A}^{0}_{-} + \int_{\rho_{0}}^{\rho} \widetilde{A}_{-}(u)\left[\widetilde{z}_{-}(u)\textnormal{e}^{u} + \textnormal{e}^{2(\rho^{*}-u)}\textnormal{O}\left(|\sigma_{0}|^2\right)\right] \textnormal{d}u, \label{Foldfixedpoint:spota}\\
    \widetilde{z}_{-}(\rho) &= \widetilde{z}^{0}_{-} + \int_{\rho_{0}}^{\rho} -\textnormal{e}^{u}\widetilde{z}_{-}^{2}(u) + c_{0}|\varepsilon_{-}(\rho_{0})|^2 \textnormal{e}^{u} + c_{3}|\widetilde{A}_{-}(u)|^2 + \textnormal{e}^{2(\rho^{*}-u)}|\sigma_{0}|^2 \textnormal{O}\left(|\widetilde{A}_{-}(u)|^{4} + |\widetilde{z}_{-}(u)| + |\varepsilon_{-}(\rho_{0})|\right) \textnormal{d}u, \nonumber
\end{align}
which, for sufficiently small variables, has a unique solution (see \cite{Sandstede1997Convergence}) of the form,
\begin{align}
    \widetilde{A}_{-}(\rho) &\approx \widetilde{A}_{-}^{0}\left[1 + \textnormal{O}\left(\sigma_{0} + \delta_{0} +  {\varepsilon^{\frac{1}{4}}}\right)\right],\qquad \qquad 
    \widetilde{z}_{-}(\rho) \approx \left(\widetilde{z}_{-}^{0} + \int_{\rho_{0}}^{\rho} c_{3}|\widetilde{A}_{-}^{0}|^{2} \textnormal{d} u\right)\left[1 + \textnormal{O}\left(\sigma_{0} + \delta_{0} +  {\varepsilon^{\frac{1}{4}}}\right)\right],\nonumber
\end{align}
\begin{align}
&A_{-}(\rho) =\widetilde{a}\delta_{0}^{-\frac{1}{2}} {\varepsilon^{\frac{1}{4}}}\textnormal{e}^{\rho/2}\textnormal{e}^{\textnormal{i}Y}(1 + \Delta_{A}),& \quad 
    &z_{-}(\rho) =-\sqrt{c_{0}} \varepsilon^{\frac{1}{2}}\left[1- (\rho-\rho_{0})\frac{c_{3}}{\sqrt{c_{0}}}\frac{\widetilde{a}^{2}}{\delta_{0}}\right]\textnormal{e}^{\rho}(1+\Delta_{A}),&\nonumber\\
&A_{-}(\rho^{*}) =\frac{\widetilde{a}}{\sqrt{\sigma_{0}\delta_{0}}} {\varepsilon^{\frac{1}{4}}}\textnormal{e}^{\textnormal{i}Y}(1 + \Delta_{A}),& \quad 
    &z_{-}(\rho^{*}) =-\sqrt{c_{0}}\frac{ \varepsilon^{\frac{1}{2}}}{\sigma_{0}}\left[1+\left|\log\left( {\varepsilon}\right)\right|\frac{c_{3}}{2\sqrt{c_{0}}}\frac{\widetilde{a}^{2}}{\delta_{0}}\right](1+\Delta_{A}),& \label{Foldinitial:spot;rho0}
\end{align}
with remainder terms
\begin{align}
    \Delta_{A} := \textnormal{O}(\sigma_{0} + \delta_{0} +  {\varepsilon^{\frac{1}{4}}}).\nonumber
\end{align}
To balance the logarithmic growth of the $\widetilde{a}^{2}$ term in $z_{-}(\rho_{0})$, we define
\begin{align}
\widetilde{a}:=\sqrt{\frac{\delta_{0}}{\left|\log\left( {\varepsilon}\right)\right|}}a,\nonumber
\end{align}
so that,
\begin{align}
&A_{-}(\rho^{*}) =\frac{a}{\sqrt{\sigma_{0}}}\left|\log\left( {\varepsilon}\right)\right|^{-\frac{1}{2}} {\varepsilon^{\frac{1}{4}}}\textnormal{e}^{\textnormal{i}Y}(1 + \Delta_{A}),& \quad 
    &z_{-}(\rho^{*}) =-\sqrt{c_{0}}\frac{ \varepsilon^{\frac{1}{2}}}{\sigma_{0}}\left[1+\frac{c_{3}}{2\sqrt{c_{0}}}a^{2}\right](1+\Delta_{A}),& \nonumber
\end{align}
Now, we return to the original $(A,B)$ coordinates by converting \eqref{Foldinitial:spot;rho0} into $(A_{1}, z_{1})$ coordinates and inverting the coordinate change \eqref{scale:1}. That is, we write $A = \sigma_{-}A_{-}$ and $B = \sigma_{-}^{2}A_{-}z_{-}$, and evaluate at $\sigma=\frac{1}{r}$; then, we obtain 
\begin{align}
    A(r_{0}) &=a\varepsilon^{\frac{1}{4}}|\log (\varepsilon)|^{-\frac{1}{2}} r_{0}^{-\frac{1}{2}}\textnormal{e}^{\textnormal{i}Y}(1 + \Delta_{A}), \qquad 
    B(r_{0}) =-\sqrt{c_{0}} a\varepsilon^{\frac{3}{4}} |\log(\varepsilon)|^{-\frac{1}{2}} r_{0}^{-\frac{1}{2}} \left[1+ \frac{c_{3}}{2\sqrt{c_{0}}}a^{2}\right]\textnormal{e}^{\textnormal{i}Y}(1+\Delta_{A}). \label{Foldinitial:spot;AB}
\end{align}
\subsubsection{Matching Core and Far-Field Solutions: Spot A for $\nu\approx 0$}\label{s:fold;match}
Following the same process as in the case when $\nu\neq 0$, we begin by substituting the initial conditions \eqref{Foldinitial:spot;AB} into the far-field parametrisation \eqref{match:farf:folia} following the normal-form transformation \eqref{ab:tilde} and equating these terms with the transformed core parametrisation \eqref{match:core;normal} for each respective coordinate $(\widetilde{a}, \widetilde{b}, a_{1}, a_{-1}, a_{2}, a_{-2}, \dots)(r_{0})$. This is equivalent to finding the zeros of the functional 
\begin{align}
    G:&\left(d_{1}, d_{2}, \widetilde{Y},  \mathbf{c}_{1}, \vec{\mathbf{d}}; a, \varepsilon\right) \mapsto \left( G^{C}_{1}, G^{C}_{2}, G^{H}_{1}, G^{H}_{-1}, G^{H}_{2}, G^{H}_{-2}, \dots \right), \label{FoldG:defn}
\end{align}
where 
\begin{align}
    G^{C}_{1} &= [-\textnormal{i}+\textnormal{O}(r_{0}^{-1})]r_{0} d_{1} + [1+\textnormal{O}(r_{0}^{-1})]d_{2} - \textnormal{e}^{\textnormal{i}\widetilde{Y}} a \varepsilon^{\frac{1}{4}}|\log(\varepsilon)|^{-\frac{1}{2}}(1 + \Delta_{A}) + \mathscr{R}^{C}_{1},\label{FoldG:SpotA;C1}\\
    G^{C}_{2} &= \left[-\textnormal{i} +\textnormal{O}(r_{0}^{-1})\right] d_{1} - \left[\nu + \textnormal{O}(r_{0}^{-\frac{1}{2}})\right] d_{2}^{2} +\textnormal{e}^{\textnormal{i}\widetilde{Y}} a \sqrt{c_{0}} \varepsilon^{\frac{3}{4}} |\log(\varepsilon)|^{-\frac{1}{2}}\left[1+\frac{c_{3}}{2\sqrt{c_{0}}}a^{2}\right](1+\Delta_{A}) + \mathscr{R}^{C}_{2},\label{FoldG:SpotA;C2}\\
    G^{H}_{n} &= \textnormal{e}^{\lambda_{n}r_{0}}r_{0}^{-\frac{1}{2}}\left[\frac{1}{\sqrt{\pi}} + \textnormal{O}(r_{0}^{-1})\right]c_{1,n} + \mathscr{R}^{H}_{1},\label{FoldG:SpotA;Hn}\\
    G^{H}_{-n} &= - \vec{d}_{n}\textnormal{e}^{-\lambda_{n}r_{0}} + \mathscr{R}^{H}_{2},\label{FoldG:SpotA;-Hn}
\end{align}
where $\widetilde{Y}$, $\mathscr{R}^{C}_{j}$, and $\mathscr{R}^{H}_{j}$ are defined as
\begin{align}
    \widetilde{Y} :&= Y + \frac{\pi}{4} + \textnormal{O}(r_{0}^{-2}) + \textnormal{O}_{r_{0}}(|\varepsilon| + |d_{1}| + |d_{2}| + |\mathbf{c}_{1}|_{1}),\label{Fold:Y;tilde}\\
    \mathscr{R}^{C}_{1} :&= \textnormal{O}_{r_{0}}\left(\left[|d_{1}| + |d_{2}| + |\mathbf{c}_{1}|_{1}\right]\left[|\varepsilon| + |d_{1}| + |d_{2}| + |\mathbf{c}_{1}|_{1}\right] + |\vec{\mathbf{d}}|_{1}\left[\varepsilon^{\frac{1}{4}}|\log\varepsilon|^{-\frac{1}{2}} + |\vec{\mathbf{d}}|_{1}\right]\right),\label{Fold:R;C1}\\
    \mathscr{R}^{C}_{2} :&= \textnormal{O}_{r_{0}}\left(\left[|d_{1}| + |d_{2}| + |\mathbf{c}_{1}|_{1}\right]\left[|\varepsilon| + |d_{1}| + |\mathbf{c}_{1}|_{1}\right] + |d_{2}|^{3} + |\vec{\mathbf{d}}|_{1}\left[\varepsilon^{\frac{1}{4}}|\log\varepsilon|^{-\frac{1}{2}} + |\vec{\mathbf{d}}|_{1}\right]\right),\label{Fold:R;C2}\\
    \mathscr{R}^{H}_{1} :&= \textnormal{O}_{r_{0}}\left(\left[|d_{1}| + |d_{2}| + |\mathbf{c}_{1}|_{1}\right]\left[|\varepsilon| + |d_{1}| + |d_{2}| + |\mathbf{c}_{1}|_{1}\right] + \left[\varepsilon^{\frac{1}{4}}|\log\varepsilon|^{-\frac{1}{2}} + |\vec{\mathbf{d}}|_{1}\right]^{2}\right),\label{Fold:R;H1}\\
    \mathscr{R}^{H}_{2} :&= \textnormal{O}_{r_{0}}\left(\left[|d_{1}| + |d_{2}| + |\mathbf{c}_{1}|_{1}\right]\left[|\varepsilon| + |d_{1}| + |d_{2}| + |\mathbf{c}_{1}|_{1}\right] + \varepsilon^{\frac{1}{4}}|\log\varepsilon|^{-\frac{1}{2}}\left[\varepsilon^{\frac{1}{4}}|\log\varepsilon|^{-\frac{1}{2}} + |\vec{\mathbf{d}}|_{1}\right]\right).\label{Fold:R;H2}
\end{align}
By introducing the scaling $(d_{1}, d_{2}^{2}) = \varepsilon^{\frac{1}{2}}|\log(\varepsilon)|^{-1}(\widetilde{d}_{1}, \widetilde{d}_{2}^{2})$ so that terms in \eqref{FoldG:SpotA;C1} and \eqref{FoldG:SpotA;C2} scale with the same order in $\varepsilon$, we find that \eqref{FoldG:SpotA;Hn}, \eqref{FoldG:SpotA;-Hn} implies that
\begin{align}
    c_{1,n} = \textnormal{O}_{r_{0}}\left(|\varepsilon|^{\frac{1}{2}}|\log\varepsilon|^{-1}\right), \qquad \vec{d}_{n}  =\textnormal{O}_{r_{0}}\left(|\varepsilon|^{\frac{1}{2}}|\log\varepsilon|^{-1}\right), \qquad \forall n\in\mathbb{N},\nonumber
\end{align}
via the same reasoning as the case when $\nu\neq0$. Returning to \eqref{FoldG:SpotA;C1} and \eqref{FoldG:SpotA;C2}, we have
\begin{align}
    \widetilde{G}_{1} &= [1+\textnormal{O}(r_{0}^{-1})]\widetilde{d}_{2} - \textnormal{e}^{\textnormal{i}\widetilde{Y}} a (1 + \Delta_{A}) + \textnormal{O}_{r_{0}}(|\varepsilon|^{\frac{1}{4}}|\log\varepsilon|^{-\frac{1}{2}}),\label{FoldMatch:SpotA;C}\\
    \widetilde{G}_{2} &= \left[-\textnormal{i} +\textnormal{O}(r_{0}^{-1})\right] \widetilde{d}_{1} - \left[\nu + \textnormal{O}(r_{0}^{-\frac{1}{2}})\right]\widetilde{d}_{2}^{2}\nonumber\\
    & \qquad \qquad +\textnormal{e}^{\textnormal{i}\widetilde{Y}} a \sqrt{c_{0}} \varepsilon^{\frac{1}{4}} |\log(\varepsilon)|^{\frac{1}{2}}\left[1+\frac{c_{3}}{2\sqrt{c_{0}}}a^{2}\right](1+\Delta_{A}) + \textnormal{O}_{r_{0}}(|\varepsilon|^{\frac{1}{4}}|\log\varepsilon|^{-\frac{1}{2}}).\nonumber
\end{align}
To leading order, we obtain the system
\begin{align}
    \widetilde{G}_{1} &= \widetilde{d}_{2} - \textnormal{e}^{\textnormal{i}\widetilde{Y}} a ,\nonumber\\
    \widetilde{G}_{2} &= -\textnormal{i}  \widetilde{d}_{1} - \nu\widetilde{d}_{2}^{2} +\textnormal{e}^{\textnormal{i}\widetilde{Y}} a \sqrt{c_{0}} \varepsilon^{\frac{1}{4}} |\log(\varepsilon)|^{\frac{1}{2}}\left[1+\frac{c_{3}}{2\sqrt{c_{0}}}a^{2}\right] .\nonumber
\end{align}
We separate into real and imaginary parts; then,
\begin{align}
    \widetilde{d}_{2} &= a \cos(\widetilde{Y}),\nonumber\\ 
    0 &= -a \sin(\widetilde{Y}), \nonumber\\
    \nu \widetilde{d}_{2}^{2} &= a \sqrt{c_{0}}\varepsilon^{\frac{1}{4}} |\log(\varepsilon)|^{\frac{1}{2}} \left[1+\frac{c_{3}}{2\sqrt{c_{0}}}a^{2}\right]\cos(\widetilde{Y}), \nonumber\\
    \widetilde{d}_{1} &= a \sqrt{c_{0}} \varepsilon^{\frac{1}{4}} |\log(\varepsilon)|^{\frac{1}{2}}\left[1+\frac{c_{3}}{2\sqrt{c_{0}}}a^{2}\right]\sin(\widetilde{Y}).\nonumber
\end{align}
We define $\nu= \varepsilon^{\frac{1}{4}}|\log(\varepsilon)|^{\frac{1}{2}}\widetilde{\nu}$; solving the second and fourth equations imply,
\begin{align}
    {\widetilde{Y}= 0 \; \textnormal{or}\; \pi, \qquad \textnormal{and}}\qquad \widetilde{d}_{1}=0,\nonumber
\end{align}
{where} the equations
\begin{align}
    \widetilde{d}_{2}&= a \cos(\widetilde{Y}), \qquad\qquad \qquad \widetilde{\nu} \widetilde{d}_{2}^{2} = a \sqrt{c_{0}} \left[1+\frac{c_{3}}{2\sqrt{c_{0}}}a^{2}\right]\cos(\widetilde{Y}),\nonumber 
\end{align} 
remain to be solved. The choices of $\widetilde{Y}=0$ and $\widetilde{Y}=\pi$ are equivalent to changing the sign of $a$, leaving $\widetilde{d}_{2}$ invariant, and so we will take $\widetilde{Y}=0$ for now. Then, we find a solution at $(\widetilde{d}_{1}, \widetilde{d}_{2}, \widetilde{Y}) = (0, a, 0),$ where $a$ satisfies,
\begin{align}
    c_{3}a^{2} - 2\widetilde{\nu}a+2\sqrt{c_{0}}=0.\label{Fold:Bif;Quad}
\end{align}
Setting $\widetilde{\nu}=0$, we see that $c_{3}<0$ implies
\begin{align}
    &a =\pm \sqrt{\frac{2\sqrt{c_{0}}}{-c_{3}}},& \quad
    &d_{2} = \pm \sqrt{\frac{2\sqrt{c_{0}}}{-c_{3}}}\frac{\varepsilon^{\frac{1}{4}}}{\sqrt{|\log(\varepsilon)|}},&\nonumber
\end{align}
whereas $c_{3}>0$ implies that the solution undergoes a fold bifurcation related to the double root of \eqref{Fold:Bif;Quad}
\begin{align}
    &\widetilde{\nu} =\pm \sqrt{2c_{3}\sqrt{c_{0}}},& \quad &\nu = \pm \sqrt{2c_{3}\sqrt{c_{0}}}\varepsilon^{\frac{1}{4}}\sqrt{|\log(\varepsilon)|},&\nonumber\\
    &a =\frac{\widetilde{\nu}}{c_{3}} = \textnormal{sgn}\;(\nu) \sqrt{\frac{2\sqrt{c_{0}}}{c_{3}}},&\quad
    &d_{2} = \textnormal{sgn}\;(\nu) \sqrt{\frac{2\sqrt{c_{0}}}{c_{3}}}\frac{\varepsilon^{\frac{1}{4}}}{\sqrt{|\log(\varepsilon)|}}.&\nonumber
\end{align}
Hence, for $(\varepsilon, \nu)$ close to zero, spot A solutions undergo a fold bifurcation along the curve $\nu=\pm \varepsilon^{\frac{1}{4}}|\log(\varepsilon)|^{\frac{1}{2}}\left[c+\textnormal{o}(1)\right]$, with constant $c:=\sqrt{2c_{3}\sqrt{c_{0}}}>0$, and the amplitude of the free surface is given by
\begin{align}
    \eta_{A}(r) = \textnormal{sgn}\;(\nu)\; \varepsilon^{\frac{1}{4}}|\log(\varepsilon)|^{-\frac{1}{2}}\frac{2}{m}\sqrt{\frac{2\sqrt{c}_{0}}{c_{3}}} \sqrt{\frac{k\pi}{2}} J_{0}(kr) + \textnormal{o}(1), \qquad \qquad 0\leq r\leq r_{0}.\nonumber
\end{align}
}

\subsection{Spot B}\label{s:spotB}
\begin{figure}
    \centering
    \includegraphics[height=6cm]{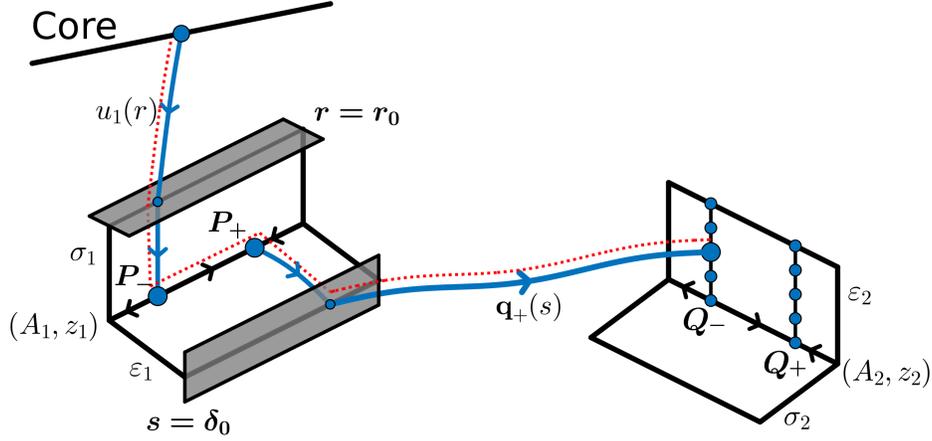}
    \caption{The path followed by the spot B solution; it begins by following $u_{1}(r)$ to the equilibrium $P_{-}$, transporting to $P_{+}$ close to $A_{1}, \sigma_{1}, \varepsilon_{1}=0$, and then tracking the connecting orbit $\mathbf{q}_{+}(s)$, as defined in \eqref{soln:chart}, to the equilibrium $ Q_{-}$.}
    \label{fig:GBUSpotB}
\end{figure}

The construction of the secondary spot, termed spot B in \cite{mccalla2013spots}, proceeds as follows; we follow the centre-stable manifold $\mathcal{W}^{cs}( Q_{-})$ near the heteroclinic orbits defined in \eqref{soln:chart;comp} for $Y=0, \pi$ backward in `spatial time' $\rho$ past the equilibria $P_{+}$ and $P_{-}$, as shown in Figure \ref{fig:GBUSpotB}. This will define the initial conditions required for solutions to decay to zero as $r\to\infty$, parametrising the centre coordinates of the far-field manifold at $r=r_{0}$. These centre coordinates are then substituted into the foliation parametrisation defined in \eqref{match:farf:folia} and the final matching between the core and far-field manifolds will take place at $r=r_{0}$.

\begin{figure}[t]
    \centering
    \includegraphics[height=7cm]{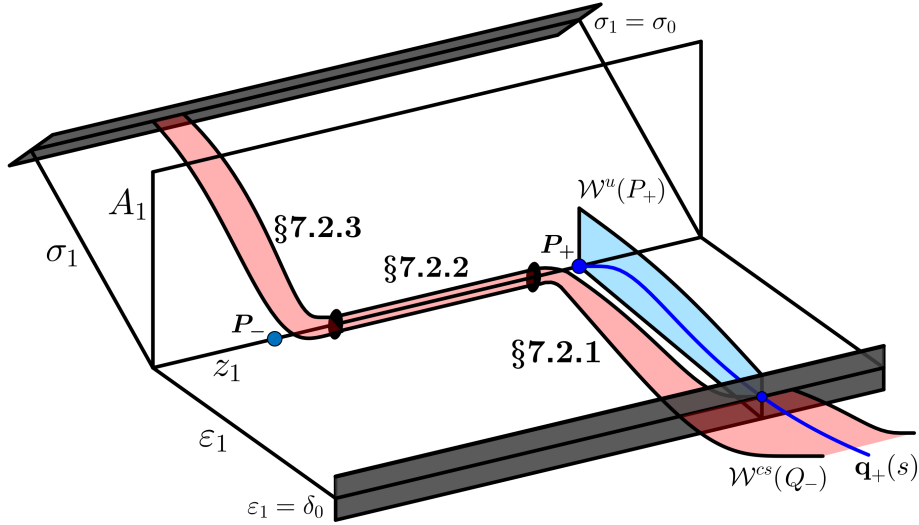}
    \caption{A schematic illustration of the various steps involved in tracking the centre-stable manifold $\mathcal{W}^{cs}( Q_{-})$ in the transition chart. We begin close to the orbit $\mathbf{q}_{+}(s)$ on the section $\Sigma_{0}$; solutions are tracked close to the equilibrium $P_{+}$ in \change{\S\ref{subs:spot;B,1}}. As solutions pass $P_{+}$, they are tracked close to the $z_{1}$-axis to the equilibrium $P_{-}$ in \change{\S\ref{subs:spot;B,2}}, until they are tracked close to $P_{-}$, back to the point $\sigma_{1}=\sigma_{0}$, or equivalently $r=r_{0}$, in \change{\S\ref{subs:spot;B,3}}.}
    \label{fig:SpotB-full}
\end{figure}

The construction of spot B in this paper is almost identical to the two-dimensional Swift-Hohenberg equation, which is completed explicitly in \cite{mccalla2013spots}. Therefore, we will provide an idea of how to construct spot B, and then state the results which have been proven previously in \cite{mccalla2013spots}. The journey in backward $\rho$ taken by the spot B solution through the transition chart can be divided as follows: Section \ref{subs:spot;B,1}, where the solution moves close to $P_{+}$; Section \ref{subs:spot;B,2}, where the solution travels between $P_{+}$ and $P_{-}$; and Section \ref{subs:spot;B,3}, where the solution moves close to $P_{-}$. The decomposition of this journey is illustrated in Figure \ref{fig:SpotB-full}.

{\paragraph{{Formal Scaling Analysis}: Spot B}
For this formal analysis, we have set all parameters equal to 1, other than the coordinates we wish to match. In the rescaling chart, ring solutions follow close to the connecting orbit $\mathbf{q}_{+}(s)$ and so, at the point $\varepsilon_{-}=1$, the manifold $\mathcal{W}^{cs}( Q_{-})$ can be written as 
\begin{align}
    \mathcal{W}^{cs}( Q_{-}): A_{+}(0)= \textnormal{e}^{\textnormal{i}Y},\qquad z_{+}(0)=-\widetilde{a}, \qquad \sigma_{+}(0)=\varepsilon^{\frac{1}{2}}, \qquad \varepsilon_{+}(0)=1.\nonumber
\end{align}
Here, we restrict $\widetilde{a}>0$, and recall that ring solutions lie close to the linear flow of the equilibrium $P_{+}$. That is, they will be subject to the equations
\begin{align}
    \frac{\textnormal{d}}{\textnormal{d}\rho}A_{+}=\frac{3}{2}A_{+}, \qquad \frac{\textnormal{d}}{\textnormal{d}\rho}z_{+}=-z_{+}, \qquad \frac{\textnormal{d}}{\textnormal{d}\rho}\sigma_{+}=-\sigma_{+}, \qquad \frac{\textnormal{d}}{\textnormal{d}\rho}\varepsilon_{+}=\varepsilon_{+}.\nonumber
\end{align}
Hence, evolving our solutions with respect to the linear flow near $P_{+}$, we find
\begin{align}
    A_{+}(\rho)= \textnormal{e}^{\textnormal{i}Y}\textnormal{e}^{\frac{3\rho}{2}},\qquad z_{+}(\rho)=-\widetilde{a}\textnormal{e}^{-\rho}, \qquad \sigma_{+}(\rho)=\varepsilon^{\frac{1}{2}}\textnormal{e}^{-\rho}, \qquad \varepsilon_{+}(\rho)=\textnormal{e}^{\rho},\nonumber
\end{align}
which we evolve backwards in $\rho$, until the point $\rho_{0}:=\log\left(2 \widetilde{a}\right)$, such that $z_{+}(\rho_{0})=-\frac{1}{2}$. Then, at $\rho=\rho_{0}$, our solutions take the form,
\begin{align}
    A_{+}(\rho_{0})= \left(2\widetilde{a}\right)^{\frac{3}{2}}\textnormal{e}^{\textnormal{i}Y},\qquad z_{+}(\rho_{0})=-\frac{1}{2}, \qquad \sigma_{+}(\rho_{0})=\frac{\varepsilon^{\frac{1}{2}}}{2\widetilde{a}}, \qquad \varepsilon_{+}(\rho_{0})=2\widetilde{a}.\nonumber
\end{align}
We take these solutions as initial conditions for the flow near to the equilibrium $P_{-}$; thus, converting to $A_{-}$, $z_{-}$, we have
\begin{align}
    A_{-}(0)= \left(2\widetilde{a}\right)^{\frac{3}{2}}\textnormal{e}^{\textnormal{i}Y},\qquad z_{-}(0)=\frac{1}{2}, \qquad \sigma_{-}(0)=\frac{\varepsilon^{\frac{1}{2}}}{2\widetilde{a}}, \qquad \varepsilon_{-}(0)=2\widetilde{a},\nonumber
\end{align}
which evolve subject to the linear flow
\begin{align}
    \frac{\textnormal{d}}{\textnormal{d}\rho}A_{-}=\frac{1}{2}A_{-}, \qquad \frac{\textnormal{d}}{\textnormal{d}\rho}z_{-}=z_{-}, \qquad \frac{\textnormal{d}}{\textnormal{d}\rho}\sigma_{-}=-\sigma_{-}, \qquad \frac{\textnormal{d}}{\textnormal{d}\rho}\varepsilon_{-}=\varepsilon_{-}.\nonumber
\end{align}
Then, we find
\begin{align}
    A_{-}(\rho)= \left(2\widetilde{a}\right)^{\frac{3}{2}}\textnormal{e}^{\textnormal{i}Y} \textnormal{e}^{\frac{\rho}{2}},\qquad z_{-}(\rho)=\frac{1}{2}\textnormal{e}^{\rho}, \qquad \sigma_{-}(\rho)=\frac{\varepsilon^{\frac{1}{2}}}{2\widetilde{a}}\textnormal{e}^{-\rho}, \qquad \varepsilon_{-}(\rho)=2\widetilde{a}\textnormal{e}^{\rho},\nonumber
\end{align}
which we evolve backwards in $\rho$, until the point $\rho_{1}:=\log\left(\frac{\varepsilon^{\frac{1}{2}}}{2 \widetilde{a}}\right)$, such that $\sigma_{-}(\rho_{1})=1$. Then, at $\rho=\rho_{1}$, our solutions take the form,
\begin{align}
    A_{-}(\rho_{1})= \left(2\widetilde{a}\right) \varepsilon^{\frac{1}{4}}\textnormal{e}^{\textnormal{i}Y} ,\qquad z_{-}(\rho_{1})=\frac{\varepsilon^{\frac{1}{2}}}{4 \widetilde{a}}, \qquad \sigma_{-}(\rho_{1})=1, \qquad \varepsilon_{-}(\rho_{1})=\varepsilon^{\frac{1}{2}}.\nonumber
\end{align}
We convert these solutions into $A,B$ coordinates via the transformation $A =\sigma_{-}A_{-}$, $B=\sigma_{-}^{2}A_{-}z_{-}$; so,
\begin{align}
    A(r_{0})= \left(2\widetilde{a}\right) \varepsilon^{\frac{1}{4}}\textnormal{e}^{\textnormal{i}Y},\qquad B(r_{0})= \frac{1}{2}\varepsilon^{\frac{3}{4}}\textnormal{e}^{\textnormal{i}Y},\nonumber
\end{align}
Our core parametrisation can be written as 
\begin{align}
    \widetilde{\mathcal{W}}^{cu}_{-}: A(r_{0})= -\textnormal{i}d_{1} + d_{2},\qquad B(r_{0})=-\textnormal{i}d_{1} - d_{2}^{2},\nonumber
\end{align}
where we have set all other parameters equal to 1, and absorbed the complex phase of \eqref{match:core;normal} into the definition of $Y$. Matching our far-field solutions to the core parametrisation and taking real and imaginary parts, we find
\begin{align}
    d_{1} &= -\left(2\widetilde{a}\right) \varepsilon^{\frac{1}{4}}\sin(Y),\qquad
    d_{2} = \left(2\widetilde{a}\right) \varepsilon^{\frac{1}{4}}\cos(Y),\qquad 
    d_{1} = -\frac{1}{2}\varepsilon^{\frac{3}{4}}\sin(Y),\qquad
    d_{2}^{2} = -\frac{1}{2}\varepsilon^{\frac{3}{4}}\cos(Y),\nonumber
\end{align}
and we see that 
\begin{align}
    Y= \pi, \qquad d_{1} = 0, \qquad d_{2}=-2\widetilde{a}\varepsilon^{\frac{1}{4}}, \qquad \widetilde{a} = \frac{\varepsilon^{\frac{1}{8}}}{2\sqrt{2}}.\nonumber
\end{align}
Hence, we expect to find $\widetilde{a}=\textnormal{O}\left(\varepsilon^{\frac{1}{8}}\right)$, $d_{2} = \textnormal{O}\left(\varepsilon^{\frac{3}{8}}\right)$, and $d_{1}$ to be of higher order in $\varepsilon$.
}

\subsubsection{Dynamics near $P_{+}$}\label{subs:spot;B,1}

We begin by tracking the centre-stable manifold $\mathcal{W}^{cs}( Q_{-})$ for $ {\varepsilon}>0$ backwards in $\rho$ as it passes near to the equilibrium $P_{+}$. We consider the transition variables $(A_{+},z_{+},\sigma_{+},\varepsilon_{+})$ such that $P_{+}$ corresponds to the origin. It is convenient to perform a change of coordinates in order to remove some of the higher order terms from the equations \eqref{amp:a1;b1,+}. Therefore, we apply the following Lemma,
\begin{Lemma}{\cite{mccalla2013spots}}
There is a smooth coordinate change of the form
\begin{align}
    \widehat{z}_{+} = z_{+} + h_{+}(A_{+}, \sigma_{+}, \varepsilon_{+}), \qquad h_{+}(A_{+}, \sigma_{+}, \varepsilon_{+}) = \textnormal{O}(|A_{+}|^2 + |\sigma_{+}|^2 + |\varepsilon_{+}|^2)\label{coord:spotB:+}
\end{align}
that transforms \eqref{amp:a1;b1,+} near the origin into
\begin{align}
    &\frac{\textnormal{d}}{\textnormal{d} \rho}A_{+} = A_{+}\left[3/2 + \textnormal{O}(|A_{+}| + |\widehat{z}_{+}| + |\sigma_{+}| + |\varepsilon_{+}|)\right],
    &\qquad
    &\frac{\textnormal{d}}{\textnormal{d} \rho}\sigma_{+} = -\sigma_{+},
    &\nonumber\\
    &\frac{\textnormal{d}}{\textnormal{d} \rho}\widehat{z}_{+} = -\widehat{z}_{+}\left[1 + \textnormal{O}(|A_{+}| + |\widehat{z}_{+}| + |\sigma_{+}| + |\varepsilon_{+}|)\right],
    &\qquad
    &\frac{\textnormal{d}}{\textnormal{d} \rho}\varepsilon_{+} =\varepsilon_{+}.
    &\label{amp:SpotB;+transf}
\end{align}
\end{Lemma}
For a proof of this result, see \cite[Lemma 3.1]{mccalla2013spots}. To find an expression for the manifold $\mathcal{W}^{cs}( Q_{-})$ near $P_{+}$, we recall that the solution \eqref{soln:chart} forms a transverse intersection of $\mathcal{W}^{cs}( Q_{-})$ and $\mathcal{W}^{u}(P_{+})$. Then, 
\begin{Lemma}{\cite{mccalla2013spots}}
For each sufficiently small $\delta_{0}>0$, there are constants $a_{0}, \varepsilon_{0}>0$ such that the following is true. Define $\Sigma_{0}:=\{\varepsilon_{+}=\delta_{0}\}$; then
\begin{align}
\mathcal{W}^{cs}( Q_{-})\cap\Sigma_{0} &= \left\{ (A_{+}, \widehat{z}_{+}, \sigma_{+}) = \left(-\textnormal{e}^{\textnormal{i}Y}[\eta_{0}(\delta_{0}) + \textnormal{O}(\widetilde{a})] + \textnormal{O}( {\varepsilon}), -\widetilde{a} + \textnormal{O}( {\varepsilon}), \frac{ \varepsilon^{\frac{1}{2}}}{\delta_{0}}\right) \,:\; |\widetilde{a}|< a_{0}, \;  {\varepsilon}<\varepsilon_{0}\right\}.\nonumber
\end{align}
where $\eta_{0}(\delta_{0}) = q_{0}\delta_{0}^{3/2} (1+\textnormal{O}(\delta_{0}))$ is smooth, $q_{0}>0$ is the constant defined in \eqref{q:defn}, and $Y\in\mathbb{R}$ is arbitrary.
\end{Lemma}
This result is proven in \cite[Lemma 3.2]{mccalla2013spots}. We begin at {$\rho=\rho_{0}$ for initial data in $\Sigma_{0}$ and need to track solutions until $\rho=\rho^{*}$, where 
\begin{align}
    \rho_{0} = {\log}\; \frac{\delta_{0}}{ \varepsilon^{\frac{1}{2}}}, \qquad \rho^{*} = {\log}\; \frac{1}{\sigma_{0}},\nonumber
\end{align}
and we consequently only solve \eqref{amp:SpotB;+transf} for $\rho_{*}\leq\rho\leq\rho_{0}$}. We choose a constant $\delta_{1}>0$ and track part of the centre-stable manifold $\mathcal{W}^{cs}_{+}( Q_{-})$ backwards in $\rho$ under the evolution of \eqref{amp:SpotB;+transf} from $\Sigma_{0}$ to $\Re \widehat{z}_{+} = -\delta_{1}$. Following {our earlier formal analysis}, we take $\widetilde{a}=\textnormal{O}( {\varepsilon^{\frac{1}{8}}})$. {For the rigorous proofs found in \cite{mccalla2013spots}, $\widetilde{a}$ is defined to be $\widetilde{a}:=\varepsilon^{\frac{1}{8}}a$, where $a\in(\varepsilon^{\kappa},\varepsilon^{-\kappa})$ for some $0<\kappa\ll1$. This means that $a$ is neither too small nor too large to change the asymptotic behaviour of our solutions as $\varepsilon\to0$ and so, for our formal approach, we will just take $a$ to be a fixed constant.} Then, we have the following;
\begin{Lemma}{\cite{mccalla2013spots}}
For fixed constants $0 < \delta_{0},\delta_{1}\ll 1$, there is an $\varepsilon_{0}>0$ such that solutions of \eqref{amp:SpotB;+transf} associated with initial data of the form
\begin{align}
    (A_{+}, \widehat{z}_{+}, \sigma_{+}, \varepsilon_{+})(\rho_{0}) = \left(-\textnormal{e}^{\textnormal{i}Y}\eta_{0}(\delta_{0}) + \textnormal{O}( {\varepsilon^{\frac{1}{8}}}), -\frac{a {\varepsilon^{\frac{1}{8}}}}{\delta_{0}} + \textnormal{O}( {\varepsilon}), \frac{ \varepsilon^{\frac{1}{2}}}{\delta_{0}}, \delta_{0}\right),\nonumber
\end{align}
in $\mathcal{W}^{cs}( Q_{-})\cap\Sigma_{0}$ with $ {\varepsilon}\in(0,\varepsilon_{0})$ {evaluated at the point $\rho=\rho_{1}$, where}
{\begin{align}
    \rho_{1} = {\log}\; \frac{a  {\varepsilon^{-\frac{3}{8}}}}{\delta_{1}}\geq \rho^{*},\nonumber
\end{align}}
{is given by}
\begin{align}
    &A_{+}(\rho_{1}) = -\left(\frac{a }{\delta_{0}\delta_{1}}\right)^{\frac{3}{2}} {\varepsilon^{\frac{3}{16}}}\textnormal{e}^{\textnormal{i}Y} \eta_{0}(\delta_{0})(1+\textnormal{O}(\delta_{0} + \delta_{1} +  {\varepsilon^{\frac{1}{8}}})),
    &\qquad 
    &\sigma_{+}(\rho_{1}) = \frac{\delta_{1}  {\varepsilon^{\frac{3}{8}}}}{a},
    &\nonumber\\
    &\widehat{z}_{+}(\rho_{1}) = -\delta_{1}(1+\textnormal{O}(\delta_{0} + \delta_{1} +  {\varepsilon^{\frac{1}{8}}})),
    &\qquad
    &\varepsilon_{+}(\rho_{1}) = \frac{a  {\varepsilon^{\frac{1}{8}}}}{\delta_{1}}.
    &\nonumber
\end{align}
\end{Lemma}
This is proven by constructing fixed-point equations and utilising the particular structure of the nonlinearities; see \cite[Lemma 3.3]{mccalla2013spots}. After inverting the coordinate change \eqref{coord:spotB:+} and returning to the original transition chart variables, we obtain
\begin{align}
    &A_{1}^{0}:= A_{1}(\rho_{1}) = -a^{\frac{3}{2}} {\varepsilon^{\frac{3}{16}}}\eta_{1} \textnormal{e}^{\textnormal{i}Y},
    &\qquad
    &\sigma_{1}^{0}:= \sigma_{1}(\rho_{1}) =\frac{\delta_{1} {\varepsilon^{\frac{3}{8}}}}{a},
    &\nonumber\\ 
    &z_{1}^{0}:= z_{1}(\rho_{1}) = \frac{1}{2} - \delta_{1}(1+ \textnormal{O}(\delta_{0} + \delta_{1} +  {\varepsilon^{\frac{1}{8}}})),
    &\qquad
    &\varepsilon_{1}^{0}:= \varepsilon_{1}(\rho_{1}) = \frac{a {\varepsilon^{\frac{1}{8}}}}{\delta_{1}},
    &\label{initial:spotB;rho1}
    \end{align}
with
\begin{align}
    \eta_{1}:= \eta_{0}(\delta_{0})(\delta_{0}\delta_{1})^{-\frac{3}{2}}(1+ \textnormal{O}(\delta_{0} + \delta_{1} +  {\varepsilon^{\frac{1}{8}}})) = q_{0}\delta_{1}^{-\frac{3}{2}}(1+ \textnormal{O}(\delta_{0} + \delta_{1} +  {\varepsilon^{\frac{1}{8}}}))>0. \label{defn:eta1}
\end{align}
Next we track this manifold to a neighbourhood of the equilibrium $P_{-}$.
\subsubsection{Dynamics between $P_{+}$ and $P_{-}$}\label{subs:spot;B,2}
We now track the centre-stable manifold $\mathcal{W}^{cs}( Q_{-})$ as it passes between the two equilibria $P_{+}$ and $P_{-}$. We note that in \eqref{initial:spotB;rho1}, $|z_{1}^{0}| \gg |A_{1}^{0}|, |\sigma_{1}^{0}|, |\varepsilon_{1}^{0}|$, and so we track solutions close to the $z_{1}$-axis, i.e. when $(A_{1}, z_{1}, \sigma_{1}, \varepsilon_{1}) = (0, z_{1}, 0, 0)$. To this end, we will introduce a small constant $\delta_{2}>0$ and integrate the transition-chart system \eqref{amp:a1;b1} given by
\begin{align}
    &\frac{\textnormal{d}}{\textnormal{d} \rho}A_{1} = A_{1}\left[1+z_{1} + \textnormal{O}(|\sigma_{1}|^{2})\right],
    &\qquad 
    &\frac{\textnormal{d}}{\textnormal{d} \rho}\sigma_{1} = -\sigma_{1},
    &\nonumber\\
    &\frac{\textnormal{d}}{\textnormal{d} \rho}z_{1} = c_{0}\varepsilon_{1}^{2} + \frac{1}{4} + c_{3}|A_{1}|^2 - z_{1}^{2} + |\sigma_{1}|^{2}\textnormal{O}\left(|A_{1}|^{4} + \left|z_{1} + \frac{1}{2}\right| + |\varepsilon_{1}|^{2}\right),
    &\qquad
    &\frac{\textnormal{d}}{\textnormal{d} \rho}\varepsilon_{1} =\varepsilon_{1},
    &\nonumber
\end{align}
with initial conditions given by \eqref{initial:spotB;rho1} backward in time until $z_{1}$ is approximately $-\frac{1}{2}+\delta_{2}$. That is, we set
{\begin{align}
    \rho_{2}:= {\log}\; \frac{a {\varepsilon^{-\frac{3}{8}}}\delta_{2}}{(1-\delta_{1})(1-\delta_{2})},\nonumber
\end{align}
and integrate \eqref{amp:a1;b1} from $\rho=\rho_{1}$ to $\rho=\rho_{2}$}. We initially set $(A_{1},\sigma_{1}, \varepsilon_{1})=0$ so that \eqref{amp:a1;b1} with the initial condition \eqref{initial:spotB;rho1} for $z_{1}$ becomes the complex differential equation
\begin{align}
    \frac{\textnormal{d}}{\textnormal{d}\rho}z_{1} = -z_{1}^{2} + \frac{1}{4}, \qquad z_{1}({\rho_{1}}) = z_{1}^{0} = \frac{1}{2} - \delta_{1}(1+ \textnormal{O}(\delta_{0} + \delta_{1} +  {\varepsilon^{\frac{1}{8}}})),\label{p-p+:zsyst}
\end{align}
whose solution $z_{1}^{*}(\rho)$ evaluated at $\rho=\rho_{2}$ is given by
\begin{align}
    z_{1}^{*}(\rho_{2}) = -\frac{1}{2} + \delta_{2}(1+ \textnormal{O}(\delta_{0} + \delta_{1} +  {\varepsilon^{\frac{1}{8}}})).\nonumber
\end{align}
Next, we expand the time-$\rho_{2}$ map of \eqref{amp:a1;b1} with initial condition $(A_{1}^{0}, z_{1}^{0}, \sigma_{1}^{0}, \varepsilon_{1}^{0})$ at $\rho=\change{\rho_{1}}$ around $(0,z_{1}^{0} ,0,0)$ and obtain
\begin{align}
    \begin{pmatrix} A_{1}(\rho_{2})\\ z_{1}(\rho_{2})\\ \sigma_{1}(\rho_{2})\\ \varepsilon_{1}(\rho_{2})\end{pmatrix} = \begin{pmatrix} A_{1}^{0} \eta_{2} (1 + \textnormal{O}(|A_{1}^{0}| + |\sigma_{1}^{0}| + |\varepsilon_{1}^{0}|))\\ z_{1}^{*}(\rho_{2}) + \textnormal{O}(|A_{1}^{0}| + |\sigma_{1}^{0}| + |\varepsilon_{1}^{0}|)\\
    \frac{ {\varepsilon^{\frac{3}{8}}}}{a \delta_{2}}(1+ \textnormal{O}(\delta_{1} + \delta_{2}))\\
    a\delta_{2} {\varepsilon^{\frac{1}{8}}}(1+ \textnormal{O}(\delta_{1} + \delta_{2}))\end{pmatrix},\nonumber
\end{align}
where $\eta_{2}$ is a constant given by $\eta_{2} = a_{1}(\rho_{2})$, and $a_{1}$ is the solution to the linear equation
\begin{align}
    \frac{\textnormal{d}}{\textnormal{d} \rho}a_{1} = (1+ z_{1}^{*}(\rho))a_{1}, \qquad a_{1}({\rho_{1}}) = 1.\nonumber
\end{align}
This equation can be solved explicitly, and we find
\begin{align}
    \eta_{2} = \delta_{1}^{\frac{3}{2}}\delta_{2}^{\frac{1}{2}}(1+ \textnormal{O}(\delta_{1} + \delta_{2})).\nonumber
\end{align}
Substituting the initial conditions \eqref{initial:spotB;rho1}, we obtain
\begin{align}
    \begin{pmatrix} A_{1}(\rho_{2})\\ z_{1}(\rho_{2})\\ \sigma_{1}(\rho_{2})\\ \varepsilon_{1}(\rho_{2})\end{pmatrix} = \begin{pmatrix} -a^{\frac{3}{2}} {\varepsilon^{\frac{3}{16}}} \eta_{3} \textnormal{e}^{\textnormal{i}Y}\\ -\frac{1}{2} + \delta_{2}(1 + \textnormal{O}(\delta_{0} + \delta_{1} +  {\varepsilon^{\frac{1}{8}}}))\\
    \frac{ {\varepsilon^{\frac{3}{8}}}}{a \delta_{2}}(1+ \textnormal{O}(\delta_{1} + \delta_{2}))\\
    a\delta_{2} {\varepsilon^{\frac{1}{8}}}(1+ \textnormal{O}(\delta_{1} + \delta_{2}))\end{pmatrix},\label{initial:spotB;rho2}
\end{align}
where $\eta_{3}$ is given by
\begin{align}
    \eta_{3}:= \eta_{1}\eta_{2}(1+\textnormal{O}( {\varepsilon^{\frac{1}{8}}})) \overset{\eqref{defn:eta1}}{=} q_{0} \delta_{2}^{\frac{1}{2}}(1 + \textnormal{O}(\delta_{0} + \delta_{1} + \delta_{2} +  {\varepsilon^{\frac{1}{8}}})). \nonumber
\end{align}
\subsubsection{Dynamics near $P_{-}$}\label{subs:spot;B,3}
The final step in tracking the spot B solution through the transition chart is to solve the system \eqref{amp:a1;b1} with initial conditions given by \eqref{initial:spotB;rho2} for {$\rho_{3}\leq\rho\leq\rho_{2}$,}
\begin{align}
    {\rho_{3} = \rho^{*} = {\log}\; \frac{1}{\sigma_{0}},} \label{defn:rho3}
\end{align}
near the equilibrium $P_{-}$. Using the variable $z_{-} = z_{1} + 1/2$, we therefore need to solve the system \eqref{amp:a1;b1,-} given by
\begin{align}
    &\frac{\textnormal{d}}{\textnormal{d} \rho}A_{-} = A_{-}\left[1/2 +z_{-} + \textnormal{O}(|\sigma_{-}|^{2})\right],
    &\qquad
    &\frac{\textnormal{d}}{\textnormal{d} \rho}\sigma_{-} = -\sigma_{-},
    &\nonumber\\
    &\frac{\textnormal{d}}{\textnormal{d} \rho}z_{-} = c_{0}\varepsilon_{-}^{2} + z_{-} - z_{-}^{2} + c_{3}|A_{-}|^2 + |\sigma_{-}|^{2}\textnormal{O}\left(|A_{-}|^{4} + \left|z_{-}\right| + |\varepsilon_{-}|\right),
    &\qquad
    &\frac{\textnormal{d}}{\textnormal{d} \rho}\varepsilon_{-} =\varepsilon_{-},
    &\nonumber
\end{align}
with initial conditions
\begin{align}
    &A_{-}({\rho_{2}})  = -a^{\frac{3}{2}} {\varepsilon^{\frac{3}{16}}} \eta_{3} \textnormal{e}^{\textnormal{i}Y},
    &\qquad
    &\sigma_{-}({\rho_{2}}) = \frac{ {\varepsilon^{\frac{3}{8}}}}{a \delta_{2}}(1+ \textnormal{O}(\delta_{1} + \delta_{2})),
    &\nonumber\\
    &z_{-}({\rho_{2}}) = \delta_{2}(1 + \textnormal{O}(\delta_{0} + \delta_{1} +  {\varepsilon^{\frac{1}{8}}})),
    &\qquad
    &\varepsilon_{-}({\rho_{2}}) = a\delta_{2} {\varepsilon^{\frac{1}{8}}}(1+ \textnormal{O}(\delta_{1} + \delta_{2})),
    &\label{initial:spotB;rho2,-}
\end{align}
from $\rho={\rho_{2}}$ to $\rho=\rho_{3}$. Using the same techniques as before, we obtain the following result 
\begin{Lemma}{\cite{mccalla2013spots}}
For all fixed small constants $\sigma_{0}, \delta_{j} >0$ with $j=0,1,2$, there is an $\varepsilon_{0}$ such that the solution of \eqref{amp:a1;b1,-} with initial condition \eqref{initial:spotB;rho2,-}, evaluated at $\rho=\rho_{3}$ with $\rho_{3}$ from \eqref{defn:rho3}, is given by
\begin{align}
    &A_{-}(\rho_{3}) =-\frac{a {\varepsilon^{\frac{3}{8}}}q_{0} \textnormal{e}^{\textnormal{i}Y}}{\sqrt{\sigma_{0}}}\left(1 + \textnormal{O}(\sigma_{0} + \delta_{0} + \delta_{1} + \delta_{2} +  {\varepsilon^{\frac{1}{8}}})\right),
    &\qquad
    &\sigma_{-}(\rho_{3}) = \sigma_{0} = \frac{1}{r_{0}} ,
    &\nonumber\\
    &z_{-}(\rho_{3}) = \frac{ {\varepsilon^{\frac{3}{8}}}}{a\sigma_{0}}\left(1 + \textnormal{O}(\sigma_{0} + \delta_{0} + \delta_{1} + \delta_{2} +  {\varepsilon^{\frac{1}{8}}})\right),
    &\qquad
    &\varepsilon_{-}(\rho_{3}) = \frac{ \varepsilon^{\frac{1}{2}}}{\sigma_{0}} =  \varepsilon^{\frac{1}{2}} r_{0},
    &\nonumber
\end{align}
uniformly in $ {\varepsilon}\in(0,\varepsilon_{0})$, where $q_{0}>0$ is the constant given in \eqref{q:defn}.
\label{Result:spotB;P-}\end{Lemma}
This result is proven in \cite[Lemma 3.4]{mccalla2013spots}. Converting back into $(A,B)$ coordinates, by $A=\sigma_{-}A_{-}$ and $B=\sigma_{-}^{2}A_{-}z_{-}$ {and evaluating at $\sigma=\frac{1}{r}$} , we obtain
\begin{align}
    A(r_{0}) &= -\frac{a {\varepsilon^{\frac{3}{8}}}q_{0} \textnormal{e}^{\textnormal{i}Y}}{\sqrt{r_{0}}} (1 + \Delta_{B}), \qquad B(r_{0}) = -\frac{ {\varepsilon^{\frac{3}{4}}}q_{0} \textnormal{e}^{\textnormal{i}Y}}{\sqrt{r_{0}}} (1 + \Delta_{B}), \label{initial:SpotB;AB}
\end{align}
where
\begin{align}
    \Delta_{B} := \textnormal{O}(\sigma_{0} + \delta_{0} + \delta_{1} + \delta_{2} +  {\varepsilon^{\frac{1}{8}}}).\nonumber
\end{align}
\subsubsection{Matching Core and Far-Field Solutions: Spot B}\label{subs:spot;B,match}
We are now fully equipped to find a nontrivial solution contained in the intersection of the core manifold $\widetilde{\mathcal{W}}^{cu}_{-}\change{(\varepsilon)}$ and the far-field manifold $\change{\widetilde{\mathcal{W}}^{s}_{+}}$. We begin by substituting the initial conditions \eqref{initial:SpotB;AB} into the far-field parametrisation \eqref{match:farf:folia} following the normal-form transformation \eqref{ab:tilde}
\begin{align}
    \widetilde{a}(r_{0}) &= A(r_{0}) +\textnormal{O}_{r_{0}}\left(|\vec{\mathbf{d}}|_{1}\left[|\varepsilon| + |A| + |B| + |\vec{\mathbf{d}}|_{1}\right]\right),\nonumber\\
    &= -\textnormal{e}^{\textnormal{i}Y}a q_{0} r_{0}^{-\frac{1}{2}}\varepsilon^{\frac{3}{8}}(1 + \Delta_{B}) +\textnormal{O}_{r_{0}}\left(|\vec{\mathbf{d}}|_{1}\left[{|\varepsilon|^{\frac{3}{8}}} + |\vec{\mathbf{d}}|_{1}\right]\right),\nonumber\\
    \widetilde{b}(r_{0}) &= B(r_{0}) +\textnormal{O}_{r_{0}}\left(|\vec{\mathbf{d}}|_{1}\left[|\varepsilon| + |A| + |B| + |\vec{\mathbf{d}}|_{1}\right]\right),\nonumber\\
    &= -\textnormal{e}^{\textnormal{i}Y}q_{0}r_{0}^{-\frac{1}{2}} \varepsilon^{\frac{3}{4}} (1 + \Delta_{B}) +\textnormal{O}_{r_{0}}\left(|\vec{\mathbf{d}}|_{1}\left[{|\varepsilon|^{\frac{3}{8}} }+ |\vec{\mathbf{d}}|_{1}\right]\right),\nonumber\\
    a_{n}(r_{0}) &= \textnormal{O}_{r_{0}}\left(\left[|A| + |B| + |\vec{\mathbf{d}}|_{1}\right]\left[|\varepsilon| + |A| + |B| + |\vec{\mathbf{d}}|_{1}\right]\right),\nonumber\\
    &= \textnormal{O}_{r_{0}}\left(\left[{|\varepsilon|^{\frac{3}{8}} }+ |\vec{\mathbf{d}}|_{1}\right]^{2}\right),\nonumber\\
    a_{-n}(r_{0}) &= \vec{d}_{n}\textnormal{e}^{-\lambda_{n}r_{0}}  + \textnormal{O}_{r_{0}}\left(\left[|A| + |B| \right]\left[|\varepsilon| + |A| + |B| + |\vec{\mathbf{d}}|_{1}\right]\right) ,\nonumber\\
     &= \vec{d}_{n}\textnormal{e}^{-\lambda_{n}r_{0}}  + \textnormal{O}_{r_{0}}\left({|\varepsilon|^{\frac{3}{8}}} \left[{ |\varepsilon|^{\frac{3}{8}} } + |\vec{\mathbf{d}}|_{1}\right]\right). \nonumber
\end{align}
We also recall the transformed core parametrisation \eqref{match:core;normal}
\begin{align}
    \widetilde{a}(r_{0}) &= \textnormal{e}^{-\textnormal{i}\left(\frac{\pi}{4} + \textnormal{O}(r_{0}^{-2}) + \textnormal{O}_{r_{0}}(|\varepsilon| + |\mathbf{d}|_{1})\right)}r_{0}^{-\frac{1}{2}}\left([-\textnormal{i}+\textnormal{O}(r_{0}^{-1})]r_{0} d_{1} + [1+\textnormal{O}(r_{0}^{-1})]d_{2}\right) + \textnormal{O}_{r_{0}}\left(|\mathbf{d}|_{1}\left[|\varepsilon| + |\mathbf{d}|_{1}\right]\right),\nonumber\\
    \widetilde{b}(r_{0}) &= \textnormal{e}^{-\textnormal{i}\left(\frac{\pi}{4} + \textnormal{O}(r_{0}^{-2}) + \textnormal{O}_{r_{0}}(|\varepsilon| + |\mathbf{d}|_{1})\right)}r_{0}^{-\frac{1}{2}}\left(\left[-\textnormal{i} +\textnormal{O}(r_{0}^{-1})\right]d_{1} - \left[\nu + \textnormal{O}(r_{0}^{-\frac{1}{2}})\right] d_{2}^{2}\right)+ \textnormal{O}_{r_{0}}\left(|\mathbf{d}|_{1}\left[|\varepsilon| + |\mathbf{d}_{2}|_{1}\right] + |d_{2}|^{3}\right),\nonumber\\
    a_{n}(r_{0}) &= \textnormal{e}^{\lambda_{n}r_{0}}r_{0}^{-\frac{1}{2}}\left[\frac{1}{\sqrt{\pi}} + \textnormal{O}(r_{0}^{-1})\right]\widetilde{c}_{1,n} + \textnormal{O}_{r_{0}}\left(|\mathbf{d}|_{1}\left[|\varepsilon| + |\mathbf{d}|_{1}\right]\right),\nonumber\\
    a_{-n}(r_{0}) &= \textnormal{O}_{r_{0}}\left(|\mathbf{d}|_{1}\left[|\varepsilon| + |\mathbf{d}|_{1}\right]\right),\nonumber
\end{align}
where $|\mathbf{d}|_{1} = |d_{1}| + |d_{2}| + |\mathbf{c}_{1}|_{1}$, and $|\mathbf{d}_{2}|_{1} = |d_{1}| + |\mathbf{c}_{1}|_{1}$. Setting these parametrisations equivalent to each other for each respective coordinate $(\widetilde{a}, \widetilde{b}, a_{1}, a_{-1}, a_{2}, a_{-2}, \dots)(r_{0})$ is the same as finding the zeros of the functional 
\begin{align}
    G:&\left(d_{1}, d_{2}, \widetilde{Y},  \mathbf{c}_{1}, \vec{\mathbf{d}}; a, \varepsilon\right) \mapsto \left( G^{C}_{1}, G^{C}_{2}, G^{H}_{1}, G^{H}_{-1}, G^{H}_{2}, G^{H}_{-2}, \dots \right), \nonumber
\end{align}
where 
\begin{align}
    G^{C}_{1} &= [-\textnormal{i}+\textnormal{O}(r_{0}^{-1})]r_{0} d_{1} + [1+\textnormal{O}(r_{0}^{-1})]d_{2} + \textnormal{e}^{\textnormal{i}\widetilde{Y}} a q_{0} \varepsilon^{\frac{3}{8}}(1 + \Delta_{B}) + \mathscr{R}^{C}_{1},\label{G:SpotB;C1}\\
    G^{C}_{2} &= \left[-\textnormal{i} +\textnormal{O}(r_{0}^{-1})\right] d_{1} - \left[\nu + \textnormal{O}(r_{0}^{-\frac{1}{2}})\right] d_{2}^{2} +\textnormal{e}^{\textnormal{i}\widetilde{Y}} q_{0} \varepsilon^{\frac{3}{4}} (1+\Delta_{B}) + \mathscr{R}^{C}_{2},\label{G:SpotB;C2}\\
    G^{H}_{n} &= \textnormal{e}^{\lambda_{n}r_{0}}r_{0}^{-\frac{1}{2}}\left[\frac{1}{\sqrt{\pi}} + \textnormal{O}(r_{0}^{-1})\right]c_{1,n} + \mathscr{R}^{H}_{1},\label{G:SpotB;Hn}\\
    G^{H}_{-n} &= - \vec{d}_{n}\textnormal{e}^{-\lambda_{n}r_{0}} + \mathscr{R}^{H}_{2},\label{G:SpotB;-Hn}
\end{align}
and we have defined
\begin{align}
    \widetilde{Y} :&= Y + \frac{\pi}{4} + \textnormal{O}(r_{0}^{-2}) + \textnormal{O}_{r_{0}}(|\varepsilon| + |d_{1}| + |d_{2}| + |\mathbf{c}_{1}|_{1}),\nonumber\\
    \mathscr{R}^{C}_{1} :&= \textnormal{O}_{r_{0}}\left(\left[|d_{1}| + |d_{2}| + |\mathbf{c}_{1}|_{1}\right]\left[|\varepsilon| + |d_{1}| + |d_{2}| + |\mathbf{c}_{1}|_{1}\right] + |\vec{\mathbf{d}}|_{1}\left[{|\varepsilon|^{\frac{3}{8}}} + |\vec{\mathbf{d}}|_{1}\right]\right),\nonumber\\
    \mathscr{R}^{C}_{2} :&= \textnormal{O}_{r_{0}}\left(\left[|d_{1}| + |d_{2}| + |\mathbf{c}_{1}|_{1}\right]\left[|\varepsilon| + |d_{1}| + |\mathbf{c}_{1}|_{1}\right] + |d_{2}|^{3} + |\vec{\mathbf{d}}|_{1}\left[{|\varepsilon|^{\frac{3}{8}}} + |\vec{\mathbf{d}}|_{1}\right]\right),\nonumber\\
    \mathscr{R}^{H}_{1} :&= \textnormal{O}_{r_{0}}\left(\left[|d_{1}| + |d_{2}| + |\mathbf{c}_{1}|_{1}\right]\left[|\varepsilon| + |d_{1}| + |d_{2}| + |\mathbf{c}_{1}|_{1}\right] + \left[{|\varepsilon|^{\frac{3}{8}}} + |\vec{\mathbf{d}}|_{1}\right]^{2}\right),\nonumber\\
    \mathscr{R}^{H}_{2} :&= \textnormal{O}_{r_{0}}\left(\left[|d_{1}| + |d_{2}| + |\mathbf{c}_{1}|_{1}\right]\left[|\varepsilon| + |d_{1}| + |d_{2}| + |\mathbf{c}_{1}|_{1}\right] + {|\varepsilon|^{\frac{3}{8}}}\left[{|\varepsilon|^{\frac{3}{8}}} + |\vec{\mathbf{d}}|_{1}\right]\right).\nonumber
\end{align}
We introduce the scaling $(d_{1}, d_{2}) = (\varepsilon^{\frac{3}{4}} \widetilde{d}_{1}, \varepsilon^{\frac{3}{8}} \widetilde{d}_{2})$ so that terms in \eqref{G:SpotB;C1} and \eqref{G:SpotB;C2} scale with the same order in $\varepsilon$. Initially setting $\varepsilon=0$, we investigate \eqref{G:SpotB;Hn} and \eqref{G:SpotB;-Hn}
\begin{align}
    F_{n} &= \textnormal{e}^{\lambda_{n}r_{0}}r_{0}^{-\frac{1}{2}}\left[\frac{1}{\sqrt{\pi}} + \textnormal{O}(r_{0}^{-1})\right]c_{1,n} + \textnormal{O}_{r_{0}}\left( |\mathbf{c}_{1}|^{2}_{1} + |\vec{\mathbf{d}}|_{1}^{2}\right),\nonumber\\
    F_{-n} &= - \vec{d}_{n}\textnormal{e}^{-\lambda_{n}r_{0}} + \textnormal{O}_{r_{0}}\left( |\mathbf{c}_{1}|^{2}_{1}\right).\nonumber
\end{align}
and, defining the functional $F:(\widetilde{\mathbf{c}}_{1}, \vec{\mathbf{d}})\mapsto \left(F_{j}\right)_{j\in\mathbb{Z}\backslash\{0\}}$, it clear that $F(\mathbf{0},\mathbf{0})=\mathbf{0}$, where ${\mathbf{0} = (0,\dots)}$. Furthermore, the {Jacobian} $DF(\mathbf{0},\mathbf{0})$ is invertible, and so we can solve \eqref{G:SpotB;Hn} and \eqref{G:SpotB;-Hn} for all values of $n\in\mathbb{N}$ uniquely for sufficiently small $0<\varepsilon\ll1$. Matching orders of $\varepsilon$ in \eqref{G:SpotB;Hn} and \eqref{G:SpotB;-Hn}, we find that 
\begin{align}
    c_{1,n} = \textnormal{O}_{r_{0}}\left({|\varepsilon|^{\frac{3}{4}}}\right), \qquad \vec{d}_{n}  =\textnormal{O}_{r_{0}}\left({|\varepsilon|^{\frac{3}{4}}}\right), \qquad \forall n\in\mathbb{N}.\nonumber
\end{align}
Returning to \eqref{G:SpotB;C1} and \eqref{G:SpotB;C2}, we have
\begin{align}
    \widetilde{G}_{1} &= [1+\textnormal{O}(r_{0}^{-1})]\widetilde{d}_{2} + \textnormal{e}^{\textnormal{i}\widetilde{Y}} a q_{0}(1 + \Delta_{B}) + [-\textnormal{i}+\textnormal{O}(r_{0}^{-1})]\varepsilon^{\frac{3}{8}}r_{0}\widetilde{d}_{1} + \varepsilon^{-\frac{3}{8}}\mathscr{R}^{C}_{1},\label{Match:SpotB;C}\\
    \widetilde{G}_{2} &= \left[-\textnormal{i} +\textnormal{O}(r_{0}^{-1})\right] \widetilde{d}_{1} - \left[\nu + \textnormal{O}(r_{0}^{-\frac{1}{2}})\right]\widetilde{d}_{2}^{2} +\textnormal{e}^{\textnormal{i}\widetilde{Y}} q_{0}(1+\Delta_{B}) + \varepsilon^{-\frac{3}{4}}\mathscr{R}^{C}_{2},\nonumber\\
    \intertext{where}
    \mathscr{R}^{C}_{1} &= |\varepsilon|^{\frac{3}{4}}\textnormal{O}_{r_{0}}\left(\left[|\varepsilon|^{\frac{3}{8}}|\widetilde{d}_{1}| + |\widetilde{d}_{2}| + {|\varepsilon|^{\frac{3}{8}}}\right]^{2} + {|\varepsilon|^{\frac{3}{8}}}\right),\nonumber\\
    \mathscr{R}^{C}_{2} &=  |\varepsilon|^{\frac{3}{4}}\textnormal{O}_{r_{0}}\left(\left[|\varepsilon|^{\frac{3}{8}}|\widetilde{d}_{1}| + |\widetilde{d}_{2}| + {|\varepsilon|^{\frac{3}{8}}}\right]\left[|\varepsilon|^{\frac{3}{8}}|\widetilde{d}_{1}| + {|\varepsilon|^{\frac{3}{8}}}\right] + |\varepsilon|^{\frac{3}{8}}|\widetilde{d}_{2}|^{3} + {|\varepsilon|^{\frac{3}{8}}}\right).\nonumber
\end{align}
Initially setting $\varepsilon=0$, we obtain the system
\begin{align}
    \widetilde{G}_{1} &= [1+\Delta_{B}]\widetilde{d}_{2} + \textnormal{e}^{\textnormal{i}\widetilde{Y}} a q_{0} (1 + \Delta_{B}),\label{Match:SpotB;C,eps0}\\
    \widetilde{G}_{2} &= \left[-\textnormal{i} +\Delta_{B}\right] \widetilde{d}_{1} - \left[\nu + \Delta_{B}\right]\widetilde{d}_{2}^{2} +\textnormal{e}^{\textnormal{i}\widetilde{Y}} q_{0}(1+\Delta_{B}).\nonumber
\end{align}
We formally set $\Delta_{B}=0$ and separate \eqref{Match:SpotB;C,eps0} into real and imaginary parts: this is equivalent to finding zeros of the functional
\begin{align}
    \widetilde{G}(\widetilde{d}_{1}, \widetilde{d}_{2}, \widetilde{Y}, a) = \begin{pmatrix} \widetilde{d}_{2} + a q_{0}\cos(\widetilde{Y}) \\ a q_{0} \sin(\widetilde{Y}) \\ -\nu \widetilde{d}_{2}^{2} + q_{0}\cos(\widetilde{Y}) \\ - \widetilde{d}_{1} + q_{0}\sin(\widetilde{Y})
    \end{pmatrix}.\nonumber
\end{align}
The vector
\begin{align}
    \left(\widetilde{d}_{1}^{*}, \widetilde{d}_{2}^{*}, \widetilde{Y}^{*}, a^{*}\right) = \left(0, \change{-\textnormal{sgn}\;(\nu)\sqrt{\frac{q_{0}}{|\nu|}}, \arccos[\textnormal{sgn}\;(\nu)], \frac{1}{\sqrt{q_{0} |\nu|}}}\right),\nonumber
\end{align}
is a root of $\widetilde{G}$ with Jacobian
\begin{align}
    D\widetilde{G}\left(\widetilde{d}_{1}^{*}, \widetilde{d}_{2}^{*}, \widetilde{Y}^{*}, a^{*}\right) = \begin{pmatrix} 0 & 1 & 0 & {\textnormal{sgn}\;(\nu)\;q_{0}} \\ 0 & 0 & \change{\sqrt{\frac{q_{0}}{|\nu|}}} & 0 \\ 0 & \change{-2\textnormal{sgn}\;(\nu)\;\sqrt{q_{0} |\nu|}} & 0 & 0 \\ -1 & 0 & {\textnormal{sgn}\;(\nu)\;q_{0}} & 0\end{pmatrix}.\nonumber
\end{align}
Since $q_{0}>0$, the Jacobian is invertible, and we can therefore solve \eqref{Match:SpotB;C,eps0} uniquely for all sufficiently small $\Delta$, that is, for $r_{0}$ large enough and $\delta_{0}$ small enough, and subsequently \eqref{Match:SpotB;C} for all $0<\varepsilon\ll1$. Reversing the scaling for $d$, we find that 
\begin{align}
    d_{1} &= \varepsilon^{\frac{3}{4}}\textnormal{O}(r_{0}^{-1 } + \delta_{0} + \delta_{1} + \delta_{2} + \varepsilon^{\frac{1}{8}}),\label{SpotB:d1}\\
    d_{2} &= \change{-\varepsilon^{\frac{3}{8}}\textnormal{sgn}\;\left(\nu\right)\sqrt{\frac{q_{0}}{|\nu|}}}\left(1 + \textnormal{O}(r_{0}^{-1} + \delta_{0} + \delta_{1} + \delta_{2} + \varepsilon^{\frac{1}{8}})\right).\label{SpotB:d2}
\end{align}
Hence, we have found a spot B solution. We recall that our solution $\mathbf{u}(r,y)$ takes the form
\begin{align}
    \mathbf{u}(r,y) &= a(r)\mathbf{e}(y) + b(r)\mathbf{f}(y) + \overline{a}(r)\overline{\mathbf{e}}(y) + \overline{b}(r)\overline{\mathbf{f}}(y) + \Sigma_{n=1}^{\infty} \left\{a_{n}(r)\mathbf{e}_{n}(y) + a_{-n}(r)\mathbf{e}_{-n}(y)\right\}, \nonumber
\end{align}
and 
\begin{align}
    \begin{pmatrix}a \\ b \end{pmatrix}(r) & = \sum_{i=1}^{4} \widetilde{d}_{i}\mathbf{V}_{i}(r), \qquad \qquad \begin{pmatrix}a_{n} \\ a_{-n} \end{pmatrix}(r)  = \sum_{i=1}^{2} \widetilde{c}_{i,n}\mathbf{W}_{i,n}(r).\nonumber
\end{align}
Substituting \eqref{SpotB:d1} and \eqref{SpotB:d2} into this form, we can write the spot B solution $\mathbf{u}_{B}$ as
\begin{align}
    \change{\mathbf{u}_{B}}(r,y) &= \change{-\varepsilon^{\frac{3}{8}}\textnormal{sgn}\;\left(\nu\right)\sqrt{\frac{q_{0}}{|\nu|}}}\sqrt{\frac{k\pi}{2}}\left[ J_{0}(k r) (\mathbf{e} + \overline{\mathbf{e}}) + \textnormal{i}J_{1}(k r)(\mathbf{e} - \overline{\mathbf{e}})\right] + \textnormal{O}(\varepsilon^{\frac{1}{2}}),\nonumber
\end{align}
for all $r\in(0,r_{0})$, where $\mathbf{u}_{B}(r, \change{y})$ decays to zero exponentially as $r\to\infty$. In particular, the height of the free surface $\eta_{B}(r)$ has the form
\begin{align}
   \change{\eta_{B}}(r) &= \change{-\varepsilon^{\frac{3}{8}}\textnormal{sgn}\;\left(\nu\right)\sqrt{\frac{4 q_{0}}{m^{2}|\nu|}}}\sqrt{\frac{k\pi}{2}} J_{0}(k r) + \textnormal{O}(\varepsilon^{\frac{1}{2}}), \nonumber
\end{align}
for all $r\in(0, r_{0})$. Following Lemma \ref{Result:spotB;P-}, we can define write the linear flow near $P_{-}$ as
\begin{align}
    A_{-}(\rho) = -a q_{0} \varepsilon^{\frac{3}{8}} \textnormal{e}^{\textnormal{i}Y} \textnormal{e}^{\frac{\rho}{2}}, \qquad z_{-}(\rho) = a^{-1} \varepsilon^{\frac{3}{8}} \textnormal{e}^{\rho},\nonumber
\end{align}
where $\rho={\log}\; r$, for $r\in\left[r_{0}, \frac{a \varepsilon^{-\frac{3}{8}}\delta_{2}}{(1-\delta_{1})(1-\delta_{2})}\right]$. For the flow between $P_{-}$ and $P_{+}$ near the $z_{1}$-axis we can explicitly solve \eqref{p-p+:zsyst} to find
\begin{align}
    A_{-}(\rho) = -q_{0} \varepsilon^{\frac{3}{4}} \textnormal{e}^{\textnormal{i}Y} \textnormal{e}^{\frac{3\rho}{2}}(1-\delta_{1}) z_{-}(\rho)^{-1}, \qquad z_{-}(\rho) = \left( 1 + \frac{a \varepsilon^{-\frac{3}{8}}}{(1-\delta_{1})} \textnormal{e}^{-\rho}\right)^{-1},\nonumber
\end{align}
where $\rho={\log}\; r$, for $r\in\left[\frac{a \varepsilon^{-\frac{3}{8}}\delta_{2}}{(1-\delta_{1})(1-\delta_{2})}, \frac{a\varepsilon^{-\frac{3}{8}}}{\delta_{1}}\right]$. Near the equilibrium $P_{+}$, we can define solutions by the linear flow, such that we find, to leading order,
\begin{align}
    A_{-}(\rho) = -q_{0} \varepsilon^{\frac{3}{8}} \textnormal{e}^{\textnormal{i}Y} \textnormal{e}^{3\frac{\rho}{2}}, \qquad z_{-}(\rho) = 1- a \varepsilon^{-\frac{3}{8}} \textnormal{e}^{-\rho},\nonumber
\end{align}
where $\rho={\log}\; r$, for $r\in\left[ \frac{a\varepsilon^{-\frac{3}{8}}}{\delta_{1}}, \delta_{0}\varepsilon^{-\frac{1}{2}}\right]$. Finally, for $(A_{2},z_{2})(s)$ in \eqref{soln:chart;comp}, we can apply the transformation \eqref{scale:rel}, such that
\begin{align}
    A_{-}(r) = - r q(\varepsilon^{\frac{1}{2}} r) \varepsilon^{\frac{1}{2}} \textnormal{e}^{\textnormal{i}Y}, \qquad z_{-}(r) = - r \frac{p(\varepsilon^{\frac{1}{2}} r)}{q(\varepsilon^{\frac{1}{2}} r)} \varepsilon^{\frac{1}{2}},\nonumber
\end{align}
where $p(\varepsilon^{\frac{1}{2}} r)$ is defined as,
\begin{align}
 p\left(s\right):= q'(s) + \frac{1}{2 s}q(s),\label{pr:defn}
\end{align}
and so, from \eqref{q:defn}, we see that
\begin{align}
    p(s) &= \left\{\begin{array}{cc}
         q_{0}s^{-1/2} + \textnormal{O}(s^{1/2}),& s\to0,  \\
         -\sqrt{c_{0}}(q_{+} + \textnormal{O}(\textnormal{e}^{-\sqrt{c_{0}}s}))\frac{\textnormal{e}^{-\sqrt{c_{0}}s}}{\sqrt{s}},& s\to\infty.
    \end{array}\right. \nonumber
\end{align}
Inverting the transformations \eqref{scale:2}, \eqref{z:def}, and \eqref{A0:A;trans}, we find the free surface profile in the far-field in each transition region, which completes the result for spot B.
\subsection{Rings}\label{s:rings}
\begin{figure}
    \centering
    \includegraphics[height=6cm]{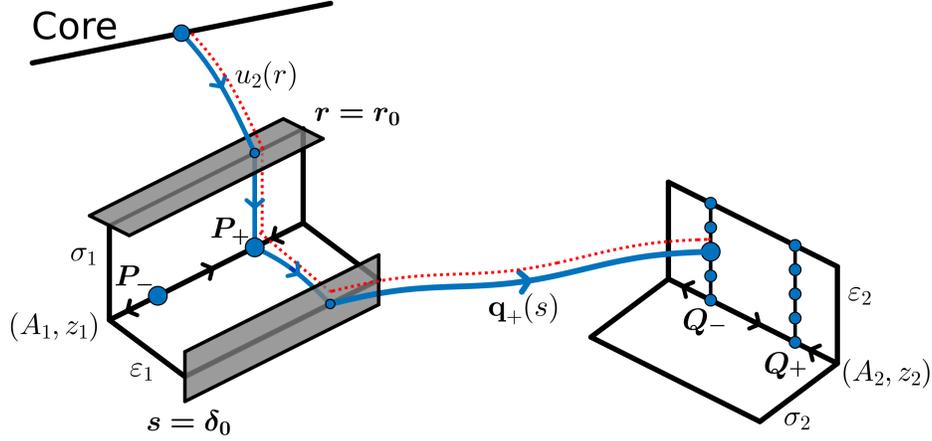}
    \caption{The path followed by the rings solution; it begins by following $u_{2}(r)$ to the equilibrium $P_{+}$, and then tracking the connecting orbit $\mathbf{q}_{+}(s)$, as defined in \eqref{soln:chart}, to the equilibrium $ Q_{-}$.}
    \label{fig:GBURings}
\end{figure}
The construction of up-ring and down-ring solutions proceeds as follows; we follow the centre-stable manifold $\mathcal{W}^{cs}( Q_{-})$ near the heteroclinic orbits defined in \eqref{soln:chart;comp} for $Y=0, \pi$ backward in `spatial time' $\rho$ past the equilibria $P_{+}$, as shown in Figure \ref{fig:GBURings}. This will define the initial conditions required for solutions to decay to zero as $r\to\infty$, parametrising the centre coordinates of the far-field manifold at $r=r_{0}$. These centre coordinates are then substituted into the foliation parametrisation defined in \eqref{match:farf:folia} and the final matching between the core and far-field manifolds will take place at $r=r_{0}$.
{\paragraph{{Formal Scaling Analysis}: Rings}
For this formal analysis, we have set all parameters equal to 1, other than the coordinates we wish to match. In the rescaling chart, ring solutions follow close to the connecting orbit $\mathbf{q}_{+}(s)$ and so, at the point $\varepsilon_{-}=1$, the manifold $\mathcal{W}^{cs}( Q_{-})$ can be written as 
\begin{align}
    \mathcal{W}^{cs}( Q_{-}): A_{+}(0)= \textnormal{e}^{\textnormal{i}Y},\qquad z_{+}(0)=-\widetilde{a}, \qquad \sigma_{+}(0)=\varepsilon^{\frac{1}{2}}, \qquad \varepsilon_{+}(0)=1.\nonumber
\end{align}
We recall that ring solutions lie close to the linear flow of the equilibrium $P_{+}$, and so will be subject to the equations
\begin{align}
    \frac{\textnormal{d}}{\textnormal{d}\rho}A_{+}=\frac{3}{2}A_{+}, \qquad \frac{\textnormal{d}}{\textnormal{d}\rho}z_{+}=-z_{+}, \qquad \frac{\textnormal{d}}{\textnormal{d}\rho}\sigma_{+}=-\sigma_{+}, \qquad \frac{\textnormal{d}}{\textnormal{d}\rho}\varepsilon_{+}=\varepsilon_{+}.\nonumber
\end{align}
Hence, evolving our solutions with respect to the linear flow near $P_{+}$, we find
\begin{align}
    A_{+}(\rho)= \textnormal{e}^{\textnormal{i}Y}\textnormal{e}^{\frac{3\rho}{2}},\qquad z_{+}(\rho)=-\widetilde{a}\textnormal{e}^{-\rho}, \qquad \sigma_{+}(\rho)=\varepsilon^{\frac{1}{2}}\textnormal{e}^{-\rho}, \qquad \varepsilon_{+}(\rho)=\textnormal{e}^{\rho},\nonumber
\end{align}
which we evolve backwards in $\rho$, until the point $\rho^{*}:=\log\left(\varepsilon^{\frac{1}{2}}\right)$, such that $\sigma_{+}(\rho^{*})=1$. Then, at $\rho=\rho^{*}$, our solutions take the form,
\begin{align}
    A_{+}(\rho^{*})= \varepsilon^{\frac{3}{4}}\textnormal{e}^{\textnormal{i}Y},\qquad z_{+}(\rho^{*})=-\widetilde{a}\varepsilon^{-\frac{1}{2}}, \qquad \sigma_{+}(\rho^{*})=1, \qquad \varepsilon_{+}(\rho^{*})=\varepsilon^{\frac{1}{2}},\nonumber
\end{align}
which we convert back into $A,B$ coordinates via the transformation $A =\sigma_{+}A_{+}$, $B=\sigma_{+}^{2}A_{+}\left[1 + z_{+}\right]$; so,
\begin{align}
    A(r_{0})= \varepsilon^{\frac{3}{4}}\textnormal{e}^{\textnormal{i}Y},\qquad B(r_{0})=\varepsilon^{\frac{3}{4}}\left[1-\widetilde{a}\varepsilon^{-\frac{1}{2}}\right]\textnormal{e}^{\textnormal{i}Y},\nonumber
\end{align}
Our core parametrisation can be written as 
\begin{align}
    \widetilde{\mathcal{W}}^{cu}_{-}: A(r_{0})= -\textnormal{i}d_{1} + d_{2},\qquad B(r_{0})=-\textnormal{i}d_{1} - d_{2}^{2},\nonumber
\end{align}
where we have set all other parameters equal to 1, and absorbed the complex phase of \eqref{match:core;normal} into the definition of $Y$. Matching our far-field solutions to the core parametrisation and taking real and imaginary parts, we find
\begin{align}
    &d_{1} = -\varepsilon^{\frac{3}{4}}\sin(Y),&\qquad
    &d_{1} = -\varepsilon^{\frac{3}{4}}\left[1-\widetilde{a}\varepsilon^{-\frac{1}{2}}\right]\sin(Y),&\nonumber\\
    &d_{2} = \varepsilon^{\frac{3}{4}}\cos(Y),&\qquad 
    &d_{2}^{2} = -\varepsilon^{\frac{3}{4}}\left[1-\widetilde{a}\varepsilon^{-\frac{1}{2}}\right]\cos(Y),&\nonumber
\end{align}
and so, we see that $d_{2}=0$, and
\begin{align}
    Y= \frac{(2\pm1)\pi}{2}, \qquad d_{1}=\pm\varepsilon^{\frac{3}{4}}, \qquad \widetilde{a}\varepsilon^{-\frac{1}{2}} =0.\nonumber
\end{align}
Hence, we expect to find $\widetilde{a}=\textnormal{O}\left(\varepsilon^{\frac{1}{2}}\right)$, $d_{1} = \textnormal{O}\left(\varepsilon^{\frac{3}{4}}\right)$, and $d_{2}$ to be of higher order in $\varepsilon$.
\subsubsection{Tracking Solutions: Rings}\label{subs:rings;track}}
We wish to track the centre-stable manifold $\mathcal{W}^{cs}( Q_{-})$ for $ {\varepsilon}>0$ in backward time towards the equilibrium $P_{+}$. We consider the equation in transition variables $(A_{+},z_{+},\sigma_{+},\varepsilon_{+})$, centred around the point $P_{+}$. The resulting system is \eqref{amp:a1;b1,+},
\begin{align}
    &\frac{\textnormal{d}}{\textnormal{d} \rho}A_{+} = A_{+}\left[3/2 +z_{+} + \textnormal{O}(|\sigma_{+}|^{2})\right],
    &\qquad
    &\frac{\textnormal{d}}{\textnormal{d} \rho}\sigma_{+} = -\sigma_{+},
    &\nonumber\\
    &\frac{\textnormal{d}}{\textnormal{d} \rho}z_{+} = c_{0}\varepsilon_{+}^{2} - z_{+} - z_{+}^{2} + c_{3}|A_{+}|^2 + \textnormal{O}\left(|\sigma_{+}|^{2}\right),
    &\qquad
    &\frac{\textnormal{d}}{\textnormal{d} \rho}\varepsilon_{+} =\varepsilon_{+}.
    &\nonumber
    \end{align}
As proven in \cite[Lemma 3.1]{mccalla2013spots}, there exists a coordinate change
\begin{align}
    \widehat{z}_{+} = z_{+} + h_{+}(A_{+}, \sigma_{+}, \varepsilon_{+}), \qquad h_{+}(A_{+}, \sigma_{+}, \varepsilon_{+}) = \textnormal{O}(|A_{+}|^2 + |\sigma_{+}|^2 + |\varepsilon_{+}|^2),\label{coord:rings}
\end{align}
that transforms \eqref{amp:a1;b1,+} near the origin into
\begin{align}
    &\frac{\textnormal{d}}{\textnormal{d} \rho}A_{+} = A_{+}\left[3/2 + \textnormal{O}(|A_{+}| + |\widehat{z}_{+}| + |\sigma_{+}| + |\varepsilon_{+}|)\right],
    &\qquad 
    &\frac{\textnormal{d}}{\textnormal{d} \rho}\sigma_{+} = -\sigma_{+},
    &\nonumber\\
    &\frac{\textnormal{d}}{\textnormal{d} \rho}\widehat{z}_{+} = -\widehat{z}_{+}\left[1 + \textnormal{O}(|A_{+}| + |\widehat{z}_{+}| + |\sigma_{+}| + |\varepsilon_{+}|)\right],
    &\qquad
    &\frac{\textnormal{d}}{\textnormal{d} \rho}\varepsilon_{+} =\varepsilon_{+}.
    &\label{amp:rings;transf}
    \end{align}
We define the section $\Sigma_{0}:=\left\{(A_{+}, \widehat{z}_{+}, \sigma_{+}, \varepsilon_{+}) \,: \;\varepsilon_{+}=\delta_{0}\right\}$ and investigate the parametrisation of $\mathcal{W}^{cs}( Q_{-})\cap \Sigma_{0}$. We know that there exists a heteroclinic orbit \eqref{soln:chart} that forms a transverse intersection between $\mathcal{W}^{u}(P_{+})$ and $\mathcal{W}^{cs}( Q_{-})$ when $ {\varepsilon}=0$. Therefore, defining $\widetilde{a}$ as the perturbation from \eqref{soln:chart} in the $\widehat{z}_{+}$ direction, we can paramatrise $\mathcal{W}^{cs}( Q_{-})\cap\Sigma_{0}$, when $\sigma_{+}=0$, as
\begin{align}
\left.\mathcal{W}^{cs}( Q_{-})\cap\Sigma_{0}\right|_{\sigma_{+}=0} &= \left\{ (A_{+}, \widehat{z}_{+}, \sigma_{+}) = (\textnormal{e}^{\textnormal{i}Y}(\eta_{0}(\delta_{0}) + \textnormal{O}(\widetilde{a})), -\widetilde{a}, 0) \,:\; |\widetilde{a}|< a_{0}\right\},\nonumber
\end{align}
for each small $\delta_{0}>0$ and constant $a_{0}$, where $\eta_{0}(\delta_{0}) = q_{0}\delta_{0}^{3/2}(1+\textnormal{O}(\delta_{0}))$, $q_{0}$ is defined in \eqref{q:defn}, and we have used the $S^{1}$-symmetry in the normal-form to introduce the parameter $Y\in\mathbb{R}$. \\

We know that the solution \eqref{soln:chart} lies in the intersection of $\mathcal{W}^{u}(P_{+})$ and $\mathcal{W}^{cs}( Q_{-}[0])$, and that $\varepsilon_{2}=\sigma_{+}=0$ is an invariant subspace. Therefore, for each small $\delta_{0}>0$, there are constants $a_{0}, \varepsilon_{0}>0$ such that
\begin{align}
\mathcal{W}^{cs}( Q_{-})\cap\Sigma_{0} &= \left\{ (A_{+}, \widehat{z}_{+}, \sigma_{+}) = \left(\textnormal{e}^{\textnormal{i}Y}(\eta_{0}(\delta_{0}) + \textnormal{O}(\widetilde{a})) + \textnormal{O}( {\varepsilon}), -\widetilde{a} + \textnormal{O}( {\varepsilon}), \frac{ \varepsilon^{\frac{1}{2}}}{\delta_{0}}\right) \,:\; |\widetilde{a}|< a_{0}, \;  {\varepsilon}<\varepsilon_{0}\right\}.\nonumber
\end{align}
 We note that adding the parameter $ {\varepsilon}$ yields additional $\textnormal{O}( {\varepsilon})$ terms, using that the remainder terms in the rescaling chart \eqref{amp:a2;z2,eps} are of order $\textnormal{O}(|\varepsilon_{2}|^{2})$. We set $\widetilde{a}= \varepsilon^{\frac{1}{2}} a$ and obtain the initial data
{\begin{align}
    A_{+}(\rho_{0}) = \textnormal{e}^{\textnormal{i}Y}\eta_{0}(\delta_{0})  + \textnormal{O}( \varepsilon^{\frac{1}{2}}), \qquad \widehat{z}_{+}(\rho_{0}) = -a \varepsilon^{\frac{1}{2}} + \textnormal{O}( {\varepsilon}), \qquad \sigma_{+}(\rho_{0}) = \frac{ \varepsilon^{\frac{1}{2}}}{\delta_{0}}, \qquad \varepsilon_{+}(\rho_{0}) = \delta_{0},\nonumber
\end{align}}
for which we need to solve \eqref{amp:rings;transf},
\begin{align}
    &\frac{\textnormal{d}}{\textnormal{d} \rho}A_{+} = A_{+}\left[3/2 + \textnormal{O}(|A_{+}| + |\widehat{z}_{+}| + |\sigma_{+}| + |\varepsilon_{+}|)\right],
    &\qquad 
    &\frac{\textnormal{d}}{\textnormal{d} \rho}\sigma_{+} = -\sigma_{+},
    &\nonumber\\
    &\frac{\textnormal{d}}{\textnormal{d} \rho}\widehat{z}_{+} = -\widehat{z}_{+}\left[1 + \textnormal{O}(|A_{+}| + |\widehat{z}_{+}| + |\sigma_{+}| + |\varepsilon_{+}|)\right],
    &\qquad
    &\frac{\textnormal{d}}{\textnormal{d} \rho}\varepsilon_{+} =\varepsilon_{+}.
    &\nonumber
    \end{align}
for{ $\rho^{*}\leq\rho\leq\rho_{0}$, where }
{\begin{align}
    \rho_{0} = {\log}\;\,\frac{\delta_{0}}{ \varepsilon^{\frac{1}{2}}},\qquad  \rho^{*} = {\log}\;\,\frac{1}{\sigma_{0}}.\nonumber
\end{align}}
In particular, we find
{\begin{align}
    \varepsilon_{+}(\rho) = \varepsilon_{+}(\rho_{0})\textnormal{e}^{\rho-\rho_{0}}, \qquad \varepsilon_{+}(\rho) =  \varepsilon^{\frac{1}{2}} \textnormal{e}^{\rho}, \qquad \sigma_{+}(\rho) = \sigma_{0}\textnormal{e}^{\rho^{*}-\rho},\nonumber
\end{align}}
for {$\rho^{*}\leq\rho\leq\rho_{0}$}. Next, we need to solve
\begin{align}
     &\frac{\textnormal{d}}{\textnormal{d} \rho}A_{+} = A_{+}\left[3/2 + \textnormal{O}(|A_{+}| + |\widehat{z}_{+}| + {\sigma_{0}\textnormal{e}^{\rho^{*}-\rho} +  \varepsilon^{\frac{1}{2}}\textnormal{e}^{\rho})}\right],& \qquad \; &A_{+}({\rho_{0}}) = \textnormal{e}^{\textnormal{i}Y}\eta_{0}(\delta_{0}) + \textnormal{O}( \varepsilon^{\frac{1}{2}}),&\nonumber\\
    &\frac{\textnormal{d}}{\textnormal{d} \rho}\widehat{z}_{+} = -\widehat{z}_{+}\left[1 + \textnormal{O}(|A_{+}| + |\widehat{z}_{+}| + {\sigma_{0}\textnormal{e}^{\rho^{*}-\rho} +  \varepsilon^{\frac{1}{2}}\textnormal{e}^{\rho})}\right],& \qquad& \widehat{z}_{+}({\rho_{0}}) = -a \varepsilon^{\frac{1}{2}} + \textnormal{O}( {\varepsilon}),& \nonumber
\end{align}
which, using the structure of the nonlinearity (see \cite{Sandstede1997Convergence}), has a unique solution given by \begin{align}
     A_{+}(\rho) &=\textnormal{e}^{\textnormal{i}Y}\eta_{0}(\delta_{0}){ {\varepsilon^{\frac{3}{4}}}\delta_{0}^{-\frac{3}{2}}}\textnormal{e}^{3\rho/2}(1 + \Delta_{R}), \qquad 
    \widehat{z}_{+}(\rho) =-a{\delta_{0}}\textnormal{e}^{-\rho}(1+\Delta_{R}),\label{soln:rings}
\end{align}
with
\begin{align}
    \Delta_{R} := \textnormal{O}(\sigma_{0} + \delta_{0} +  {\varepsilon^{\frac{1}{4}}}),\nonumber
\end{align}
for $\delta_{0}$, $\sigma_{0}$, $ {\varepsilon}$ small. Inverting the coordinate transformation \eqref{coord:rings}, we obtain
\begin{align}
    A_{+}({\rho^{*}}) &=q_{0} { {\varepsilon^{\frac{3}{4}}} \sigma_{0}^{-\frac{3}{2}}} \textnormal{e}^{\textnormal{i}Y}(1 + \Delta_{R}), \qquad 
    z_{+}({\rho^{*}}) =-a\delta_{0}\sigma_{0}(1+\Delta_{R}), \qquad\sigma_{+}({\rho^{*}}) = \sigma_{0}, \qquad \varepsilon_{+}({\rho^{*}}) = \frac{ \varepsilon^{\frac{1}{2}}}{\sigma_{0}},\label{initial:rings;rho0}
\end{align}
Now, we return to the original $(A,B)$ coordinates by converting \eqref{initial:rings;rho0} into $(A_{1}, z_{1})$ coordinates and inverting the coordinate change \eqref{scale:1}. That is, we write $A = \sigma_{+}A_{+}$ and $B = \sigma_{+}^{2}A_{+}[1 + z_{+}]$, and {evaluate at $\sigma=\frac{1}{r}$; then,} we obtain 
\begin{align}
    A({r_{0}}) &=q_{0}\sqrt{r_{0}}  {\varepsilon^{\frac{3}{4}}}\textnormal{e}^{\textnormal{i}Y}(1 + \Delta_{R}), \qquad 
    B({r_{0}}) =\frac{q_{0}}{\sqrt{r_{0}}} {\varepsilon^{\frac{3}{4}}}\textnormal{e}^{\textnormal{i}Y}(1+\Delta_{R}). \label{initial:rings;AB}
\end{align}
\subsubsection{Matching Core and Far-Field Solutions: Rings}\label{subs:rings;match}
We are now fully equipped to find a nontrivial solution contained in the intersection of the core manifold $\widetilde{\mathcal{W}}^{cu}_{-}\change{(\varepsilon)}$ and the far-field manifold $\change{\widetilde{\mathcal{W}}^{s}_{+}}$. We begin by substituting the initial conditions \eqref{initial:rings;AB} into the far-field parametrisation \eqref{match:farf:folia} following the normal-form transformation \eqref{ab:tilde}
\begin{align}
    \widetilde{a}(r_{0}) &= A(r_{0}) +\textnormal{O}_{r_{0}}\left(|\vec{\mathbf{d}}|_{1}\left[|\varepsilon| + |A| + |B| + |\vec{\mathbf{d}}|_{1}\right]\right),\nonumber\\
    &= \textnormal{e}^{\textnormal{i}Y}q_{0} r_{0}^{\frac{1}{2}} \varepsilon^{\frac{3}{4}}(1 + \Delta_{R}) +\textnormal{O}_{r_{0}}\left(|\vec{\mathbf{d}}|_{1}\left[|\varepsilon|^{\frac{3}{4}} + |\vec{\mathbf{d}}|_{1}\right]\right),\nonumber\\
    \widetilde{b}(r_{0}) &= B(r_{0}) +\textnormal{O}_{r_{0}}\left(|\vec{\mathbf{d}}|_{1}\left[|\varepsilon| + |A| + |B| + |\vec{\mathbf{d}}|_{1}\right]\right),\nonumber\\
    &= \textnormal{e}^{\textnormal{i}Y} q_{0} r_{0}^{-\frac{1}{2}}\varepsilon^{\frac{3}{4}}(1+\Delta_{R}) +\textnormal{O}_{r_{0}}\left(|\vec{\mathbf{d}}|_{1}\left[|\varepsilon|^{\frac{3}{4}} + |\vec{\mathbf{d}}|_{1}\right]\right),\nonumber\\
    a_{n}(r_{0}) &= \textnormal{O}_{r_{0}}\left(\left[|A| + |B| + |\vec{\mathbf{d}}|_{1}\right]\left[|\varepsilon| + |A| + |B| + |\vec{\mathbf{d}}|_{1}\right]\right),\nonumber\\
    &= \textnormal{O}_{r_{0}}\left(\left[|\varepsilon|^{\frac{3}{4}} + |\vec{\mathbf{d}}|_{1}\right]^{2}\right),\nonumber\\
    a_{-n}(r_{0}) &= \vec{d}_{n}\textnormal{e}^{-\lambda_{n}r_{0}}  + \textnormal{O}_{r_{0}}\left(\left[|A| + |B| \right]\left[|\varepsilon| + |A| + |B| + |\vec{\mathbf{d}}|_{1}\right]\right) ,\nonumber\\
     &= \vec{d}_{n}\textnormal{e}^{-\lambda_{n}r_{0}}  + \textnormal{O}_{r_{0}}\left(|\varepsilon|^{\frac{3}{4}}\left[ |\varepsilon|^{\frac{3}{4}} + |\vec{\mathbf{d}}|_{1}\right]\right). \nonumber
\end{align}
We also recall the transformed core parametrisation \eqref{match:core;normal}
\begin{align}
    \widetilde{a}(r_{0}) &= \textnormal{e}^{-\textnormal{i}\left(\frac{\pi}{4} + \textnormal{O}(r_{0}^{-2}) + \textnormal{O}_{r_{0}}(|\varepsilon| + |\mathbf{d}|_{1})\right)}r_{0}^{-\frac{1}{2}}\left([-\textnormal{i}+\textnormal{O}(r_{0}^{-1})]r_{0} d_{1} + [1+\textnormal{O}(r_{0}^{-1})]d_{2}\right) + \textnormal{O}_{r_{0}}\left(|\mathbf{d}|_{1}\left[|\varepsilon| + |\mathbf{d}|_{1}\right]\right),\nonumber\\
    \widetilde{b}(r_{0}) &= \textnormal{e}^{-\textnormal{i}\left(\frac{\pi}{4} + \textnormal{O}(r_{0}^{-2}) + \textnormal{O}_{r_{0}}(|\varepsilon| + |\mathbf{d}|_{1})\right)}r_{0}^{-\frac{1}{2}}\left(\left[-\textnormal{i} +\textnormal{O}(r_{0}^{-1})\right]d_{1} - \left[\nu + \textnormal{O}(r_{0}^{-\frac{1}{2}})\right] d_{2}^{2}\right)+ \textnormal{O}_{r_{0}}\left(|\mathbf{d}|_{1}\left[|\varepsilon| + |\mathbf{d}_{2}|_{1}\right] + |d_{2}|^{3}\right),\nonumber\\
    a_{n}(r_{0}) &= \textnormal{e}^{\lambda_{n}r_{0}}r_{0}^{-\frac{1}{2}}\left[\frac{1}{\sqrt{\pi}} + \textnormal{O}(r_{0}^{-1})\right]\widetilde{c}_{1,n} + \textnormal{O}_{r_{0}}\left(|\mathbf{d}|_{1}\left[|\varepsilon| + |\mathbf{d}|_{1}\right]\right),\nonumber\\
    a_{-n}(r_{0}) &= \textnormal{O}_{r_{0}}\left(|\mathbf{d}|_{1}\left[|\varepsilon| + |\mathbf{d}|_{1}\right]\right),\nonumber
\end{align}
where $|\mathbf{d}|_{1} = |d_{1}| + |d_{2}| + |\mathbf{c}_{1}|_{1}$, and $|\mathbf{d}_{2}|_{1} = |d_{1}| + |\mathbf{c}_{1}|_{1}$. Setting these parametrisations equivalent to each other for each respective coordinate $(\widetilde{a}, \widetilde{b}, a_{1}, a_{-1}, a_{2}, a_{-2}, \dots)(r_{0})$ is the same as finding the zeros of the functional 
\begin{align}
    G:&\left(d_{1}, d_{2}, \widetilde{Y},  \mathbf{c}_{1}, \vec{\mathbf{d}}; a, \varepsilon\right) \mapsto \left( G^{C}_{1}, G^{C}_{2}, G^{H}_{1}, G^{H}_{-1}, G^{H}_{2}, G^{H}_{-2}, \dots \right), \nonumber
\end{align}
where 
\begin{align}
    G^{C}_{1} &= [-\textnormal{i}+\textnormal{O}(r_{0}^{-1})] d_{1} + [1+\textnormal{O}(r_{0}^{-1})]r_{0}^{-1} d_{2} - \textnormal{e}^{\textnormal{i}\widetilde{Y}}q_{0} \varepsilon^{\frac{3}{4}}(1 + \Delta_{R}) + \mathscr{R}^{C}_{1},\label{G:Rings;C1}\\
    G^{C}_{2} &= \left[-\textnormal{i} +\textnormal{O}(r_{0}^{-1})\right] d_{1} - \left[\nu + \textnormal{O}(r_{0}^{-\frac{1}{2}})\right] d_{2}^{2} - \textnormal{e}^{\textnormal{i}\widetilde{Y}} q_{0} \varepsilon^{\frac{3}{4}}(1+\Delta_{R}) + \mathscr{R}^{C}_{2},\label{G:Rings;C2}\\
    G^{H}_{n} &= \textnormal{e}^{\lambda_{n}r_{0}}r_{0}^{-\frac{1}{2}}\left[\frac{1}{\sqrt{\pi}} + \textnormal{O}(r_{0}^{-1})\right]c_{1,n} + \mathscr{R}^{H}_{1},\label{G:Rings;Hn}\\
    G^{H}_{-n} &= - \vec{d}_{n}\textnormal{e}^{-\lambda_{n}r_{0}} + \mathscr{R}^{H}_{2},\label{G:Rings;-Hn}
\end{align}
and we have defined
\begin{align}
    \widetilde{Y} :&= Y + \frac{\pi}{4} + \textnormal{O}(r_{0}^{-2}) + \textnormal{O}_{r_{0}}(|\varepsilon| + |d_{1}| + |d_{2}| + |\mathbf{c}_{1}|_{1}),\nonumber\\
    \mathscr{R}^{C}_{1} :&= \textnormal{O}_{r_{0}}\left(\left[|d_{1}| + |d_{2}| + |\mathbf{c}_{1}|_{1}\right]\left[|\varepsilon| + |d_{1}| + |d_{2}| + |\mathbf{c}_{1}|_{1}\right] + |\vec{\mathbf{d}}|_{1}\left[|\varepsilon|^{\frac{3}{4}} + |\vec{\mathbf{d}}|_{1}\right]\right),\nonumber\\
    \mathscr{R}^{C}_{2} :&= \textnormal{O}_{r_{0}}\left(\left[|d_{1}| + |d_{2}| + |\mathbf{c}_{1}|_{1}\right]\left[|\varepsilon| + |d_{1}| + |\mathbf{c}_{1}|_{1}\right] + |d_{2}|^{3} + |\vec{\mathbf{d}}|_{1}\left[|\varepsilon|^{\frac{3}{4}} + |\vec{\mathbf{d}}|_{1}\right]\right),\nonumber\\
    \mathscr{R}^{H}_{1} :&= \textnormal{O}_{r_{0}}\left(\left[|d_{1}| + |d_{2}| + |\mathbf{c}_{1}|_{1}\right]\left[|\varepsilon| + |d_{1}| + |d_{2}| + |\mathbf{c}_{1}|_{1}\right] + \left[|\varepsilon|^{\frac{3}{4}} + |\vec{\mathbf{d}}|_{1}\right]^{2}\right),\\\
    \mathscr{R}^{H}_{2} :&= \textnormal{O}_{r_{0}}\left(\left[|d_{1}| + |d_{2}| + |\mathbf{c}_{1}|_{1}\right]\left[|\varepsilon| + |d_{1}| + |d_{2}| + |\mathbf{c}_{1}|_{1}\right] + |\varepsilon|^{\frac{3}{4}}\left[|\varepsilon|^{\frac{3}{4}} + |\vec{\mathbf{d}}|_{1}\right]\right).\nonumber
\end{align}
We introduce the scaling $(d_{1}, d_{2}) = (\varepsilon^{\frac{3}{4}}\widetilde{d}_{1}, \varepsilon^{\frac{3}{4}}r_{0} \widetilde{d}_{2})$ so that terms in \eqref{G:Rings;C1} and \eqref{G:Rings;C2} scale with the same order in $\varepsilon$. Initially setting $\varepsilon=0$, we investigate \eqref{G:Rings;Hn} and \eqref{G:Rings;-Hn}
\begin{align}
    F_{n} &= \textnormal{e}^{\lambda_{n}r_{0}}r_{0}^{-\frac{1}{2}}\left[\frac{1}{\sqrt{\pi}} + \textnormal{O}(r_{0}^{-1})\right]c_{1,n} + \textnormal{O}_{r_{0}}\left( |\mathbf{c}_{1}|^{2}_{1} + |\vec{\mathbf{d}}|_{1}^{2}\right),\nonumber\\
    F_{-n} &= - \vec{d}_{n}\textnormal{e}^{-\lambda_{n}r_{0}} + \textnormal{O}_{r_{0}}\left( |\mathbf{c}_{1}|^{2}_{1}\right).\nonumber
\end{align}
and, defining the functional $F:(\widetilde{\mathbf{c}}_{1}, \vec{\mathbf{d}})\mapsto \left(F_{j}\right)_{j\in\mathbb{Z}\backslash\{0\}}$, it clear that $F(\mathbf{0},\mathbf{0})=\mathbf{0}$, where ${\mathbf{0} = (0,\dots)}$. Furthermore, the {Jacobian} $DF(\mathbf{0},\mathbf{0})$ is invertible, and so we can solve \eqref{G:Rings;Hn} and \eqref{G:Rings;-Hn} for all values of $n\in\mathbb{N}$ uniquely for sufficiently small $0<\varepsilon\ll1$. Matching orders of $\varepsilon$ in \eqref{G:Rings;Hn} and \eqref{G:Rings;-Hn}, we find that 
\begin{align}
    c_{1,n} = \textnormal{O}_{r_{0}}\left(|\varepsilon|^{\frac{3}{2}}\right), \qquad \vec{d}_{n}  =\textnormal{O}_{r_{0}}\left(|\varepsilon|^{\frac{3}{2}}\right), \qquad \forall n\in\mathbb{N}.\nonumber
\end{align}
Returning to \eqref{G:Rings;C1} and \eqref{G:Rings;C2}, we have
\begin{align}
    \widetilde{G}_{1} &= [-\textnormal{i}+\textnormal{O}(r_{0}^{-1})]\widetilde{d}_{1} + [1+\textnormal{O}(r_{0}^{-1})]\widetilde{d}_{2} - \textnormal{e}^{\textnormal{i}\widetilde{Y}} q_{0}(1 + \Delta_{R}) + \varepsilon^{-\frac{3}{4}}\mathscr{R}^{C}_{1},\label{Match:Rings;C}\\
    \widetilde{G}_{2} &= \left[-\textnormal{i} +\textnormal{O}(r_{0}^{-1})\right] \widetilde{d}_{1} -\textnormal{e}^{\textnormal{i}\widetilde{Y}} q_{0} (1+\Delta_{R}) - \left[\nu + \textnormal{O}(r_{0}^{-\frac{1}{2}})\right]\varepsilon^{\frac{3}{4}}r_{0}^{2}\widetilde{d}_{2}^{2} + \varepsilon^{-\frac{3}{4}}\mathscr{R}^{C}_{2},\nonumber\\
    \intertext{where}
    \mathscr{R}^{C}_{1} &= |\varepsilon|^{\frac{3}{2}}\textnormal{O}_{r_{0}}\left(\left[|\widetilde{d}_{1}| + |\widetilde{d}_{2}| + |\varepsilon|^{\frac{3}{4}}\right]\left[|\widetilde{d}_{1}| +  |\widetilde{d}_{2}| + |\varepsilon|^{\frac{1}{4}}\right] + |\varepsilon|^{\frac{3}{4}}\right),\nonumber\\
    \mathscr{R}^{C}_{2} &=  |\varepsilon|^{\frac{3}{2}}\textnormal{O}_{r_{0}}\left(\left[|\widetilde{d}_{1}| + |\widetilde{d}_{2}| + |\varepsilon|^{\frac{3}{4}}\right]\left[|\widetilde{d}_{1}| +  |\varepsilon|^{\frac{1}{4}}\right] + |\varepsilon|^{\frac{3}{4}}|\widetilde{d}_{2}|^{3} + |\varepsilon|^{\frac{3}{4}}\right).\nonumber
\end{align}
Initially setting $\varepsilon=0$, we obtain the system
\begin{align}
    \widetilde{G}_{1} &= \left[-\textnormal{i} +\Delta_{R}\right] \widetilde{d}_{1} + [1+\Delta_{R}]\widetilde{d}_{2} - \textnormal{e}^{\textnormal{i}\widetilde{Y}} q_{0} (1 + \Delta_{R}),\label{Match:Rings;C,eps0}\\
    \widetilde{G}_{2} &= \left[-\textnormal{i} +\Delta_{R}\right] \widetilde{d}_{1} -\textnormal{e}^{\textnormal{i}\widetilde{Y}} q_{0} (1+\Delta_{R}).\nonumber
\end{align}
We formally set $\Delta_{R}=0$ and separate \eqref{Match:Rings;C,eps0} into real and imaginary parts: this is equivalent to finding zeros of the functional
\begin{align}
    \widetilde{G}(\widetilde{d}_{1}, \widetilde{d}_{2}, \widetilde{Y}) = \begin{pmatrix} \widetilde{d}_{2} - q_{0} \cos(\widetilde{Y}) \\ -\widetilde{d}_{1} -q_{0} \sin(\widetilde{Y}) \\ -q_{0}\cos(\widetilde{Y}) \end{pmatrix}.\nonumber
\end{align}
It is apparent that the vectors
\begin{align}
    \mathbf{U}^{\pm}:\left(\widetilde{d}_{1}, \widetilde{d}_{2}, \widetilde{Y}\right) = \left(\pm q_{0},\, 0,\, \frac{(2\pm1)\pi}{2} \right),\nonumber
\end{align}
are roots of $\widetilde{G}$ with Jacobian
\begin{align}
    D\widetilde{G}\left[\mathbf{U}^{\pm}\right] = \begin{pmatrix} 0 & 1 & \mp q_{0} \\ -1 & 0 & 0  \\ 0 & 0 & \mp q_{0} \end{pmatrix}.\nonumber
\end{align}
Since $q_{0}>0$, the Jacobian is invertible, and we can therefore solve \eqref{Match:Rings;C,eps0} uniquely for all sufficiently small $\Delta_{R}$, that is, for $r_{0}$ large enough and $\delta_{0}$ small enough, and subsequently \eqref{Match:Rings;C} for all $0<\varepsilon\ll1$. Reversing the scaling for $d$, we find that 
\begin{align}
    d_{1} &= \pm\varepsilon^{\frac{3}{4}} q_{0}\left(1 + \textnormal{O}(r_{0}^{-1} + \delta_{0} + \varepsilon^{\frac{1}{4}})\right),\label{Rings:d1}\\
    d_{2} &= \textnormal{O}(\varepsilon).\label{Rings:d2}
\end{align}
Hence, we have found the ring solutions. We recall that our solution $\mathbf{u}(r,y)$ takes the form
\begin{align}
    \mathbf{u}(r,y) &= a(r)\mathbf{e}(y) + b(r)\mathbf{f}(y) + \overline{a}(r)\overline{\mathbf{e}}(y) + \overline{b}(r)\overline{\mathbf{f}}(y) + \Sigma_{n=1}^{\infty} \left\{a_{n}(r)\mathbf{e}_{n}(y) + a_{-n}(r)\mathbf{e}_{-n}(y)\right\}, \nonumber
\end{align}
and 
\begin{align}
    \begin{pmatrix}a \\ b \end{pmatrix}(r) & = \sum_{i=1}^{4} \widetilde{d}_{i}\mathbf{V}_{i}(r), \qquad \qquad \begin{pmatrix}a_{n} \\ a_{-n} \end{pmatrix}(r)  = \sum_{i=1}^{2} \widetilde{c}_{i,n}\mathbf{W}_{i,n}(r).\nonumber
\end{align}
Substituting \eqref{Rings:d1} and \eqref{Rings:d2} into this form, we can write the ring solution $\mathbf{u}_{R}$ as
\begin{align}
    \mathbf{u}_{R} &= \pm\varepsilon^{\frac{3}{4}}q_{0}\sqrt{\frac{k\pi}{2}}\left[ r J_{1}(k r) (\mathbf{e} + \overline{\mathbf{e}}) + \textnormal{i}[k J_{1}(k r) - r J_{0}(k r)](\mathbf{e} - \overline{\mathbf{e}}) + J_{1}(k r)(\mathbf{f} + \overline{\mathbf{f}}) -\textnormal{i}J_{0}(k r)(\mathbf{f} - \overline{\mathbf{f}})\right] + \textnormal{O}(\varepsilon),\nonumber
\end{align}
for all $r\in(0,r_{0})$, where $\mathbf{u}_{R}(r, \change{y})$ decays to zero exponentially as $r\to\infty$. In particular, the height of the free surface $\eta_{R}(r)$ has the form
\begin{align}
    \eta^{\pm}_{R}(r) &= \pm\varepsilon^{\frac{3}{4}}\frac{2 q_{0}}{m}\sqrt{\frac{k\pi}{2}}\left[ r J_{1}(k r) + \widetilde{b}_{D} J_{0}(k r)\right] + \textnormal{O}(\varepsilon), \nonumber
\end{align}
for all $r\in[0, r_{0}]$, where $\widetilde{b}_{D}:=b_{D} - D\tanh(k D)+ k^{-1}$ and $b_{D}$ is defined in \eqref{bD:defn} in Appendix \ref{app:basis}.  Similarly, for $(A_{+},z_{+})(\rho)$ in \eqref{soln:rings}, where $\rho={\log}\; r $, we can invert the coordinate transformations \eqref{scale:1}, \eqref{z:def}, and \eqref{A0:A;trans} to find the free surface profile in the transition chart,
\begin{align}
    \eta^{\pm}_{R}(r) &= \pm\varepsilon^{\frac{3}{4}}\frac{2 q_{0}}{m}\left[ r^{\frac{1}{2}} \sin\left(k r-\frac{\pi}{4}\right) + \widetilde{b}_{D} r^{-\frac{1}{2}}\cos\left(k r - \frac{\pi}{4}\right)\right] + \textnormal{O}(\varepsilon), \nonumber
\end{align}
for $r\in[r_{0}, \delta_{0}\varepsilon^{-\frac{1}{2}}]$. Finally, for $(A_{2},z_{2})(s)$ in \eqref{soln:chart;comp}, we can invert transformations \eqref{scale:2}, \eqref{z:def}, and \eqref{A0:A;trans} to find the free surface profile in the far-field,
\begin{align}
    \eta^{\pm}_{R}(r) &= \pm\varepsilon^{\frac{3}{4}}\frac{2}{m}\left[ \varepsilon^{-\frac{1}{4}}q(\varepsilon^{\frac{1}{2}}r) \sin\left(k r - \frac{\pi}{4}\right) + \widetilde{b}_{D}\, p(\varepsilon^{\frac{1}{2}} r) \cos\left(k r - \frac{\pi}{4}\right)\right] + \textnormal{O}(\varepsilon), \nonumber
\end{align}
for $r\in[\delta_{0}\varepsilon^{-\frac{1}{2}}, \infty)$, where $p(\varepsilon^{\frac{1}{2}} r)$ is defined in \eqref{pr:defn}. This completes the result for rings.
\section{Appendix}
\subsection{Spectrum of L as \textit{r} approaches infinity}\label{app:spect}
\begin{figure}[t!]
    \centering
    \includegraphics[height=7cm]{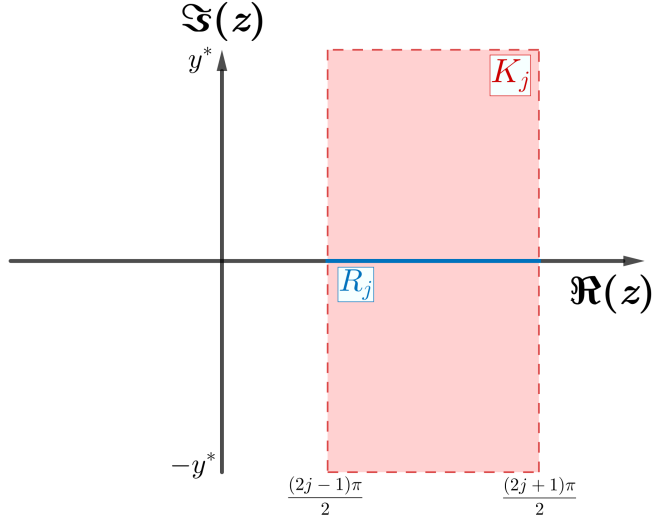}
    \caption{The {set $\{\varsigma\in\mathbb{C}: \Delta_{2}(\varsigma)=0\}$} can be decomposed into {the regions $K_{j}$, where} $K_{j}=\left\{z=x+\textnormal{i}y \, :\; x\in\left(\frac{(2j-1)\pi}{2}, \frac{(2j+1)\pi}{2}\right), y\in(-y^{*},y^{*})\right\}$ for all $j\in\mathbb{Z}$. The boundary $\partial K_{j}$ is a closed contour for all $j\in\mathbb{Z}$.}
    \label{fig:Rouche}
\end{figure}

We fix $(\mathcal{M}, \widetilde{\Upsilon}_{0}) = (\mathcal{M}_{H}, \widetilde{\Upsilon}_{H})\left(kD\right)$ as defined in \eqref{M:Upsilon}. Then, $\lambda = \pm\textnormal{i} k$ are the only purely imaginary eigenvalues of the linear differential operator $\mathbf{L}_{\infty}:=\lim_{r\to\infty}\mathbf{L}(r)$, where $\mathbf{L}(r)$ is defined in \eqref{L:lin}. We recall, from dispersion relation \eqref{disp}, a complex number $\lambda\in\mathbb{C}$ is an eigenvalue of $\mathbf{L}_{\infty}$ if and only if $\Delta_{0}(\lambda D)=0$. It is clear that, for any $\varsigma\in\mathbb{C}$, 
\begin{align}
    \Delta_{0}(\varsigma) &= 0 \qquad \iff \qquad \sin(\varsigma)=0 \quad \textnormal{or}\quad \Delta_{1}(\varsigma) = 0\nonumber\\
    \intertext{where}
    \Delta_{1}(\varsigma) &:= \mathcal{M}\;\varsigma \sin(\varsigma) - \left(\varsigma^2 - \widetilde{\Upsilon}_{0}\right)\cos(\varsigma).\label{disp:1}
\end{align}
Then, for any eigenvalue $\lambda\in\mathbb{C}$, 
\begin{align}
    \lambda D \in\left\{n\pi\right\}_{n\in\mathbb{Z}\backslash\{0\}}\cup\,{\{\varsigma\in\mathbb{C}: \Delta_{1}(\varsigma)=0\}}.\nonumber
\end{align}
Before we can explicitly find { $\{\varsigma\in\mathbb{C}: \Delta_{1}(\varsigma)=0\}$}, we first decompose $\varsigma=x+\textnormal{i}y$, where $x,y\in\mathbb{R}$, and prove some ancillary results:
\begin{Lemma} Let $j\in\mathbb{Z}$, and $y\in\mathbb{R}$. Then, $\Delta_{1}\left(\frac{(2j+1)\pi}{2}+\textnormal{i}y\right)\neq0$.
\label{app:dispcomp;x}\end{Lemma}
\begin{Proof}
Assume there exists some $y\in\mathbb{R}$, such that $\Delta_{1}\left(\frac{(2j+1)\pi}{2} + \textnormal{i}y\right)=0$. By taking real and imaginary parts of \eqref{disp:1}, one can find
\begin{align}
    \Re[\Delta_{1}]&=0, &\quad \implies \qquad\, y&=\widetilde{y}, & & \textnormal{where $\widetilde{y}$ satisfies }\quad 2\widetilde{y}\tanh(\widetilde{y}) = \mathcal{M},&\nonumber\\
    \Im[\Delta_{1}]&=0, &\quad \implies \quad 4\widetilde{\Upsilon}_{0} &= (2j+1)^2\pi^2 + 4\widetilde{y}^2,& &\,&\nonumber\\
     \,& \, & \, &= (2j+1)^2\pi^2 + \mathcal{M}^2\coth^2(\widetilde{y}).& &\,& \nonumber
\end{align}
However, by taking the explicit values $(\mathcal{M}, \widetilde{\Upsilon}_{0}) = (\mathcal{M}_{H},\widetilde{\Upsilon}_{H})$ seen in \eqref{M:Upsilon}, one can see that $4\widetilde{\Upsilon}_{0}<\mathcal{M}^2$. Thus, by contradiction, $\Delta_{1}\left(\frac{(2j+1)\pi}{2} + \textnormal{i}y\right)\neq0$.
\end{Proof}
Following Lemma \ref{app:dispcomp;x}, we see that $\Delta_{1}(x+\textnormal{i}y)=0$ only if $x\neq \frac{(2j+1)\pi}{2}$ for all $j\in\mathbb{Z}$. Then, we define the complex function
\begin{align}
    \Delta_{2}(\varsigma):=\mathcal{M}\;\varsigma \tan(\varsigma) - \left(\varsigma^2 - \widetilde{\Upsilon}_{0}\right), \qquad \textnormal{such that} \; \,{\{\varsigma\in\mathbb{C}: \Delta_{2}(\varsigma)=0\} = \{\varsigma\in\mathbb{C}: \Delta_{1}(\varsigma)=0\}}.\label{disp:2}
\end{align}
Hence, we will proceed by analysing the zeros of $\Delta_{2}(\varsigma)$. 
\begin{Lemma} We define the region $R_{j}=\left(\frac{(2j-1)\pi}{2},\frac{(2j+1)\pi}{2}\right)$, for each $j\in\mathbb{Z}$. Then, there exists a value $x_{j}\in R_{j}$ for all $j\in\mathbb{Z}\backslash\{0\}$, such that $\Delta_{2}(x_{j})=0$ and $x_{-j}=-x_{j}$. Furthermore, $\Delta_{2}(x)\neq0$ for all $x\in R_{0}$.
\label{app:disp;real}\end{Lemma}
\begin{Proof} As $\Delta_{2}(\varsigma)$ is an even function, we restrict $\varsigma=x\in R_{j}$ for $j\in\mathbb{N}$; any results for $x_{j}\in R_{j}$ are mirrored for $x_{-j}=-x_{j}\in R_{-j}$. Since, for each $j\in\mathbb{N}$, $\Delta_{2}(x) \to \pm\infty$ as $x\to \frac{(2j\pm1)\pi}{2}$, we apply a limiting version of the Intermediate Value Theorem: one can find constants $a_{j}, b_{j} \in R_{j}$ such that $a_{j}<b_{j}$, $\Delta_{2}(a_{j})<0$, and $\Delta_{2}(b_{j})>0$. Then, by the IVT, there exists a point $x_{j}\in [a_{j},b_{j}]$ such that $\Delta_{2}(x_{j})=0$.  Therefore, we have proven that $\Delta_{2}(\varsigma)$ has at least one root $\varsigma=x_{j}\in R_{j}$, for each $j\in\mathbb{N}$. Using the even properties of $\Delta_{2}(\varsigma)$, there is an $x_{j}\in R_{j}$ for $j\in\mathbb{Z}\backslash\{0\}$, where $x_{-j}=-x_{j}$.

{For $x\in R_{0}$, we note that $\Delta_{2}^{'}(0) =0$ and $\Delta_{2}^{''}(0) =2(\mathcal{M} - 1) > 0$, and so $x=0\in R_{0}$ is a local minimum. Then, computing the second derivative of $\Delta_{2}(x)$ for all $x\in R_{0}$, we see that
\begin{align}
    \Delta_{2}^{''}(x) &=2\mathcal{M}\left[x\tan(x)\left(1+ \tan^{2}(x)\right) + 1 + \tan^{2}(x)\right] - 2,\nonumber\\ 
    &=2\mathcal{M}\left[x\tan(x)\left(1+ \tan^{2}(x)\right) + \tan^{2}(x)\right] + 2(\mathcal{M}-1).\nonumber
\end{align}
For $|x|<\frac{\pi}{2}$, $x\tan(x)\geq0$, and so $\Delta_{2}^{''}(x)>0$ for all $x\in R_{0}$, i.e. $\Delta_{2}(x)$ is convex. Therefore, $\Delta_{2}(0)$ is also a global minimum; hence, $\Delta_{2}(x) \geq \Delta_{2}(0) = \widetilde{\Upsilon}_{0} > 0$. Consequently, $\Delta_{2}(x)>0$ for all $x \in R_{0}$.}
\end{Proof}
\begin{Lemma} There exists some fixed $y^{*}$ such that for all ${|y|}>y^{*}$, $\Delta_{2}\left(x+\textnormal{i}y\right)\neq0$.
\label{app:dispcomp;y}\end{Lemma}
\begin{Proof}
{ We recall the following inequalities,
\begin{align}
    |\sin(\varsigma)|\leq \frac{1}{2}(\textnormal{e}^{\Im(\varsigma)} + \textnormal{e}^{-\Im(\varsigma)}), \qquad |\cos(\varsigma)|\geq \frac{1}{2}|\textnormal{e}^{\Im(\varsigma)} - \textnormal{e}^{-\Im(\varsigma)}|,\nonumber
\end{align}
and so, taking $\varsigma=x+\textnormal{i}y$, we can write that there exists some $y_{0}$ such that, for $|y|>y_{0}$,
\begin{align}
    |\sin(x+\textnormal{i}y)|\leq \textnormal{e}^{|y|}, \qquad |\cos(x+\textnormal{i}y)|\geq \frac{1}{4}\textnormal{e}^{|y|}.\nonumber
\end{align}
We also note that $|\varsigma^{2} - \widetilde{\Upsilon}_{0}|\geq |\varsigma|^{2} - \widetilde{\Upsilon}_{0}>4\mathcal{M}|\varsigma|$ for $|\varsigma|>r_{0}$, where $r_{0}$ is some large fixed value. Then, for $|y|> y^{*}:=\max\{y_{0},r_{0}\}$, we see that
\begin{align}
    |\mathcal{M}\varsigma \sin(\varsigma)|\leq \mathcal{M}\sqrt{x^{2}+y^{2}} \textnormal{e}^{|y|}, \qquad |(\varsigma^{2}-\widetilde{\Upsilon}_{0})\cos(\varsigma)|> \mathcal{M}\sqrt{x^{2}+y^{2}} \textnormal{e}^{|y|},\nonumber
\end{align}
so that $\Delta_{1}(\varsigma)\neq0$, and consequently $\Delta_{2}(\varsigma)\neq0$.
}
\end{Proof}
We define the complex regions $K_{j} := R_{j}\times\textnormal{i}\left[-y^{*},y^{*}\right]$, where $R_{j}$ is defined in Lemma \ref{app:disp;real}, and $y^{*}$  is some large real number as seen in Lemma \ref{app:dispcomp;y}; see Figure \ref{fig:Rouche}. Thanks to Lemmas \ref{app:dispcomp;x} and \ref{app:dispcomp;y}, we can {write $\{\varsigma\in\mathbb{C}: \Delta_{2}(\varsigma)=0\}\subseteq \bigcup_{j\in\mathbb{Z}} K_{j}$}. From now on, we will investigate zeros of $\Delta_{2}(\varsigma)$ with $\varsigma$ restricted to $K_{j}$ for some arbitrary $j\in\mathbb{Z}$.
\begin{Lemma} Define $\Delta_{3}(\varsigma):\mathbb{C}\to\mathbb{C}, \; \varsigma\mapsto -\varsigma\left[\mathcal{M}\tan(\varsigma)-\varsigma\right]$. Then, $\Delta_{3}(\varsigma)$ has four zeros in $K_{0}$ and one zero in $K_{j}$, for all $j\in\mathbb{Z}\backslash\{0\}$.
\label{app:dispcomp;zeros}\end{Lemma}
\begin{Proof} We define $\varsigma = x+\textnormal{i}y$ and begin with the case when $y=0$. Since $\mathcal{M}\tan(x)$ is a monotonically-increasing surjective function for $x\in R_{j}$, $\Delta_{3}(x)$ has a repeated zero at $x=0\in K_{0}$, and a unique non-trivial zero $x=x_{j}\in K_{j}$, for all $j\in\mathbb{Z}\backslash\{0\}$.

\noindent For $y\neq0$, { any solution $\varsigma=x+\textnormal{i} y$ such that $\Delta_{3}(\varsigma)=0$ must also satisfy 
\begin{align}
    \left[y \sin(2x) + x\sinh(2y)\right] \Re[\Delta_{3}(\varsigma)] + \left[y \sinh(2y) - x\sin(2x)\right]\Im[\Delta_{3}(\varsigma)] &= 0,\nonumber\\
    \implies \left(x^{2} + y^{2}\right)\left[ y\sin(2x) - x\sinh(2y) \right] &= 0.\nonumber
\end{align}
  When $y\neq0$, this is only satisfied when $x=0$.} Furthermore, $\Im[\Delta_{3}(\textnormal{i}y)]=0$ implies that $\mathcal{M}\tanh(y)-y=0$. This has two non-trivial zeros at $y=\pm y_{0}\in K_{0}$. Hence, the Lemma is proven.
\end{Proof}

Finally, we can explicitly find $\{\varsigma\in\mathbb{C}: \Delta_{2}(\varsigma)=0\}$. We will make use of the symmetric version of Rouch\'e's theorem (See Glicksberg \cite{Glicksberg1976Rouche}), and so we state it here:
\begin{Theorem}{\cite{Glicksberg1976Rouche}}
Let $K\subset\mathbb{C}$ be a bounded region with continuous boundary $\partial K$. Two holomorphic functions $f,g:\mathbb{C}\to\mathbb{C}$ have the same number of zeros in $K$, if the strict inequality
\begin{align}
\left|f(\varsigma)+g(\varsigma)\right|<\left|f(\varsigma)\right| + \left|g(\varsigma)\right|,\nonumber
\end{align}
holds for all values of $\varsigma$ on the boundary $\partial K$,
\label{Glicksberg}\end{Theorem}
\begin{Proposition}
$\Delta_{1}(\varsigma)$ has precisely two purely imaginary zeros $\varsigma=\pm\textnormal{i}kD$, each with multiplicity 2, and an infinite set of purely real zeros $\varsigma=\pm\widetilde{\lambda}_{j}D$, $j\in\mathbb{N}$, with multiplicity 1.
\end{Proposition}
\begin{Proof} By Lemmas \ref{app:dispcomp;x} \& \ref{app:dispcomp;y}, {$\,\{\varsigma\in\mathbb{C}: \Delta_{2}(\varsigma)=0\}\subseteq\bigcup_{j\in\mathbb{Z}}K_{j}$}. We will use Theorem \ref{Glicksberg} where $K = K_{j}$ for some $j\in\mathbb{Z}$, with $f(\varsigma)=\Delta_{3}(\varsigma)$ and $g(\varsigma)=\Delta_{2}(\varsigma)$. By applying the triangle inequality and verifying that $\widetilde{\Upsilon}_{0} \neq \left|\Delta_{2}(\varsigma)\right| + \left|\Delta_{3}(\varsigma)\right|$ for all $\varsigma\in\partial K_{j}$, we see that $\Delta_{3}(\varsigma)$ and $\Delta_{2}(\varsigma)$ have the same number of zeros in each domain $K_{j}$. Then, by Lemma \ref{app:dispcomp;zeros}, $\Delta_{2}(\varsigma)$ has four zeros in $K_{0}$ and one zero in $K_{j}$ for each $j\in\mathbb{Z}\backslash\{0\}$.

However, for $(\mathcal{M},\widetilde{\Upsilon}_{0}) = (\mathcal{M}_{H},\widetilde{\Upsilon}_{H})$, we have already four zeros in $K_{0}$, namely $\varsigma=\pm \textnormal{i}k D$ with double algebraic multiplicity. We define the unique root $K_{j}$ for each $j\in\mathbb{Z}\backslash\{0\}$, namely $\widetilde{\lambda}_{j}D := x_{j}$ as found in Lemma \ref{app:disp;real}. By \eqref{disp:2}, the result is proven.
\end{Proof}

Therefore, a complex number $\lambda$ is an eigenvalue of $\mathbf{L}_{\infty}$ if and only if
\begin{align}
    \lambda \in \left\{\pm\textnormal{i}k\right\} \cup\left\{\pm\lambda_{n}\right\}_{n\in\mathbb{N}},\qquad    \textnormal{where}\quad\lambda_{2n-1}:=\min\left\{\widetilde{\lambda}_{n}, \frac{n\pi}{D}\right\}, \qquad \lambda_{2n} :=\max\left\{\widetilde{\lambda}_{n}, \frac{n\pi}{D}\right\},\nonumber
\end{align}
for each $n\in\mathbb{N}$. 
\subsection{Basis of \textit{y}-dependent eigenmodes}\label{app:basis}
In Section \ref{s:spect}, we introduce a `spectral' decomposition to reduce the problem to a system of amplitude equations. To do this, we find eigenmodes of the linear operator $\mathbf{L}_{\infty}$. Explicitly, for the imaginary eigenvalues $\lambda=\pm \textnormal{i}k$ with double algebraic multiplicity, we have two eigenmodes $\mathbf{e}(y)$ and $\mathbf{f}(y):= \widetilde{\mathbf{f}}(y)+\textnormal{i}b_{D} \mathbf{e}(y)$ of the form,
\begin{align}
    &\mathbf{e}(y) := \frac{1}{m}\begin{pmatrix} M_{0} \frac{\cosh(k(D+y))}{\cosh(kD)}\\ -\mu M_{0} \frac{\cosh(k(D-y))}{\cosh(kD)}\\ 1\\ \textnormal{i} k M_{0} \frac{\cosh(k(D+y))}{\cosh(kD)}\\ -\textnormal{i} k \mu M_{0} \frac{\cosh(k(D-y))}{\cosh(kD)}\\ \textnormal{i} k \end{pmatrix},& \qquad &\widetilde{\mathbf{f}}(y) := \frac{1}{k m}\begin{pmatrix} -\textnormal{i} M_{0} \frac{k(D+y)\sinh(k(D+y)) - \cosh(k(D+y))}{\cosh(kD)}\\ \textnormal{i}\mu M_{0} \frac{k(D-y)\sinh(k(D-y)) - \cosh(k(D-y))}{\cosh(kD)}\\ -\textnormal{i} [k D \tanh(k D) - 1]\\ M_{0} \frac{k^2(D+y)\sinh(k(D+y))}{\cosh(kD)}\\ -\mu M_{0} \frac{k^2(D-y)\sinh(k(D-y))}{\cosh(kD)}\\ k^2 D \tanh(k D) \end{pmatrix},& \label{eigvec;c:defn}\\
    \intertext{where}
    &M_{0} := \frac{\mu-1}{\mu+1},&  \qquad   &m := \sqrt{[kD\tanh(kD)-1]\mathcal{M} \,\textnormal{sech}^{2}(kD) + 1},& \label{const:defn}
\end{align}
and $b_{D}$ is a constant chosen such that $\Omega(\mathbf{f},\overline{\mathbf{f}}) = 0$; in particular, $b_{D}$ takes the form,
\begin{align}
    b_{D}=\frac{k D \cosh^2(k D ) + \frac{4}{3}k^{3} D^{3} -\sinh(k D )\cosh(k D )-2 k^{2} D^{2} \tanh(k D )- k^{3} D^{3}\;\textnormal{sech}^{2}(k D )}{k^2[\sinh(k D )\cosh(k D ) + k D ]\left[1+[k D \tanh(k D) - 1]\mathcal{M}\;\textnormal{sech}^{2}(k D)\right]}.\label{bD:defn}
\end{align}
Similarly, for the infinite series of eigenvalues $\{\pm\lambda_{n}\}_{n\in\mathbb{N}}$ (as found in Section \ref{app:spect}), there are eigenmodes
\begin{align}
    \mathbf{e}_{\pm n} &:=\frac{1}{c_{1,j}} \begin{pmatrix} M_{0} \cos(\lambda_{\pm j}(D+y))\\ -\mu M_{0} \cos({\lambda}_{\pm j}(D-y))\\ \cos({\lambda}_{\pm j}D)\\ {\lambda}_{\pm j} M_{0} \cos({\lambda}_{\pm j}(D+y))\\ -{\lambda}_{\pm j} \mu M_{0} \cos({\lambda}_{\pm j}(D-y))\\ {\lambda}_{\pm j}\cos({\lambda}_{\pm j}D) \end{pmatrix}, \nonumber
\end{align}
    if $\lambda_{n} = \widetilde{\lambda}_{j}$ for some $j\in\mathbb{N}$, and 
    \begin{align}
    \mathbf{e}_{\pm n} &:=\frac{1}{c_{2,j}} \begin{pmatrix} \cos(\frac{j \pi}{D}(D+y))\\ \cos(\frac{j \pi}{D}(D-y))\\ 0\\ \pm\frac{j \pi}{D} \cos(\frac{j \pi}{D}(D+y))\\ \pm\frac{j \pi}{D} \cos(\frac{j \pi}{D}(D-y))\\ 0 \end{pmatrix}, \nonumber 
\end{align}
if $\lambda_{n} = \frac{j \pi}{D}$, for some $j\in\mathbb{N}$. Here, $c_{1,n}$ and $c_{2,n}$ are defined such that $\Omega(\mathbf{e}_{-i}, \mathbf{e}_{j}) = \delta_{i,j}$ for any $i,j\in\mathbb{N}$.
\subsection{Core Problem: Quadratic Term}\label{app:core;quad}
In order to determine the existence of spot A solutions, we need to find the coefficient of the $d_{2}^{2}\mathbf{V}_{3}(r)$ term in the \eqref{fixed:a} expansion for $b(r)$. To do this, we Taylor expand $\widetilde{d}_{3}$ and isolate the $d_{2}^{2}$ term. We recall from \eqref{fixed:a} that 
\begin{align}
    \widetilde{d}_{3}(d_{1}, d_{2}, \mathbf{c}_{1}) = \int_{0}^{r}  \langle \mathbf{V}_{3}^{*}, \mathbf{F}\rangle \textnormal{d} s, \qquad \textnormal{and} \qquad \mathbf{F}= (-\Omega( \mathcal{F}, \overline{\mathbf{f}}), \Omega( \mathcal{F}, \overline{\mathbf{e}})).\nonumber
\end{align}
We set $\varepsilon=0$ and recall that $\mathbf{V}^{*}_{3}(r) = \sqrt{\frac{k \pi}{2}}\Big(0, \;-r\left[J_{0}(k r) + \textnormal{i} J_{1}(k r) \right]\Big)$, and so we can write
\begin{align}
    \widetilde{d}_{3} = \sqrt{\frac{k\pi}{2}}\int_{0}^{r}  \left\{ s J_{1}(k s) \,\Omega(\mathcal{F}, \Im[\mathbf{e}]) - s J_{0}(k s)\,\Omega(\mathcal{F}, \Re[\mathbf{e}]) \right\} \textnormal{d} s. \label{int:besssel}
\end{align}
The nonlinearity is computed up to quadratic order; we take the quadratic nonlinearity of \eqref{Ham:psi;auto}-\eqref{Ham:gamma;auto}, and use the definition of $\mathcal{F}$ in \eqref{full:syst;linbd} to write
\begin{align}
    \mathcal{F}(\mathbf{u},\varepsilon, r) = \left[\textnormal{d}\mathcal{G}\left[\mathcal{G}^{-1}(\mathbf{u})\right]- \mathbb{1}\right]\mathbf{L}(r)\mathbf{u} + \mathscr{F}_{2}(\mathbf{u},\varepsilon,r) + \textnormal{O}(|\mathbf{u}|^{3}),\nonumber
\end{align}
where the subscript $\mathscr{F}_{2}$ denotes elements of $\mathscr{F}$ with quadratic order. Computing the {Jacobian} of $\mathcal{G}$ (see \cite{Hill2021Thesis}), we can write the quadratic part of the nonlinearity {(when $\varepsilon=0$) as}
    \begin{align}
    \mathcal{F}_{2} :&= \begin{pmatrix} 
    \int_{0}^{y}\left\{-\left(\frac{\alpha^{-}_{y}\eta + \psi^{-}_{y}\gamma}{D}\right) + \left(\frac{D+t}{D}\right)\left[\psi^{-}_{yy}\gamma - \alpha^{-}\Lambda_{0} + \frac{2}{r}\alpha^{-}\gamma\right]\right\} \textnormal{d} t
    \\
    \int_{0}^{y} \left\{\left(\frac{\alpha^{+}_{y}\eta + \psi^{+}_{y}{\gamma}}{D}\right) + \left(\frac{D-t}{D}\right)\left[\psi^{+}_{yy}{\gamma} - {\alpha^{+}}\Lambda_{0} + \frac{2}{r}{\alpha^{+}}{\gamma}\right]\right\}\textnormal{d} t + (\mu-1)\Gamma
    \\
    -\Gamma
    \\
    \left(\frac{\psi^{-}_{yy}{\eta} - {\alpha^{-}}{\gamma}}{D}\right) +  \left(\frac{D+y}{D}\right)\left[\alpha^{-}_{y}{\gamma} + \psi^{-}_{y}\Lambda_{0}\right]
    \\
    - \left(\frac{\psi^{+}_{yy}{\eta} - \alpha^{+}{\gamma}}{D}\right) +  \left(\frac{D-y}{D}\right)\left[\alpha^{+}_{y}{\gamma} + \psi^{+}_{y}\Lambda_{0}\right] -(\mu-1)\widetilde{\Gamma}
    \\
    - \widehat{\Gamma}
    \end{pmatrix} + 
    \begin{pmatrix} 
    \Lambda_{1}\\
    \Lambda_{1}\\
    0\\
    \Lambda_{2}\\
    \Lambda_{2}\\
    0
    \end{pmatrix},\label{F2:expl}
\end{align}
    where
\begin{align}
    {\Lambda_{0} :}&{= \Upsilon_{0}\eta - \mu(\mu-1) \psi^{-}_{y}|_{y=0},}\nonumber\\
    \Gamma :&= \mu\int_{-D}^{0} \left(\frac{D+y}{D}\right)\left\{\alpha^{-}\psi^{-}_{y}\right\}\textnormal{d}y + \int_{0}^{D}\left(\frac{D-y}{D}\right)\left\{\alpha^{+}\psi^{+}_{y}\right\}\textnormal{d}y,\nonumber\\
    \widetilde{\Gamma} :&= \mu \int_{-D}^{0} \left(\frac{D+y}{D}\right)\left\{\alpha^{-}\alpha^{-}_{y} - \psi^{-}_{y}\psi^{-}_{yy}\right\}\textnormal{d}y + \int_{0}^{D} \left(\frac{D-y}{D}\right)\left\{\alpha^{+}\alpha^{+}_{y} - \psi^{+}_{y}\psi^{+}_{yy}\right\}\textnormal{d}y, \nonumber\\
     {\widehat{\Gamma}:}&{= \frac{\mu}{2D}\int_{-D}^{0} \left[(\alpha^{-})^2 - (\psi^{-}_{y})^2\right]\textnormal{d} y - \frac{1}{2D}\int_{0}^{D} \left[(\alpha^{+})^2 - (\psi^{+}_{y})^2\right]\textnormal{d} y,}\nonumber
\end{align}
and $\Lambda_{1}(y)$, $\Lambda_{2}(y)$ are functions of $y$, where $\Lambda_{i}(-y) = \Lambda_{i}(y)$. The second vector in \eqref{F2:expl}, containing only terms that are even in $y$, will disappear in the upcoming projection. As defined in \eqref{eigvec;c:defn}, 
\begin{align}
    \Re(\mathbf{e}) = \frac{1}{m}\begin{pmatrix} M_{0} \frac{\cosh(k(D+y))}{\cosh(kD)}\\ -\mu M_{0} \frac{\cosh(k(D-y))}{\cosh(kD)}\\ 1\\ 0 \\ 0 \\ 0 \end{pmatrix}, \quad \Im(\mathbf{e}) = \frac{k}{m}\begin{pmatrix} 0\\ 0\\ 0\\  M_{0} \frac{\cosh(k(D+y))}{\cosh(kD)}\\ -  \mu M_{0} \frac{\cosh(k(D-y))}{\cosh(kD)}\\ 1 \end{pmatrix},\nonumber
\end{align}
and we recall the symplectic two-form $\Omega$
\begin{align}
    \Omega\left(\mathbf{u}_{1}, \mathbf{u}_{2}\right) = \mu\int_{-D}^{0} \left[\psi^{-}_{1}\alpha^{-}_{2} - \alpha^{-}_{1}\psi^{-}_{2}\right]\;\textnormal{d}y + \int_{0}^{D} \left[\psi^{+}_{1}\alpha^{+}_{2} - \alpha^{+}_{1}\psi^{+}_{2}\right]\;\textnormal{d}y - \left[\eta_{1}\gamma_{2} - \gamma_{1}\eta_{2}\right]. \nonumber
\end{align}
Denoting $\mathcal{F}_{2}(\mathbf{u},\varepsilon,r) =(f_{1},f_{2},f_{3},f_{4},f_{5},f_{6})^\intercal$, we can write
\begin{align}
    \Omega(\mathcal{F}_{2}, \Re[\mathbf{e}]) &= \frac{1}{m}\left[ \mu M_{0}\left[\int_{0}^{D} \frac{\cosh(k(D-y))}{\cosh(kD)} f_{5}(y)\;\textnormal{d} y - \int_{-D}^{0} \frac{\cosh(k(D+y))}{\cosh(kD)} f_{4}(y)\;\textnormal{d} y \right]+ f_{6}\right], \nonumber\\
    \Omega(\mathcal{F}_{2}, \Im[\mathbf{e}]) &= -\frac{k}{m}\left[ \mu M_{0}\left[\int_{0}^{D} \frac{\cosh(k(D-y))}{\cosh(kD)} f_{2}(y)\;\textnormal{d} y - \int_{-D}^{0} \frac{\cosh(k(D+y))}{\cosh(kD)} f_{1}(y)\;\textnormal{d} y \right]+ f_{3}\right]. \nonumber
\end{align}
We define $j(y) := \frac{\cosh(k(D+y))}{\cosh(kD)}$ and perform a change of integration variables such that we can write
\begin{align}
    \Omega(\mathcal{F}_{2}, \Re[\mathbf{e}]) &= \frac{1}{m}\left[ \frac{\mathcal{M}}{(\mu-1)D}\int_{-D}^{0} j(y) \left[f_{5}(-y) - f_{4}(y)\right]\;\textnormal{d} y + f_{6}\right], \nonumber\\
    \Omega(\mathcal{F}_{2}, \Im[\mathbf{e}]) &= -\frac{k}{m}\left[ \frac{\mathcal{M}}{(\mu-1)D}\int_{-D}^{0} j(y) \left[f_{2}(-y)- f_{1}(y)\right]\;\textnormal{d} y + f_{3}\right], \nonumber
\end{align}
where $\mathcal{M}=\mu(\mu-1)D M_{0}$. In order to isolate the coefficient of $d_{2}^{2}$, we set $d_{1} = c_{1,n} = 0$, for all $n\in\mathbb{N}$. Then, 
\begin{align}
    \mathbf{u} &= a\mathbf{e} + \overline{a}\,\overline{\mathbf{e}}, \qquad \qquad \textnormal{where} \quad a= d_{2} \sqrt{\frac{k\pi}{2}}(J_{0}(k r) + \textnormal{i} J_{1}(k r)),\nonumber
\end{align}
or, more specifically,
{\begin{align}
   \psi^{-}(y) &=d_{2}\frac{\left(2\pi k\right)^{\frac{1}{2}}}{m} \eta M_{0} j(y) \qquad \psi^{+}(y) = -\mu d_{2}\frac{\left(2\pi k\right)^{\frac{1}{2}}}{m}\eta M_{0} j(-y), \nonumber\\ \alpha^{-}(y) &= d_{2}\frac{\left(2\pi k\right)^{\frac{1}{2}}}{m}\gamma  M_{0} j(y) \qquad \alpha^{+}(y) = -\mu d_{2}\frac{\left(2\pi k\right)^{\frac{1}{2}}}{m}\gamma M_{0} j(-y), \nonumber\end{align}}
  where
  {\begin{align} \eta &=  J_{0}(k r), \qquad
   \gamma = -k J_{1}(k r).\nonumber
\end{align}}
Then, {one can verify that}
{\begin{align}
    &f_{2}(-y) = \mu f_{1}(y) + (\mu-1)\Gamma, \qquad \qquad f_{5}(-y) = \mu f_{4}(y) - (\mu-1)\widetilde{\Gamma},& & &\nonumber\\
    \intertext{and,}
    &f_{1}(y) = d^{2}_{2}\frac{4\pi k}{m^{2}}M_{0}\left( \gamma\eta \left[-\frac{2 j(t)}{D} + \left(\frac{D+t}{D}\right)j'(t)\right]^{y}_{{t=}0} + \frac{\gamma^{2}}{k^{2} r}\left[-\frac{j(t)}{D} + \left(\frac{D+t}{D}\right)j'(t)\right]^{y}_{{t=}0}\right),& \quad &f_{3} = - \Gamma ,&\nonumber\\
    &f_{4}(y) = d^{2}_{2}\frac{2\pi k}{m^{2}} M_{0}\left[\gamma^{2} - k^{2} \eta^{2}\right]\left[-\frac{j(y)}{D} + \left(\frac{D+y}{D}\right)j'(y)\right],&\quad &f_{6} = - \widehat{\Gamma},&\nonumber
    \end{align}
    where
    \begin{align}
    \Gamma &= -d^{2}_{2}\frac{2\pi k}{m^{2}} M_{0}\gamma\eta \frac{\mathcal{M}}{D}\int_{-D}^{0} \left(\frac{D+y}{D}\right)j(y) j'(y) \textnormal{d}y,\nonumber\\
    \widetilde{\Gamma} &= -d^{2}_{2}\frac{2\pi k}{m^{2}}M_{0}\left[\gamma^{2} - k^{2} \eta^{2}\right]\frac{\mathcal{M}}{D}\int_{-D}^{0} \left(\frac{D+y}{D}\right)j(y) j'(y) \textnormal{d}y,\nonumber\\
    \widehat{\Gamma} &= d^{2}_{2}\frac{2\pi k}{m^{2}}\frac{M_{0}}{2D}\left[\frac{\mathcal{M}}{D}\int_{-D}^{0} \left\{\eta^{2} j'(y)^{2} - \gamma^{2} j(y)^{2}\right\}\textnormal{d} y\right].\nonumber
\end{align}
Based on numerical analysis, and for notational simplicity, we will introduce the constant
\begin{align}
C:= d_{2}^{2}M_{0}\sqrt{\frac{k\pi}{2}}\frac{k\mathcal{M}}{D m^{3}}, \nonumber
\end{align}
so that the nonlinear terms $\Omega(\mathcal{F}_{2}, \Re[\mathbf{e}])$ and $\Omega(\mathcal{F}_{2}, \Im[\mathbf{e}])$ can be written as
\begin{align}
    \Omega(\mathcal{F}_{2}, \Re[\mathbf{e}]) &= \frac{1}{m}\left[ \frac{\mathcal{M}}{D}\int_{-D}^{0} j(y) \left[f_{4}(y) - \widetilde{\Gamma}\right]\;\textnormal{d} y - \widehat{\Gamma}\right], \nonumber\\
     &= 2\pi C \sqrt{\frac{2}{k\pi}}\left[\gamma^{2} - k^{2}\eta^{2}\right]\left[-\frac{\mathcal{I}_{1}}{D} + \mathcal{I}_{3}\right] - \frac{\mathcal{M}}{D m}\mathcal{I}_{4}\widetilde{\Gamma} - \frac{1}{m}\widehat{\Gamma}, \nonumber\\
    \Omega(\mathcal{F}_{2}, \Im[\mathbf{e}]) &= -\frac{k}{m}\left[ \frac{\mathcal{M}}{D}\int_{-D}^{0} j(y) \left[f_{1}(y)+ \Gamma\right]\;\textnormal{d} y - \Gamma\right], \nonumber\\
    &= -4 k \pi C \sqrt{\frac{2}{k\pi}} \left[\left(2\gamma\eta + \frac{\gamma^{2}}{k^{2} r}\right)\left(\frac{\mathcal{I}_{4} - \mathcal{I}_{1}}{D}\right) + \left(\gamma\eta + \frac{\gamma^{2}}{k^{2} r}\right)\left(\mathcal{I}_{3} - k\tanh(kD) \mathcal{I}_{4}\right)\right]\nonumber\\
    & \qquad -\frac{k}{m} \left(\frac{\mathcal{M}}{D}\mathcal{I}_{4}-1\right)\Gamma, \nonumber
\end{align}
where we have introduced the integrals,
\begin{align}
    &\mathcal{I}_{1} = \int_{-D}^{0} \left[j(y)\right]^{2}\textnormal{d}y,& \qquad  &\mathcal{I}_{2} = \int_{-D}^{0} \left[j'(y)\right]^{2}\textnormal{d}y,& 
    \nonumber\\ 
    &\mathcal{I}_{3} = \int_{-D}^{0} \left(\frac{D+y}{D}\right)j(y)j'(y)\textnormal{d}y,& \qquad &\mathcal{I}_{4} = \int_{-D}^{0} j(y)\textnormal{d}y.& \nonumber
\end{align}
With these integrals, we can write
\begin{align}
    \Gamma &= -2\pi m C \sqrt{\frac{2}{k\pi}}\gamma\eta \mathcal{I}_{3},\qquad
    \widetilde{\Gamma} = -2\pi m C \sqrt{\frac{2}{k\pi}}\left[\gamma^{2} - k^{2} \eta^{2}\right]\mathcal{I}_{3},\qquad 
    \widehat{\Gamma} = 2\pi m C \sqrt{\frac{2}{k\pi}}\frac{1}{2D}\left[\eta^{2} \mathcal{I}_{2} - \gamma^{2}\mathcal{I}_{1}\right],\nonumber
\end{align}
and so,
\begin{align}
    \Omega(\mathcal{F}_{2}, \Re[\mathbf{e}]) &= 2\pi C \sqrt{\frac{2}{k\pi}}\left(\left[\gamma^{2} - k^{2}\eta^{2}\right]\left[-\frac{\mathcal{I}_{1}}{D} + \mathcal{I}_{3} + \left(\frac{\mathcal{M}}{D} \mathcal{I}_{4}\right)\mathcal{I}_{3} \right] - \frac{1}{2D}\left[\eta^{2} \mathcal{I}_{2} - \gamma^{2}\mathcal{I}_{1}\right]\right), \nonumber\\
    \Omega(\mathcal{F}_{2}, \Im[\mathbf{e}]) &= 2\pi k C \sqrt{\frac{2}{k\pi}}\left[\gamma \eta \left(\frac{\mathcal{M}}{D}\mathcal{I}_{4}-1\right)\mathcal{I}_{3}- 2\left(2\gamma\eta + \frac{\gamma^{2}}{k^{2} r}\right)\left(\frac{\mathcal{I}_{4} - \mathcal{I}_{1}}{D}\right)\right.\nonumber\\
    & \qquad \qquad \qquad  \qquad \left. -2\left(\gamma\eta + \frac{\gamma^{2}}{k^{2} r}\right)\left(\mathcal{I}_{3} - k\tanh(kD) \mathcal{I}_{4}\right)\right].\nonumber
\end{align}}
We note that 
\begin{align}
    &\int_{0}^{r} s J_{0}(k s) k\eta^{2}(s)\textnormal{d}s = {\frac{2}{\pi k \sqrt{3}}+ \textnormal{O}(r^{-\frac{1}{2}}),}& \qquad &\int_{0}^{r} s J_{0}(k s) \frac{\gamma^{2}(s)}{k}\textnormal{d}s = {\frac{1}{\pi k \sqrt{3}}+ \textnormal{O}(r^{-\frac{1}{2}}),}&\nonumber\\
    &\int_{0}^{r} s J_{1}(k s) \eta(s)\gamma(s)\textnormal{d}s = {-\frac{1}{\pi k \sqrt{3}}+ \textnormal{O}(r^{-\frac{1}{2}}),}& \qquad &\int_{0}^{r}  J_{1}(k s) \frac{\gamma^{2}(s)}{k^{2}}\textnormal{d}s = {\frac{3}{2 \pi k\sqrt{3}}+ \textnormal{O}(r^{-\frac{1}{2}}),}&\nonumber
\end{align}
for large values of $r$, where the $\textnormal{O}(r^{-\frac{1}{2}})$ term is an estimate found in \cite{mccalla2013spots}, improving the $\textnormal{o}(1)$ term seen in \cite{watson1944bessel,lloyd2009localized}. {Then, after formally setting the $\textnormal{O}(r^{-\frac{1}{2}})$ terms to zero, we can apply the above integrals to find, 
\begin{align}
    \int_{0}^{r} s J_{0}(k s)\Omega(\mathcal{F}_{2}, \Re[\mathbf{e}]) \textnormal{d} s &= C \sqrt{\frac{2}{3 k\pi}}\left(2\left[\frac{\mathcal{I}_{1}}{D} - \mathcal{I}_{3} - \left(\frac{\mathcal{M}}{D} \mathcal{I}_{4}\right)\mathcal{I}_{3} \right] - \frac{1}{k^{2}D}\left[2\mathcal{I}_{2} - k^{2}\mathcal{I}_{1}\right]\right), \nonumber\\
    \int_{0}^{r} s J_{1}(k s)\Omega(\mathcal{F}_{2}, \Im[\mathbf{e}]) \textnormal{d} s &= C \sqrt{\frac{2}{3 k\pi}}\left[-2 \left(\frac{\mathcal{M}}{D}\mathcal{I}_{4}-1\right)\mathcal{I}_{3}+ 2\left(\frac{\mathcal{I}_{4} - \mathcal{I}_{1}}{D}\right) -2\left(\mathcal{I}_{3} - k\tanh(kD) \mathcal{I}_{4}\right)\right],\nonumber
\end{align}
and so,
\begin{align}
    \widetilde{d}_{3}&= \sqrt{\frac{k\pi}{2}}\int_{0}^{r} \left\{s J_{1}(k s)\Omega(\mathcal{F}_{2}, \Im[\mathbf{e}]) - s J_{0}(k s)\Omega(\mathcal{F}_{2}, \Re[\mathbf{e}])\right\}\textnormal{d} s,\nonumber\\
    &= C \frac{1}{\sqrt{3}}\left[\frac{2}{k^{2}D}\mathcal{I}_{2} + \left(\frac{2\mathcal{I}_{4} - 5\mathcal{I}_{1}}{D}\right) + 2\mathcal{I}_{3} + 2k\tanh(kD)\mathcal{I}_{4}\right] .\nonumber
\end{align}
Hence, we see that
\begin{align}
    \widetilde{d}_{3}(r_{0}) &= \left[\nu + \textnormal{O}\left(r_{0}^{-\frac{1}{2}}\right)\right]d_{2}^{2} + \textnormal{O}_{r_{0}}\left(|\mathbf{d}|_{1}\left[|\varepsilon| + |\mathbf{d}_{2}|_{1}\right] + |d_{2}|^{3}\right),\nonumber
\end{align}
where
\begin{align}
    \nu &:= M_{0}\sqrt{\frac{k\pi}{2}}\frac{k\mathcal{M}}{D m^{3}} \frac{1}{\sqrt{3}}\left[\frac{2}{k^{2}D}\mathcal{I}_{2} + \left(\frac{2\mathcal{I}_{4} - 5\mathcal{I}_{1}}{D}\right) + 2\mathcal{I}_{3} + 2k\tanh(kD)\mathcal{I}_{4}\right] .\label{nu:defn;int}
\end{align}
Computing the integrals $\mathcal{I}_{j}$ for $j=1,2,3,4$, we find,
\begin{align}
    &\mathcal{I}_{1} = \frac{\tanh(kD) + kD \;\textnormal{sech}^{2}(kD)}{2k},& \qquad  &\mathcal{I}_{2} =  k\tanh(kD) - k^{2}\mathcal{I}_{1},& 
    \nonumber\\ 
    &\mathcal{I}_{3} = \frac{1}{2} - \frac{1}{2D}\mathcal{I}_{1} ,& \qquad &\mathcal{I}_{4} = \frac{\tanh(kD)}{k}& \nonumber
\end{align}
Hence, we have found that,
\begin{align}
    \nu &=M_{0}\sqrt{\frac{3 k \pi}{2}}\frac{k \mathcal{M}}{ D m^{3}}\left[1- 2\;\textnormal{sech}^{2}(k D)\right].\label{app:nu;defn}
\end{align}
Then, $\nu>0$ for $k D>x_{c}$ and $\nu<0$ for $k D<x_{c}$, where $x_{c}:=\textnormal{arccosh}(\sqrt{2}) = \log(\sqrt{2} +1)$. }

\bibliographystyle{abbrv}


\begin{thebibliography}{10}

\bibitem{abramowitz1972handbook}
M.~Abramowitz and I.~Stegun.
\newblock {\em Handbook of Mathematical Functions with Formulas, Graphs, and
  Mathematical Tables}.
\newblock Dover, New York, 1972.

\bibitem{aulbach2003foliation}
B.~Aulbach and T.~Wanner.
\newblock Invariant foliations for {C}arath\'eodory type differential equations
  in {B}anach spaces.
\newblock In {\em Advances in stability theory at the end of the 20th century},
  volume~13 of {\em Stability Control Theory Methods Appl.}, pages 1--14.
  Taylor \& Francis, London, 2003.

\bibitem{beck2019exponential}
M.~Beck, G.~Cox, C.~Jones, Y.~Latushkin, and A.~Sukhtayev.
\newblock Exponential dichotomies for elliptic {PDE} on radial domains, 2019.
\newblock arXiv preprint.

\bibitem{blyth2014solitary}
M.~Blyth and E.~P\u{a}r\u{a}u.
\newblock Solitary waves on a ferrofluid jet.
\newblock {\em J. Fluid Mech.}, 750:401--420, 2014.

\bibitem{Bohlius2007Adjoint}
S.~Bohlius, H.~Pleiner, and H.~Brand.
\newblock Solution of the adjoint problem for instabilities with a deformable
  surface: {R}osensweig and {M}arangoni instability.
\newblock {\em Phys. Fluids}, 19(9):094103, 2007.

\bibitem{bohlius2011amplitude}
S.~Bohlius, H.~Pleiner, and H.~Brand.
\newblock The amplitude equation for the {R}osensweig instability in magnetic
  fluids and gels.
\newblock {\em Progr. Theor. Exp. Phys.}, 125(1):1--46, 2011.

\bibitem{buffoni1996plethora}
B.~Buffoni, M.~Groves, and J.~Toland.
\newblock A plethora of solitary gravity-capillary water waves with nearly
  critical {B}ond and {F}roude numbers.
\newblock {\em Philos. Trans. Roy. Soc. London Ser. A}, 354(1707):575--607,
  1996.

\bibitem{cao2014formation}
Y.~Cao and Z.~Ding.
\newblock Formation of hexagonal pattern of ferrofluid in magnetic field.
\newblock {\em J. Magn. Magn. Mater.}, 355:93--99, 2014.

\bibitem{castillo2019extended}
C.~Castillo-Pinto, M.~Clerc, and G.~Gonz{\'a}lez-Cort{\'e}s.
\newblock Extended stable equilibrium invaded by an unstable state.
\newblock {\em Sci. Rep.}, 9(1):1--8, 2019.

\bibitem{chen2019center}
R.~Chen, S.~Walsh, and M.~Wheeler.
\newblock Center manifolds without a phase space for quasilinear problems in
  elasticity, biology, and hydrodynamics, 2019.
\newblock arXiv preprint.

\bibitem{chow1990smooth}
S.-N. Chow, X.-B. Lin, and K.~Lu.
\newblock Smooth invariant foliations in infinite-dimensional spaces.
\newblock {\em J. Differential Equations}, 94(2):266--291, 1991.

\bibitem{cowley1967interfacial}
M.~Cowley and R.~Rosensweig.
\newblock The interfacial stability of a ferromagnetic fluid.
\newblock {\em J. Fluid Mech.}, 30(4):671--688, 1967.

\bibitem{faye2013localized}
G.~Faye, J.~Rankin, and D.~Lloyd.
\newblock Localized radial bumps of a neural field equation on the {E}uclidean
  plane and the {P}oincar\'e disc.
\newblock {\em Nonlinearity}, 26(2):437--478, 2013.

\bibitem{friedrichs2001pattern}
R.~Friedrichs and A.~Engel.
\newblock Pattern and wave number selection in magnetic fluids.
\newblock {\em Phys. Rev. E}, 64(2):021406, 2001.

\bibitem{gailitis1977formation}
A.~Gail{\=i}tis.
\newblock Formation of the hexagonal pattern on the surface of a ferromagnetic
  fluid in an applied magnetic field.
\newblock {\em J. Fluid Mech.}, 82(3):401--413, 1977.

\bibitem{Glicksberg1976Rouche}
I.~Glicksberg.
\newblock A remark on {R}ouch\'{e}'s theorem.
\newblock {\em Amer. Math. Monthly}, 83(3):186--187, 1976.

\bibitem{Groves2019Complete}
M.~Groves.
\newblock personal communication, 2019.

\bibitem{Groves2018periodic}
M.~Groves and J.~Horn.
\newblock Small-amplitude static periodic patterns at a fluid-ferrofluid
  interface.
\newblock {\em Proc. R. Soc. A.}, 474(2216):20180038, 2018.

\bibitem{groves2017pattern}
M.~Groves, D.~Lloyd, and A.~Stylianou.
\newblock Pattern formation on the free surface of a ferrofluid: {S}patial
  dynamics and homoclinic bifurcation.
\newblock {\em Phys. D}, 350:1--12, 2017.

\bibitem{Groves2018Jet}
M.~Groves and D.~Nilsson.
\newblock Spatial dynamics methods for solitary waves on a ferrofluid jet.
\newblock {\em J. Math. Fluid Mech.}, 20(4):1427--1458, 2018.

\bibitem{Haragus2011Bifurcation}
M.~Haragus and G.~Iooss.
\newblock {\em Local bifurcations, center manifolds, and normal forms in
  infinite-dimensional dynamical systems}.
\newblock Universitext. Springer-Verlag London, Ltd., London; EDP Sciences, Les
  Ulis, 2011.

\bibitem{Hill2021Thesis}
D.~Hill.
\newblock {\em Localised radial patterns on the free surface of a ferrofluid}.
\newblock {PhD} thesis, University of Surrey, Guildford, In Preparation.

\bibitem{Horn2015Masters}
J.~Horn.
\newblock {\em Bifurcation theory for static periodic patterns at a
  fluid-ferrofluid interface}.
\newblock Master's thesis, Universit{\"a}t des Saarlandes, Saarbr{\"u}cken,
  2015.

\bibitem{Iooss1992Small}
G.~Iooss and K.~Kirchg\"{a}ssner.
\newblock Water waves for small surface tension: an approach via normal form.
\newblock {\em Proc. Roy. Soc. Edinburgh Sect. A}, 122(3-4):267--299, 1992.

\bibitem{Knieling2007Growth}
H.~Knieling, R.~Richter, I.~Rehberg, G.~Matthies, and A.~Lange.
\newblock Growth of surface undulations at the {R}osensweig instability.
\newblock {\em Phys. Rev. E}, 76:066301, Dec 2007.

\bibitem{knobloch2008spatially}
E.~Knobloch.
\newblock Spatially localized structures in dissipative systems: {O}pen
  problems.
\newblock {\em Nonlinearity}, 21(4):T45--T60, 2008.

\bibitem{KOPELL1981Target}
N.~Kopell and L.~N. Howard.
\newblock Target pattern and spiral solutions to reaction-diffusion equations
  with more than one space dimension.
\newblock {\em Advances in Applied Mathematics}, 2(4):417 -- 449, 1981.

\bibitem{lavrova2008numerical}
O.~Lavrova, G.~Matthies, and L.~Tobiska.
\newblock Numerical study of soliton-like surface configurations on a magnetic
  fluid layer in the {R}osensweig instability.
\newblock {\em Commun. Nonlinear Sci. Numer. Simul.}, 13(7):1302--1310, 2008.

\bibitem{lavrova2016modeling}
O.~Lavrova, V.~Polevikov, and L.~Tobiska.
\newblock Modeling and simulation of magnetic particles diffusion in a
  ferrofluid layer.
\newblock {\em Magnetohydrodynamics}, 52(4):417--430, 2016.

\bibitem{lloyd2015homoclinic}
D.~Lloyd, C.~Gollwitzer, I.~Rehberg, and R.~Richter.
\newblock Homoclinic snaking near the surface instability of a polarisable
  fluid.
\newblock {\em J. Fluid Mech.}, 783:283--305, 2015.

\bibitem{lloyd2009localized}
D.~Lloyd and B.~Sandstede.
\newblock Localized radial solutions of the {S}wift-{H}ohenberg equation.
\newblock {\em Nonlinearity}, 22(2):485--524, 2009.

\bibitem{Mccalla2010snaking}
S.~McCalla and B.~Sandstede.
\newblock Snaking of radial solutions of the multi-dimensional
  {S}wift–{H}ohenberg equation: A numerical study.
\newblock {\em Physica D}, 239(16):1581 -- 1592, 2010.

\bibitem{mccalla2013spots}
S.~McCalla and B.~Sandstede.
\newblock Spots in the {S}wift-{H}ohenberg equation.
\newblock {\em SIAM J. Appl. Dyn. Syst.}, 12(2):831--877, 2013.

\bibitem{mcquighan2014oscillons}
K.~McQuighan and B.~Sandstede.
\newblock Oscillons in the planar {G}inzburg-{L}andau equation with {$2:1$}
  forcing.
\newblock {\em Nonlinearity}, 27(12):3073--3116, 2014.

\bibitem{Mielke1986Reduction}
A.~Mielke.
\newblock A reduction principle for nonautonomous systems in
  infinite-dimensional spaces.
\newblock {\em J. Differential Equations}, 65(1):68--88, 1986.

\bibitem{mielke2006hamiltonian}
A.~Mielke.
\newblock {\em Hamiltonian and {L}agrangian flows on center manifolds}, volume
  1489.
\newblock Springer-Verlag, Berlin, 1991.

\bibitem{Pomeau1985Axisymmetric}
Y.~Pomeau, S.~Zaleski, and P.~Manneville.
\newblock Axisymmetric cellular structures revisited.
\newblock {\em Z. angew. Math. Phys.}, 36(3):367-- 394, 1985.

\bibitem{Reimann2003Oscillatory}
B.~Reimann, R.~Richter, I.~Rehberg, and A.~Lange.
\newblock Oscillatory decay at the {R}osensweig instability: Experiment and
  theory.
\newblock {\em Phys. Rev. E}, 68:036220, Sep 2003.

\bibitem{Richter2011Mountains}
R.~Richter.
\newblock Mag(net)ic liquid mountains.
\newblock {\em Europhys. News}, 42(3):17--19, 2011.

\bibitem{richter2005two}
R.~Richter and I.~Barashenkov.
\newblock Two-dimensional solitons on the surface of magnetic fluids.
\newblock {\em Phys. Rev. Lett.}, 94(18):184503, 2005.

\bibitem{rosensweig1987magnetic}
R.~Rosensweig.
\newblock Magnetic fluids.
\newblock {\em Annu. Rev. Fluid Mech.}, 19(1):437--461, 1987.

\bibitem{rosensweig2013ferrohydrodynamics}
R.~Rosensweig.
\newblock {\em Ferrohydrodynamics}.
\newblock Dover Publications, New York, 2013.

\bibitem{Sandstede1997Convergence}
B.~Sandstede.
\newblock Convergence estimates for the numerical approximation of homoclinic
  solutions.
\newblock {\em IMA J. Numer. Anal.}, 17(3):437--462, 1997.

\bibitem{scheel2003radially}
A.~Scheel.
\newblock Radially symmetric patterns of reaction-diffusion systems.
\newblock {\em Mem. Amer. Math. Soc.}, 165(786):viii+86, 2003.

\bibitem{silber1988pattern}
M.~Silber and E.~Knobloch.
\newblock Pattern selection in ferrofluids.
\newblock {\em Phys. D}, 30(1-2):83--98, 1988.

\bibitem{spyropoulos2019spike}
A.~N. Spyropoulos, A.~G. Papathanasiou, and A.~G. Boudouvis.
\newblock The 2-3-4 spike competition in the {R}osensweig instability.
\newblock {\em J. Fluid Mech.}, 870:389–404, 2019.

\bibitem{torres2014recent}
I.~Torres-Diaz and C.~Rinaldi.
\newblock Recent progress in ferrofluids research: novel applications of
  magnetically controllable and tunable fluids.
\newblock {\em Soft Matter}, 10(43):8584--8602, 2014.

\bibitem{twombly1983bifurcating}
E.~Twombly and J.~Thomas.
\newblock Bifurcating instability of the free surface of a ferrofluid.
\newblock {\em SIAM J. Math. Anal.}, 14(4):736--766, 1983.

\bibitem{vandenberg2015Rigorous}
J.~van~den Berg, C.~Groothedde, and J.~Williams.
\newblock Rigorous computation of a radially symmetric localized solution in a
  {G}inzburg-{L}andau problem.
\newblock {\em SIAM J. Appl. Dyn. Syst.}, 14(1):423--447, 2015.

\bibitem{vanderbauwhede1992center}
A.~Vanderbauwhede and G.~Iooss.
\newblock Center manifold theory in infinite dimensions.
\newblock In {\em Dynamics reported: expositions in dynamical systems},
  volume~1, pages 125--163. Springer, Berlin, 1992.

\bibitem{Walter1998ODEs}
W.~Walter.
\newblock {\em Ordinary differential equations}, volume 182 of {\em Graduate
  Texts in Mathematics}.
\newblock Springer-Verlag, New York, 1998.
\newblock Translated from the sixth German (1996) edition by Russell Thompson,
  Readings in Mathematics.

\bibitem{watson1944bessel}
G.~Watson.
\newblock {\em A treatise on the theory of {B}essel functions}.
\newblock Cambridge University Press, Cambridge, 1944.

\bibitem{Zaitsev1970Nature}
V.~{Zaitsev} and M.~{Shliomis}.
\newblock Nature of the instability of the interface between two liquids in a
  constant field.
\newblock {\em Sov. Phys. Dokl.}, 14:1001, Apr 1970.

\end{thebibliography}
\end{document}